\documentclass[11pt,letterpaper]{article}
\usepackage{fontspec}
\usepackage{amsmath,amssymb,amsthm,mathtools,mathrsfs}
\usepackage{tikz}
\usepackage{graphicx,array,booktabs,enumitem,needspace,longtable,float}
\usepackage[margin=1in,headheight=0pt,headsep=0pt,footskip=30pt]{geometry}
\usepackage[titles]{tocloft}

\usepackage[noadjust]{cite}
\usepackage[hidelinks]{hyperref}
\hypersetup{pdftitle={A Simplified Proof of Finite-time Blowup of Navier-Stokes with Smooth Forcing},pdfauthor={Zhen Lei and Xiao Ren},pdfsubject={Axisymmetric background, oscillatory stress realization and localization}}
\numberwithin{equation}{section}
\newtheorem{theorem}{Theorem}[section]
\newtheorem{lemma}[theorem]{Lemma}
\newtheorem{proposition}[theorem]{Proposition}
\newtheorem{corollary}[theorem]{Corollary}
\theoremstyle{definition}
\newtheorem{definition}[theorem]{Definition}
\theoremstyle{remark}
\newtheorem{remark}[theorem]{Remark}
\newcommand{\NS}{\operatorname{NS}}
\newcommand{\dd}{\,\mathrm d}
\newcommand{\RR}{\mathbb R}
\newcommand{\supp}{\operatorname{supp}}
\newcommand{\Div}{\operatorname{div}}
\newcommand{\curl}{\operatorname{curl}}
\setlist[enumerate]{leftmargin=2em,itemsep=4pt,topsep=5pt}
\setlist[itemize]{leftmargin=1.8em,itemsep=3pt,topsep=4pt}
\hypersetup{bookmarksdepth=3}
\makeatletter
\renewcommand{\@seccntformat}[1]{%
  \@ifundefined{#1@cntformat}
    {\csname the#1\endcsname\quad}
    {\csname #1@cntformat\endcsname}}
\newcommand{\section@cntformat}{Section~\thesection.\quad}
\makeatother
\newcommand{\pulsebridge}[1]{\ifcsname pulsebridge@#1\endcsname
  \csname pulsebridge@#1\endcsname\par\fi}
\newcommand{\newpulsebridge}[2]{\expandafter\def\csname pulsebridge@#1\endcsname{#2}}
\newpulsebridge{2}{For a zero-mean pulse \(w\), the mean tangential force contains \(+\mathcal D\langle w^r(w^\theta,w^z)\rangle\). To cancel the background force, we prescribe \(\langle w^r(w^\theta,w^z)\rangle=\mathsf T_B\) at leading order.}

\newpulsebridge{2.1.4}{A label reaching \(Z=\pm1\) needs an extension for Fourier truncation. We extend only in \(Z\), keeping \(\lambda,R\) and the auxiliary variables fixed. We shorten the extension as \(k\) increases so that the label supports remain disjoint.}
\newpulsebridge{2.2}{We decompose \(\mathsf T_B\) into its components parallel and perpendicular to the shear. The strict cone in \eqref{intro:uniform-cone} bounds the ratio of these components. Using the native clocks in \eqref{pulse:eq:2.22}, we construct two phases. The directions of their averaged fluxes bound a cone containing the prescribed stress.}
\newpulsebridge{2.2.4}{The same scalar envelope estimates the amplitude on both sides of the clock midpoint and the propagator from any earlier clock time to a later one.}
\newpulsebridge{2.2.7}{The averaged flux depends on the mass of the pulse and its direction during the part of the clock interval where this mass is concentrated. We bound both quantities before applying the strict margin in \eqref{intro:uniform-cone}.}
\newpulsebridge{3}{After the primary pulses cancel the leading stress, we prepare the remaining force for the source equation. An auxiliary velocity removes the leading oscillation in its mean tangential component. A pure swirl on the mixed window of Section~\ref{mid:section} cancels three integral defects while preserving both velocity moments. A pressure correction then removes the angular mean of the radial force.}
\newpulsebridge{3.2}{If the primitive vanishes at the inner edge, its value at the outer edge is the full radial integral. A compensator on the fixed middle interval cancels this trace. After auxiliary evaluation the primitive includes a torus translation between the two radii.}
\newpulsebridge{3.3}{The radial pressure correction leaves one integral defect, and the two mixed tangential fluxes give two further defects. We use an axisymmetric swirl to change these three quantities while retaining zero angular momentum. This last constraint is the row with weight \(r^2\) in the moment system.}
\newpulsebridge{3.4}{We first take the curl, then substitute the divergence-free velocity in the physical equation. We thus obtain the terms in the force coming from derivatives of the amplitude, the phase and the cylindrical frame.}
\newpulsebridge{3.6}{We solve the nonlinear four-row moment system including the quadratic terms. We fix this corrected velocity before linearizing the force in Section~\ref{pulse:sec:4}.}
\newpulsebridge{4}{The initialized force has a compatible angular mean and a nonzero oscillatory part. For a prescribed trial source we construct velocity and pressure corrections. Their force residual consists of a linear term and a quadratic term, estimated below in the same source coordinates.}
\newpulsebridge{4.1}{The trial force may contain every carrier harmonic. Projection onto the plane perpendicular to the phase gradient and the additional viscous damping at higher harmonics give a bound uniform in the harmonic number.}
\newpulsebridge{4.3}{After angular averaging and removal of the Haar mean, the leading time derivative has a periodic inverse on the auxiliary torus. A curl correction enforces incompressibility for its radial and axial velocity components.}
\newpulsebridge{4.4}{We specify a unit for each source component and estimate the successive corrections in these units. We apply the corrections in this order to preserve the radial moments.}
\newpulsebridge{5}{The force estimates give a fixed gain in \(\varepsilon_k\) and a fixed loss of derivatives. We choose the Fourier cutoff so that a positive power of \(\varepsilon_k\) remains after differentiation. We then compare solutions on neighboring windows, where the Fourier truncations are taken on different covers.}
\newpulsebridge{5.2}{On an overlap we differentiate the physical sum before restricting it to an output row. The restricted force therefore includes the terms obtained by differentiating the partition functions.}
\newpulsebridge{5.4}{We retain the central component of each window solution. To estimate the force on a fixed band, we compare the contributing components with those of one window solution in coefficient coordinates, then apply physical evaluation.}

\makeatletter
\def\PulseSuffix#1.#2\@nil{#2}
\newcommand{\PulseResult}[1]{\phantomsection\edef\@currentlabel{\thesection.\PulseSuffix#1\@nil}%
  \label{pulse:thm:#1}\@currentlabel}
\newcommand{\NumberedPulseResult}[1]{\refstepcounter{theorem}\label{pulse:thm:#1}\thetheorem}
\makeatother

\begin{document}
\begin{center}
{\Large\bfseries A Simplified Proof of Finite-time Blowup of\\
Navier-Stokes with Smooth Forcing\par}
\vspace{1em}
\begingroup\renewcommand{\thefootnote}{\fnsymbol{footnote}}
{\large Zhen Lei\footnotemark[1]\qquad Xiao Ren\footnotemark[2]\par}
\footnotetext[1]{\raggedright Center for Applied Mathematics \& School of Mathematical
Sciences, Fudan University, Shanghai 200433, P.~R.~China.
Email: zlei@fudan.edu.cn.}
\footnotetext[2]{\raggedright Center for Applied Mathematics, Fudan University,
Shanghai 200433, P.~R.~China. Email: xren@fudan.edu.cn.}
\endgroup

\vspace{.6em}27 November 2026
\end{center}
\begin{abstract}
We give a self-contained, simplified proof of OpenAI's finite-time blowup
result for the three-dimensional incompressible Navier--Stokes equations
with smooth forcing \cite{OpenAI-NS}. The solution starts from rest, has
bounded kinetic energy and finite dissipation, and develops a Type~II
singularity at \((t,x)=(1,0)\). The force is smooth and compactly supported,
and all its space--time derivatives vanish to arbitrary order there.
Our quantitative decay estimate gives a common constant and radius at each
derivative order and decay rate. For fixed geometry and margins, both depend
only on finite-order bounds on the normalized background fields and weighted stress.
A restart near blowup followed by rescaling gives the same Type~II rate with
smooth nonzero initial data and a uniformly small force.

We obtain \(E_B=0\) by a finite assembly that replaces the successive
corrections and their infinite sum. Our new construction uses an exact mixed exterior flow to balance axial and angular energy
and eliminates the main-pulse, steepening and waiting stages.
An exact power law and two radially analytic connections join the unforced
core to this exterior. They match all five
moments and give strictly admissible stress.
Two pulse families match this stress in the mean. After initialization, our new formulation replaces the successive
pulse corrections of the original construction by a fixed-point problem for the source
coefficients on each window of neighboring scale bands. At sufficiently fine scales,
Fourier truncation makes the quadratic source map contractive, and its fixed point
determines the velocity and pressure corrections on that window.
Comparing source coefficients on shared bands near the window centers controls
the assembly error. Since only the pulse residual contributes, this bound and
the cutoff and tail estimates give the local force estimate.
\end{abstract}
\tableofcontents

\section{Introduction}\label{intro:section}

In September 2026, OpenAI announced finite-time blowup for the
three-dimensional incompressible Navier--Stokes equations with smooth
forcing \cite{OpenAI-NS},
\begin{equation}\label{intro:NS}
 \partial_tu+(u\cdot\nabla)u+\nabla p=\Delta u+f,
 \qquad \nabla\cdot u=0.
\end{equation}
The solution starts from rest and develops unbounded velocity with
uniformly bounded kinetic energy. The construction combines a singular
axisymmetric background with oscillatory perturbations that turn its
singular stress into smooth forcing.

Our earlier papers \cite{Part-I,Part-II} give a two-part exposition of
this construction. Here we give a self-contained proof and simplify both
the background and the pulse realization. We use two ideas of OpenAI
\cite{OpenAI-NS}. Cancelling five moment defects at fixed axis pressure
keeps profile changes localized. The instability and cone conditions let
us match the stress by a positive combination of two averaged pulse fluxes.

\subsection{Background and main results}

In cylindrical coordinates $(r,\theta,z)$, introduce the similarity
variables and leading normalization
\begin{equation}\label{intro:coordinates}\csname ltx@label\endcsname{intro:lambda}
\begin{alignedat}{4}
 R&=\frac{r^2}{2\lambda^2},\qquad&
 Z&=\frac{z}{\lambda^{1-\delta}},\qquad&
 1-t&=\lambda^2(1-Z^2),\qquad& 0&<\delta<\frac1{200},\\
 u_{(0)}^\theta&=\frac{U^\theta}{\lambda^{1+\delta}},&
 u_{(0)}^z&=\frac{U^z}{\lambda^{1+\delta}},&
 u_{(0)}^r&=\frac{U^r}{\lambda},&
 p_{(0)}&=\frac{P}{\lambda^{2+2\delta}}.
\end{alignedat}
\end{equation}
Denote
\[
 \NS(u,p)=\partial_tu+(u\cdot\nabla)u+\nabla p-\Delta u.
\]
For $\mathsf T=(\mathsf T^\theta,\mathsf T^z)$, define
\[
 \mathcal D\mathsf T
 =\left(\partial_r+\frac2r\right)\mathsf T^\theta e_\theta
  +\left(\partial_r+\frac1r\right)\mathsf T^z e_z.
\]
For the complete background define
\begin{equation}\label{intro:physical-shear}
\begin{alignedat}{2}
 U_B^\theta&=\lambda^{1+\delta}u_B^\theta,\qquad\qquad&
 U_B^z&=\lambda^{1+\delta}u_B^z,\\
 F_B&=\frac{U_B^\theta}{\sqrt{2R}},\qquad\qquad&
 \mathcal T_B&=\lambda^{2+\delta}\mathsf T_B,\\
 \mathcal S_B&=(2R\partial_RF_B,
       \sqrt{2R}\partial_RU_B^z),\qquad\qquad&
 \mathcal S_B^\perp&=(-\mathcal S_B^z,\mathcal S_B^\theta),\\
 \kappa_B&=-\frac{|\mathcal S_B|^2}{F_B\mathcal S_B^\theta}.
\end{alignedat}
\end{equation}
Complete profiles depend on $(\lambda,R,Z)$; leading profiles depend
on $(R,Z)$.
Let $\mathscr D^{\le m}$ denote compositions
of at most $m$ operators from
$\{\lambda\partial_\lambda,R\partial_R,\partial_Z\}$.
On a fixed stress annulus $R_{\rm in}<R<R_e$ set
\begin{equation}\label{intro:edge-weight}
\begin{aligned}
 d_e&=\min\{1,\log(R/R_{\rm in}),\log(R_e/R)\},\\
 \omega&=\exp[-c_a/\log^2(R/R_{\rm in})-c_h/\log^2(R_e/R)],
 \qquad c_a,c_h>0.
\end{aligned}
\end{equation}
The quantitative flat-edge bounds are
\begin{equation}\label{intro:stress-bounds}
 |\mathcal T_B|\ge c\omega,\qquad
 |\mathscr D^{\le m}\mathcal T_B|
       \le C_m\omega d_e^{-q_m},
\end{equation}
for fixed positive $c$ and finite $C_m,q_m$.

\begin{definition}[Uniform admissibility]
We require fixed positive upper and lower bounds for $F_B$ and
$-\mathcal S_B^\theta$ on the closed annulus. On the open annulus
$R_{\rm in}<R<R_e$, we require a fixed $c_{\rm cone}\in(0,1)$ such that
\begin{equation}\label{intro:uniform-cone}
 \kappa_B\ge2+c_{\rm cone},\qquad
 -\mathcal T_B\cdot\frac{\mathcal S_B}{|\mathcal S_B|}
                 \ge c_{\rm cone}|\mathcal T_B|,\qquad
 \sqrt{\frac{\kappa_B-2}{2}}
 \left|\frac{\mathcal T_B\cdot\mathcal S_B^\perp}
                  {\mathcal T_B\cdot\mathcal S_B}\right|
                  \le1-c_{\rm cone}.
\end{equation}
The first inequality in \eqref{intro:uniform-cone}
is the uniform instability condition and the other two are the cone conditions.
\end{definition}

The bounds \eqref{intro:stress-bounds} give smooth extension by zero at
both stress edges. We write the background equations with remainder \(E_B\) as
\begin{equation}\label{intro:prior-equation}
 \Div u_B=0,\qquad \NS(u_B,p_B)=-\mathcal D\mathsf T_B+E_B.
\end{equation}
Admissibility of a leading profile uses the same instability and cone conditions
with $F=U^\theta/\sqrt{2R}$,
$\mathcal S=(2R\partial_R F,\sqrt{2R}\partial_R U^z)$, and
$\mathcal T=\mathcal I+\mathcal S$.

\begin{theorem}[Exact background with an intermediate power law]
\label{intro:leading}\label{intro:background}\label{fin:background}
Let $P_0(Z)=-3P_*^2/(1+Z^2)^2$.
For sufficiently large pressure amplitude $P_*$, choose the remaining
parameters as in the proof of Theorem~\ref{thm:connection}. Then there exist
$0<\delta<1/200$, $0<R_{\rm in}<R_{ac}<5<R_b<R_{\rm mid}<e^\ell R_{\rm mid}<e^{2\ell}R_{\rm mid}=R_e$,
$h_\infty>0$ and $\epsilon_{\mathrm{em}}>0$ such that the following hold.

\emph{(i) The leading background.}
There is a leading pair $(U^\theta,U^z)$, smooth for $R>0$, $|Z|\le1$.
Both $U^\theta/\sqrt{2R}$ and $U^z$ are smooth up to the axis.
The pair agrees with the inner profile on $[0,R_{ac}]$ and with the
exact mixed exterior on $[R_e,\infty)$. It is
real analytic in $R$ on each of $(R_{ac},R_{\rm mid})$, $(R_{\rm mid},e^\ell R_{\rm mid})$ and
$(e^\ell R_{\rm mid},R_e)$, with all derivatives matching at the joins.
On $[R_{\rm mid},e^\ell R_{\rm mid}]$ it equals
\((1-\epsilon_0)h_\infty R^{-(1+\delta)/2}
(1,\sqrt{(1+\epsilon_{\rm em}^2)/2})\) for some constant \(0<\epsilon_0<1\).
Its axis pressure is $P_0$.

The stress $\mathcal T=\mathcal I+\mathcal S$ is zero for $R\le R_{\rm in}$
and $R\ge R_e$, nonzero and strictly admissible in between, and flat
at both edges. For $R\ge R_e$, all five moments integrated from the
axis by \eqref{pr:moments} agree with \eqref{ue:right-moments}. In particular,
\begin{equation}\label{intro:terminal}
 M^{z\theta}(R,Z)=-\frac{\epsilon_{\mathrm{em}}^2}{2}
                \int_R^\infty(U_e^\theta)^2\,d\rho,
 \qquad M^p=P_e-P_0.
\end{equation}
The primitive $M_e^z$ includes the additive term $Z$.
\par\smallskip
\noindent\emph{(ii) Exact physical assembly.}
There exist smooth real axisymmetric fields
$(u_B,p_B,\mathsf T_B)$ on
$0<\lambda\le\lambda_*$, $R\ge0$, $|Z|\le1$, for some
$\lambda_*>0$, satisfying
\begin{equation}\label{intro:exact}
 \Div u_B=0,\qquad \NS(u_B,p_B)=-\mathcal D\mathsf T_B.
\end{equation}
They equal the exact unforced core for $R\le R_{\rm in}=4/\Lambda$ and
the exact mixed exterior \eqref{ue:physical} for $R\ge R_e$.
For some fixed $j_0>0$, their axis data are
\begin{equation}\label{intro:axis-data}
 \lambda^{2+\delta}\left.\frac{u_B^\theta}{r}\right|_{r=0}=F_0(Z),
 \quad \lambda^{1+\delta}u_B^z(t,0,z)=4Z+j_0,
 \quad \lambda^{2+2\delta}p_B(t,0,z)=P_0(Z).
\end{equation}
Every fixed $\mathscr D$ derivative of the normalized velocity
$(U_B^\theta,U_B^z,\lambda u_B^r)$ is bounded. For every fixed $R_0>0$,
the difference from $(U^\theta,U^z,U^r)$ on $0\le R\le R_0$ is
$O_{m,R_0}(\lambda^{2\delta}(1+|\log\lambda|))$, including all
normalized derivatives of order at most $m$.

The stress is supported on $[R_{\rm in},R_e]$, extends smoothly by zero at
both edges, and satisfies \eqref{intro:stress-bounds}. The background
is uniformly admissible, with the cone bounds \eqref{intro:uniform-cone}.
The physical moment corrections are supported on a fixed interval \(I_2\).
On a disjoint interval \(I_3\Subset(R_{ac},R_{\rm mid})\), the tangential
velocities retain their leading values and the four-row moment matrix is uniformly invertible
(Lemma~\ref{fin:mixed-inverse}).
\end{theorem}

The exact axis trace is
\begin{equation}\label{intro:blowup-trace}
 u_B^z(t,0,0)=j_0(1-t)^{-(1+\delta)/2}.
\end{equation}

\begin{corollary}[Blowup with smooth forcing]\label{intro:blowup}
There exist $0<t_0<1$, a compact set $K\subset\mathbb R^3$, real fields
$u,p\in C^\infty(\mathbb R^3\times[0,1))$, and
$f\in C_c^\infty(\mathbb R^3\times(0,\infty);\mathbb R^3)$ such that
\eqref{intro:NS} holds, $u=p=0$ for $0\le t\le t_0$, and
$\operatorname{supp}u(t)\cup\operatorname{supp}p(t)\subset K$.
Moreover,
\begin{equation}\label{intro:energy}
 \sup_{0\le t<1}\|u(t)\|_{L^2}
       +\|\nabla u\|_{L^2(\mathbb R^3\times(0,1))}<\infty,
 \qquad \|u(t)\|_{L^\infty}\ge c(1-t)^{-(1+\delta)/2}
\end{equation}
for all $t$ sufficiently close to $1$, with $c>0$.
For each \(\varepsilon>0\), we also obtain a solution satisfying
\eqref{intro:energy} with smooth, compactly supported, nonzero initial
velocity and a smooth force such that
\[
 \|f\|_{L^{4/3}(0,1;L^2(\mathbb R^3))}^{2}
 +\|f\|_{L^\infty((0,1)\times\mathbb R^3)}<\varepsilon .
\]
\end{corollary}

For the force in Corollary~\ref{intro:blowup} we also obtain a quantitative
local estimate. Let \(Q_\rho=(1-\rho^2,1)\times B_\rho(0)\). For every
fixed integer \(m\ge0\) and every \(N>0\), there are \(C_{m,N}\) and
\(\rho_{m,N}>0\) such that, for \(0<\rho\le\rho_{m,N}\),
\begin{equation}\label{intro:parabolic}
 \max_{2a+|\beta|\le m}
 \rho^{3+2a+|\beta|}\|\partial_t^a\partial_x^\beta f\|_{L^\infty(Q_\rho)}
 \le C_{m,N}\rho^{N+3}.
\end{equation}
With the geometry and margins fixed as in the proof of
Proposition~\ref{pulse:thm:6.1}, finite-order bounds on the normalized
background fields and weighted stress give uniform \(C_{m,N}\) and \(\rho_{m,N}\).

\subsection{Differences from the original construction}\label{intro:count}

We prescribe the axis pressure and obtain \(E_B=0\) by a finite assembly,
replacing the successive corrections and infinite sum in \cite[Section~8]{Part-I}.
An exact power law and two radially analytic connections match all five moments
with strictly admissible stress (Figure~\ref{fig:short-connection}).

We balance axial and angular energy in the mixed exterior by taking
\[
 U_e^z=\sqrt{\frac{1+\epsilon_{\rm em}^2}{2}}\,U_e^\theta,\qquad
 \int_R^\infty\!\left((U_e^z)^2-\tfrac12(U_e^\theta)^2\right)d\rho
   =\frac{\epsilon_{\rm em}^2}{2}\int_R^\infty (U_e^\theta)^2\,d\rho .
\]
At \(\epsilon_{\rm em}=0\) the energies cancel. Varying \(\epsilon_{\rm em}\) sets
the tail energy at fixed angular pressure and removes the main-pulse, steepening
and waiting stages of \cite[Section~4]{Part-I}. In the auxiliary left profile,
the negative angular contribution to the pole energy over a log interval of length
\(1/\mu\) lets us choose the exterior mixing parameter (Section~\ref{np:section}).

Two pulse families match the stress in the mean. In the original construction,
each correction cycle updates the velocity and pressure and produces the residual
for the next cycle \cite[Proposition~9.6]{OpenAI-NS}. After initialization, we
organize the nonzero-mode, stress, auxiliary, moment and pressure corrections
into a map from a trial source to velocity and pressure increments. The corresponding
fixed-point map in source coordinates has the quadratic form
\(H_k(g)=f_k+E_kg+B_k(g,g)\).

The estimates for this map lose a fixed number of derivatives. We truncate the
Fourier series in all active coefficient variables, including the slow, auxiliary
and carrier variables. Bernstein's inequality then converts the derivative loss
into a power of the cutoff frequency. We choose the cutoff so that a positive
power of the small parameter remains. At sufficiently fine scales, the Fourier-truncated
map is therefore a contraction on each window of neighboring scale bands.
For each band, we take the source coefficients from the window centered on that
band. These coefficients determine a single physical correction, and comparison
on overlaps controls the resulting error in the force.

OpenAI proves all-order decay of the momentum residual for its constructed flow
\cite[Theorem~3.1(iii), (3.4)]{OpenAI-NS}. Our earlier Part~II also proves such
decay \cite[Theorem~1.1(ii)]{Part-II}. Here we use the comparison on overlapping
windows to quantify how the estimate depends on the background.
At fixed geometry and margins, all backgrounds with the same finite-order bounds
on normalized background fields and weighted stress share the constant \(C_{m,N}\) and
radius \(\rho_{m,N}>0\) in \eqref{intro:parabolic}.

\subsection{Related results}\label{intro:literature}

Leray \cite{Leray} constructed global finite-energy weak solutions.
Caffarelli, Kohn and Nirenberg \cite{CKN} proved partial regularity. For a
Morrey condition on the force, see \cite{KukavicaMorrey}. For the Prodi--Serrin
criterion and its $L_t^\infty L_x^3$ endpoint, see \cite{Prodi1959,Serrin62,Serrin63,ESS}.
Backward self-similar profiles are trivial in $L^3(\RR^3)$ \cite{NRS}
or under local energy estimates \cite{Tsai}.
Axisymmetric criteria control the velocity or its swirl
\cite{CSTY08,CSTY09,KNSS,Lei-Zhang11,Lei-Zhang17,Wei,CFZ,OzanskiPalasek}.
Constantin, Ignatova and Vicol \cite{CIV} prove regularity for a shrinking
axisymmetric core under anisotropic bounds on the angular mean and spatially
analytic forcing uniformly bounded in $C_x^2$.
Albritton, Bru\'e and Colombo \cite{ABC,ABCGluing} proved nonuniqueness from
rest with a common force singular at the initial time. For forced singularities
in related equations, see \cite{CM2023,CMIPM,CMZ}. For shear waves and
instability criteria, see \cite{CraikCriminale,SinghSridhar,FV91}, and
for centrifugal instability in vortices and jets, see
\cite{LeibovichStewartson,BillantGallaire2005,BillantGallaire2013}.
For Type~I concentration and possible Type~II growth, see
\cite{BarkerPrange2020,SereginTypeII}. Alp\"oge and Buckmaster \cite{AB-NS}
announce unforced Type~I blowup from smooth data through a backward discretely
self-similar profile. Here the axis trace \eqref{intro:blowup-trace} gives
Type~II blowup with smooth forcing.

There are also many important related results. See, for instance,
\cite{Fefferman,KatoFujita,Kato,FurioliEtAl,Heywood,Solonnikov,Temam,Lin,Lei-Ren-regularity,TaoLocalization}
on regularity and localization,
\cite{JS-forward,Tsai-DSS,BradshawTsaiDSS,ChaeWolfForward,ChaeWolfBack,BradshawTsaiRot,PineauVicol,AlbrittonBarker,SereginSverakTypeI,HouPotential}
on self-similarity and concentration, and
\cite{JS-nonunique,GS-numerical,BV,CL,DLS09,DLS12,DS2017,BCV2022,BVSurvey,IonescuJiaPalasek,CoiculescuPalasek,HouWangYang,BuckmasterCoiculescuLi,VishikI,VishikII,VishikBook}
on instability and nonuniqueness. For singularities in related fluid
equations and models, see
\cite{LiSinai,Tao2016,AlpogeBuckmasterEuler,AlpogeBuckmasterBoussinesq,AlpogeBuckmasterCoiculescu,Elgindi,ElgindiGhoulMasmoudi,DrivasElgindi,ABCEuler,ChenHouI,ChenHouII,ChenHouHuangDG,ChenHouHuangHL,MerleEtAlI,MerleEtAlII,BuckmasterCaoGomez,CaoGomezShiStaff,ChenShkollerVicol}.
For numerical profiles and computer-assisted proofs, see
\cite{WangLaiGomezBuckmaster,WangEtAlDiscovery,GomezSurvey}.

\subsection{Organization}
Section~\ref{pr:section} gives the notation, moment identities and
correction tools. Sections~\ref{profiles:section}--\ref{mid:section}
construct the exact exterior, inner core and intermediate connection.
Section~\ref{pulse:sec:2} constructs the primary pulses, and
Sections~\ref{pulse:sec:3}--\ref{pulse:sec:5} solve the compatible source
equations on finite scale windows. Section~\ref{fin:blowup} proves
the physical force bounds and localizes the fields to obtain
Corollary~\ref{intro:blowup}.
\section{Preliminaries}
\label{pr:section}

We recover the leading radial velocity and pressure from the tangential velocities.
The five radial moments below determine the stress and the conditions for its admissibility.
Notation introduced for a particular proof or calculation applies only there.

\subsection{The equations and coordinates}
For an axisymmetric velocity $u=u^re_r+u^\theta e_\theta+u^ze_z$,
write $\Delta=\partial_r^2+r^{-1}\partial_r+\partial_z^2$. The three
components of $\NS(u,p)$ are
\begin{equation}\label{pr:NS}
\begin{aligned}
\NS^r&=\partial_t u^r+u^r\partial_r u^r+u^z\partial_z u^r-(u^\theta)^2/r+\partial_r p
                 -(\Delta-r^{-2})u^r,\\
\NS^\theta&=\partial_t u^\theta+u^r\partial_r u^\theta+u^z\partial_z u^\theta+u^ru^\theta/r
                 -(\Delta-r^{-2})u^\theta,\\
\NS^z&=\partial_t u^z+u^r\partial_r u^z+u^z\partial_z u^z+\partial_z p-\Delta u^z.
\end{aligned}
\end{equation}
Incompressibility is $\partial_r(ru^r)+\partial_z(ru^z)=0$.

Fix $0<\delta<1/200$ and define
\begin{equation}\label{pr:coordinates}
R=\frac{r^2}{2\lambda^2},\qquad Z=\frac z{\lambda^{1-\delta}},\qquad
1-t=\lambda^2(1-Z^2),\qquad L=1-Z^2,\quad L_\delta=1-\delta Z^2.
\end{equation}
For fixed $t<1$ and $z$, the left side of $\lambda^2-z^2\lambda^{2\delta}=1-t$
increases from zero to infinity on $\lambda>|z|^{1/(1-\delta)}$. Its derivative is $2\lambda(1-\delta Z^2)>0$ since $|Z|<1$.
Thus the solution is unique and, by the implicit function theorem, smooth in $(t,z)$.
For each $t<1$, this identifies $(r,z)\in[0,\infty)\times\RR$ with
$(R,Z)\in[0,\infty)\times(-1,1)$.
At $t=1$, $z\ne0$, one has $\lambda=|z|^{1/(1-\delta)}$,
with $Z=1$ for $z>0$ and $Z=-1$ for $z<0$. The implicit function theorem still applies.
For fixed $r>0,z=0$, $R\to\infty$ as $t\uparrow1$.

Define the leading profiles
\begin{equation}\label{pr:physical}
\begin{aligned}
 u^\theta_{(0)}&=\frac{U^\theta}{\lambda^{1+\delta}},\quad &
 u^z_{(0)}&=\frac{U^z}{\lambda^{1+\delta}},\quad &
 u^r_{(0)}&=\frac{U^r}{\lambda},\\
 p_{(0)}&=\frac{P}{\lambda^{2+2\delta}},\quad &
 F&=\frac{U^\theta}{\sqrt{2R}}.
\end{aligned}
\end{equation}
The leading profiles depend on $(R,Z)$. We omit their subscript $(0)$.
Smoothness at the axis means that $u^r/r,u^\theta/r,u^z,p$ are smooth
functions of $(r^2,t,z)$. In leading variables this requires
$U^r/\sqrt R,F,U^z,P$ to be smooth in $(R,Z)$ at $R=0$.
For the inertial primitives in \eqref{pr:inertia-primitives}, the corresponding
regularity conditions are smoothness of $\mathcal I^\theta/R$ and
$\mathcal I^z/\sqrt R$. Indeed, $Q^\theta/r$ and $Q^z$ in
\eqref{pr:conservative-inertia} are smooth functions of $(r^2,t,z)$.
The support and admissibility conditions on the stress are imposed separately.

Differentiating \eqref{pr:coordinates} gives
\begin{equation}\label{pr:derivatives}
\begin{aligned}
\partial_t\lambda&=-\frac1{2\lambda L_\delta},\quad&
\partial_z\lambda&=\frac{\lambda^\delta Z}{L_\delta},\quad&
\partial_t Z&=\frac{(1-\delta)Z}{2\lambda^2L_\delta},\quad&
\partial_z Z&=\frac{L}{\lambda^{1-\delta}L_\delta},\\
\partial_t R&=\frac{R}{\lambda^2L_\delta},\quad&
\partial_z R&=-\frac{2RZ}{\lambda^{1-\delta}L_\delta},\quad&
\partial_r R&=\frac{\sqrt{2R}}\lambda.
\end{aligned}
\end{equation}
For $G=G(R,Z)$ and any real $\beta$, the chain rule gives
\begin{equation}\label{pr:chain-rules}
\begin{aligned}
\partial_t(\lambda^\beta G)
 &=\frac{\lambda^{\beta-2}}{L_\delta}
   \left(-\frac\beta2G+\frac{1-\delta}2Z\partial_Z G+R\partial_R G\right),\\
\partial_z(\lambda^\beta G)
 &=\frac{\lambda^{\beta-1+\delta}}{L_\delta}
   (\beta ZG+L\partial_Z G-2ZR\partial_R G),\\
\partial_r(\lambda^\beta G)&=\lambda^{\beta-1}\sqrt{2R}\,\partial_R G.
\end{aligned}
\end{equation}
Fields depending on $\lambda$ require the additional terms in \eqref{fin:chain-rules}.

\subsection{Five moments and reconstruction}
For fixed $Z$, define the five cumulative moments by
\begin{equation}\label{pr:moments}
\begin{aligned}
M^\theta(R,Z)&=\int_0^R\sqrt{2\rho}\,U^\theta(\rho,Z)\dd\rho,\quad&
M^z(R,Z)&=\int_0^R U^z(\rho,Z)\dd\rho,\\
M^{\theta z}(R,Z)&=\int_0^R\sqrt{2\rho}\,U^\theta U^z(\rho,Z)\dd\rho,\quad&
M^{z\theta}(R,Z)&=\int_0^R\bigl((U^z)^2-\tfrac12(U^\theta)^2\bigr)\dd\rho,\\
M^p(R,Z)&=\int_0^R F(\rho,Z)^2\dd\rho.
\end{aligned}
\end{equation}
On a subinterval, use the same integrals with prescribed incoming values,
or integrate backward from prescribed right-endpoint values.
At leading order the radial equation balances $\partial_r p$ with
$(u^\theta)^2/r$. Thus, with $P_0(Z)=P(0,Z)$,
\begin{equation}\label{pr:pressure}
\partial_R P=F^2,\qquad P=P_0+M^p.
\end{equation}
At fixed $(t,z)$ we have $r\,dr=\lambda^2\,dR$. Thus the leading
streamfunction satisfies $r\Psi_{(0)}=\int_0^r s u_{(0)}^z\,ds
=\lambda^{1-\delta}M^z(R,Z)$. Incompressibility and regularity at the
axis give $ru_{(0)}^r=-\partial_z(r\Psi_{(0)})$. Differentiating at
fixed $(t,r)$ and using $\partial_R M^z=U^z$ gives
\begin{equation}\label{pr:radial}
U^r=\frac{2ZRU^z-(1-\delta)ZM^z-L\partial_Z M^z}
              {L_\delta\sqrt{2R}}.
\end{equation}
If $P(\infty,Z)=0$, then
$P_0=-M^p(\infty,Z)$ and $P=-\int_R^\infty F^2\dd\rho$.

\subsection{Inertia, shear, and the actual stress}
For a physical pair $\mathsf T=(\mathsf T^\theta,\mathsf T^z)$ set
\begin{equation}\label{pr:D}
\mathcal D\mathsf T=
 (\partial_r+2/r)\mathsf T^\theta e_\theta
 +(\partial_r+1/r)\mathsf T^ze_z.
\end{equation}
The associated tensor is
\begin{equation}\label{pr:tensor}
\mathbb T=\mathsf T^\theta(e_r\otimes e_\theta+e_\theta\otimes e_r)
 +\mathsf T^z(e_r\otimes e_z+e_z\otimes e_r)
 +r\partial_z\mathsf T^z e_\theta\otimes e_\theta.
\end{equation}
Since $\mathbb T_{\theta\theta}/r=\partial_z\mathsf T^z$, the radial terms
cancel and $\Div\mathbb T=\mathcal D\mathsf T$.

By incompressibility, the angular and axial inertial terms, including the
axial pressure gradient, have the conservation form
\begin{equation}\label{pr:conservative-inertia}
\begin{aligned}
Q^\theta
 &=\partial_tu^\theta+\partial_z(u^zu^\theta)
                    +(\partial_r+2/r)(u^ru^\theta),\\
Q^z
 &=\partial_tu^z+\partial_z\bigl((u^z)^2+p\bigr)
                    +(\partial_r+1/r)(u^ru^z).
\end{aligned}
\end{equation}
Define the inertial primitive $\mathcal I=(\mathcal I^\theta,\mathcal I^z)$ by
\begin{equation}\label{pr:inertia-primitives}
\lambda^{-2-\delta}\mathcal I^\theta
 =-\frac1{r^2}\int_0^r s^2Q^\theta(t,s,z)\dd s,\qquad
\lambda^{-2-\delta}\mathcal I^z
 =-\frac1r\int_0^r sQ^z(t,s,z)\dd s.
\end{equation}

Integration by parts in the angular primitive gives
\begin{equation}\label{pr:angular-integrations}
\begin{aligned}
\int_0^R\sqrt{2\rho}\,\rho \partial_\rho U^\theta\dd\rho
 &=R\sqrt{2R}\,U^\theta-\tfrac32M^\theta,\\
\int_0^R\sqrt{2\rho}\,\rho\partial_\rho(U^zU^\theta)\dd\rho
 &=R\sqrt{2R}\,U^zU^\theta-\tfrac32M^{\theta z}.
\end{aligned}
\end{equation}
Thus \eqref{pr:chain-rules} and \eqref{pr:inertia-primitives} give
\begin{equation}\label{pr:inertia-theta}
\begin{aligned}
\mathcal I^\theta={}&
 \frac{(1-\delta/2)M^\theta-\frac{1-\delta}{2}Z\partial_Z M^\theta
              -R\sqrt{2R}U^\theta}{2L_\delta R}\\
&+\frac{(2\delta-1)ZM^{\theta z}-L\partial_Z M^{\theta z}
               +2R\sqrt{2R}ZU^zU^\theta}{2L_\delta R}
 -U^rU^\theta.
\end{aligned}
\end{equation}
For the axial primitive, $\int_0^R\rho \partial_\rho U^z\dd\rho=RU^z-M^z$.
Integration by parts in \eqref{pr:pressure} gives
\begin{equation}\label{pr:pressure-integral}
 \int_0^R P\dd\rho=RP-\tfrac12\int_0^R(U^\theta)^2\dd\rho,
 \qquad
 \int_0^R\bigl((U^z)^2+P\bigr)\dd\rho=M^{z\theta}+RP.
\end{equation}
Hence
\begin{equation}\label{pr:inertia-z}
\begin{aligned}
\mathcal I^z={}&\frac{\frac{1-\delta}{2}(M^z-Z\partial_Z M^z)-RU^z
                 +2\delta ZM^{z\theta}-L\partial_Z M^{z\theta}}
                  {L_\delta\sqrt{2R}}\\
&+\frac{R\{2(1+\delta)ZP-L\partial_Z P+2Z(U^z)^2\}}
                  {L_\delta\sqrt{2R}}-U^rU^z.
\end{aligned}
\end{equation}
Eliminating $U^r$ by \eqref{pr:radial} yields
\begin{equation}\label{pr:moment-inertia}
\begin{aligned}
\mathcal I^\theta={}&\frac{U^\theta[-R+(1-\delta)ZM^z+L\partial_Z M^z]}
                         {L_\delta\sqrt{2R}}\\
&+\frac{(1-\delta/2)M^\theta-\frac{1-\delta}{2}Z\partial_Z M^\theta
        -L\partial_Z M^{\theta z}+(2\delta-1)ZM^{\theta z}}{2L_\delta R},\\
\mathcal I^z={}&\frac{1}{L_\delta\sqrt{2R}}
 \bigl\{[-R+(1-\delta)ZM^z+L\partial_Z M^z]U^z
       +\tfrac{1-\delta}{2}(M^z-Z\partial_Z M^z)\\
&\hspace{3em}+2\delta ZM^{z\theta}-L\partial_Z M^{z\theta}
       +R[2(1+\delta)Z-L\partial_Z]P\bigr\}.
\end{aligned}
\end{equation}
Define the radial shear and normalized stress by
\begin{equation}\label{pr:shear}
\mathcal S=(\mathcal S^\theta,\mathcal S^z)=(2R\partial_R F,\sqrt{2R}\,\partial_R U^z),\qquad
\mathcal T=\mathcal I+\mathcal S,\qquad \mathsf T_{(0)}=\lambda^{-2-\delta}\mathcal T.
\end{equation}
The physical radial derivatives of the leading velocities satisfy
\begin{equation}\label{pr:radial-viscosity}
 \bigl(\partial_r u^\theta-u^\theta/r,\ \partial_r u^z\bigr)
 =\lambda^{-2-\delta}\bigl(2R\partial_R F,\sqrt{2R}\,\partial_R U^z\bigr)
 =\lambda^{-2-\delta}\mathcal S .
\end{equation}
Here $2R\partial_R F=\sqrt{2R}\,\partial_R U^\theta-U^\theta/\sqrt{2R}$.
Since $(\partial_r+2/r)(\partial_r-1/r)=\partial_r^2+r^{-1}\partial_r-r^{-2}$
and $(\partial_r+1/r)\partial_r=\partial_r^2+r^{-1}\partial_r$, we obtain
\begin{equation}\label{pr:leading-residual}
\begin{aligned}
\NS^\theta(u_{(0)},p_{(0)})
 &=- (\partial_r+2/r)\mathsf T_{(0)}^\theta-\partial_z^2u_{(0)}^\theta,\\
\NS^z(u_{(0)},p_{(0)})
 &=- (\partial_r+1/r)\mathsf T_{(0)}^z-\partial_z^2u_{(0)}^z,\\
\NS^r(u_{(0)},p_{(0)})&=O(\lambda^{-3}).
\end{aligned}
\end{equation}
On compact $(R,Z)$ sets, radial viscosity in both tangential equations has scale
$\lambda^{-3-\delta}$, while axial viscosity is $O(\lambda^{-3+\delta})$ and gains $\lambda^{2\delta}$.

Differentiating the two integrals in \eqref{pr:inertia-primitives} gives
\begin{equation}\label{pr:radial-stress}
R\partial_R \mathcal I^\theta+\mathcal I^\theta=N^\theta,\qquad
R\partial_R \mathcal I^z+\tfrac12 \mathcal I^z=N^z,
\qquad N^j=-\lambda^{3+\delta}\sqrt{R/2}\,Q^j.
\end{equation}
The identity $\mathcal T=0$ is equivalent to
\begin{equation}\label{pr:stress-free}
2R\partial_R F=-\mathcal I^\theta,\qquad \sqrt{2R}\,\partial_R U^z=-\mathcal I^z,
\end{equation}
together with \eqref{pr:moments}, \eqref{pr:pressure}, and
\eqref{pr:radial}.
Applying $R\partial_R+1$ to the angular equation and
$R\partial_R+1/2$ to the axial equation gives
\begin{equation}\label{pr:stress-free-differential}
2R(R\partial_R^2 F+2\partial_R F)=-N^\theta,\qquad
\sqrt{2R}(R\partial_R^2 U^z+\partial_R U^z)=-N^z.
\end{equation}
\subsection{Instability, cone conditions, and admissibility}
\begin{figure}[H]
\centering
\begin{minipage}[t]{0.485\textwidth}
\centering
\includegraphics[width=\linewidth]{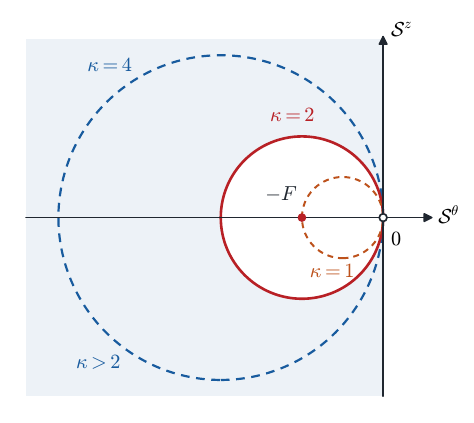}
\caption{Circles of constant $\kappa$ for fixed $F>0$. The shaded
region has $\mathcal S^\theta<0$ and $\kappa>2$.}
\label{fig:shear-circle}
\end{minipage}\hfill
\begin{minipage}[t]{0.485\textwidth}
\centering
\includegraphics[width=\linewidth]{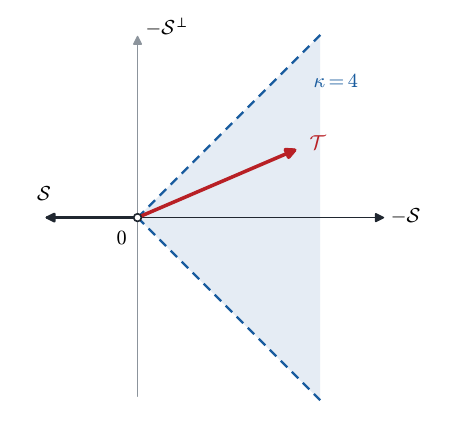}
\caption{Open stress cone for $\kappa=4$, with axis $-\mathcal S$ and half-angle $\pi/4$.}
\label{fig:stress-cone}
\end{minipage}
\end{figure}

On the annulus where the pulses will be supported, assume $F>0$ and
$\mathcal S^\theta<0$, and define
$\mathcal S^\perp=(-\mathcal S^z,\mathcal S^\theta)$ and $\kappa=-|\mathcal S|^2/(F\mathcal S^\theta)$.
Strict admissibility means
\begin{equation}\label{pr:cone}
\kappa>2,\qquad \mathcal T\cdot \mathcal S<0,\qquad
(\kappa-2)(\mathcal T\cdot \mathcal S^\perp)^2<2(\mathcal T\cdot \mathcal S)^2.
\end{equation}
Figure~\ref{fig:shear-circle} shows the \emph{instability condition} $\kappa>2$.
The direction and angle inequalities give the open stress cone, shown in
Figure~\ref{fig:stress-cone} for $\kappa=4$.
We impose them where $\mathcal T\ne0$. Define relaxed admissibility by
\begin{equation}\label{pr:relaxed}
\begin{cases}
\mathcal T\cdot\mathcal S<0,\quad
(\kappa-2)(\mathcal T\cdot\mathcal S^\perp)^2<2(\mathcal T\cdot\mathcal S)^2,
&\kappa>2,\\
-\mathcal T\cdot \mathcal S>-F\mathcal S^\theta(2-\kappa),&\kappa\le2.
\end{cases}
\end{equation}
Since $F\mathcal S^\theta<0$,
\[
\kappa>2\quad\Longleftrightarrow\quad
(\mathcal S^\theta+F)^2+(\mathcal S^z)^2>F^2.
\]
For $\mathcal S=F(-a,b)$, one has $a>0$, $\kappa=(a^2+b^2)/a$.
For $\kappa\le2$, the second branch of \eqref{pr:relaxed} reduces to
$a\mathcal I^\theta-b\mathcal I^z>2aF$, which also gives
$\mathcal T\cdot\mathcal S<0$.

\subsection{Local matching and flat cutoffs}
\begin{lemma}[Moment matching]\label{pr:matching}
Suppose two leading pairs have the same axis pressure and the same
$U^\theta,U^z$ on $[R_1,R_2]$, and their five cumulative moments
at $R_1$ agree as functions of $Z$. Then their moments, pressure,
radial velocity, inertia, shear, and stress agree throughout the interval.
\end{lemma}
\begin{proof}
Subtract the five definitions \eqref{pr:moments} and integrate from
$R_1$. The differences vanish, including their $Z$ derivatives.
Equations \eqref{pr:pressure}, \eqref{pr:radial},
\eqref{pr:inertia-theta}--\eqref{pr:shear} give the conclusion.
\end{proof}

For flat cutoffs we use
\begin{equation}\label{pr:cutoffs}
 \sigma(s)=\begin{cases}e^{-1/s^2},&s>0,\\0,&s\le0,\end{cases}
 \qquad \sigma_{\rm norm}(s)=\frac{\sigma(s)}{\sigma(s)+\sigma(1-s)}
\end{equation}
and their translates and dilates. The Gaussian integral formula
\eqref{fin:gaussian-smooth-quotient} and the tail bound \eqref{fin:gaussian-primitive}
are proved in Appendix~\ref{app:tools}.

\subsection{Admissibility tests and radial equations}\label{oc:calculus-section}
Here $y=\log R$. Later we specify its reference radius, as in $y=\log(R/R_{\rm in})$. Set
\begin{equation}\label{oc:normalized}
\begin{aligned}
 a&=1-2\partial_y\log U^\theta,\quad & b&=2\partial_yU^z/U^\theta,\quad &
 \zeta&=\partial_Z\log U^\theta,\\
 \mathcal{I}^{z\theta}&=\mathcal I^z/\mathcal I^\theta,\quad &
 \widetilde I^\theta&=\frac{L_\delta \mathcal I^\theta}{\sqrt{R/2}U^\theta},\quad &
 \widetilde I^z&=\frac{L_\delta \mathcal I^z}{\sqrt{R/2}(U^\theta)^2},\\
 \Pi&=\frac{2(1+\delta)ZP-L\partial_Z P}{(U^\theta)^2},\quad &
 W&=1-\frac{(1-\delta)ZM^z+L\partial_ZM^z}{R}.
\end{aligned}
\end{equation}
In each coordinate $\partial_y=R\partial_R$. We have
$\mathcal I^\theta/F=R\widetilde I^\theta/L_\delta$,
$\mathcal S=F(-a,b)$ and $\kappa=a+b^2/a$. Where $\mathcal I^\theta\ne0$,
we have $\mathcal{I}^{z\theta}=U^\theta\widetilde I^z/\widetilde I^\theta$
and $\mathcal T/F=(\mathcal I^\theta/F-a,\mathcal I^z/F+b)$.

\begin{lemma}[Uniform admissibility test]\label{oc:cone-test}
For this lemma and its proof only, set
\[
 w=\mathcal{I}^{z\theta},\qquad d=a-bw,\qquad s^2=a^2+b^2,\qquad
 \gamma=2-2bw-b^2/a-(a-2)w^2 .
\]
Suppose $a>0$, $\widetilde I^\theta>0$, and
\begin{equation}\label{oc:cone-test-conditions}
 a-bw>0,\qquad 2bw+b^2/a+(a-2)w^2<2.
\end{equation}
Then \eqref{pr:relaxed} holds for sufficiently large $\mathcal I^\theta/F$, and
\eqref{pr:cone} holds if $\kappa>2$.
The threshold is uniform if $a,b^2,a-bw$ are bounded above and both
inequalities in \eqref{oc:cone-test-conditions} have positive uniform margins.
For $b=0$ and $0<a\leq2$, relaxed admissibility is exactly $\mathcal I^\theta/F>2$.
In particular, a sufficient condition for \eqref{pr:relaxed}, and for
\eqref{pr:cone} when $\kappa>2$, is
\[
 1\le a\le4,\qquad \mathcal I^\theta/F\ge16,\qquad b^2/a+bw+(a-2)_+w^2\le\tfrac12 .
\]
\end{lemma}
\begin{proof}
Since $\mathcal S/F=(-a,b)$ and $\mathcal I^z=w\mathcal I^\theta$, we have
\begin{equation}\label{oc:detail-cone-algebra}
 \frac{\mathcal T\cdot \mathcal S}{F^2}=s^2-d\frac{\mathcal I^\theta}{F},\qquad
 \frac{\mathcal T^\theta \mathcal S^z-\mathcal T^z\mathcal S^\theta}{F^2}=\frac{\mathcal I^\theta}{F}(aw+b),\qquad
 a\kappa=s^2.
\end{equation}
For $\kappa\le2$, \eqref{pr:relaxed} reduces to $\mathcal I^\theta/F>2a/d$. For $\kappa>2$,
the direction condition in \eqref{pr:cone} is $\mathcal I^\theta/F>s^2/d$. Dividing the angle
condition by $F^2(\mathcal I^\theta)^2$ gives
\[
 2\left(d-s^2\frac{F}{\mathcal I^\theta}\right)^2
                  -(\kappa-2)(aw+b)^2>0.
\]
The identity
\[
 2d^2-(\kappa-2)(aw+b)^2
 =(a^2+b^2)\{2-2bw-b^2/a-(a-2)w^2\}=s^2\gamma
\]
gives $s^2[\gamma-4d(F/\mathcal I^\theta)+2s^2(F/\mathcal I^\theta)^2]$ for the left-hand side.
If $\mathcal I^\theta/F>\max\{2a/d,s^2/d,8d/\gamma\}$, both direction bounds hold
and this expression is at least $s^2\gamma/2$. The assumed margins make this threshold uniform.
Under the last sufficient condition, we have $bw\le1/2$, so
$d\ge a/2$ and $s^2\le4(a+b^2/a)\le4d+2\le8d$. Moreover,
$4\gamma-d\ge6-a-3bw\ge1/2$.
Thus $\mathcal I^\theta/F>\max\{2a/d,s^2/d\}$ and
$\gamma-4d(F/\mathcal I^\theta)\ge1/8$, which proves the claim.
\end{proof}

Using \eqref{pr:radial} and $\partial_y M^z=RU^z$, we obtain
\begin{equation}\label{oc:detail-W-identity}
 W=1-2ZU^z+L_\delta\sqrt{2/R}\,U^r,\qquad
 \partial_y W=1-W-(1-\delta)ZU^z-L\partial_Z U^z.
\end{equation}
By \eqref{pr:moment-inertia}, we have
\begin{align*}
 R\sqrt{2R}U^\theta(\widetilde I^\theta+W)
 ={}&(1-\delta/2)M^\theta-\tfrac{1-\delta}{2}Z\partial_Z M^\theta\\
 &-L\partial_Z M^{\theta z}+(2\delta-1)ZM^{\theta z},\\
 (U^\theta)^2\widetilde I^z={}&-WU^z+\frac{1-\delta}{2R}(M^z-Z\partial_Z M^z)\\
 &+\frac{2\delta ZM^{z\theta}-L\partial_Z M^{z\theta}}{R}
 +(U^\theta)^2\Pi.
\end{align*}
Differentiate at fixed $Z$, using $\partial_y M^\theta=R\sqrt{2R}U^\theta$,
$\partial_y M^{\theta z}=R\sqrt{2R}U^\theta U^z$ and
$\partial_y M^{z\theta}=R[(U^z)^2-(U^\theta)^2/2]$. The logarithmic derivatives are
\[
 \partial_y\log\!\bigl(R\sqrt{2R}\,U^\theta\bigr)=2-\frac a2,\qquad
 \partial_y\log\!\left(\frac1{R(U^\theta)^2}\right)=a-2.
\]
The undifferentiated moments cancel, and pressure balance gives
\[
 \partial_y\Pi-(a-1)\Pi=(1+\delta)Z-L\zeta,
 \qquad 2\partial_y P=(U^\theta)^2.
\]
Together with \eqref{oc:detail-W-identity}, this gives the radial equations
\begin{equation}\label{oc:radial-system}
\begin{gathered}
 \partial_y\widetilde I^\theta+(2-a/2)\widetilde I^\theta=-W(1-a/2)-\tfrac\delta2(1-2ZU^z)
 -\left(\tfrac{1-\delta}{2}Z+LU^z\right)\zeta,\\
 \partial_y\widetilde I^z+(2-a)\widetilde I^z=Z+\Pi-\frac{1}{(U^\theta)^2}
 \left[W \partial_y U^z+\tfrac{1+\delta}{2}(1-2ZU^z)U^z
 +\left(\tfrac{1-\delta}{2}Z+LU^z\right)\partial_Z U^z\right].
\end{gathered}
\end{equation}

\subsection{A finite correction of the five moments}
\label{rp:local-correction-section}

Differentiating the five moment densities with respect to $(U^\theta,U^z)$
gives the vector weights
\begin{equation}\label{rr:rows}
\begin{aligned}
 \omega_\theta&=(\sqrt{2R},0),\quad&
 \omega_z&=(0,1),\quad&
 \omega_{\theta z}&=(\sqrt{2R}U^z,\sqrt{2R}U^\theta),\\
 \omega_{z\theta}&=(-U^\theta,2U^z),\quad&
 \omega_p&=(U^\theta/R,0).
\end{aligned}
\end{equation}
\begin{lemma}[Finite moment correction]\label{rp:local-correction}
Let $m\ge1$ and let $\chi\ge0$ be a smooth cutoff supported in
$J\Subset(R_-,R_+)$. Evaluate the five weights in \eqref{rr:rows} at the
uncorrected profile. For each $Z\in[-1,1]$, assume they are linearly
independent as functions of $R$ on $\{\chi>0\}$.
Their Gram matrix is
\begin{equation}\label{rp:moment-matrix}
 A_{j\ell}(Z)=\int_J\chi(R)\,\omega_j(R,Z)\cdot\omega_\ell(R,Z)\,dR .
\end{equation}
For each increment $M^j(R_+,Z)-M^j(R_-,Z)$, let $(m_{\rm def})_j(Z)$
be its target value minus its current value, in the order \eqref{rr:rows}.
If $\|m_{\rm def}\|_{C_Z^m}$ is sufficiently small, all five targets are attained by
\begin{equation}\label{rp:local-velocity-correction}
 (\Delta U^\theta,\Delta U^z)
                  =\chi(R)\sum_{\ell=1}^5c_\ell(Z)\omega_\ell(R,Z),
 \qquad \|c\|_{C_Z^m}\le C\|m_{\rm def}\|_{C_Z^m}.
\end{equation}
The small coefficient vector $c$ is unique and is smooth for smooth data.
The pressure and energy increments may remain fixed
while $M^\theta,M^z,M^{\theta z}$ change.
\end{lemma}
\begin{proof}
Expanding the five moment increments gives the exact equation
\begin{equation}\label{rp:quadratic-moment-map}
 Ac+\mathcal Q(c,c)=m_{\rm def},
\end{equation}
where only the mixed, energy, and pressure components have quadratic
terms. Row independence gives
\begin{equation}\label{rp:local-gram}
 v^\top A(Z)v=\int_J\chi(R)\left|\sum_{j=1}^5v_j\omega_j(R,Z)\right|^2\,dR>0
 \qquad (0\ne v\in\mathbb R^5).
\end{equation}
Thus $A^{-1}$ is smooth on $[-1,1]$. Since $C_Z^m$ is an algebra,
Lemma~\ref{pr:quadratic} applied in $C_Z^m$ gives the unique small
solution and \eqref{rp:local-velocity-correction}. For smooth data the
linearization in $c$ remains invertible, so the solution is smooth in $Z$.
\end{proof}

With $P_0$ and the right moment values fixed, backward integration preserves
all fields to the right of $J$. The correction preserves the endpoint jets of
$U^\theta,U^z$ and, if the axial moment is fixed, those of $U^r$.
On compact subintervals of $R>0$, the fields, shear, and inertia change by
$O(\|m_{\rm def}\|_{C_Z^1})$ in $C^0$, so uniform strict-admissibility margins persist.

\section{Exact mixed exterior, inner core and its admissible exit}\label{profiles:section}
\subsection{Exact mixed exterior}
\label{ue:section}

Let $0<\delta\le1/20$ and $\epsilon_{\mathrm{em}}\ge0$. With
$L=1-Z^2$ and $L_\delta=1-\delta Z^2$, the exact heat profile of amplitude $h_\infty>0$ is
\begin{equation}\label{ue:heat-factor}
 U_e^\theta(R,Z)=h_\infty R^{-(1+\delta)/2}H_e(2L/R),\qquad
 H_e(\xi)=\frac{1}{\Gamma(\delta/2)}\int_0^\infty
 e^{-v}v^{\delta/2-1}(1+\xi v)^{-(1+\delta)/2}\,dv .
\end{equation}
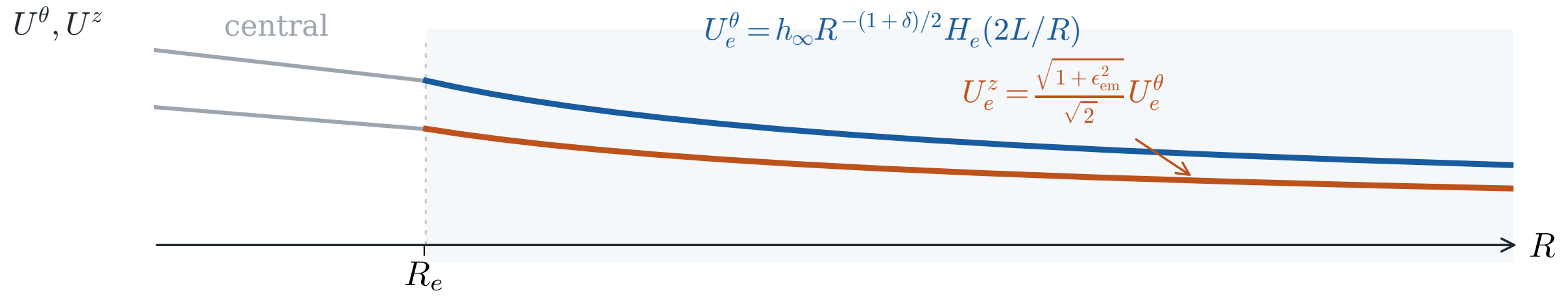
\begin{figure}[htbp]
\centering
\begin{tikzpicture}
  \node[inner sep=0,anchor=south west] (panel) at (0,0)
    {\includegraphics[width=\linewidth]{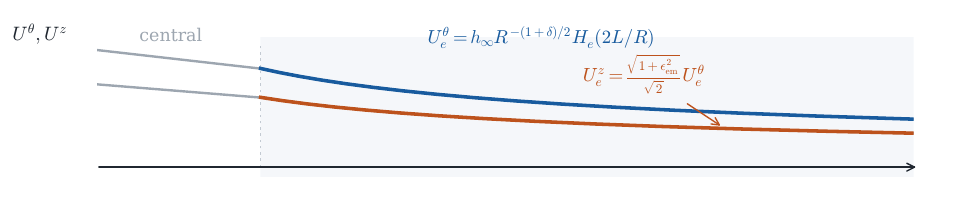}};
  \begin{scope}[x={(panel.south east)},y={(panel.north west)}]
    \draw[thin] (.267,.191)--(.267,.16);
    \node[anchor=north,inner sep=1pt,font=\small] at (.267,.155) {$R_e$};
    \node[anchor=west,inner sep=1pt,font=\small] at (.946,.191) {$R$};
  \end{scope}
\end{tikzpicture}
\caption{Exact mixed exterior. The parameter $\epsilon_{\mathrm{em}}$
changes the axial velocity and tail energy at fixed angular velocity
and pressure. The exterior satisfies \(\NS(u_e,p_e)=0\) and
\(\mathcal T_e=0\).}
\label{fig:exterior}
\end{figure}

\subsubsection{The exact exterior and its five moment primitives}

\begin{proposition}[Free axial amplitude]\label{ue:exact}
Let
\begin{equation}\label{ue:physical}
\begin{aligned}
u_e^r&=-1/r,
 &\qquad u_e^\theta&=h_\infty(r^2/2)^{-(1+\delta)/2}
                   H_e\bigl(4(1-t)/r^2\bigr),\\
 u_e^z&=\frac{\sqrt{1+\epsilon_{\mathrm{em}}^2}}{\sqrt2}u_e^\theta,
 &\qquad p_e&=-\frac1{2r^2}-\int_r^\infty
                    \frac{(u_e^\theta(t,s))^2}{s}\,ds.
\end{aligned}
\end{equation}
Then $\NS(u_e,p_e)=0$ and $\Div u_e=0$ for $r>0$, $t<1$.
The leading profiles are \eqref{ue:heat-factor}, $U_e^z=\frac{\sqrt{1+\epsilon_{\mathrm{em}}^2}}{\sqrt2}U_e^\theta$,
$U_e^r=-1/\sqrt{2R}$, and
\begin{equation}\label{ue:pressure}
 P_e(R,Z)=-\frac12\int_R^\infty(U_e^\theta(\rho,Z))^2\frac{d\rho}{\rho}.
\end{equation}
The five exterior primitives, normalized at infinity as below, give
$\mathcal T_e=0$:
\begin{equation}\label{ue:right-moments}
\begin{aligned}
 M_e^\theta={}&\frac{\sqrt2h_\infty}{1-\delta/2}R^{1-\delta/2}
 -\int_R^\infty\sqrt{2\rho}
       [U_e^\theta-h_\infty\rho^{-(1+\delta)/2}]\,d\rho,\\
 M_e^z={}&\frac{\sqrt{2(1+\epsilon_{\mathrm{em}}^2)} h_\infty}{1-\delta}R^{(1-\delta)/2}
 -\frac{\sqrt{1+\epsilon_{\mathrm{em}}^2}}{\sqrt2}\int_R^\infty
       [U_e^\theta-h_\infty\rho^{-(1+\delta)/2}]\,d\rho+Z,\\
 M_e^{\theta z}={}&\frac{\sqrt{1+\epsilon_{\mathrm{em}}^2} h_\infty^2}{1/2-\delta}R^{1/2-\delta}
 -\sqrt{1+\epsilon_{\mathrm{em}}^2}\int_R^\infty\sqrt\rho
       [(U_e^\theta)^2-h_\infty^2\rho^{-1-\delta}]\,d\rho,\\
 M_e^{z\theta}={}&-\epsilon_{\mathrm{em}}^2\mathcal J^\theta(R,Z),\qquad
 \mathcal J^\theta(R,Z):=\frac12\int_R^\infty(U_e^\theta)^2\,d\rho,\qquad M_e^p=P_e-P_0.
\end{aligned}
\end{equation}
\end{proposition}

\begin{proof}
Differentiating \eqref{ue:heat-factor} gives
\[
 \xi^2\partial_\xi^2H_e+
 \bigl[1+(\tfrac32+\delta)\xi\bigr]\partial_\xi H_e+
 \tfrac{\delta(1+\delta)}4 H_e
 =\frac{1+\delta}{2\Gamma(\delta/2)}
  \int_0^\infty\partial_v\!\left[e^{-v}v^{\delta/2}
                (1+\xi v)^{-(3+\delta)/2}\right]\,dv=0 .
\]
The boundary terms vanish since $\delta/2>0$. Both tangential equations reduce
to $\partial_t v=\partial_r^2v+2(\partial_r v)/r$. The field is divergence
free, and the radial vector Laplacian of $-1/r$ vanishes. The radial equation
follows from $u_e^r\partial_ru_e^r=-r^{-3}$ and
$\partial_rp_e=r^{-3}+(u_e^\theta)^2/r$.

The subtractions in \eqref{ue:right-moments} fix the constants at infinity.
The integrals converge after any fixed number of $Z$ derivatives, and
differentiation in $R$ gives the five densities in \eqref{pr:moments}.
The even functions
$M_e^z-Z,M_e^\theta,M_e^{\theta z},M_e^{z\theta},P_e$
are homogeneous in $(R,L)$ of degrees
\[
 (1-\delta)/2,\quad 1-\delta/2,\quad
 1/2-\delta,\quad-\delta,\quad-1-\delta.
\]
Hence
\[
 (1-\delta)ZM_e^z+L\partial_Z M_e^z=2ZRU_e^z+L_\delta.
\]
By \eqref{pr:radial}, $U_e^r=-1/\sqrt{2R}$.

In \eqref{pr:inertia-theta}--\eqref{pr:inertia-z}, homogeneity of
$M_e^{\theta z}$ cancels the angular quadratic terms. Moreover,
\[
\begin{aligned}
 2\delta ZM_e^{z\theta}-L\partial_Z M_e^{z\theta}
   &=-2ZR\bigl((U_e^z)^2-(U_e^\theta)^2/2\bigr),\\
 R\{2(1+\delta)ZP_e-L\partial_Z P_e+2Z(U_e^z)^2\}
   &=2ZR\bigl((U_e^z)^2-(U_e^\theta)^2/2\bigr).
\end{aligned}
\]
The common heat equation $-\partial_L v=2R\partial_R^2 v+3\partial_R v$, together with the
finite-part normalizations at infinity, gives
\[
 \partial_L M_e^\theta=-2R\sqrt{2R}\,\partial_R U_e^\theta,\qquad
 \partial_L(M_e^z-Z)=-2R\partial_R U_e^z-U_e^z.
\]
Thus, with $F_e=U_e^\theta/\sqrt{2R}$,
$\mathcal I_e^\theta=-2R\partial_R F_e$ and
$\mathcal I_e^z=-\sqrt{2R}\partial_R U_e^z$.
Hence $\mathcal T_e=\mathcal I_e+\mathcal S_e=0$.
\end{proof}

\Needspace{8\baselineskip}
The energy moment singles out $\epsilon_{\mathrm{em}}^2/\delta$ since
$\mathcal J^\theta(R,Z)$ grows like $1/\delta$ at fixed $R$.

\begin{lemma}[The finite-energy coordinate]\label{ue:regularized-tail}
For $R\ge1$ and $|Z|\le1$, the quotient below extends smoothly to
$\delta=0$ and satisfies
\[
 \lim_{\delta\downarrow0}
 \frac{\delta\,\mathcal J^\theta(R,Z)}{R\bigl(U_e^\theta(R,Z)\bigr)^2}=\frac12 .
\]
The quotient and all its $Z,\delta$ derivatives are uniformly bounded on
$R\ge1$, $|Z|\le1$, $0\le\delta\le1/20$, with bounds depending on the order.
At $Z=\pm1$ and $\delta>0$ we have
\[
 \epsilon_{\mathrm{em}}^2\mathcal J^\theta(R,\pm1)
     =\frac{\epsilon_{\mathrm{em}}^2}{2\delta}
            R(U_e^\theta(R,\pm1))^2.
\]
With $\epsilon_{\mathrm{em}}^2/\delta$ as the energy parameter, all
five exterior primitives extend smoothly to $\delta=0$ on fixed
compact radial intervals.  In particular,
\[
 M_e^\theta(R,Z)\big|_{\delta=0}=\sqrt2h_\infty(R+1-Z^2).
\]
For each $Z\in[-1,1]$, the exterior velocities $U_e^\theta,U_e^z,U_e^r$
are real analytic in $R>0$. All their partial derivatives in $R,Z,\delta$
extend continuously to $Z=\pm1$ and $\delta=0$, also at the corners.
Uniformly in $Z$,
\begin{equation}\label{rr:exterior-input}
 U_e^z=\frac{\sqrt{1+\epsilon_{\mathrm{em}}^2}}{\sqrt2}U_e^\theta,
 \qquad U_e^\theta=h_\infty R^{-(1+\delta)/2}
                    [1+O(\delta/R)]\quad(R\to\infty).
\end{equation}
\end{lemma}
\begin{proof}
Substitute $R\rho$ for the radial variable in $\mathcal J^\theta$ and use
$\delta\int_1^\infty \rho^{-1-\delta}\,d\rho=1$ for $\delta>0$. This gives
\[
 \frac{\delta \mathcal J^\theta(R,Z)}{R(U_e^\theta(R,Z))^2}
 =\frac{1+\delta\int_1^\infty \rho^{-1-\delta}
                  [H_e(2L/(R\rho))^2-1]\,d\rho}{2H_e(2L/R)^2}.
\]
Inserting the Gamma integral for $(1+\xi v)^{-(1+\delta)/2}$ in
\eqref{ue:heat-factor} and exchanging the order of integration
(the integrand is nonnegative), we obtain
\[
 H_e(\xi)=\frac1{\Gamma((1+\delta)/2)}
 \int_0^\infty e^{-v}v^{(\delta-1)/2}(1+\xi v)^{-\delta/2}\,dv .
\]
On bounded $\xi$ intervals in $[0,\infty)$, every derivative in $\xi,\delta$
is dominated under the integral by
$Ce^{-v/2}(1+v^{-1/2})(1+|\log v|^m)$ for some $C,m$, also at $\delta=0$.
This proves smoothness up to both poles and $\delta=0$. For $|Z|<1$ the
integral is holomorphic in $1/R$ on $\operatorname{Re}(1/R)>0$, using the
principal power. At $Z=\pm1$ the Gamma integral gives $H_e(0)=1$, so
$U_e^\theta(R,\pm1)=h_\infty R^{-(1+\delta)/2}$. Integrating its square
gives the identity at the poles.

For $R,\rho\ge1$, the bracket $H_e(2L/(R\rho))^2-1$ and all its
$Z,\delta$ derivatives are $O(1/(R\rho))$, since $H_e(0)=1$ for every
$\delta$. Differentiating $\rho^{-1-\delta}$ adds only powers of $\log\rho$.
The bound $C\rho^{-2}(1+\log\rho)^m$ is integrable, and the Gamma integral
over $1\le v\le2$ bounds $H_e(2L/R)$ away from zero uniformly.
We may therefore differentiate the quotient under its integral. At
$\delta=0$ the Gamma formula gives $H_e(\xi)=1$ for all $\xi\ge0$,
so the quotient has numerator $1$ and denominator $2$.

Taylor expansion at $\xi=0$, with $\partial_\xi H_e(0)=-\delta(1+\delta)/4$, gives
\[
 H_e(2L/\rho)
  =1-\frac{\delta(1+\delta)L}{2\rho}
          +O_{C_Z^m}(\delta/\rho^2).
\]
This gives \eqref{rr:exterior-input}. Inserting the expansion into
\eqref{ue:right-moments} gives
\[
 \begin{aligned}
 M_e^\theta={}&\frac{\sqrt2h_\infty}{1-\delta/2}R^{1-\delta/2}
       +\sqrt2h_\infty(1+\delta)L R^{-\delta/2}\\
 &-\sqrt2h_\infty\int_R^\infty\rho^{-\delta/2}
   \left[H_e(2L/\rho)-1+\frac{\delta(1+\delta)L}{2\rho}\right]d\rho.
 \end{aligned}
\]
The last integral vanishes at $\delta=0$. Its differentiated integrand
is bounded for $\rho\ge1$ by $C\rho^{-2}(1+\log\rho)^m$.
For $M_e^z,M_e^{\theta z},P_e$ the same Taylor estimate gives the integrable
bound $C\rho^{-3/2}(1+\log\rho)^m$. Finally,
$M_e^{z\theta}=-(\epsilon_{\mathrm{em}}^2/\delta)(\delta\mathcal J^\theta)$.
\end{proof}

\subsubsection{The free coefficient and the moment conditions}

\begin{lemma}[Tail energy and fixed pressure]\label{ue:budget}
Let $R\ge R_e>1$ and $|Z|\le1$. Then, uniformly in $Z$,
\begin{equation}\label{ue:tail-bounds}
\begin{aligned}
 \frac{h_\infty^2}{2\delta}R^{-\delta}
       \left(1-\frac{\delta^2L}{R}\right)
 &\le \mathcal J^\theta(R,Z)
 \le\frac{h_\infty^2}{2\delta}R^{-\delta},\\
 \frac{h_\infty^2}{2(1+\delta)}R^{-1-\delta}
 \left(1-\frac{\delta(1+\delta)^2L}{(2+\delta)R}\right)
 &\le-P_e(R,Z)
 \le\frac{h_\infty^2}{2(1+\delta)}R^{-1-\delta}.
\end{aligned}
\end{equation}
In particular,
\begin{equation}\label{ue:tail-scales}
 \mathcal J^\theta=\frac{R(U_e^\theta)^2}{2\delta}\,[1+O(\delta/R)],\qquad
 -P_e=\frac{(U_e^\theta)^2}{2(1+\delta)}\,[1+O(\delta/R)].
\end{equation}
For fixed $R,Z,h_\infty$, the range of $M_e^{z\theta}$ over $\epsilon_{\mathrm{em}}\ge0$
is $(-\infty,0]$. Moreover,
\begin{equation}\label{ue:forward-tail}
 \int_R^\infty\bigl((U_e^z)^2-(U_e^\theta)^2/2\bigr)\,d\rho
       =\epsilon_{\mathrm{em}}^2\mathcal J^\theta(R,Z).
\end{equation}
The integral is positive for $\epsilon_{\mathrm{em}}>0$; the pressure
tail is independent of $\epsilon_{\mathrm{em}}$.
\end{lemma}

\begin{proof}
For $\xi,v\ge0$, convexity gives
$(1+\xi v)^{-(1+\delta)/2}\ge1-(1+\delta)\xi v/2$.
Integrating in \eqref{ue:heat-factor} and using
$\Gamma(1+\delta/2)/\Gamma(\delta/2)=\delta/2$ gives
$1\ge H_e(\xi)\ge1-\xi\delta(1+\delta)/4$.
Thus $0\le1-H_e(2L/\rho)^2\le\delta(1+\delta)L/\rho$.
Integration against $h_\infty^2\rho^{-1-\delta}/2$ and
$h_\infty^2\rho^{-2-\delta}/2$ gives \eqref{ue:tail-bounds}.
Using $H_e(2L/R)=1+O(\delta/R)$ gives \eqref{ue:tail-scales}.
Equations \eqref{ue:pressure},
\eqref{ue:right-moments}, and \eqref{rr:exterior-input} give the remaining claims.
\end{proof}

Suppose a middle profile joins pressure $P_a(Z)<0$ at $R_{\rm in}$ to
the exterior at $R_e$. Pressure and energy matching require
\begin{equation}\label{ue:global-budgets}
\begin{aligned}
 \int_{R_{\rm in}}^{R_e}\frac{(U^\theta)^2}{2R}\,dR-P_e(R_e,Z)&=-P_a(Z),\\
 M_a^{z\theta}(Z)+\int_{R_{\rm in}}^{R_e}\bigl((U^z)^2-(U^\theta)^2/2\bigr)\,dR
               +\epsilon_{\mathrm{em}}^2\mathcal J^\theta(R_e,Z)&=0.
\end{aligned}
\end{equation}
\begin{remark}[Role of the free axial amplitude]
The pressure condition requires $-P_e(R_e,Z)\le-P_a(Z)$. As
Figure~\ref{fig:exterior} illustrates, varying $\epsilon_{\mathrm{em}}$ fixes
$U_e^\theta$ and hence $P_e,M_e^\theta$. In \eqref{ue:right-moments} the
factor $\sqrt{1+\epsilon_{\mathrm{em}}^2}$ multiplies $M_e^z-Z$ and
$M_e^{\theta z}$ relative to their values at zero amplitude.
For a fixed middle profile the energy condition gives
\[
 \epsilon_{\mathrm{em}}^2=
 \frac{-M_a^{z\theta}(Z)-\int_{R_{\rm in}}^{R_e}
       \bigl((U^z)^2-(U^\theta)^2/2\bigr)\,dR}{\mathcal J^\theta(R_e,Z)}.
\]
A constant amplitude satisfies this equation precisely when the ratio is
nonnegative and independent of $Z$. The other three endpoint moments
must still match \eqref{ue:right-moments}.
\end{remark}
\subsection{Exact inner core and its admissible exit}
\label{ic:section}

\begin{figure}[htbp]
\centering
\includegraphics[width=\textwidth]{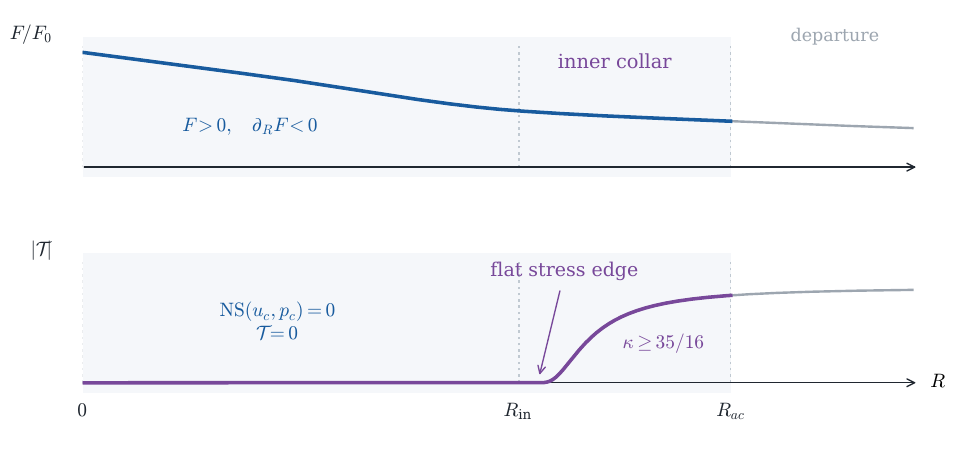}
\caption{Exact core and exit collar. The upper curve shows the normalized swirl
$F/F_0$. Below, the stress is zero for $R\le R_{\rm in}$, vanishes to every order
at $R_{\rm in}$, and is nonzero for $R_{\rm in}<R\le R_{ac}$.}
\label{fig:inner}
\end{figure}

The exact core supplies the singular axial trace. Figure~\ref{fig:inner} shows
the collar from the stress edge $R_{\rm in}$ to the connection data at $R_{ac}$.
We prescribe the axis pressure $P_0$ and recover the pressure elsewhere as $P_0+M^p$.

\subsubsection{Prescribed axis pressure}
\label{ap:section}

For the parameter $P_*\ge1$ of Section~\ref{np:section}, prescribe
\begin{equation}\label{ap:datum}
 P_0(Z)=-\frac{3P_*^2}{(1+Z^2)^2}.
\end{equation}

\begin{lemma}[Pressure shape]\label{ap:shape}
Let $P_0$ be given by \eqref{ap:datum}. Then $P_0$ is even and
holomorphic on the strip $|\operatorname{Im}Z|<1/4$. For real $|Z|\le1$, one has
\begin{equation}\label{ap:bounds}
\begin{aligned}
 \frac{5P_*^2}{2(1+Z^2)^2}&\le -P_0(Z)
 \le\frac{7P_*^2}{2(1+Z^2)^2},\\
 Z\partial_ZP_0(Z)&=\frac{12P_*^2Z^2}{(1+Z^2)^3}
 \ge\frac{10P_*^2Z^2}{(1+Z^2)^3}.
\end{aligned}
\end{equation}
The derivative bounds, with $C_m$ independent of $P_*$, are
\begin{align}
 \|\partial_Z^mP_0\|_{L^\infty([-1,1])}&\le C_mP_*^2
       \qquad (m\ge0),\label{ap:derivative-bounds}\\
 \|P_0\|_{C_Z^1}&\le10P_*^2.\label{ic:pressure-data}
\end{align}
\end{lemma}
\begin{proof}
Direct differentiation gives \eqref{ap:bounds}. Since $2|Z|\le1+Z^2$,
we have $|P_0|\le3P_*^2$ and $|\partial_ZP_0|\le6P_*^2$.
On this strip, $\operatorname{Re}(1+Z^2)
=1+(\operatorname{Re}Z)^2-(\operatorname{Im}Z)^2>15/16$.
Thus $P_0$ is holomorphic there and
Cauchy's estimate gives \eqref{ap:derivative-bounds}.
\end{proof}

\subsubsection{Axis data and the inner construction}

Choose $0<j_0\le1/20$ at most one eighth of the prescribed axial tolerance,
and choose $P_*\ge1$ with $\mu$
as in Section~\ref{np:section}. Choose
$\Lambda\ge\max\{500,j_0^{-2}\}$ sufficiently large for the required
derivative orders and set $R_{\rm in}=4/\Lambda$. Choose $C_*$ next and
$\delta>0$ last. The leading estimates are uniform for $0\le\delta\le1/200$.
Put $U_0^z=4Z+j_0$ and define
\begin{equation}\label{ic:axis-data}
\begin{aligned}
&\begin{aligned}
 H_0&=\tfrac{1-\delta}{2}Z+LU_0^z,
 &\quad\sigma_0&=j_0/500,
 &\quad\chi_c&=\frac{H_0^2}{H_0^2+\sigma_0^2},
 \end{aligned}\\
 &\begin{aligned}
 G(Z)&=\int_{Z_0}^Z\frac{L_\delta(w)H_0(w)}{H_0(w)^2+\sigma_0^2}\dd w,
 &\quad\beta&=\frac{3-\delta/2+j_0Z}{L_\delta},
 &\quad F_0&=C_*^{-1}e^{-\Lambda G},
 \end{aligned}\\
 &\begin{aligned}
 g&=-\tfrac{1+\delta}{2}(1-2ZU_0^z)U_0^z-4H_0
          +2(1+\delta)ZP_0-L\partial_Z P_0.
 \end{aligned}
\end{aligned}
\end{equation}
Let $Z_0$ be the zero of $H_0$. Since
$\partial_Z(H_0/L)=4+(1-\delta)(1+Z^2)/(2L^2)>0$ on $(-1,1)$, it is unique.
Moreover $-j_0/4<Z_0<-j_0/5$, $G\ge0$ on the real interval, and
$29/10<\beta<31/10$.

The denominators are nonzero on the compact real rectangle, and the
simple zero $Z_0$ depends analytically on $\delta$. Thus $P_0,G$,
$L_\delta^{-1}$ and $(H_0^2+\sigma_0^2)^{-1}$ are jointly holomorphic
near $[-1,1]\times[0,1/200]$. Since $G\ge0$ on the real rectangle,
shrink this neighborhood so that $\Re G\ge-3$. Write $\Omega$ for
its $Z$-domain and choose $r_\Omega>0$ with the closed
$r_\Omega$-neighborhood of $[-1,1]$ in $\Omega$.
For $\Lambda\ge1$, require
\begin{equation}\label{ic:complex-smallness}
\begin{aligned}
C_*&\ge\Lambda^2e^{3\Lambda},
 &\quad\sup_\Omega|F_0|&\le\Lambda^{-2}.
\end{aligned}
\end{equation}
Then $F_0,4Z+j_0,P_0$ have common holomorphic bounds independent
of $\Lambda,C_*$.

\subsubsection{An exact analytic Navier--Stokes core}
\label{ex:appendix}

The physical axis data corresponding to \eqref{ic:axis-data} are
\begin{equation}\label{ex:axis-data}
 \left.\frac{u_c^\theta}{r}\right|_{r=0}
   =\lambda^{-2-\delta}F_0(Z),\qquad
 u_c^z(t,0,z)=\lambda^{-1-\delta}(4Z+j_0),\qquad
 p_c(t,0,z)=\lambda^{-2-2\delta}P_0(Z).
\end{equation}
The radial series gives the physical solution at $\bar{\lambda}=\lambda^{2\delta}$. Its normalized value at
$\bar{\lambda}=0$ gives the leading exit profiles.

\begin{theorem}[Uniform exact core]\label{ex:core}\label{ex:local}
Assume $0<\delta\le1/200$ and \eqref{ic:complex-smallness}.
Then there exist $R_0,\bar{\lambda}_0>0$,
independent of $\Lambda\ge1$ and $C_*$, such that
\eqref{ex:axis-data} define an exact physical solution
$(u_c,p_c)$, regular at the axis, on
\[
\begin{aligned}
0&\le R\le R_0,\qquad |Z|\le1,
 &\quad0&<\lambda^{2\delta}\le\bar{\lambda}_0/2.
\end{aligned}
\]
Thus $0\le r\le\sqrt{2R_0}\lambda(t,z)$ with $R_0$ independent of $(t,z)$.
The solution satisfies
\[
 \nabla\cdot u_c=0,\qquad
 \partial_tu_c+(u_c\cdot\nabla)u_c-\Delta u_c+\nabla p_c=0,
\]
and is unique among solutions with these axis data that are regular at the
axis and analytic in $r^2$.
The regular normalized quantities
\begin{equation}\label{ex:normalized}
 \left(
 \lambda^{2+\delta}\frac{u_c^\theta}{r},\
 \lambda^{1+\delta}u_c^z,\
 \lambda^2\frac{u_c^r}{r},\
 \lambda^{2+2\delta}p_c
 \right)
\end{equation}
are jointly holomorphic in the independent variables $(\bar{\lambda},R,Z,\delta)$
on a fixed neighborhood of $\{0\}\times[0,R_0]\times[-1,1]\times[0,1/200]$.
On smaller neighborhoods each derivative is bounded uniformly in $\Lambda,C_*$.
Setting $\bar{\lambda}=\lambda^{2\delta}$ gives the physical solution for $\delta>0$.
At $\bar{\lambda}=0$ the leading coefficients are smooth down to $\delta=0$.
For $\Lambda>4.1/R_0$, the core contains $R\le R_{\rm in}=4/\Lambda$, equivalently,
\[
 \frac{r}{\lambda(t,z)}\le\sqrt{2R_{\rm in}}=\sqrt{\frac8\Lambda},
\]
and the collar $R_{\rm in}<R\le4.1/\Lambda$, uniformly for subsequent increases of $C_*$.
\end{theorem}

\begin{proof}
We use the temporary notation $s=r^2/2$ and $A,B,C,D$ in this proof only.
The radial Taylor coefficients $A_n,B_n,C_n,D_n$ are defined by
\begin{equation}\label{ex:series}
 \left(\frac{u_c^\theta}{r},\frac{u_c^r}{r},u_c^z,p_c\right)=(A,B,C,D)
 =\sum_{n\ge0}(A_n,B_n,C_n,D_n)(t,z)s^n.
\end{equation}
In the next three identities, $h$ denotes any smooth scalar function of $s$.
\[
\begin{gathered}
 \partial_r h(s)=r\partial_s h,\qquad
 (\partial_r^2+r^{-1}\partial_r)h(s)=2s\partial_s^2 h+2\partial_s h,\\
 (\partial_r^2+r^{-1}\partial_r-r^{-2})(rh(s))=r(2s\partial_s^2 h+4\partial_s h).
\end{gathered}
\]
Incompressibility and the two tangential equations become
\begin{equation}\label{ex:regular-system}
\begin{aligned}
2B+2s\partial_s B+\partial_z C&=0,\\
\partial_t A+C\partial_z A+2B(A+s\partial_s A)&=2s\partial_s^2 A+4\partial_s A+\partial_z^2 A,\\
\partial_t C+C\partial_z C+2sB\partial_s C+\partial_z D&=2s\partial_s^2 C+2\partial_s C+\partial_z^2 C.
\end{aligned}
\end{equation}
Coefficient matching gives
\begin{equation}\label{ex:div-rec}
 B_n=-\frac{\partial_z C_n}{2(n+1)},
\end{equation}
\begin{equation}\label{ex:theta-rec}
 2(n+1)(n+2)A_{n+1}
 =(\partial_t-\partial_z^2)A_n
 +\sum_{i+k=n}\bigl[C_i\partial_zA_k+2(k+1)B_iA_k\bigr],
\end{equation}
\begin{equation}\label{ex:z-rec}
 2(n+1)^2C_{n+1}
 =(\partial_t-\partial_z^2)C_n+\partial_zD_n
 +\sum_{i+k=n}\bigl[C_i\partial_zC_k+2kB_iC_k\bigr].
\end{equation}
The pressure equation is $\Delta p_c=-\operatorname{tr}((\nabla u_c)^2)$.
In the orthonormal cylindrical frame,
\begin{equation}\label{ex:gradient-matrix}
\nabla u_c=
\begin{pmatrix}
B+2s\partial_s B&-A&r\partial_z B\\
A+2s\partial_s A&B&r\partial_z A\\
r\partial_s C&0&\partial_z C
\end{pmatrix}.
\end{equation}
Thus
\[
\operatorname{tr}((\nabla u_c)^2)
 =(B+2s\partial_s B)^2+B^2+(\partial_z C)^2-2A(A+2s\partial_s A)+4s(\partial_z B)(\partial_s C).
\]
Equating coefficients in the pressure equation gives
\begin{equation}\label{ex:p-rec}
 \begin{aligned}
 &2(n+1)^2D_{n+1}
 =-\partial_z^2D_n
       +2\sum_{i+k=n}(1+2k)A_iA_k\\
 &-\sum_{i+k=n}\left\{
       \bigl[(1+2i)(1+2k)+1\bigr]B_iB_k
       +(\partial_zC_i)(\partial_zC_k)
       +4k(\partial_zB_i)C_k\right\}.
 \end{aligned}
\end{equation}
Equation \eqref{ex:div-rec} gives $B_n$, and the other recurrences determine
degree $n+1$ from degrees at most $n$. In particular,
$D_1=A_0^2-\tfrac34(\partial_zC_0)^2-\tfrac12\partial_z^2D_0$.

Each recurrence takes at most two $(t,z)$ derivatives, so we assign weight $2n$ at degree $n$.
Let $\mathcal U_2$ be a $(t,z)$ bidisk on which the axis data are bounded and holomorphic.
For the estimates on this bidisk, $h$ denotes a holomorphic function of $(t,z)$. Set
\[
\begin{aligned}
d(t,z)&=\min\{1,\operatorname{dist}((t,z),\partial\mathcal U_2)\},
 &\quad\|h\|_k&=\sup_{\mathcal U_2}d^k|h|.
\end{aligned}
\]
Cauchy's formula with radius $d/(k+2)$ gives, for $\xi=t,z$,
\begin{equation}\label{ex:nagumo}
\begin{aligned}
\|\partial_\xi h\|_{k+1}&\le6(k+1)\|h\|_k,
 &\quad\|fg\|_{k+l}&\le\|f\|_k\|g\|_l,\\
 \|\partial_z^2h\|_{k+2}&\le36(k+1)(k+2)\|h\|_k,
 &\|h\|_{k+1}&\le\|h\|_k .
\end{aligned}
\end{equation}
Let
$m_n=\max\{\|A_n\|_{2n},\|C_n\|_{2n},\|D_n\|_{2n}\}$.
Equation \eqref{ex:div-rec} gives
$\|B_n\|_{2n+1}\le6m_n$.
Substitution in the recurrences gives
\begin{equation}\label{ex:majorant}
 m_{n+1}\le c m_n+c\sum_{i+k=n}m_im_k,
\end{equation}
with an absolute constant $c$ in this majorant argument. The denominator $2(n+1)^2$ absorbs the
second-derivative factor $36(2n+1)(2n+2)$. For $i+k=n$, the pressure
products use $(2i+1)(2k+1)\le(n+1)^2$ and $k(i+1)\le(n+1)^2/4$.
In particular, $\|\partial_z B_i\|_{2i+2}\le72(i+1)m_i$ bounds
$4k(\partial_z B_i)C_k/[2(n+1)^2]$ by $36m_i m_k$ in the required norm.
Each tangential product costs at most one derivative with coefficient
at most $2(n+1)$.

Choose $M_0>0$ with $M_0\ge m_0$ and define
$M_{n+1}=cM_n+c\sum_{i+k=n}M_iM_k$. Then $m_n\le M_n$ by induction.
For an auxiliary complex variable $x$, summing the recurrence against $x^{n+1}$ gives
\begin{equation}\label{ex:generating}
 \mathcal M(x):=\sum_{n=0}^\infty M_nx^n
 =M_0+cx\mathcal M(x)+cx\mathcal M(x)^2 .
\end{equation}
For $x_0=[8c(1+M_0)]^{-1}$ and $|x|\le x_0$, the right side maps
$|\mathcal M|\le2M_0$ into itself with derivative at most $1/2$.
The unique holomorphic solution with $\mathcal M(0)=M_0$ has the same Taylor
recurrence, so the series in \eqref{ex:generating} converges for $|x|<x_0$.
On $d\ge d_0>0$, the terms of the three series for $A,C,D$ are bounded by
$M_n(|s|/d_0^2)^n$, and those for $B$ have the additional factor $6/d_0$.
Taking $x=|s|/d_0^2<x_0$ gives normal convergence of all four radial series.
Cauchy's estimates give termwise differentiation on smaller domains.

The tangential components of $\NS(u_c,p_c)$ vanish. The pressure equation
and incompressibility give $\nabla\cdot\NS(u_c,p_c)=0$, so
$r\,\NS^r(u_c,p_c)$ is constant in $r$. Its limit at the axis is zero,
which proves the radial equation.
The Cartesian components $(Bx_1-Ax_2,Bx_2+Ax_1)$ are analytic across the axis.
The recurrences give uniqueness.

Fix $(t_0,z_0)$ and set $\ell=\lambda(t_0,z_0)$. Use $r^2/(2\ell^2)$
for the radial series and write
\[
\begin{aligned}
t&=t_0+\ell^2\sigma,
 &\quad z&=z_0+\ell^{1-\delta}\zeta .
\end{aligned}
\]
The scale equation \eqref{pr:coordinates} for the ratio
$\Lambda_{\rm loc}(\sigma,\zeta):=\lambda(t,z)/\lambda(t_0,z_0)$ becomes
\begin{equation}\label{ex:implicit}
 \Lambda_{\rm loc}^2
 -\Lambda_{\rm loc}^{2\delta}\bigl(Z(t_0,z_0)+\zeta\bigr)^2
 =1-Z(t_0,z_0)^2-\sigma,\qquad
 \Lambda_{\rm loc}(0,0)=1 .
\end{equation}
At $(\Lambda_{\rm loc},\sigma,\zeta)=(1,0,0)$ the derivative of the
left side in $\Lambda_{\rm loc}$ is $2(1-\delta Z(t_0,z_0)^2)\ge2(1-\delta)$.
Using the power branch through $1$, the analytic implicit function theorem
and compactness give bounded holomorphic $\Lambda_{\rm loc}^{\pm1}$
and $Z=(Z(t_0,z_0)+\zeta)\Lambda_{\rm loc}^{\delta-1}$ on common neighborhoods
of $(\sigma,\zeta)=(0,0)$, $Z(t_0,z_0)\in[-1,1]$ and $\delta\in[0,1/200]$.
The normalized axis data and the factors multiplying
$(A_n,B_n,C_n,D_n)$ in \eqref{ex:series} are, respectively,
\begin{equation}\label{ex:coefficient-scales}
 \begin{gathered}
 (\Lambda_{\rm loc}^{-2-\delta}F_0(Z),\
  \Lambda_{\rm loc}^{-1-\delta}(4Z+j_0),\
  \Lambda_{\rm loc}^{-2-2\delta}P_0(Z)),\\
 (\ell^{2+\delta+2n},\ \ell^{2+2n},\
  \ell^{1+\delta+2n},\ \ell^{2+2\delta+2n}).
 \end{gathered}
\end{equation}
The axis data are uniformly bounded by \eqref{ic:complex-smallness}.
In \eqref{ex:theta-rec}--\eqref{ex:z-rec},
$\partial_t-\partial_z^2$ becomes
$\partial_\sigma-\ell^{2\delta}\partial_\zeta^2$.
The transport terms and axial pressure gradient have coefficient one.
In \eqref{ex:p-rec}, the pressure and swirl terms have scale
$\ell^{-4-2\delta-2n}$; the axial pressure derivative and non-swirl
quadratic terms have scale $\ell^{-4-2n}$, hence coefficient $\ell^{2\delta}$.
The divergence recurrence and all denominators are unchanged.

Replace $\ell^{2\delta}$ by an independent complex parameter $\bar{\lambda}$.
The coefficients are polynomial in $\bar{\lambda}$ and holomorphic in the reference
latitude and $\delta$. For $|\bar{\lambda}|\le1$, \eqref{ex:majorant} holds on common
neighborhoods of $[-1,1]$ and $[0,1/200]$, uniformly in $\ell,\Lambda,C_*$.
Normal convergence on a fixed disk in $r^2/(2\ell^2)$ gives joint holomorphy.

At $\sigma=\zeta=0$, the sums equal \eqref{ex:normalized}, since
$\Lambda_{\rm loc}=1$, $Z=z_0/\ell^{1-\delta}$, and $R=r^2/(2\ell^2)$.
Taking $\bar{\lambda}=\ell^{2\delta}$ gives the physical recurrences.
Uniqueness in \eqref{ex:div-rec}--\eqref{ex:p-rec} identifies solutions on overlaps.
By \eqref{ex:implicit}, $\lambda\asymp\ell$, giving a uniform $R$ radius;
Cauchy's estimates give the derivative bounds. The same argument at
$Z=\pm1$, $\lambda>0$, gives smooth one-sided extensions.
Thus $R_0,\bar{\lambda}_0$ depend only on $P_*,j_0$, the common neighborhoods
and the axis-data bounds, uniformly in $0\le\delta\le1/200$ and
independently of $t,z,\Lambda,C_*$.
Choose $\Lambda>4.1/R_0$ and then $C_*$ satisfying
\eqref{ic:complex-smallness}; later increases of $C_*$ preserve these radii.
\end{proof}
\subsubsection{The leading profile and bounds at its exit}
\label{ic:core-section}

In the radial equation divide by $\lambda^{-3-2\delta}\sqrt{2R}$.
The pressure and centrifugal terms give $\partial_RP-F^2$.
By \eqref{ex:normalized}, the remaining terms are $\bar{\lambda}$ times
functions regular at the axis, with $\bar{\lambda}=\lambda^{2\delta}$.
Evaluation at $\bar{\lambda}=0$ and integration from $P(0,Z)=P_0(Z)$ give
\eqref{pr:pressure}.

At $\bar{\lambda}=0$, the exact core satisfies the leading equations

\begin{align}
 2L_\delta(R\partial_R^2 F+2\partial_R F)
 &=W(F+R\partial_R F)+\tfrac\delta2(1-2ZU^z)F+H\partial_Z F,\label{ic:core-equations}\\
 2L_\delta(R\partial_R^2 U^z+\partial_R U^z)
 &=WR\partial_R U^z+\tfrac{1+\delta}{2}(1-2ZU^z)U^z+H\partial_Z U^z
       \notag\\
 &\quad+L\partial_Z P-2(1+\delta)ZP-2ZRF^2,\notag\\
 W&=1-\frac{(1-\delta)ZM^z+L\partial_Z M^z}{R},\qquad
 H=\tfrac{1-\delta}{2}Z+LU^z.\notag
\end{align}
\begin{proposition}[Bounds at the core exit]\label{ic:core-exit}
For sufficiently large $\Lambda$, and then sufficiently large $C_*$,
the leading core exists on $0\le R\le4.1/\Lambda$, has
$F_c>0$, $-\partial_R\log F_c\ge1/64$, and at $R_{\rm in}=4/\Lambda$
\begin{equation}\label{ic:kappa-exit}
 \kappa_c=-2R_{\rm in}\frac{\partial_R F_c}{F_c}
             +\frac{(\partial_R U_c^z)^2}{-F_c\partial_R F_c}\ge\frac52.
\end{equation}
The angular contribution alone exceeds $5/2$ on
$\{\chi_c\ge99/100\}$; on its complement
$-2L_\delta \partial_R U_c^z(R_{\rm in})\ge3j_0/10$. Moreover
\begin{align}
 |\partial_Z\log F_c(R_{\rm in})+\Lambda\partial_Z G|
       &\le B(\sqrt{\chi_c}+\Lambda^{-1}),\label{ic:exit-z}\\
 \|C_*F_c\|_{C^k}+\|(C_*F_c)^{-1}\|_{C^k}
        +\|U_c^z\|_{C^k}&\le B_{k,\Lambda}.
 \label{ic:core-norms}
\end{align}
The last norms use $(R/R_{\rm in},Z)\in[0,e^{1/100}]\times[-1,1]$.
These bounds and the threshold on $\Lambda$ are independent of later $C_*$.
\end{proposition}

\begin{proof}
For the comparison with the frozen linear recurrence, we rescale $R$
at fixed $Z$ and set
\[
\begin{aligned}
q(R)&=\tfrac R2\bigl(\Lambda\chi_c(Z)+\beta(Z)\bigr),
 &B(q)&=\sum_{n\ge0}\frac{(-q)^n}{n!(n+1)!},\\
 F^m&=F_0B(q),
 &U^{z,m}&=U_0^z-\frac{gR}{2L_\delta}.
\end{aligned}
\]
The axis slopes in \eqref{ic:core-equations} are
\[
 4L_\delta \partial_R F(0)
   =\bigl(W_0+\tfrac\delta2(1-2ZU_0^z)\bigr)F_0+H_0\partial_Z F_0,\qquad
 2L_\delta \partial_R U^z(0)=-g,
\]
where $W_0=1-(1-\delta)ZU_0^z-4L$. Since
$W_0+\delta(1-2ZU_0^z)/2=-L_\delta\beta$ and
$H_0\partial_Z F_0=-\Lambda L_\delta\chi_cF_0$, these become
\[
\begin{aligned}
\partial_R F(0)&=-\tfrac14(\Lambda\chi_c+\beta)F_0,
 &\partial_R U^z(0)&=-g/(2L_\delta).
\end{aligned}
\]
\paragraph{Bounds for the model series.}
For $\Lambda\ge\max(500,j_0^{-2})$ and $\Lambda R\le4.1$,
$0\le q<21/10$. Alternating-series bounds give
\begin{equation}\label{ic:bessel-bounds}
\begin{aligned}
\tfrac14&<B(q)\le1,
 &\tfrac18&<-B'(q)\le\tfrac12,\\
 -2qB'(q)/B(q)&>3\qquad (19/10\le q\le21/10).
\end{aligned}
\end{equation}
The successive term ratios for $B$ and $-B'$ are $q/((n+1)(n+2))$
and $q/((n+1)(n+3))$. For $3B+2qB'$ the ratio is
$\frac{2n+5}{2n+3}\frac{q}{(n+1)(n+2)}\le1/2$ for $n\ge1$.
Alternating truncation therefore gives
\[
\begin{aligned}
 -B'(q)&\in[1/2-q/6,\,1/2],\\
 B(q)&\ge1-q/2+q^2/12-q^3/140,\\
 3B(q)+2qB'(q)&\le3-\tfrac52q+\tfrac7{12}q^2
                                    -\tfrac1{16}q^3+\tfrac1{250}q^4.
\end{aligned}
\]
The cubic has derivative at most $-3/20$ on $[0,21/10]$
and exceeds $1/4$ at $q=21/10$. Also $B(0)=1$.
The quartic has derivative $-5/2+7q/6-q^2(3/16-2q/125)\le-1/20$
on $[19/10,21/10]$ and lies below $-1/50$ at $q=19/10$.
These estimates give \eqref{ic:bessel-bounds}.
\paragraph{Angular and axial contributions.}
In the strong sector $\chi_c\ge99/100$, the exit value
$q=2\chi_c+2\beta/\Lambda\ge1.98$ gives
$-2R_{\rm in}\partial_R F^m/F^m=-2qB'(q)/B(q)>3$.
In the weak sector $\chi_c\le99/100$, $|H_0|\le\sqrt{99}\,\sigma_0<j_0/50$.
Since $H_0=(1-\delta)Z/2+L(4Z+j_0)$, we have
\[
\begin{aligned}
 H_0&\le-(1-\delta)j_0/8<-j_0/50 &&(-1\le Z\le-j_0/4),\\
 H_0&\ge j_0+(9/2-\delta/2)Z-j_0Z^2\ge j_0/10-j_0^3/25>j_0/50
                                      &&(-j_0/5\le Z\le0),\\
 H_0&\ge3j_0/4                          &&(0\le Z\le1/2),\\
 H_0&\ge(1-\delta)/4>j_0/50              &&(1/2\le Z\le1).
\end{aligned}
\]
Thus the weak sector lies
in $-j_0/4\le Z\le-j_0/5$. On this interval
$0\le U_0^z\le j_0/5$, $H_0\le j_0/5$, and $|Z|\le1/80$.
The two possibly negative contributions to $g$ satisfy
\[
 \tfrac{1+\delta}{2}(1-2ZU_0^z)U_0^z
 \le\tfrac{1+\delta}{2}(1+j_0^2/10)j_0/5\le j_0/5,\qquad
 4H_0\le4j_0/5.
\]
The pressure derivative gives
\[
 -L\partial_Z P_0\ge\frac{10P_*^2|Z|(1-Z^2)}{(1+Z^2)^3}
                  \ge8|Z|,
\]
while $2(1+\delta)ZP_0\ge0$. Thus
\begin{equation}\label{ic:weak-sector}
 g\ge-j_0/5-4j_0/5+8|Z|\ge3j_0/5.
\end{equation}
If $F_0\le\Lambda^{-2}$, then
\[
\begin{aligned}
\frac{(\partial_R U^{z,m})^2}{-F^m\partial_R F^m}
 &\ge\frac{(2j_0/5)^2}{F_0^2(\Lambda\chi_c+\beta)}
 \ge\frac{(2j_0/5)^2\Lambda^3}{1+4/\Lambda}\\
 &\ge\frac{2j_0^2\Lambda^3}{25}\ge\frac2{25j_0^4}>3
 \qquad\hbox{on }\{\chi_c\le99/100\}.
\end{aligned}
\]
\paragraph{The rescaled model equations.}
Fix $0<h<r_\Omega/4$ before $\Lambda,C_*$ and assume $4.1/\Lambda<R_0$.
The functions $P_0$, $1/L_\delta$, $1/(H_0^2+\sigma_0^2)$,
$\chi_c$, $\beta$, and $g/L_\delta$ have common holomorphic bounds on
$\Omega$. Constants in this analytic estimate depend only on these
bounds and $h$, not on later $\Lambda,C_*$ satisfying
\eqref{ic:complex-smallness}.
Let $x=\Lambda R$, $\varepsilon=\Lambda^{-1}$,
$F^m=F_0\Phi^m$, and $U^{z,m}=U_0^z+\varepsilon\Psi^m$. Then
\begin{equation}
\begin{aligned}
 2(x\partial_{x}^2
       +2\partial_{x})\Phi^m
       +(\chi_c+\varepsilon\beta)\Phi^m&=0,\\
 2(x\partial_{x}^2
       +\partial_{x})\Psi^m&=-g/L_\delta,\\
 \Phi^m(0,Z)=1,\qquad \Psi^m(0,Z)&=0.
\end{aligned}\label{ic:detail-model}
\end{equation}
Integrating from the regular axis gives
$\partial_x\Psi^m=-g/(2L_\delta)$.

\paragraph{The analytic space and radial inversion.}
We estimate the difference from the model on $0\le x\le4.1$.
The norm below controls axial differentiation after radial inversion,
as required by the transport terms in \eqref{ic:core-equations}.
For $f=\sum_{n\ge0}f_n(Z)x^n$, with smooth coefficients
on $[-1,1]$, set
\begin{equation}
 \|f\|_h
 =\sup_{\substack{n,m\ge0\\Z\in[-1,1]}}
 \frac{20^n h^m(n+1)^2(m+1)^2
              |\partial_Z^m f_n(Z)|}
      {m!\binom{n+m}{m}}.
 \label{ic:analytic-norm}
\end{equation}
Let $\mathcal X_h$ be the space of series with finite norm.
Completeness follows by taking uniform limits of each coefficient and
its derivatives and passing to the limit in \eqref{ic:analytic-norm}.
On $|x|<5$, $|Z-Z_*|<h/2$, $Z_*\in[-1,1]$, the joint Taylor series
is bounded in absolute value by
$\|f\|_h\sum_{n,m}\binom{n+m}{m}(1/4)^n(1/2)^m=4\|f\|_h$.
The local extensions agree on overlaps. Cauchy's estimates give
\begin{equation}
 \|f\|_{C^k([0,4.1]\times[-1,1])}
       \le K_{k,h}\|f\|_h
       \qquad(k\ge0).
 \label{ic:detail-analytic-3}
\end{equation}
Here and below these $C^k$ norms use $(x,Z)$.

At fixed $Z$, let $\mathcal M_x$ denote the radial average. Multiplying by $x$ gives the radial primitive.
\[
\begin{aligned}
(\mathcal M_x f)(x,Z)&=\int_0^1 f(sx,Z)\,ds,
 &x(\mathcal M_x f)(x,Z)&=\int_0^x f(v,Z)\,dv.
\end{aligned}
\]
In this proof we write $\mathcal T_\nu$ ($\nu=1,2$) for the regular inverse of $x\partial_x^2+\nu\partial_x$
with zero value at the axis. Thus
\[
\begin{aligned}
(\mathcal T_\nu f)_{n+1}&=\frac{f_n}{(n+1)(n+\nu)},
 &(\mathcal T_\nu f)_0&=0,\\
 \mathcal T_\nu f
   &=\int_0^{x}s^{-\nu}
            \int_0^s t^{\nu-1}f(t,Z)\,dt\,ds.
\end{aligned}
\]
We claim
\begin{equation}
\begin{aligned}
 &\|fk\|_h\le K\|f\|_h\|k\|_h,\qquad
 \|\mathcal M_x f\|_h\le\|f\|_h,\qquad
 \|x\mathcal M_x f\|_h\le80\|f\|_h,\\
 &\|\mathcal T_\nu((\partial_Zf)k)\|_h
 +\|\mathcal T_\nu(fx\partial_{x}k)\|_h
 +\|\mathcal T_\nu((\partial_Zf)x
                  \partial_{x}k)\|_h
 \le K_h\|f\|_h\|k\|_h.
\end{aligned}\label{ic:analytic-algebra}
\end{equation}
Let $f_{n,m}(Z)=\partial_Z^m f_n(Z)/m!$ denote the $m$th Taylor coefficient
of $f_n$ at $Z$, and set
$w_{n,m}=\binom{n+m}{m}/((n+1)^2(m+1)^2)$.
Then $|f_{n,m}(Z)|\le
\|f\|_h20^{-n}h^{-m}w_{n,m}$.
Splitting a convolution at its midpoint gives
\[
 \sum_{i=0}^n\frac1{(i+1)^2(n-i+1)^2}
 \le\frac8{(n+1)^2}\sum_{i=0}^\infty\frac1{(i+1)^2}
 \le\frac{16}{(n+1)^2}.
\]
For $i+i'=n$, $p+q=m$, Vandermonde's identity gives
$\binom{i+p}{p}\binom{i'+q}{q}\le\binom{n+m}{m}$.
Leibniz's rule and the preceding convolution in both indices give
$\sum_{i+i'=n,\,p+q=m}w_{i,p}w_{i',q}\le260w_{n,m}$, proving the
algebra estimate.

For the mixed derivative in \eqref{ic:analytic-algebra}, the normalized output
coefficient at $(n+1,m)$ is
\[
 \frac1{(n+1)(n+\nu)}
 \sum_{\substack{i+i'=n\\p+q=m}}
             (p+1)f_{i,p+1}(Z)\,i' k_{i',q}(Z).
\]
Use
\[
\begin{aligned}
 (p+1)\binom{i+p+1}{p+1}
     &=(i+1)\binom{i+p+1}{p},\\
 \binom{i+p+1}{p}\binom{i'+q}{q}
     &\le\binom{n+m+1}{m},\qquad
       \frac{(i+1)i'}{(n+1)(n+\nu)}\le1.
\end{aligned}
\]
Its weighted absolute value is at most
\[
 \frac{20}{h}\|f\|_h\|k\|_h
 (n+2)^2(m+1)^2
 \sum_{\substack{i+i'=n\\p+q=m}}
 \frac1{(i+1)^2(i'+1)^2(p+2)^2(q+1)^2}.
\]
The convolution bounds and $(n+2)^2/(n+1)^2\le4$ give the claim.
Replacing $i'$ by $1$ gives the first differentiated-product bound;
omitting the axial shift and using $i'/((n+1)(n+\nu))\le1$ gives the second.
The average divides degree $n$ by $n+1$, and multiplication by $x$
increases its weight by at most $20(n+2)^2/(n+1)^2\le80$.
This proves the primitive bound and, for series starting at degree $k$,
$\|\mathcal T_\nu f\|_h\le80\|f\|_h/((k+1)(k+\nu))$.
For bounded holomorphic $a$ on $\Omega$, the Cauchy estimate on
disks of radius $r_\Omega/2$ and $h<r_\Omega/4$ give
\begin{equation}
 \|a\|_h
 \le\sup_{m\ge0}(m+1)^2(2h/r_\Omega)^m
                    \sup_\Omega|a|
 \le K_h\sup_\Omega|a|.
 \label{ic:detail-analytic-5}
\end{equation}
Together with \eqref{ic:analytic-algebra}, this bounds multiplication by $a$.
\paragraph{The inverse of the model operator.}
For $0\le\varepsilon\le1/500$, set
\begin{equation}
 \mathcal A_\varepsilon f
    =\tfrac12\mathcal T_2((\chi_c+\varepsilon\beta)f).
 \label{ic:detail-analytic-6}
\end{equation}
By \eqref{ic:detail-analytic-5}, multiplication by $\chi_c+\varepsilon\beta$
is uniformly bounded on $\mathcal X_h$. Increasing $K_h$ if necessary, for
series starting at degree $k$ we obtain
$\|\mathcal A_\varepsilon f\|_h
\le K_h\|f\|_h/((k+1)(k+2))$.
Iteration gives
\begin{equation}
 \|\mathcal A_\varepsilon^n\|
 \le\frac{K_h^n}{n!(n+1)!},\qquad
 (I+\mathcal A_\varepsilon)^{-1}
     =\sum_{n=0}^\infty(-\mathcal A_\varepsilon)^n.
\label{ic:detail-analytic-7}
\end{equation}
The factorial bound gives norm convergence and a uniformly bounded
inverse. Integrating \eqref{ic:detail-model} at the regular axis gives
$(I+\mathcal A_\varepsilon)\Phi^m=1$ and
$\Psi^m=-\tfrac12\mathcal T_1(g/L_\delta)$, hence
\begin{equation}\label{ic:detail-analytic-8}
\begin{aligned}
\Phi^m&=(I+\mathcal A_\varepsilon)^{-1}1,
 &\quad\Psi^m&=-\frac{g}{2L_\delta}x.
\end{aligned}
\end{equation}
Writing $\Phi^m=\sum_{n\ge0}\varphi_n x^n$, the recurrence $2(n+1)(n+2)\varphi_{n+1}
=-(\chi_c+\varepsilon\beta)\varphi_n$, $\varphi_0=1$, gives
$\Phi^m=B((\chi_c+\varepsilon\beta)x/2)$.
Thus $X^m=(\Phi^m,\Psi^m)$ is uniformly bounded in
$\mathcal X_h\times\mathcal X_h$.
\paragraph{The equation for the difference.}
By \eqref{pr:moments} and \eqref{pr:pressure},
\begin{equation}\label{ic:detail-analytic-9}
 M^z(R,Z)=\int_0^R U^z(s,Z)\,ds,\qquad
 P(R,Z)=P_0(Z)+\int_0^R F(s,Z)^2\,ds.
\end{equation}
Set
\begin{equation}\label{ic:detail-analytic-10}
\begin{aligned}
F(R,Z)&=F_0(Z)\Phi(x,Z),
 &\quad U^z(R,Z)&=U_0^z(Z)+\varepsilon\Psi(x,Z),
\end{aligned}
\end{equation}
where $\Phi(0,Z)=1$ and $\Psi(0,Z)=0$. In particular,
$M^z(\varepsilon x,Z)/(\varepsilon x)
=U_0^z+\varepsilon\mathcal M_x\Psi$, also at the axis by
continuity. With $W,H$ from \eqref{ic:core-equations},
\begin{equation}
\begin{aligned}
 M^z(\varepsilon x,Z)
     &=\varepsilon x U_0^z+\varepsilon^2 x\mathcal M_x\Psi,\\
 P(\varepsilon x,Z)&=P_0+\varepsilon\mathcal P,
       \qquad \mathcal P=x F_0^2\mathcal M_x(\Phi^2),\\
 W&=W_0-\varepsilon((1-\delta)Z\mathcal M_x\Psi
                              +L\partial_Z\mathcal M_x\Psi),\\
 H&=H_0+\varepsilon L\Psi.
\end{aligned}\label{ic:detail-analytic-11}
\end{equation}
Here $W_0=1-(1-\delta)ZU_0^z-4L$, and
$W_0+\delta(1-2ZU_0^z)/2
=-3+\delta/2-j_0Z=-L_\delta\beta$.
Substituting in \eqref{ic:core-equations} gives
\begin{equation}
\begin{aligned}
 2(x\partial_{x}^2\Phi
           +2\partial_{x}\Phi)+\chi_c\Phi
       &=\varepsilon E_\theta,\\
 2(x\partial_{x}^2\Psi
           +\partial_{x}\Psi)
       &=-\frac g{L_\delta}+\varepsilon E_z,
\end{aligned}\label{ic:normalized-core}
\end{equation}
where
\begin{equation}
\begin{aligned}
 E_\theta={}&\frac1{L_\delta}
 \left[W(\Phi+x\partial_{x}\Phi)
       +\frac\delta2(1-2Z(U_0^z+\varepsilon\Psi))\Phi
       +H\partial_Z\Phi\right]
       -\frac{LH_0}{H_0^2+\sigma_0^2}\Psi\Phi,\\
 E_z={}&\frac1{L_\delta}\left[
       Wx\partial_{x}\Psi
       +\left(\frac{1+\delta}{2}(1-4ZU_0^z)+4L\right)\Psi
       -\varepsilon(1+\delta)Z\Psi^2+H\partial_Z\Psi\right.\\
     &\hspace{4em}\left.
       +L\partial_Z\mathcal P-2(1+\delta)Z\mathcal P
       -2Zx F_0^2\Phi^2\right].
\end{aligned}\label{ic:detail-analytic-13}
\end{equation}
Here $\partial_R=\varepsilon^{-1}\partial_x$ and
$\partial_Z\log F_0=-\varepsilon^{-1}L_\delta H_0/(H_0^2+\sigma_0^2)$.
The latter identity gives the term $-\chi_c\Phi$ in the angular transport,
and the terms independent of $\varepsilon$ in the axial equation,
including pressure, give $-g/L_\delta$.
Put $X=(\Phi,\Psi)$, $Y=X-X^m$, and add
$\varepsilon\beta\Phi$ to both sides of the angular equation.
The remaining source is
\begin{equation}
\begin{aligned}
 G_\theta(X)&=E_\theta(X)+\beta\Phi\\
 &=\frac1{L_\delta}\left[
   (W-W_0-\varepsilon\delta Z\Psi)\Phi
    +Wx\partial_{x}\Phi
    +H\partial_Z\Phi\right]
       -\frac{LH_0}{H_0^2+\sigma_0^2}\Psi\Phi.
\end{aligned}\label{ic:detail-analytic-14}
\end{equation}
Subtracting the model equations and applying the inverses with zero axis value gives
\begin{equation}
 \begin{aligned}
 Y&=\varepsilon\mathcal N_{\varepsilon,C_*}(X^m+Y),\\
 \mathcal N_{\varepsilon,C_*}(X)
   &=\tfrac12\bigl((I+\mathcal A_\varepsilon)^{-1}\mathcal T_2G_\theta(X),
                         \mathcal T_1E_z(X)\bigr).
 \end{aligned}\label{ic:detail-analytic-15}
\end{equation}
Let $\mathbb B$ be the closed unit ball of pairs with zero axis value
in the sum norm. Both radial inverses preserve this subspace.
For $X\in X^m+\mathbb B$, \eqref{ic:analytic-algebra} bounds the
products of $\partial_Z\mathcal M_x\Psi$ with $\Phi$, $x\partial_x\Phi$,
and $x\partial_x\Psi$ after radial inversion. The other transport
terms use the same bounds with at most one derivative.
The $Z$-dependent multipliers commute with $\mathcal T_\nu$.
For pressure, $\sup_\Omega|F_0|\le C_*^{-1}e^{3\Lambda}
\le\varepsilon^2$ by \eqref{ic:complex-smallness}, so
$\|F_0^2\|_h\le K_h\varepsilon^4$ by \eqref{ic:detail-analytic-5}.
The formula $\mathcal P(\Phi)=x F_0^2\mathcal M_x(\Phi^2)$ gives
\[
 \|\mathcal P(\Phi_1)-\mathcal P(\Phi_2)\|_h
 \le K_h\varepsilon^4(\|\Phi_1\|_h+\|\Phi_2\|_h)\|\Phi_1-\Phi_2\|_h.
\]
Equation \eqref{ic:analytic-algebra} also bounds $\mathcal T_1\partial_Z$
of this difference by a constant times the same right-hand side.
Taking differences in the remaining products gives a uniform Lipschitz
bound on the ball.
Hence, for $X_i=X^m+Y_i$, $Y_i\in\mathbb B$,
\begin{equation}
 \|\mathcal N_{\varepsilon,C_*}(X_1)\|\le K_*,\qquad
 \|\mathcal N_{\varepsilon,C_*}(X_1)
       -\mathcal N_{\varepsilon,C_*}(X_2)\|
       \le K_*\|Y_1-Y_2\|.
 \label{ic:detail-analytic-16}
\end{equation}
Here $K_*$ depends only on the fixed analytic data, not on
$0<\varepsilon\le1/500$ or the later $C_*$.

For $\Lambda\ge\max\{500,j_0^{-2},2K_*\}$, \eqref{ic:detail-analytic-15}
is a contraction of $\mathbb B$ into its radius-$1/2$ ball, with
constant at most $1/2$. Its fixed point has real coefficients and satisfies
\begin{equation}
 \|\Phi-\Phi^m\|_h
       +\|\Psi-\Psi^m\|_h
       \le K_*\varepsilon.
 \label{ic:comparison}
\end{equation}
By \eqref{ic:detail-analytic-3}, the rescaled fixed point is analytic
near $[0,4.1/\Lambda]\times[-1,1]$ and satisfies
\eqref{ic:normalized-core}, hence \eqref{ic:core-equations}.
Since $x\mathcal M_x$ raises degree and $\mathcal M_x,x\partial_x$ preserve it,
degree $n$ on the right uses only degrees $\le n$. The divisors at degree
$n+1$ are $2(n+1)(n+2)$ and $2(n+1)^2$. Hence induction and analytic
continuation identify the unique solution on the common connected domain
with the $\bar{\lambda}=0$ coefficient of Theorem~\ref{ex:core}.

Combining \eqref{ic:comparison} with \eqref{ic:detail-analytic-3}, for every fixed $k\ge0$,
\begin{equation}
 \|\Phi-\Phi^m\|_{C^k}
       +\|\Psi-\Psi^m\|_{C^k}
       \le K_k\varepsilon.
 \label{ic:detail-analytic-18}
\end{equation}
The norms use $(x,Z)\in[0,4.1]\times[-1,1]$. Hence, at $R=x/\Lambda$,
\begin{equation}
 \left\|\frac{F-F^m}{F_0}\right\|_{C^2}
       +\Lambda\|U^z-U^{z,m}\|_{C^2}
       \le\frac{K_2}{\Lambda}.
 \label{ic:detail-analytic-19}
\end{equation}
Axial derivatives act on the whole normalized difference.

\paragraph{Comparison of radial slopes.}
For this slope estimate near $Z=Z_0$, where $\chi_c=0$, write $D=\Phi-\Phi^m$.
Since
\[
 |H_0|=\sqrt{H_0^2+\sigma_0^2}\sqrt{\chi_c}
       \le K\sqrt{\chi_c},\qquad
 \frac{|H_0|}{H_0^2+\sigma_0^2}
       \le\sigma_0^{-1}\sqrt{\chi_c},
\]
the model formula and \eqref{ic:detail-analytic-18} give
$|\partial_{x}\Phi|\le K(\chi_c+\varepsilon)$.
The same comparison bounds $\Phi,\Psi,\partial_Z\Phi,
\partial_Z\Psi$; hence \eqref{ic:detail-analytic-11} implies
$|W-W_0|\le K\varepsilon$ and
$|H|\le K(\sqrt{\chi_c}+\varepsilon)$.
Equation \eqref{ic:detail-analytic-14} now gives
$|G_\theta(X)|\le K(\sqrt{\chi_c}+\varepsilon)$.
The difference satisfies
\[
\begin{aligned}
 2(x\partial_{x}^2D
           +2\partial_{x}D)
        &=-(\chi_c+\varepsilon\beta)D+\varepsilon G_\theta(X),\\
 \partial_{x}(x^2
       \partial_{x}D)
        &=\tfrac12x[-(\chi_c+\varepsilon\beta)D+\varepsilon G_\theta(X)].
\end{aligned}
\]
Since $|D|\le K\varepsilon$ and $0\le\chi_c\le1$, integrating
from the axis gives the following estimate, valid also at $R=0$ by continuity.
\begin{equation}
 |\partial_{x}\Phi
           -\partial_{x}\Phi^m|
       \le K(\varepsilon\sqrt{\chi_c}+\varepsilon^2).
 \label{ic:relative-slope}
\end{equation}
By \eqref{ic:bessel-bounds} and $\beta>1$,
$-\partial_{x}\Phi^m
\ge(\chi_c+\varepsilon\beta)/16
\ge(\chi_c+\varepsilon)/16$.
Thus
\begin{equation}
 \left|\frac{\partial_RF}{\partial_RF^m}-1\right|
 \le K\frac{\varepsilon\sqrt{\chi_c}+\varepsilon^2}
              {\chi_c+\varepsilon}
 \le K\sqrt{\varepsilon}.
 \label{ic:detail-analytic-21}
\end{equation}
Indeed, $\varepsilon\sqrt{\chi_c}/(\chi_c+\varepsilon)
\le\sqrt\varepsilon/2$ and
$\varepsilon^2/(\chi_c+\varepsilon)\le\varepsilon$.
Increase $\Lambda$ so that the last relative error is at most
$1/2$. Together with
$-\partial_{x}\Phi^m
\le(\chi_c+\varepsilon\beta)/4$ and $\beta<4$, this proves
\begin{equation}
 \frac1{32}(\chi_c+\varepsilon)
       \le-\partial_{x}\Phi
       \le2(\chi_c+\varepsilon).
 \label{ic:detail-analytic-22}
\end{equation}
Also $1/4<\Phi^m\le1$ and \eqref{ic:detail-analytic-18} imply
$1/8\le\Phi\le2$ and $F/F^m=1+O(\varepsilon)$.
Thus $F>0$, $\partial_RF<0$, and, for $R>0$,
$U^\theta=\sqrt{2R}F>0$ and $\mathcal S^\theta=2R\partial_RF<0$.
In particular,
\[
 c_0(\Lambda\chi_c+1)
       \le-\frac{\partial_RF}{F}
       =-\Lambda\frac{\partial_{x}\Phi}{\Phi}
       \le K(\Lambda\chi_c+1).
\]

\paragraph{The exit and its continuation.}
By \eqref{ic:weak-sector}, $g\ge 2j_0/5$ on $[-j_0/4,-j_0/5]$,
which contains $\{\chi_c\le99/100\}$. The axial comparison gives
\begin{equation}\label{ic:detail-analytic-23}
\begin{aligned}
 \partial_RU^z(R,Z)&=-\frac{g(Z)}{2L_\delta(Z)}+O(\varepsilon),\\
 \partial_RU^z(R,Z)&\le-\frac{j_0}{10}\qquad (-j_0/4\le Z\le-j_0/5).
\end{aligned}
\end{equation}
Indeed, the model slope is at most $-j_0/5$ there and the error
is at most $j_0/10$ for $\Lambda\ge10K/j_0$. At $R_{\rm in}$, the ratios
of the angular and axial contributions to their respective model values are
\[
 \frac{\partial_RF}{\partial_RF^m}\frac{F^m}{F}
       =1+O(\sqrt\varepsilon),\qquad
 \left(\frac{\partial_RU^z}{\partial_RU^{z,m}}\right)^2
 \frac{F^m\partial_RF^m}{F\partial_RF}=1+O(\sqrt\varepsilon)
\]
in the strong and weak sectors, respectively. For the second ratio,
\eqref{ic:detail-analytic-23} gives a relative axial slope error at
most $5K\varepsilon/j_0$. Choose $\Lambda\ge100K/j_0$ large enough that
both ratios exceed $5/6$. The model bounds following
\eqref{ic:bessel-bounds} and \eqref{ic:weak-sector} then give
\begin{equation}
 \begin{aligned}
 -2R_{\rm in}\frac{\partial_RF}{F}&>\frac52
       &&(\chi_c\ge99/100),\\
 \frac{(\partial_RU^z)^2}{-F\partial_RF}&>\frac52
       &&(\chi_c\le99/100)
 \end{aligned}
 \qquad\text{at }R=R_{\rm in}=4/\Lambda.
 \label{ic:detail-analytic-24}
\end{equation}
Both contributions are nonnegative, so $\kappa(R_{\rm in},Z)\ge5/2$.
By compactness, $\kappa>2$ for all $|Z|\le1$ on a fixed right
neighborhood of $R_{\rm in}$ contained in $[R_{\rm in},4.1/\Lambda]$. Since $4e^{1/100}<4.1$ and
$1/8\le\Phi\le2$, \eqref{ic:detail-analytic-18} applied to
$C_*F_c=e^{-\Lambda G}\Phi$ and $U_c^z=U_0^z+\Lambda^{-1}\Psi$
gives \eqref{ic:core-norms} uniformly in $C_*$ for fixed $\Lambda$.
Equation \eqref{ic:detail-analytic-22} gives $-\partial_R\log F_c\ge1/64$.
On the weak sector, \eqref{ic:weak-sector} and \eqref{ic:detail-analytic-23}
give $-2L_\delta\partial_RU_c^z=g+O(\varepsilon)\ge3j_0/10$ for large $\Lambda$.
The axial comparison and its radial average give
\[
\begin{aligned}
 U_c^z(R_{\rm in})&=4Z+j_0-\frac{2g}{\Lambda L_\delta}
                         +O_{C_Z^2}(\Lambda^{-2}),\\
 \frac{M_c^z(R_{\rm in})}{R_{\rm in}}
       &=4Z+j_0-\frac{g}{\Lambda L_\delta}
                         +O_{C_Z^2}(\Lambda^{-2}).
\end{aligned}
\]
Finally, $\partial_Z\log F_c(R_{\rm in})+\Lambda\partial_ZG
=\partial_Z\log\Phi(4,Z)$. From $\chi_c=H_0^2/(H_0^2+\sigma_0^2)$
we have $|\partial_Z\chi_c|\le B\sqrt{\chi_c}$, so
\[
 \left|\partial_Z\log\Phi^m\right|
 =\left|\frac{B'(q)}{B(q)}\frac{x}{2}
          (\partial_Z\chi_c+\varepsilon\partial_Z\beta)\right|
 \le B(\sqrt{\chi_c}+\varepsilon)
\]
by \eqref{ic:bessel-bounds}. The $C^1$ comparison and the positive
lower bounds for $\Phi,\Phi^m$ give \eqref{ic:exit-z}.
\end{proof}

The leading equations give zero radial derivatives for
$R\mathcal T^\theta$ and $\sqrt R\,\mathcal T^z$. Both constants vanish by axis
regularity, so the leading stress vanishes throughout the core.

\subsubsection{The full analytic core remainder}
\label{ex:matching}

Let $G_c(\bar{\lambda},R,Z)$ denote any of the four normalized quantities in
\eqref{ex:normalized}, with physical parameter $\bar{\lambda}=\lambda^{2\delta}$.
For the normalized angular and axial velocities, the triangular recurrence
and Proposition~\ref{ic:core-exit} identify their values at $\bar{\lambda}=0$ with $F_c$
and $U_c^z$, respectively, on $0\le R\le4.1/\Lambda$.

\begin{proposition}[A single analytic remainder]\label{ex:remainder}
For each such $G_c$, we have
\begin{equation}\label{ex:scale-remainder}
 G_c(\bar{\lambda},R,Z)-G_c(0,R,Z)
 =\bar{\lambda}\int_0^1
       (\partial_{\bar{\lambda}}G_c)(s\bar{\lambda},R,Z)\,\mathrm ds.
\end{equation}
The integral is holomorphic for small $|\bar{\lambda}|$ and $(R,Z)$ in a fixed
neighborhood of $[0,4.1/\Lambda]\times[-1,1]$. On smaller neighborhoods,
all its $(\bar{\lambda},R,Z)$ derivatives are bounded uniformly in
$0\le\delta\le1/200$ for fixed $P_*,j_0,\Lambda,C_*$.
On the real domain, for $\delta>0$ and all integers $i,j,k\ge0$,
\begin{equation}\label{ex:remainder-factor}
 (\lambda\partial_\lambda)^i\partial_R^j\partial_Z^k
 \bigl[G_c(\lambda^{2\delta},R,Z)-G_c(0,R,Z)\bigr]
 =O_{i,j,k}(\lambda^{2\delta}).
\end{equation}
\end{proposition}

\begin{proof}
Integrating in $\bar{\lambda}$ at fixed $(R,Z,\delta)$ gives
\eqref{ex:scale-remainder}. Cauchy's formula on the common neighborhood in
Theorem~\ref{ex:core} bounds the derivatives of the integral.
For $\bar{\lambda}=\lambda^{2\delta}$, the identity
$\lambda\partial_\lambda=2\delta\bar{\lambda}\partial_{\bar{\lambda}}$
preserves the factor $\bar{\lambda}$ in \eqref{ex:scale-remainder}, giving \eqref{ex:remainder-factor}.
\end{proof}

Define $\kappa_c$ by \eqref{intro:physical-shear}. At $\bar{\lambda}=0$,
the swirl is positive, the angular shear is negative, and
$\kappa_c(0,R_{\rm in},Z)\ge5/2$ by \eqref{ic:kappa-exit}.
By Proposition~\ref{ex:remainder}, decreasing $\lambda_*$ gives
\begin{equation}\label{ex:exit}
 \kappa_c(\lambda,R_{\rm in},Z)\ge\frac94
 \qquad(0<\lambda\le\lambda_*,\ |Z|\le1).
\end{equation}
Compactness gives a fixed collar of $R_{\rm in}$, within $0<R<4.1/\Lambda$,
on which $\kappa_c>2$.

\subsubsection{Leaving the core with the actual stress}
\label{ic:exit-connection}

For a pair $(F,U^z)$ and its five moments, use the quantities in
\eqref{oc:normalized}, with $y=\log R$.
Since $U^\theta=\sqrt{2R}\,F$, we have $a=-2\partial_y\log F$ and
$\zeta=\partial_Z\log F$. Recall the axial transport coefficient
$H=(1-\delta)Z/2+LU^z$ from \eqref{ic:core-equations}.
In the following identity only, write $\mathcal Pw=2(1+\delta)Zw-L\partial_Z w$.
The differentiated moment identities \eqref{pr:inertia-theta}--\eqref{pr:moment-inertia} give
\begin{equation}\label{ic:transport}
\begin{aligned}
 \partial_y\widetilde I^\theta+(2-a/2)\widetilde I^\theta
 &=\frac{L_\delta}{R}\left[\partial_y\!\left(\frac{\mathcal I^\theta}{F}\right)
                  +(1-a/2)\frac{\mathcal I^\theta}{F}\right]\\
 &=-W(1-a/2)-\tfrac\delta2(1-2ZU^z)-H\zeta,\\
 (\partial_y+1)\bigl[(U^\theta)^2\widetilde I^z\bigr]
 &=g[U^z]+\mathcal P(P-P_0)+Z(U^\theta)^2-W\partial_y U^z,\\
 g[U^z]&=-\tfrac{1+\delta}{2}(1-2ZU^z)U^z-H\partial_Z U^z+\mathcal PP_0 .
\end{aligned}
\end{equation}
Indeed, $\mathcal I^\theta/F=R\widetilde I^\theta/L_\delta$. Multiply the axial equation in
\eqref{oc:radial-system} by $(U^\theta)^2$ and use $2\partial_y\log U^\theta=1-a$.
This gives the axial identity since $(U^\theta)^2\widetilde I^z=L_\delta\sqrt{2/R}\,\mathcal I^z$.
At the axis, $g[U_0^z]=g$ by \eqref{ic:axis-data}.

\paragraph{Constants for the exit collar.}
Define
\begin{equation}\label{ic:detail-data-sizes}
\begin{aligned}
 A&=10+\sup_{0\le\delta\le1/200}
       \bigl(\|\log(C_*F_c)\|_{C^3}+\|U_c^z\|_{C^3}\bigr),\\
 K&=10^6+100+R_{\rm in}^{-1}+C_*+P_*+A\\
  &\quad+\sup_{0\le\delta\le1/200}\left(
       \|F_c\|_{C^3}+\|1/F_c\|_{C^3}
       +\|P_0\|_{C_Z^3}+\sum_j\|M_c^j\|_{C^3}\right).
\end{aligned}
\end{equation}
The sum runs over the five moments, and the joint norms use
$(R/R_{\rm in},Z)\in[0,e^{1/100}]\times[-1,1]$.
The suprema are taken once for each fixed $C_*$; thus $A,K$ and
the cutoff and radii below are independent of $\delta$.
By \eqref{ic:core-norms}, $A=O_\Lambda(1)$ and
$C_*\le K=O_\Lambda(C_*)$.
Let $\sigma_{\rm norm}$ be given by \eqref{pr:cutoffs}, with $0\le\sigma_{\rm norm}'\le8$.

\paragraph{The core comparison on the short exit interval.}
For $y=\log(R/R_{\rm in})$, set
\begin{equation}\label{ic:comparison-profile}
\begin{aligned}
(\bar F,\bar M^j)(R,Z)
       &=(F_c,M_c^j)(R,Z),\\
 0&\le y\le\frac1{C_*K^{100}},
 &\bar P&=P_0+\bar M^p.
\end{aligned}
\end{equation}
Since $R_{\rm in}e^{1/(C_*K^{100})}<4.1/\Lambda$, the continued core is defined
throughout the collar. For this collar comparison only, write
$q=(q^\theta,q^z)=\bar{\mathcal I}/\bar F$.
The leading core has zero stress, so its shear parameter from
\eqref{ic:kappa-exit} is $\kappa_c=|q|^2/q^\theta$.

\begin{lemma}[Direction of the continued leading core]\label{ic:direction}
For sufficiently large $\Lambda$, then sufficiently large $C_*$, we have
\begin{equation}\label{ic:direction-margin}
\begin{aligned}
q^\theta&\ge R_{\rm in}/64,
 &\quad\kappa_c&\ge\frac94\qquad\left(0\le y\le\frac1{C_*K^{100}}\right).
\end{aligned}
\end{equation}
For each fixed $m$ there is a constant $B_{m,\Lambda}$, independent
of later $C_*$, such that
\[
 \|q^\theta\|_{C_Z^m}+\|1/q^\theta\|_{C_Z^m}
          +C_*^{-1}\|q^z\|_{C_Z^m}\le B_{m,\Lambda}.
\]
\end{lemma}

\begin{proof}
Since the core stress vanishes,
\[
\begin{aligned}
q^\theta&=-2R\,\partial_R\log F_c,
 &q^z&=-\sqrt{2R}\,\partial_R U_c^z/F_c.
\end{aligned}
\]
By \eqref{ic:core-norms}, all derivatives of $q^\theta$ and $-2L_\delta\partial_R U_c^z$ in $R,Z$
are bounded independently of $C_*$ at each order. Their changes
from $R_{\rm in}$ in $C_Z^m$ are at most $B_{m,\Lambda}/(C_*K^{100})$.
At $R_{\rm in}$, $q^\theta\ge R_{\rm in}/32$ and $q^\theta>5/2$ on
$\{\chi_c\ge99/100\}$. On its complement, $-2L_\delta\partial_R U_c^z\ge3j_0/10$.
Choose $C_*$ so that
\[
 \frac{B_{2,\Lambda}}{C_*K^{100}}\le\min\{R_{\rm in}/64,1/8,j_0/30\}.
\]
Thus throughout the collar
\begin{equation}\label{ic:layer-direction}
 q^\theta\ge R_{\rm in}/64,\qquad
 q^\theta\ge19/8\qquad (\chi_c\ge99/100),\qquad
 -2L_\delta\partial_R U_c^z\ge4j_0/15\qquad (\chi_c\le99/100).
\end{equation}
For larger $C_*$, $\bar F\le B_{0,\Lambda}/C_*$ gives
$q^z\ge9/8$ on the weak sector, hence $\kappa_c\ge2|q^z|\ge9/4$.
On the strong sector, $\kappa_c\ge q^\theta\ge19/8$.
Differentiating $1/q^\theta$ and the formula for $q^z$ gives the remaining bounds.
\end{proof}

\paragraph{The actual connection.}
For larger $C_*$,
\begin{equation}\label{ic:crude-bounds}
 q^\theta\ge(2K)^{-1},\qquad |q|\le2K^3,\qquad
 \|q\|_{C_Z^2}\le K^7,\qquad \kappa_c\le K^{10}.
\end{equation}
Set
\begin{equation}\label{ic:exit-cutoff}
 \chi_b(y)=1-\left(1-\frac1{C_*K^{100}}\right)\sigma_{\rm norm}(C_*K^{100}y).
\end{equation}
Define the actual admissible exit by
\begin{equation}\label{ic:strict-outlet}
 t_{ac}=\sigma_{\rm norm}^{-1}\left(\frac1{16K^{10}}\right),\qquad
 R_{ac}=R_{\rm in}\exp\left(\frac{t_{ac}}{C_*K^{100}}\right),\qquad 0<t_{ac}<1.
\end{equation}
Prescribe
\begin{equation}\label{ic:actual-bridge}
\begin{aligned}
\partial_y\log F&=-\tfrac12\chi_b q^\theta,
 &\quad\partial_y U^z&=-\chi_b\sqrt{R/2}\,F q^z,\\
 (F,U^z)|_{R_{\rm in}}&=(F_c,U_c^z)|_{R_{\rm in}},
 &\quad R_{\rm in}&<R\le R_{ac}.
\end{aligned}
\end{equation}
Integrating gives
\begin{align}
 F(y,Z)&=F_c(R_{\rm in},Z)\exp\left[-\frac12\int_0^y
                  \chi_b(t)q^\theta(t,Z)\dd t\right],\notag\\
 U^z(y,Z)&=U_c^z(R_{\rm in},Z)-\int_0^y\chi_b(t)
       \sqrt{R_{\rm in}e^t/2}\,F(t,Z)q^z(t,Z)\dd t.
 \label{ic:detail-bridge-quadratures}
\end{align}
Thus $\mathcal S=-\chi_bFq$ and $F>0$. Continue all five moments
from their core values, retain $P_0$, and define $\mathcal I$ by
\eqref{pr:inertia-theta}--\eqref{pr:inertia-z}.
For the continued core, \eqref{ic:actual-bridge} holds with $\chi_b=1$.
To estimate $e=\mathcal I/F-q$, put $\omega=1-\chi_b$ and
$\ell=\log(F/\bar F)$. Subtracting gives
\[
 \partial_y\ell=\tfrac12\omega q^\theta,\qquad
 \partial_y(U^z-U_c^z)=(\chi_be^\ell-1)\partial_y U_c^z.
\]
Integrate in $y$, using $\int_0^y\omega(s)\dd s\le\omega(y)/(C_*K^{100})$
by monotonicity and $\chi_be^\ell-1=-\omega+\chi_b(e^\ell-1)$. We obtain
\[
 \|\ell\|_{C_Z^2}+\|U^z-U_c^z\|_{C_Z^2}
          \le\frac{C\omega}{C_*K^{99}},\qquad
 \|F-\bar F\|_{C_Z^2}\le\frac{C\omega}{C_*K^{98}}.
\]
In particular, throughout the actual collar,
\[
 X(y):=\sup_{s\le y}
   \big(\|F-\bar F\|_{C_Z^1}+\|U^z-U_c^z\|_{C_Z^1}\big)
       \le\frac{C\omega(y)}{C_*K^{98}}.
\]
Let $\Delta F=F-\bar F$, $\Delta U^z=U^z-U_c^z$.
Subtracting \eqref{pr:moments}, with zero differences at $R_{\rm in}$, gives
\begin{equation}\label{ic:detail-moment-differences}
\begin{gathered}
 \Delta M^\theta=2\int_{R_{\rm in}}^R\rho\Delta F\dd\rho,\qquad
 \Delta M^z=\int_{R_{\rm in}}^R\Delta U^z\dd\rho,\\
 \Delta M^{\theta z}=2\int_{R_{\rm in}}^R
        \rho(U^z\Delta F+\bar F\Delta U^z)\dd\rho,\qquad
 \Delta M^p=\int_{R_{\rm in}}^R(F+\bar F)\Delta F\dd\rho,\\
 \Delta M^{z\theta}=\int_{R_{\rm in}}^R
       [(U^z+U_c^z)\Delta U^z-\rho(F+\bar F)\Delta F]\dd\rho.
\end{gathered}
\end{equation}
Since $R\le R_{ac}<1$, the $C_Z^1$ product rule and the common
$CK$ velocity bound give
\begin{equation}\label{ic:detail-moment-error}
 \sum_j\|\Delta M^j\|_{C_Z^1}+\|\Delta P\|_{C_Z^1}
          \le CKX(y),\qquad \Delta P=\Delta M^p.
\end{equation}
In \eqref{pr:moment-inertia}, $R^{-1}\le K$, $L_\delta^{-1}\le2$,
$F^{-1}+\bar F^{-1}\le CK$, and $|F^{-1}-\bar F^{-1}|\le CK^2X$.
The numerators of the moment fractions are $O(K^3)$ and their
differences are $O(K^2X)$ by \eqref{ic:detail-moment-error}.
Taking differences after division by $F$ also gives
$|(1-\delta)Z\Delta M^z+L\partial_Z\Delta M^z|/L_\delta\le CKX$ for the remaining angular term. Hence
\[
 |e|\le CK^6X(y)\le\frac{C\omega(y)}{C_*K^{84}}.
\]
By \eqref{ic:crude-bounds}, for sufficiently large $C_*$,
\begin{equation}\label{ic:cone-error}
 |e|\le\min\left\{\frac{q^\theta}{16(1+|q|)},
       \frac{\omega|q|}{10\sqrt{1+\kappa_c}}\right\}.
\end{equation}
Hence
\[
 \mathcal T/F=\omega q+e,\qquad \kappa=\chi_b\kappa_c,\qquad
 |\mathcal T|/F\ge\tfrac9{10}\omega|q|>0.
\]
If $\kappa\le2$, the relaxed admissibility test is
$\kappa_c+e\cdot q/q^\theta>2$, and its left side is at least
$9/4-1/16$.

If $\kappa>2$, then $-\mathcal T\cdot \mathcal S>0$ and
\begin{equation}\label{ic:actual-angle}
 \frac{|\mathcal T\cdot \mathcal S^\perp|}{-\mathcal T\cdot \mathcal S}
 \le\frac{|e|}{\omega|q|-|e|}
 \le\frac1{9\sqrt{1+\kappa_c}}
 <\sqrt{\frac2{\kappa-2}}.
\end{equation}
For $R_{\rm in}<R\le R_{ac}$, \eqref{ic:strict-outlet} gives
$\omega\kappa_c\le\sigma_{\rm norm}(t_{ac})K^{10}\le1/16$, hence
$\kappa\ge35/16>2$. Thus the collar is strictly admissible.
On $0<y\le y_0<\min\{1/2,t_{ac}\}/(C_*K^{100})$, $F|q|$ is bounded
above and below by positive constants, so uniformly for $|Z|\le1$,
\begin{equation}\label{ic:inner-flat-size}
 c\omega(y)\le |\mathcal T(R_{\rm in}e^y,Z)|\le C\omega(y),\qquad c,C>0.
\end{equation}

The differences have zero initial data, and $\omega$ is a Gaussian
$e^{-1/(C_*^2K^{200}y^2)}$ times a smooth positive factor.
Integrating the difference equations for $\ell,U^z-U_c^z$ and using
$e^\ell-1=\ell\int_0^1e^{s\ell}\dd s$ gives the same Gaussian bound
for the velocity differences and their axial derivatives. Each further
$y$ derivative costs at most a finite power of $y^{-1}$. By
\eqref{ic:detail-moment-differences} and $F,\bar F>0$ at $R_{\rm in}$,
all jets there agree with the core and, for each $i,k\ge0$,
\begin{equation}\label{ic:flat-error}
 |\partial_y^i\partial_Z^ke|
       \le C_{ik}e^{-1/(C_*^2K^{200}y^2)}y^{-m_{ik}}
       \quad\left(0<y<\frac{\min\{1/2,t_{ac}\}}{C_*K^{100}}\right).
\end{equation}
Here $m_{ik},C_{ik}$ may depend on the fixed profile.

\subsubsection{The actual exit, its pressure and five moments}

\begin{lemma}[Bounds at the actual admissible exit]\label{ic:small-defect}\label{ic:new-entrance-data}
Fix $P_*\ge1$, an integer $m\ge2$, and a positive tolerance for the axial errors.
For $P_0$ as in \eqref{ap:datum}, choose $j_0>0$ sufficiently small,
then $\Lambda$ and $C_*$ sufficiently large, in that order.
The following bounds put both axial errors below the prescribed tolerance, uniformly for $0\le\delta\le1/200$.
\begin{equation}\label{ic:axial-tolerance}
 \sup_{0\le R\le R_{ac}}\|U^z(R)-4Z\|_{C_Z^m}<2j_0,\qquad
 \left\|\frac{M^z(R_{ac})}{R_{ac}}-4Z\right\|_{C_Z^m}<2j_0.
\end{equation}
On the same interval,
\begin{equation}\label{ic:prefix-bound}
 \sup_{0\le R\le R_{ac}}
 \left(C_*\|F(R)\|_{C_Z^m}
       +C_*^{-1}\|F(R)^{-1}\|_{C_Z^m}
       +\|U^z(R)\|_{C_Z^m}\right)\le B_{m,\Lambda},
\end{equation}
and the five moments satisfy
\begin{equation}\label{ic:prefix-moment-bounds}
\begin{aligned}
 \sup_{0\le R\le R_{ac}}\bigl(&
 C_*\|M^\theta(R)\|_{C_Z^m}
 +C_*\|M^{\theta z}(R)\|_{C_Z^m}
 +\|M^z(R)\|_{C_Z^m}\\
 &+\|M^{z\theta}(R)\|_{C_Z^m}
 +C_*^2\|M^p(R)\|_{C_Z^m}\bigr)\le B_{m,\Lambda}.
\end{aligned}
\end{equation}
Here $P=P_0+M^p$ and $\|L_\delta\sqrt{2/R_{ac}}\,\mathcal I^z(R_{ac})\|_{C_Z^m}\le B_{m,\Lambda}$.
The scaled bounds \eqref{ic:prefix-bound}--\eqref{ic:prefix-moment-bounds}
and this axial bound hold at every fixed derivative order, with
order-dependent constants, once $\Lambda$ is fixed. At the exit we also have
\begin{equation}\label{ic:new-exit-linear}
\begin{aligned}
 &\|U^z(R_{ac})-(4Z+j_0)\|_{C_Z^m}
   +\|M^z(R_{ac})/R_{ac}-(4Z+j_0)\|_{C_Z^m}
       \le B_m/\Lambda,\\
 &\|M^z(R_{ac})-R_{ac}U^z(R_{ac})\|_{C_Z^m}
       \le B_m/\Lambda^2,
\end{aligned}
\end{equation}
and
\begin{equation}\label{ic:new-exit-transport}
 H(R_{ac},Z)\partial_Z\log F(R_{ac},Z)
       \le-\frac34\Lambda L_\delta\chi_c(Z)+\frac B\Lambda.
\end{equation}
With $g[\cdot]$ as in \eqref{ic:transport}, the exit satisfies
\begin{equation}\label{ic:new-exit-sectors}
\begin{aligned}
 \frac{\mathcal I^\theta(R_{ac},Z)}{F(R_{ac},Z)}&\ge R_{ac}/256,\\
 \frac{\mathcal I^\theta(R_{ac},Z)}{F(R_{ac},Z)}&>9/4 &&(\chi_c\ge99/100),\\
 L_\delta\sqrt{2/R_{ac}}\,\mathcal I^z(R_{ac},Z)&\ge j_0/5,\qquad
 g[U^z(R_{ac},Z)]\ge2j_0/5 &&(\chi_c\le99/100).
\end{aligned}
\end{equation}
The constants $B_m,B$ may depend on $P_*,j_0,m$, and $B_{m,\Lambda}$ also on
$\Lambda$. All are independent of $\delta$ and of subsequent increases in $C_*$.
At fixed $P_*,j_0,\Lambda,C_*$ the pair $(F,U^z)$ and its five moments are
smooth in $(R,Z,\delta)$ up to the fixed radius $R_{ac}$, including $\delta=0$.
\end{lemma}

\begin{proof}
Choose $0<j_0\le1/20$ at most one eighth of the prescribed tolerance.
The model $U^{z,m}=4Z+j_0-gR/(2L_\delta)$ and \eqref{ic:detail-analytic-18}
give $\|U_c^z-(4Z+j_0)\|_{C_Z^m}\le B_m/\Lambda$ throughout the continued
core. The quantities $C_*F_c$, its reciprocal, $U_c^z$, $q^\theta$ and
$q^z/C_*$ have $C_Z^m$ bounds independent of later $C_*$ by
\eqref{ic:core-norms} and Lemma~\ref{ic:direction}.
Since $\chi_b$ is independent of $Z$ and $0\le y\le1/(C_*K^{100})$,
integrating $\partial_y\ell=\omega q^\theta/2$ gives
$\|\ell\|_{C_Z^m}\le B_{m,\Lambda}/(C_*K^{100})$. Using $F=F_ce^\ell$ in
\eqref{ic:detail-bridge-quadratures} gives on the collar
\[
 C_*\|F-F_c\|_{C_Z^m}+\|U^z-U_c^z\|_{C_Z^m}\le\frac{B_{m,\Lambda}}{C_*K^{100}}.
\]
This gives \eqref{ic:prefix-bound} at every fixed order. For the chosen $m$,
increase $\Lambda$ and then $C_*$ so that the last error is at most
$1/\Lambda$ and \eqref{ic:complex-smallness}, \eqref{ic:cone-error} hold.
Then $\|U^z-(4Z+j_0)\|_{C_Z^m}\le B_m/\Lambda$ throughout the inner interval.
Radial averaging gives the same bound for $M^z/R-(4Z+j_0)$. These inequalities
give the first line of \eqref{ic:new-exit-linear} after enlarging $B_m$,
and \eqref{ic:axial-tolerance} upon taking $\Lambda$ large enough that $B_m/\Lambda<j_0$.
By \eqref{ic:detail-analytic-18}, $\|\partial_R U_c^z\|_{C_Z^m}\le B_m$.
On the collar $\partial_R U^z=\chi_b e^\ell\partial_R U_c^z$, with
$\|e^\ell-1\|_{C_Z^m}\le B_{m,\Lambda}(C_*K^{100})^{-1}$.
The second line of \eqref{ic:new-exit-linear} follows from
\[
 M^z(R_{ac})-R_{ac}U^z(R_{ac})
       =-\int_0^{R_{ac}}R\,\partial_R U^z(R)\,dR,
 \qquad R_{ac}<4.1/\Lambda.
\]

Integrating the five densities in \eqref{pr:moments} in $C_Z^m$ gives
\eqref{ic:prefix-moment-bounds} at that order, with
$P-P_0=M^p=O_{m,\Lambda}(C_*^{-2})$. In \eqref{pr:moment-inertia},
$\mathcal I^z$ uses at most one axial derivative. Since $R_{ac}\ge4/\Lambda$,
the moment bounds at order $m+1$ give the stated bound for
$L_\delta\sqrt{2/R_{ac}}\,\mathcal I^z(R_{ac})$.

Equation \eqref{ic:exit-z}, the core derivative bounds, and the
short collar give
\[
 H(R_{ac},Z)=H_0+O(\Lambda^{-1}),\qquad
 \partial_Z\log F(R_{ac},Z)
       =-\Lambda\partial_Z G+O(\sqrt{\chi_c}+\Lambda^{-1}).
\]
Since $H_0\partial_Z G=L_\delta\chi_c$ and
$|H_0|+|\partial_Z G|\le C\sqrt{\chi_c}$, multiplication gives
\[
 H(R_{ac})\partial_Z\log F(R_{ac})
       \le-\Lambda L_\delta\chi_c
          +C(\chi_c+\sqrt{\chi_c}+\Lambda^{-1}).
\]
Young's inequality absorbs the error into
$\Lambda L_\delta\chi_c/4+B/\Lambda$ for large $\Lambda$,
proving \eqref{ic:new-exit-transport}.

By definition $\mathcal I^\theta/F=q^\theta+e^\theta$. Equations \eqref{ic:layer-direction}
and \eqref{ic:cone-error} give $\mathcal I^\theta/F\ge(15/16)q^\theta>R_{ac}/256$
and, on the strong sector, $\mathcal I^\theta/F\ge19/8-1/16>9/4$.
On the weak sector $\ell\ge0$, $q^z\ge9/8$, and $|e^z|\le1/16$, so
\[
 L_\delta\sqrt{2/R_{ac}}\,\mathcal I^z(R_{ac})
 =\frac F{F_c}(-2L_\delta\partial_R U_c^z)
                 \left(1+\frac{e^z}{q^z}\right)\bigg|_{R=R_{ac}}
       \ge\frac{17}{18}\frac{4j_0}{15}>j_0/5.
\]
There $g[4Z+j_0]\ge3j_0/5$ by \eqref{ic:weak-sector}, and
\eqref{ic:new-exit-linear} gives $g[U^z(R_{ac})]\ge2j_0/5$ for large $\Lambda$.

The normalized core is holomorphic in $\delta$ by
Theorem~\ref{ex:core}. The constants in
\eqref{ic:detail-data-sizes}, the cutoff, and $R_{ac}$ are fixed
over the whole $\delta$ interval. Thus the bridge quadratures and
the five moment integrals depend smoothly on $\delta$ down to zero.
\end{proof}

\begin{theorem}[Inner profile and its actual admissible exit]
\label{ic:inner-connection}
Let $P_0$ be given by \eqref{ap:datum}. Fix the axial tolerance and the
derivative order in \eqref{ic:axial-tolerance}. Then one can choose $j_0>0$, followed by sufficiently
large $\Lambda$ and $C_*$ in that order, so that there is a leading pair with
$(F,U^z)$ smooth on $[0,R_{ac}]\times[-1,1]$,
where $R_{\rm in}=4/\Lambda$ and $R_{ac}$ is given by \eqref{ic:strict-outlet}.

On $R\le R_{\rm in}$, it equals the leading coefficient of the exact core
in Theorem~\ref{ex:core}, with axis data \eqref{ic:axis-data} and $P_0$.
Its stress vanishes in the core, is flat at $R_{\rm in}$, and is nonzero and
strictly admissible on $R_{\rm in}<R\le R_{ac}$, with $\kappa\ge35/16$.
The axis pressure is $P_0$, and all five moments continue from their core values at $R_{\rm in}$.

The exit estimates and smooth dependence on $\delta$ at fixed $R_{\rm in},R_{ac}$
are those of Lemma~\ref{ic:new-entrance-data}.
\end{theorem}
\begin{proof}
Combine Theorem~\ref{ex:core}, Proposition~\ref{ic:core-exit}, the collar
construction \eqref{ic:actual-bridge} and its stress estimates
\eqref{ic:cone-error}--\eqref{ic:flat-error}, and Lemma~\ref{ic:new-entrance-data}.
\end{proof}
\section{The intermediate connection and the exact background}\label{mid:section}
\subsection{The exact mixed power and its connection to the exterior}
\label{pow:section}

The middle power has the exterior velocity ratio
\(U_e^z/U_e^\theta=\sqrt{(1+\epsilon_{\rm em}^2)/2}\).
We reduce the power-law amplitude \(h_\infty\) by a fixed fraction \(0<\epsilon_0<1\)
so that a single cutoff joins the profiles with strictly admissible stress.
Set \(V_p=(1-\epsilon_0)h_\infty R^{-(1+\delta)/2}\),
\(R_e=e^{2\ell}R_{\rm mid}\), and \(d=\log(R_e/R)\).
With the cutoff from \eqref{pr:cutoffs}, the two profiles are
\begin{equation}\label{pow:profiles}
 U^\theta=
 \begin{cases}
 V_p,&R_{\rm mid}\le R\le e^\ell R_{\rm mid},\\
 U_e^\theta-\sigma_{\rm norm}(d/\ell)(U_e^\theta-V_p),
     &e^\ell R_{\rm mid}<R<R_e .
 \end{cases}
 \qquad U^z=\frac{U_e^z}{U_e^\theta}\,U^\theta .
\end{equation}
On both intervals the energy density is
\((U^z)^2-\tfrac12(U^\theta)^2=\tfrac12\epsilon_{\rm em}^2(U^\theta)^2\),
as in the exterior.
For \(R\ge R_e\) we use the exact mixed exterior.
As \(R\) increases across the transition, the cutoff decreases from one to zero. It is analytic on
the open transition and flat relative to its endpoint values. Hence the
velocities match to every order at both ends. By \eqref{ue:heat-factor}
and Lemma~\ref{ue:regularized-tail}, the heat factor differs from one
by \(O(\delta/R_{\rm mid})\) on this collar, also after a logarithmic
radial derivative or any fixed number of \(Z\) derivatives. Consequently
\[
 U_e^\theta-V_p\asymp\epsilon_0 U_e^\theta>0,\qquad
 |R\partial_R(U_e^\theta-V_p)|+|\partial_Z(U_e^\theta-V_p)|\le C(U_e^\theta-V_p)
\]
when \(\delta/R_{\rm mid}\ll\epsilon_0\).

\begin{figure}[!bhp]
\centering
\definecolor{powink}{RGB}{32,51,65}
\definecolor{powprep}{RGB}{53,100,133}
\definecolor{powteal}{RGB}{12,99,94}
\definecolor{powochre}{RGB}{139,83,31}
\begin{tikzpicture}[x=1cm,y=1cm,font=\small,draw=powink,text=powink,
  every node/.style={inner sep=1.2pt},line cap=round,line join=round]
  \node[anchor=west,font=\small\bfseries] at (0,1.48)
    {Radial profiles};
  \node[anchor=east,font=\footnotesize] at (15.0,1.48)
    {schematic spacing};
  \fill[powink!7]    (0,0) rectangle (2.38,0.82);
  \fill[powprep!12]  (2.38,0) rectangle (6.82,0.82);
  \fill[powteal!15]  (6.82,0) rectangle (9.70,0.82);
  \fill[powochre!15] (9.70,0) rectangle (12.48,0.82);
  \fill[powink!7]    (12.48,0) rectangle (15.0,0.82);
  \foreach \a/\b/\ink in {0/2.38/powink,2.38/6.82/powprep,
      6.82/9.70/powteal,9.70/12.48/powochre,12.48/15/powink}
    \draw[\ink,line width=1.05pt] (\a,.82)--(\b,.82);
  \node at (1.19,.40) {Inner profile};
  \node[text=powprep] at (4.60,.40) {Analytic preparation};
  \node[text=powteal] at (8.26,.40) {Exact power};
  \node[text=powochre] at (11.09,.40) {Transition};
  \node at (13.74,.40) {Mixed exterior};
  \foreach \x/\idx in {4.60/1,8.26/2,11.09/3}
    \node[font=\footnotesize] at (\x,1.05) {$M_{\idx}$};
  \draw[->] (0,0)--(15.22,0);
  \node[anchor=west,font=\footnotesize] at (15.27,0) {$R$};
  \foreach \x/\rad in {0/0,1.13/R_{\rm in},2.38/R_{ac},6.82/R_{\rm mid},
       9.70/e^\ell R_{\rm mid},12.48/R_e} {
    \draw (\x,0)--(\x,-.12);
    \node[anchor=north] at (\x,-.16) {$\rad$};}
  \draw[powprep,line width=.95pt] (3.31,-.14)--(3.91,-.14);
  \draw[powprep,line width=.95pt] (4.56,-.14)--(5.16,-.14);
  \node[anchor=north,font=\footnotesize,text=powprep] at (3.61,-.19) {$I_2$};
  \node[anchor=north,font=\footnotesize,text=powprep] at (4.86,-.19) {$I_3$};
  \begin{scope}[yshift=-1.04cm,>=latex]
    \fill[powprep!14] (5.27,-.26) rectangle (6.04,.23);
    \draw[powprep!70] (6.04,-.26)--(6.04,.23);
    \node[font=\footnotesize,anchor=north east] at (6.04,-.26) {$R_b\exp(J)$};
    \draw[powprep,line width=.9pt,->] (2.38,0)--(6.01,0);
    \draw[powochre,line width=.9pt,dash pattern=on 4pt off 2pt,->]
      (12.48,0)--(6.07,0);
    \fill[powprep] (2.38,0) circle[radius=1.35pt];
    \fill[powochre] (12.48,0) circle[radius=1.35pt];
    \node[font=\footnotesize,text=powprep,anchor=south] at (3.76,.04) {From the inner data};
    \node[font=\footnotesize,text=powochre,anchor=south] at (9.58,.04) {From the exterior data};
    \node[font=\footnotesize,anchor=north west] at (6.15,-.08)
      {$M^j_{\to}=M^j_{\leftarrow}$, \(\ j=\theta,z,\theta z,z\theta,p\)};
  \end{scope}
  \begin{scope}[yshift=-.51cm]
  \node[anchor=west,font=\small\bfseries] at (0,-1.39)
     {The transition in logarithmic radius};
  \draw[powink!22] (1.27,-2.45)--(7.33,-2.45);
  \draw[powink!36,dashed] (4.30,-3.85)--(4.30,-2.30);
  \draw[powink!36,dashed] (7.33,-3.85)--(7.33,-2.30);
  \draw[->] (1.27,-3.85)--(7.79,-3.85);
  \draw[->] (1.27,-3.85)--(1.27,-2.05);
  \node[anchor=east] at (1.16,-2.45) {$1$};
  \node[anchor=east] at (1.16,-3.85) {$0$};
  \node[anchor=south west,font=\footnotesize] at (1.32,-2.20) {$\sigma_{\rm norm}(d/\ell)$};
  \foreach \x/\tick in {1.27/R_{\rm mid},4.30/e^\ell R_{\rm mid},7.33/R_e}{
    \draw (\x,-3.85)--(\x,-3.95);
    \node[anchor=north,font=\footnotesize] at (\x,-4.00) {$\tick$};}
  \node[anchor=west] at (7.86,-3.85) {$R$};
  \draw[powteal,line width=1.15pt] (1.27,-2.45)--(4.30,-2.45);
  \draw[powochre,line width=1.15pt] plot coordinates {
       (4.300,-2.45000) (4.351,-2.45000) (4.401,-2.45000) (4.452,-2.45000) (4.502,-2.45000) 
      (4.553,-2.45000) (4.603,-2.45000) (4.654,-2.45000) (4.704,-2.45000) (4.755,-2.45000) 
      (4.805,-2.45000) (4.855,-2.45000) (4.906,-2.45000) (4.957,-2.45000) (5.007,-2.45000) 
      (5.058,-2.45000) (5.108,-2.45001) (5.159,-2.45004) (5.209,-2.45016) (5.260,-2.45056) 
      (5.310,-2.45164) (5.361,-2.45424) (5.411,-2.45990) (5.462,-2.47119) (5.512,-2.49216) 
      (5.563,-2.52866) (5.613,-2.58824) (5.663,-2.67889) (5.714,-2.80591) (5.764,-2.96725) 
      (5.815,-3.15000) (5.865,-3.33275) (5.916,-3.49409) (5.966,-3.62111) (6.017,-3.71176) 
      (6.067,-3.77134) (6.118,-3.80784) (6.168,-3.82881) (6.219,-3.84010) (6.269,-3.84576) 
      (6.320,-3.84836) (6.370,-3.84944) (6.421,-3.84984) (6.471,-3.84996) (6.522,-3.84999) 
      (6.572,-3.85000) (6.623,-3.85000) (6.673,-3.85000) (6.724,-3.85000) (6.774,-3.85000) 
      (6.825,-3.85000) (6.875,-3.85000) (6.926,-3.85000) (6.976,-3.85000) (7.027,-3.85000) 
      (7.077,-3.85000) (7.128,-3.85000) (7.178,-3.85000) (7.229,-3.85000) (7.279,-3.85000) 
      (7.330,-3.85000) };
  \draw[powink,line width=.85pt] (7.33,-3.85)--(7.59,-3.85);
  \node[anchor=west,font=\footnotesize] at (8.30,-2.47)
    {$U^\theta=U_e^\theta-\sigma_{\rm norm}(d/\ell)(U_e^\theta-V_p)$};
  \node[anchor=west,font=\footnotesize] at (8.30,-3.06)
    {$V_p=(1-\epsilon_0)h_\infty R^{-(1+\delta)/2}$};
  \node[anchor=west] at (8.30,-3.67) {$U^z/U^\theta=U_e^z/U_e^\theta$};
  \end{scope}
\end{tikzpicture}
\caption{Stress starts at \(R_{\rm in}\). At fixed \(Z\) the moments integrated from
the inner exit \(R_{ac}\) and from \(R_e\) match at \(R=R_b e^{\sup J}\).
The physical moment correction uses \(I_2\), the pulses \(I_3\).
The lower graph shows the cutoff in \eqref{pow:profiles}.}
\label{fig:short-connection}\label{fig:moment-transfer}
\end{figure}
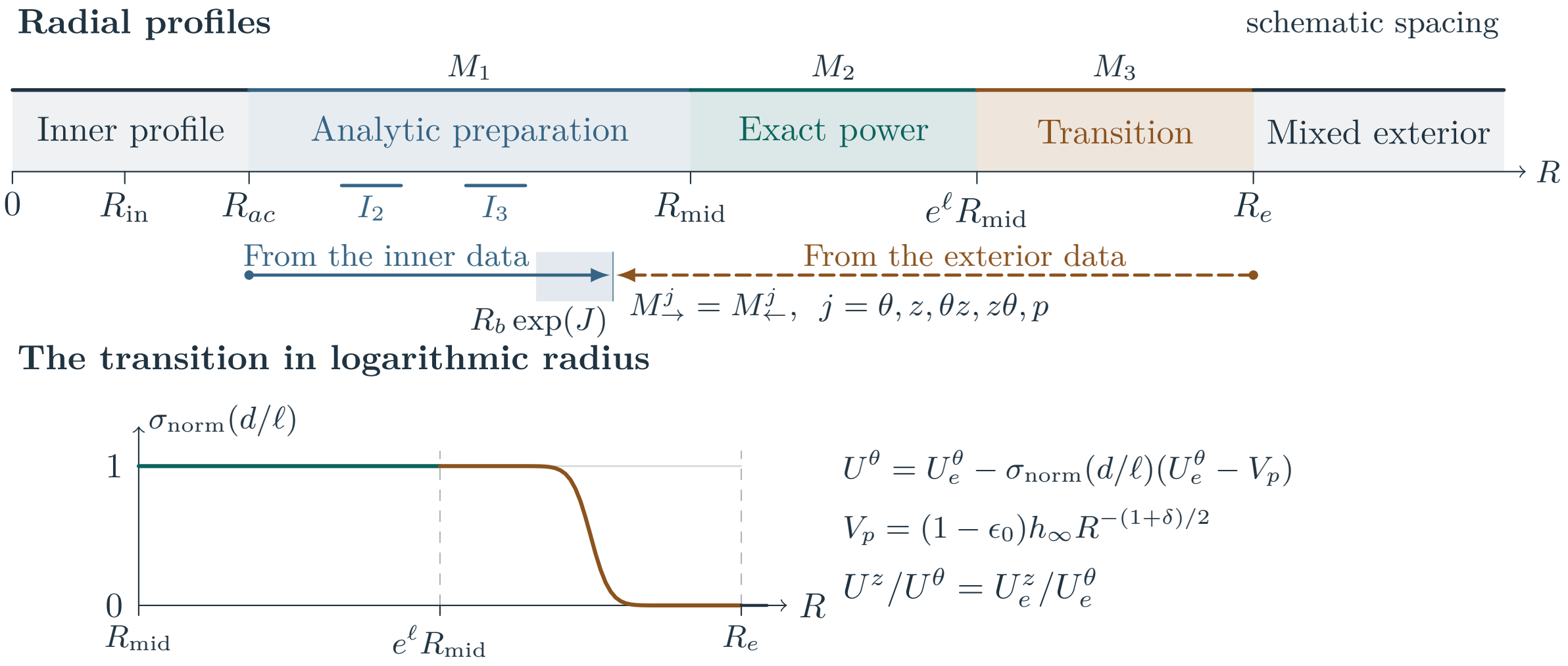

Once the exterior and the two intervals are fixed, the five terminal
values needed by the left profile are already determined. Integrating
backward from the exterior moments \eqref{ue:right-moments} gives
\begin{equation}\label{pow:inlet-moments}
 M^j(R_{\rm mid},Z)=M_e^j(R_e,Z)
       -\int_{R_{\rm mid}}^{R_e}m_j(\rho,Z)\,d\rho ,
 \qquad j=\theta,z,\theta z,z\theta,p ,
\end{equation}
where \(m_j\) are the densities in \eqref{pr:moments} evaluated on
\eqref{pow:profiles}. Any left profile ending with these five values and
the jets of the middle power therefore gives the prescribed exterior moments
at \(R_e\) when integrated from the axis (Figure~\ref{fig:moment-transfer}).

\paragraph{Strict admissibility on the right.}
We take \(R_{\rm mid}^{-1}\) and the supremum of \(U_e^\theta\) on the collar small as in
Proposition~\ref{nc:paper-Gamma}, fix \(\ell\) small, then \(\epsilon_0\)
small, and finally \(\delta>0\) sufficiently small. In this estimate write
\(A=U_e^\theta-U^\theta=\sigma_{\rm norm}(d/\ell)(U_e^\theta-V_p)\).
The interval from \(R\) to \(R_e\) has length at most \(C\ell R\), and
the cutoff is nonincreasing. By \eqref{pr:moments}, the five moment differences
from the exterior are therefore bounded, also after one \(Z\) derivative, by
\[
 C\ell A(R,Z)\bigl(R^{3/2},\,R,\,R^{3/2}U_e^\theta,\,RU_e^\theta,\,U_e^\theta\bigr)
\]
in the order of \eqref{pow:inlet-moments}. For the exterior moments,
\[
 b_0=1-\frac{(1-\delta)ZM_e^z+L\partial_ZM_e^z}{R}
     =1-2ZU_e^z-L_\delta/R>1/2
\]
by Proposition~\ref{ue:exact}, including the additive \(Z\) in \(M_e^z\).
Since the exterior stress vanishes, \eqref{pr:moment-inertia} and
\eqref{pr:shear} give
\begin{equation}\label{pow:stress-comparison}
 L_\delta\sqrt{2R}\,\mathcal T
  =\left[b_0RA+\frac{2L_\delta}{\ell}\sigma_{\rm norm}'(d/\ell)(U_e^\theta-V_p)\right]\,
      \frac{(U_e^\theta,U_e^z)}{U_e^\theta}
       +O\bigl((\ell+R_{\rm mid}^{-1})RA\bigr).
\end{equation}
The first term comes from the velocity change. The five bounds above
control the moment terms, including \(R[2(1+\delta)Z-L\partial_Z]\Delta P\),
by \(C\ell RA\). The bound on \(R\partial_R(U_e^\theta-V_p)\) controls
the remaining shear terms by \(CA\), with \(C\) absolute.

Both leading terms point along \((U_e^\theta,U_e^z)\), since \(U_e^\theta>V_p\) and \(\sigma_{\rm norm}'\ge0\).
The normalized exterior shear is \((-a_e,-(a_e-1)U_e^z/U_e^\theta)\), with
\(a_e=2+O(\delta)\). At \(a_e=2\),
\[
 \kappa=\frac{9+\epsilon_{\rm em}^2}{4},\qquad
 -\frac{(U_e^\theta,U_e^z)\cdot\mathcal S_e}{U_e^\theta F_e}
   =\frac{5+\epsilon_{\rm em}^2}{2},\qquad
 \frac{|(U_e^\theta,U_e^z)\cdot\mathcal S_e^\perp|}{U_e^\theta F_e}
   =\frac{U_e^z}{U_e^\theta}.
\]
Thus the exterior velocity direction lies strictly inside the cone
\eqref{pr:cone}, uniformly for \(1/\sqrt2\le U_e^z/U_e^\theta\le1\).
With \(\ell\) fixed, taking \(\epsilon_0\) and then \(\delta\) small keeps
the new normalized shear close to the exterior shear. The relative error
above is small compared with \(b_0RA\), so the stress is nonzero and
strictly admissible on \([R_{\rm mid},R_e)\).

Near \(R_e\), the cutoff \(\sigma_{\rm norm}(d/\ell)\) divided by
\(e^{-\ell^2/d^2}\) is smooth and positive up to \(d=0\).
The tangential velocity differences and the moment density differences
are the cutoff times smooth functions. By \eqref{fin:gaussian-smooth-quotient},
backward integration of the moment densities gains a factor \(d^3\).
Equation \eqref{pr:radial} gives the same cutoff factor for \(U^r-U_e^r\).
Only differentiation of the cutoff in the shear produces the factor \(d^{-3}\).
Hence
\begin{equation}\label{pow:edge-data}
 \begin{gathered}
 (U^\theta-U_e^\theta,U^z-U_e^z,U^r-U_e^r)=\sigma_{\rm norm}(d/\ell)\,\mathbf B(d,Z),
 \qquad M^j-M_e^j=\sigma_{\rm norm}(d/\ell)d^3 B_j(d,Z),\\
 \mathcal T=\sigma_{\rm norm}(d/\ell)d^{-3}\mathbf A(d,Z),\qquad
 \mathbf A(0,Z)=2\ell^2\sqrt{2/R_e}\,(U_e^\theta(R_e,Z)-V_p(R_e))
                 (1,U_e^z/U_e^\theta).
 \end{gathered}
\end{equation}
The coefficients are smooth at \(R_e\), and \(\mathbf A^\theta(0,Z)>0\).
Hence we have the Gaussian derivative bounds and \(\mathcal T/|\mathcal T|\)
extends smoothly to \(R_e\).
\subsection{From the actual inner exit to the mixed power}
\label{conn:section}\label{cl:section}\label{cp:section}
The inner exit fixes both velocities and all five moments at \(R_{ac}\),
with axis pressure \(P_0(Z)=-3P_*^2/(1+Z^2)^2\).
To reach \eqref{pow:inlet-moments} at \(R_{\rm mid}\), we first match pressure and
energy, then vary the shear to make the stress strictly admissible and use five
independent variations to correct the remaining moment errors.
At fixed $Z$, write $\Delta$ for the increment from $R_-$ to $R_+$.
Since $F>0$ decreases, $\partial_RP=F^2$ and Cauchy--Schwarz give
\[
\begin{gathered}
 (R_+-R_-)F(R_+)^2<\Delta P<(R_+-R_-)F(R_-)^2,\qquad
 (\Delta M^\theta)^2\le\tfrac43(R_+^3-R_-^3)\Delta P,\\
 \Delta M^{z\theta}+R_+\Delta P\ge (\Delta M^z)^2/(R_+-R_-) .
\end{gathered}
\]
Indeed, $\Delta M^\theta=2\int_{R_-}^{R_+}RF\,dR$ and
$\Delta M^{z\theta}+R_+\Delta P\ge\int_{R_-}^{R_+}(U^z)^2\,dR$.
\subsubsection{From the actual inner exit to \texorpdfstring{\(R=5\)}{R = 5}}
\label{en:section}

We start at \(R_{ac}\) with the velocities and five moments of the inner exit. By varying the shear on short intervals, we prescribe the outlet derivatives at \(R=5\) while keeping the velocities and moments close to their comparison values.

\begin{lemma}[A small change from the inner exit]\label{en:passage}
Let $0<R_1<R_2$, and let a smooth inner profile be given on a left collar
of $R=R_1$, with its five moments and $F(R_1,Z)>0$. For comparison, take
$F^0(R,Z)=F(R_1,Z)$ and $U^{z,0}(R,Z)=U^z(R_1,Z)$ on $[R_1,R_2]$.
Integrate their five moment densities from the given values at $R_1$.
The superscript $0$ on $\mathcal I^{\theta,0},\mathcal I^{z,0}$ means that
we use these comparison velocities and moments in \eqref{pr:moment-inertia}.
Suppose the inner shear is strictly admissible at $R_1$, and
\begin{equation}\label{en:comparison}
 \frac{\mathcal I^{\theta,0}}{F^0}>0,\qquad
 \frac{|\mathcal I^0|^2}{F^0\mathcal I^{\theta,0}}>2
 \quad(R_1\le R\le R_2),\qquad
 \frac{\mathcal I^{\theta,0}}{F^0}>3\quad\hbox{near }R_2.
\end{equation}
Fix $a_*\in(0,2)$. Then there is a smooth pair on $[R_1,R_2]$ satisfying
relaxed admissibility, with the same derivatives at $R_1$ and with
$a=a_*$, $b=0$ near $R_2$. It is strictly admissible near $R_1$.
Its velocities and five moments can be made arbitrarily close to
the comparison values in any fixed $C_Z^m$ norm, uniformly in $R$.
In particular, $\mathcal I^\theta/F>3$ at $R_2$.
If the data depend smoothly on a parameter in a compact set and the
assumed inequalities have uniform margins, the pair can also be
chosen to depend smoothly on that parameter.
\end{lemma}

\begin{proof}
For a fixed state $F,\mathcal I$, put $t=-b/a$.  Then
$(a,b)=\kappa(1,-t)/(1+t^2)$.  When $\kappa>2$, admissibility is
\begin{equation}\label{en:geometry}
 \begin{aligned}
 \mathcal I^\theta+t\mathcal I^z&>\kappa F,\\
 (\kappa-2)(\mathcal I^z-t\mathcal I^\theta)^2
   &<2(\mathcal I^\theta+t\mathcal I^z-\kappa F)^2.
 \end{aligned}
\end{equation}
For $0<\kappa\le2$, it is
$\mathcal I^\theta+t\mathcal I^z>2F$.
Decreasing $\kappa$ at fixed $t$ preserves these conditions.
After decreasing $\kappa$ to $1$, one can vary $t$ within the interval
$\{t:\mathcal I^\theta+t\mathcal I^z>2F\}$.
By \eqref{en:comparison}, this interval contains
$t=\mathcal I^{z,0}/\mathcal I^{\theta,0}$ throughout $[R_1,R_2]$,
and also contains $t=0$ near $R_2$.

Choose short layers at the two endpoints. In this proof, $\epsilon_1>0$ controls the shear between the layers.
Away from the layers prescribe
\begin{equation}\label{en:small-shear}
 a=\epsilon_1\frac{\mathcal I^{\theta,0}}{F^0},\qquad
 b=-\epsilon_1\frac{\mathcal I^{z,0}}{F^0},
\end{equation}
with $\epsilon_1$ small enough that $\kappa<1$.
In the left layer, extend the inner profile smoothly, decrease $\kappa$,
vary $t$ as above, and reduce the shear to \eqref{en:small-shear}.

In the right layer, vary $t$ to zero and
then set $\kappa=a_*$. Use smooth steps constant near their endpoints,
so all derivatives match. These choices are smooth in $Z$.

Recover the velocities by
\begin{equation}\label{en:reconstruction}
 \begin{aligned}
 F(R,Z)&=F(R_1,Z)\exp\left[-\frac12\int_{R_1}^R
                                a(\rho,Z)\frac{d\rho}{\rho}\right],\\
 U^z(R,Z)&=U^z(R_1,Z)+\frac12\int_{R_1}^R
                 b(\rho,Z)\sqrt{2\rho}F(\rho,Z)\frac{d\rho}{\rho}.
 \end{aligned}
\end{equation}
As the layer widths and $\epsilon_1$ tend to zero, the shear integrals and
their fixed $Z$ derivatives tend to zero. Hence the velocities,
moments, and inertia approach their comparison values; the inertia
estimate uses one extra $Z$ derivative.

On the part of the path with \(\kappa>2\), the margins in
\eqref{en:geometry} persist down to \(\kappa=2\). On the other part,
the relaxed inequality divided by \(a\) is \(\mathcal I^\theta+
t\mathcal I^z>2F\), with a margin independent of \(\epsilon_1\).
Choose the layer widths and then \(\epsilon_1\) sufficiently small to
preserve these margins and \(\mathcal I^\theta(R_2)/F(R_2)>3\).
For a compact parameter family, choose the layers, $\epsilon_1$, and smooth
steps once, using the uniform bounds and margins. Then
\eqref{en:reconstruction} gives smooth parameter dependence.
\end{proof}

\paragraph{Estimates at the inner exit.}
Apply the lemma with $R_1=R_{ac}$, $R_2=5$, and $a_*=4/5$, with $j_0$ small and $\Lambda$
chosen for at least four $Z$ derivatives. Theorem~\ref{ic:inner-connection}
gives the scaled bounds \eqref{ic:prefix-bound}--\eqref{ic:prefix-moment-bounds}
and the exit data \eqref{ic:new-exit-linear}--\eqref{ic:new-exit-sectors}.
It remains to check \eqref{en:comparison} on $[R_{ac},5]$.
For the constant comparison velocities, differentiating \eqref{pr:moment-inertia} gives
\begin{equation}\label{en:frozen-angular}
 \begin{aligned}
 \partial_R\left(R\frac{\mathcal I^{\theta,0}}{F^0}\right)
 ={}&\frac R{L_\delta}\left[-1-\frac\delta2+ZU^z(R_{ac})
      +L\partial_ZU^z(R_{ac})
      -\left(\frac{1-\delta}{2}Z+LU^z(R_{ac})\right)
                             \partial_Z\log F(R_{ac})\right]\\
 &+\frac{[(1-\delta)Z+L\partial_Z]
                 [M^z(R_{ac})-R_{ac}U^z(R_{ac})]}{L_\delta}.
 \end{aligned}
\end{equation}
By \eqref{ic:new-exit-linear}, $ZU^z(R_{ac})+L\partial_ZU^z(R_{ac})=4+j_0Z+O(\Lambda^{-1})$.
We use \eqref{ic:new-exit-transport} for the term involving $\partial_Z\log F(R_{ac})$.
Since $R\ge R_{ac}>4/\Lambda$, the $O(\Lambda^{-2})$ moment term
is $O(R/\Lambda)$. Thus, for large $\Lambda$,
\begin{equation}\label{en:frozen-angular-lower}
 \partial_R\left(R\frac{\mathcal I^{\theta,0}}{F^0}\right)
       \ge2R+\frac34\Lambda\chi_c(Z)R,
 \qquad \frac{\mathcal I^{\theta,0}(R_{ac})}{F^0(R_{ac})}\ge R_{ac}/256.
\end{equation}
Here $\chi_c$ is defined in \eqref{ic:axis-data}, and the initial
bound follows from \eqref{ic:new-exit-sectors}.
On $\{\chi_c\ge99/100\}$ the initial angular ratio exceeds $9/4$.
At a first crossing of $9/4$, its $\log R$ derivative
would be at least $2R+\frac34\Lambda\chi_cR-9/4>0$, since
$\Lambda R\ge\Lambda R_{ac}>4$.  Thus it remains greater than $9/4$.
For every $Z$, integration gives $\mathcal I^{\theta,0}/F^0\ge R-255R_{ac}/256>0$.
This lower bound exceeds $3$ near $R=5$.

On $\{\chi_c\le99/100\}$, the axial identity \eqref{ic:transport}
for the same comparison profile gives
\begin{equation}\label{en:frozen-axial}
 \begin{aligned}
 \partial_R\bigl(L_\delta\sqrt{2R}\,\mathcal I^{z,0}\bigr)
 ={}&-\frac{1+\delta}{2}(1-2ZU^z(R_{ac}))U^z(R_{ac})
      -\left(\frac{1-\delta}{2}Z+LU^z(R_{ac})\right)\partial_ZU^z(R_{ac})\\
 &+[2(1+\delta)Z-L\partial_Z]
          [P_0+M^p(R_{ac})-R_{ac}F(R_{ac})^2]\\
 &+2R[(2+\delta)ZF(R_{ac})^2-LF(R_{ac})\partial_ZF(R_{ac})].
 \end{aligned}
\end{equation}
The terms with $P_0$ form $g[U^z(R_{ac})]\ge2j_0/5$ by
\eqref{ic:transport} and \eqref{ic:new-exit-sectors}.
At $R=R_{ac}$ the comparison and exit profiles have the same velocities and
moments. Thus \eqref{pr:moment-inertia} and \eqref{ic:new-exit-sectors} give
$L_\delta\sqrt{2R_{ac}}\,\mathcal I^{z,0}(R_{ac},Z)\ge R_{ac} j_0/5$ on this set.
The remaining terms are $O_\Lambda(C_*^{-2})$ on $[R_{ac},5]$.
For large $C_*$ the right side of \eqref{en:frozen-axial} is at least
$j_0/5$. Integrating from $R_{ac}$ gives
$L_\delta\sqrt{2R}\,\mathcal I^{z,0}(R,Z)\ge R j_0/5$ throughout $[R_{ac},5]$.
Since $F(R_{ac})=O_\Lambda(C_*^{-1})$ by \eqref{ic:prefix-bound}, increasing
$C_*$ gives $\mathcal I^{z,0}/F^0\ge j_0\sqrt{R/2}/(5L_\delta F(R_{ac}))>9/8$.
Since $\mathcal I^{\theta,0}>0$, one has
\[
 \frac{|\mathcal I^0|^2}{F^0\mathcal I^{\theta,0}}
       \ge\frac{2|\mathcal I^{z,0}|}{F^0}\ge\frac94.
\]
Thus Lemma~\ref{en:passage} gives $a=4/5$, $b=0$ and
$\mathcal I^\theta/F>3$ at $R=5$. Choose its changes small enough to
preserve the scaled bounds there. For the comparison,
$M^{z,0}(5)=5U^z(R_{ac})+M^z(R_{ac})-R_{ac}U^z(R_{ac})$. Hence \eqref{ic:new-exit-linear},
with $B_4/\Lambda$ below $j_0$ and the odd-part tolerance, makes the odd
parts of $U^z(5),M^z(5)/5$ close to $4Z$ and the even parts $O(j_0)$ in $C_Z^4$.
For large $\Lambda$ and small changes, the transport bound is below
$1/20$ at $\delta=0$. Uniform choices on the compact $\delta$ interval
give smooth dependence down to zero; hence the bound remains below
$1/10$ for sufficiently small $\delta>0$. These are the hypotheses of
Lemma~\ref{np:preparation}.
\begin{figure}[!htbp]
\centering
\includegraphics[width=\textwidth]{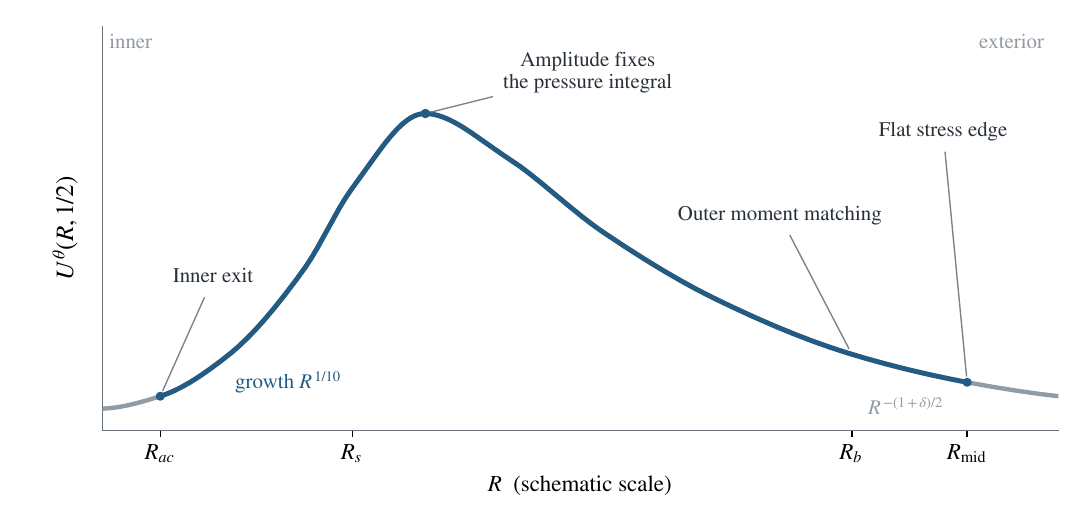}
\par
\includegraphics[width=\textwidth]{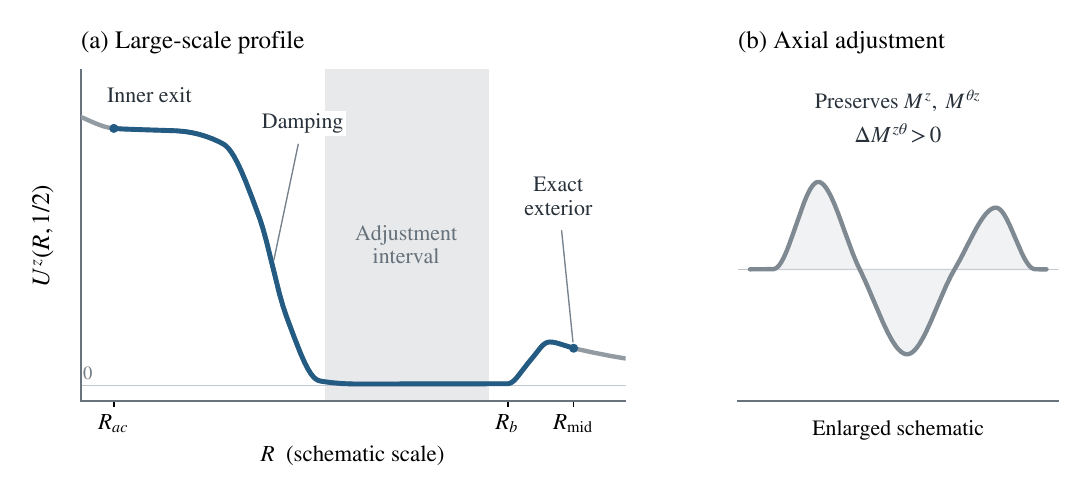}
\caption{Auxiliary preparation at \(Z=1/2\) before insertion of the power law.
At this stage the exterior begins at \(R_{\rm mid}\). Above, the angular
amplitude fixes the pressure integral. Below, (a) marks the adjustment interval
on the axial profile, and (b) enlarges the axial correction. Its two linear
moments vanish. The integral of its square gives the prescribed increment
in \(M^{z\theta}\).}
\label{fig:one-theta}\label{fig:one-axial}
\end{figure}
\subsubsection{Preparing the pressure and energy up to \texorpdfstring{\(R_b\)}{Rb}}\label{np:section}

Here \(R_1\) denotes the entrance radius for the preparation. We take \(R_1=5\),
the right endpoint of Lemma~\ref{en:passage}. We first match the pressure
with the angular profile, then correct the energy at fixed axial and mixed linear moments.

\begin{lemma}[Preparation of the actual incoming data]\label{np:preparation}
Let $R_1\ge5$. Fix the central family of
Proposition~\ref{nc:paper-Gamma}. \emph{Incoming data.} Suppose the incoming smooth pair has
$a=4/5$, $b=0$ near $R_1$, with its five moments integrated from the
core. A subscript $1$ on a field or moment denotes its value at $R=R_1$,
with $Z$ left free. Write $U_1^z=U^z_{1o}+U^z_{1e}$, where the two terms are respectively odd and even in $Z$. Assume
\begin{gather*}
 \mathcal I_1^\theta/F_1>3,\qquad
 \left(Z/2+(1-Z^2)U_1^z\right)\partial_Z\log U_1^\theta\le1/10,\\
 \|U^z_{1o}-4Z\|_{C_Z^4}+\|(M_1^z)_o/R_1-4Z\|_{C_Z^4}\le\epsilon_{\rm odd},
 \qquad \|U^z_{1e}\|_{C_Z^4}+\|(M_1^z)_e/R_1\|_{C_Z^4}\le Cj_0.
\end{gather*}
Let $j_0,\epsilon_{\rm odd}>0$ be sufficiently small. Suppose
$C_*U_1^\theta$ has positive upper and lower bounds, and the
$C_Z^4$ norms of $C_*U_1^\theta$ and
\[
 C_*M_1^\theta,\quad C_*M_1^{\theta z},\quad M_1^z,\quad
 M_1^{z\theta},\quad C_*^2M_1^p
\]
are bounded independently of large $C_*$, with the other core
parameters fixed. Assume smooth extension to $\delta=0$. Set
$P_0=-3P_*^2/(1+Z^2)^2$ and $P(R_1)=P_0+M_1^p$.

\emph{Construction.} For every sufficiently small $\tau>0$, one can choose $\mu,P_*$
with $\mu P_*^4$ fixed and small, then $C_*$, and finally $\delta>0$,
as in the \hyperref[np:closure]{parameter choices below}, to obtain a smooth pair on $[R_1,R_b]$
matching all incoming derivatives, with $U^\theta>0$, $a>0$, and
relaxed admissibility for the forward moments. On the reserve
immediately preceding the terminal transition, $U^z=0$, the angular
velocity is independent of $Z$, and $2<a<4$ with strict cone margins.
The short terminal transition remains strictly admissible and matches
both central inlet velocities to all orders at $R_b$.

\emph{Terminal conditions.} Appending the central family and exact exterior, with some
$\epsilon_{\mathrm{em}}>0$, gives
\[
 M_1^p+\frac12\int_{R_1}^{\infty}(U^\theta)^2\frac{dR}{R}=-P_0,\qquad
 M_1^{z\theta}+
 \int_{R_1}^{\infty}\left[(U^z)^2-(U^\theta)^2/2\right]dR=0.
\]
Let $U_b^\theta=U^\theta(R_b,1)$. At $R_b$, the forward values of
$M^\theta,M^z,M^{\theta z}$ differ from the central backward values by
less than $\tau$ in $C_Z^3$, in the first three units of \eqref{nc:paper-moment-units}.
\end{lemma}

\begin{proof}
First take $\delta=0$.

\paragraph*{Amplification and prescribed pressure.}\label{np:amplification}

Use $y=\log(R/R_1)$.
Keep $U^z$ independent of $R$ and choose $0<a<9/5$.
Then $m=M^z/R$ solves $\partial_y m=U^z-m$, so
$1-Zm-(1-Z^2)\partial_Zm\le-2$. The angular equation is
\[
 \begin{aligned}
 \partial_y(\mathcal I^\theta/F)+(1-a/2)\mathcal I^\theta/F
 &=-R\bigl[1-Zm-(1-Z^2)\partial_Zm\bigr](1-a/2)\\
 &\quad-R\left(Z/2+(1-Z^2)U^z\right)\partial_Z\log U^\theta .
 \end{aligned}
\]
At $\mathcal I^\theta/F=3$ we have
$\partial_y(\mathcal I^\theta/F)\ge(2R-3)(1-a/2)-R/10\ge(R-3)/10>0$.
Thus $\mathcal I^\theta/F$ stays above $3$, which gives relaxed admissibility since $b=0$ and $a<2$.
Interpolate the logarithmic angular shape slowly to a positive
$q(Z)$ near $2/(1+Z^2)$, with $q(1)=1$, keeping $a$ in a compact
subinterval of $(0,9/5)$. The transport expression is a convex
combination of its endpoint values. At $q=2/(1+Z^2)$, it is
nonpositive for $j_0=0$ and $O(j_0^2)$ for small $j_0$.
Small shape perturbations preserve the transport bound.

Use the nondecreasing step $\sigma_{\rm norm}$ from \eqref{pr:cutoffs},
which is zero for $s\le0$ and one for $s\ge1$. Set
\[
 U^\theta=U_1^\theta(Z)e^{y/10}
 \exp\left\{\sigma_{\rm norm}(y/L_g)
    \log\frac{q(Z)U_1^\theta(1)}{U_1^\theta(Z)}\right\}.
\]
Choose the fixed $L_g$ so that
\[
 \frac{2\|\sigma_{\rm norm}'\|_\infty}{L_g}
 \sup_{q\text{ in the chosen neighborhood}}
 \left\|\log\frac{qU_1^\theta(1)}{U_1^\theta}\right\|_\infty
 <1/10.
\]
Then $7/10<a<9/10$. End the growth at $y_{\mathrm{g},+}>L_g+2$.
After $L_g$, the angular shape is $q(Z)$ and $c=(a-2)/2=-3/5$.

From $R_1$ to the end of the transition, set $U^z(R,Z)=U^z(R_1,Z)$.
Put $d_0=\partial_Z U^z_{1o}(0)\in[3,5]$.
Throughout the transition and damping, let $d_1=\partial_ZU^z(R,0)$ and
$d_2=\partial_Z(M^z/R)(R,0)$. The comparison pair uses the actual
initial $d_2(y_{\mathrm{g},+})$ and has $U^z=d_1 Z$, $U^\theta=2A/(1+Z^2)$,
$\partial_y\log A=-1/2-c$. Here $A=U^\theta_{\rm model}(R,1)$.
For this pair the normalized angular moment is independent of $Z$, and the
normalized mixed moment is linear in $Z$. Thus their values at $Z=1$ give
\[
 \frac{(M^\theta_{\rm model},M^{\theta z}_{\rm model})(R,Z)}
      {R\sqrt{2R}\,U^\theta_{\rm model}(R,Z)}
 =\frac{(M^\theta_{\rm model}(R,1),\,Z M^{\theta z}_{\rm model}(R,1))}
       {R\sqrt{2R}\,A}.
\]
We have $M^z_{\rm model}=R d_2 Z$. By \eqref{pr:moments},
\[
 \partial_y d_2=d_1-d_2,\qquad
 (\partial_y+1-c)
 \frac{(M^\theta_{\rm model},M^{\theta z}_{\rm model})(R,1)}{R\sqrt{2R}\,A}
 =(1,d_1).
\]
The first equation also holds for the actual moment. At $Z=0$ the normalized angular inertia is
$\widetilde I_0^\theta=-1+d_2+(M^\theta_{\rm model}-M^{\theta z}_{\rm model})(R,1)/(R\sqrt{2R}\,A)$.
With primes denoting $\partial_y$, \eqref{pr:moment-inertia} gives
\begin{equation}\label{new:angular-inertia}
\begin{aligned}
 (\widetilde I_0^\theta)'+(1-c)\widetilde I_0^\theta&=c(1-d_2),\\
 \widetilde I_{\rm model}^\theta:=\frac{\mathcal I^\theta_{\mathrm{model}}}{RF}
 &=\widetilde I_0^\theta+\frac{2Z^2}{1+Z^2}
 \frac{M^\theta_{\rm model}(R,1)/2+(1-Z^2)M^{\theta z}_{\rm model}(R,1)}{R\sqrt{2R}\,A}.
\end{aligned}
\end{equation}
On the transition $d_1=d_0$. Thus
$d_2(y_{\mathrm{g},+})-d_0=(\partial_Z M_1^z(0)/R_1-d_0)e^{-y_{\mathrm{g},+}}$, of absolute value at most $1/4$
for large $y_{\mathrm{g},+}$. Set the two normalized moments at $Z=1$ to their limiting
growth values $(5/8,5d_0/8)$ at $y=y_{\mathrm{g},+}$. Then
$\widetilde I_0^\theta(y_{\mathrm{g},+})=\tfrac38(d_0-1)+d_2(y_{\mathrm{g},+})-d_0\ge\tfrac12$.

Fix the logarithmic length $L_{\rm tr}>0$ of the transition before $\mu$ and the shooting parameters.
It runs from $R_1e^{y_{\mathrm{g},+}}$ to $R_s=R_1e^{y_{\mathrm{g},+}+L_{\rm tr}}$.
On $[y_{\mathrm{g},+},y_{\mathrm{g},+}+L_{\rm tr}]$, interpolate $c$ smoothly
from $-3/5$ to $\mu(1-d_2)$.
Take the cutoff flat at both endpoints.
For small $\mu$ and large $y_{\mathrm{g},+}$, one has $d_2>1$ and $-3/5\le c\le0$.
Here $c(1-d_2)\ge0$ and both normalized moments at $Z=1$ stay positive. Hence
\eqref{new:angular-inertia} gives throughout this interval
\[
 \widetilde I_{\rm model}^\theta\ge \widetilde I_0^\theta\ge \widetilde I_0^\theta(y_{\mathrm{g},+})
       e^{-\int_{y_{\mathrm{g},+}}^y(1-c(s))\,ds}\ge\tfrac12e^{-8L_{\rm tr}/5}.
\]

For the actual pair, $M^\theta/(R\sqrt{2R}U^\theta)$ and
$M^{\theta z}/(R\sqrt{2R}U^\theta)$ satisfy the same equations at each $Z$,
with right-hand sides $1$ and $U_1^z(Z)$, respectively.
Also $\partial_y(M^z/R)=U_1^z-M^z/R$. The growth errors decay exponentially.
Since $\|U_1^z-d_0Z\|_{C_Z^2}\le C(j_0+\epsilon_{\rm odd})$, variation of constants in
these equations and \eqref{pr:moment-inertia} give
\[
 \sup_{y_{\mathrm{g},+}\le y\le y_{\mathrm{g},+}+L_{\rm tr}}
 \left\|\widetilde I^\theta
                  -\widetilde I_{\rm model}^\theta\right\|_{C_Z^0}
 \le C\bigl(j_0+\epsilon_{\rm odd}+\xi+e^{-y_{\mathrm{g},+}}\bigr).
\]
Here $\xi$ bounds the normalized prefix and pressure-shape errors, and
$C$ depends only on the fixed transition and $C_Z^2$ bounds for
$U_1^z,q,q^{-1}$, uniformly for $d_0\in[3,5]$. Require both $C(j_0+\epsilon_{\rm odd})$
and $C(\xi+e^{-y_{\mathrm{g},+}})$ to be at most $e^{-8L_{\rm tr}/5}/8$.
Then $\widetilde I^\theta\ge
e^{-8L_{\rm tr}/5}/4$. For $R_1e^{y_{\mathrm{g},+}}\ge16e^{8L_{\rm tr}/5}$
we have $R\widetilde I^\theta\ge4$, giving \eqref{pr:relaxed} since $b=0$ and $0<a\le2$.

Let $A_s=A(R_s)$. From here through the shape change, measure
$y=\log(R/R_s)$ from $R_s$. End the damping interval at $y=\mu^{-1/2}$.
There $e^{-y/4}=O_N(\mu^N)$ for every fixed $N>0$ and $\mu y=\sqrt\mu$.
Choose both damping factors equal to one near $y=0$.
The flat transition then matches every radial derivative.
Since $d_2-d_0$ and its $y_{\mathrm{g},+}$ derivative are $O(e^{-y_{\mathrm{g},+}})$, integration
of $A'/A=-1/2-c$ on the fixed transition gives
\[
 A_s\asymp U_1^\theta(1)e^{y_{\mathrm{g},+}/10},\qquad
 \partial_{y_{\mathrm{g},+}}\log A_s=1/10+O(e^{-y_{\mathrm{g},+}}).
\]
The incoming axial-moment error is $O(e^{-y_{\mathrm{g},+}})$.
Define $\mathcal E_s=R_sA_s^2$.

Remove the $Z$ dependence after $y=\mu^{-1/2}$ as in the \hyperref[np:separated]{separated interval below}. Define
\[
 \mathcal P(Z;y_{\mathrm{g},+},q)=M_1^p(Z)
    +\frac12\int_{R_1}^{\infty}
             (U^\theta(R,Z;y_{\mathrm{g},+},q))^2\frac{dR}{R}.
\]
At $Z=1$ extend the pure growth law backwards to $y=-\infty$ and denote
the resulting pole amplitude through damping by $A(y)$, with $A(0)=A_s$.
The axial profile fixes $d_2$ and $A$ independently of $q$. Thus
\[
 \mathcal P(Z)=C_p(y_{\mathrm{g},+})A_s^2q(Z)^2+\mathcal R_p(Z),\qquad
 C_p(y_{\mathrm{g},+})=\frac12\int_{-\infty}^{\mu^{-1/2}}\!\bigl(A(y)/A_s\bigr)^2\,dy .
\]
The growth law, fixed transition and damping bound below give $0<c<C_p<C$.
At fixed $R/R_s$, $d_2'=d_1-d_2$ gives
$\partial_{y_{\mathrm{g},+}}d_2(y)=O(e^{-y_{\mathrm{g},+}}e^{-(y+L_{\rm tr})})$ from the start
of the transition. The cutoff defining $c$ is fixed there, and $c=\mu(1-d_2)$
during damping. Hence $\int_{-\infty}^{\mu^{-1/2}}|\partial_{y_{\mathrm{g},+}}c|\,dy
=O(\mu e^{-y_{\mathrm{g},+}})$. Integrating $\partial_y\log(A/A_s)=-1/2-c$
from $0$ gives $|\partial_{y_{\mathrm{g},+}}C_p|\le C\mu e^{-y_{\mathrm{g},+}}C_p$.
At fixed $y_{\mathrm{g},+}$ and $Z$, the pressure depends on the shape only through $q(Z)$.
We write $\partial_q$ for differentiation in this scalar value.
In a fixed $C_Z^4$ neighborhood of $2/(1+Z^2)$ we have
\[
 \|\mathcal R_p\|+
 \|\partial_q\mathcal R_p\|+
 \|\partial_{y_{\mathrm{g},+}}\mathcal R_p\|
 \le C_{\Lambda}C_*^{-2}
       +C A_s^2e^{-c/\sqrt\mu}.
\]
The fixed initial shape interval and the extrapolated growth prefix cost
$O_\Lambda(C_*^{-2})$. After $y=\mu^{-1/2}$ the angular amplitude is
$O(A_s e^{-c/\sqrt\mu})$, with the differentiated tail integrals estimated
\hyperref[np:closure]{below}. The shape factors $q^{\sigma_{\rm norm}},q^{1-\sigma_{\rm norm}}$ have
bounded derivatives on this positive range. Since $q(1)=1$ and
$\partial_{y_{\mathrm{g},+}}\log A_s=1/10+O(e^{-y_{\mathrm{g},+}})$, we obtain
\[
 \partial_{y_{\mathrm{g},+}}\mathcal P(1)\ge cA_s^2,
 \qquad \partial_q\mathcal P(Z)\ge cA_s^2.
\]
At $Z=1$ the target is $3P_*^2/4$. Choose fixed constants
$0<\alpha_-<\alpha_+$ with $C\alpha_-^2<3/8$ and
$c\alpha_+^2>3/2$, using the bounds on $C_p$ above.
Since $A_s\asymp U_1^\theta(1)e^{y_{\mathrm{g},+}/10}$ and $U_1^\theta(1)\asymp_\Lambda C_*^{-1}$,
the trial values $A_s=\alpha_\pm P_*$ are reached at radii tending to infinity
with $C_*$. The remainder bound then brackets $3P_*^2/4$ for small $\mu$,
and monotonicity gives a unique $y_{\mathrm{g},+}$ on this branch.

Fix this $y_{\mathrm{g},+}$. At the trial shape $2/(1+Z^2)$, subtracting $4/(1+Z^2)^2$ times the pole equation cancels the leading pressure term and gives
\[
 \|\mathcal P+P_0\|_{C_Z^4}
 \le C A_s^2\varepsilon_p,\qquad
 \varepsilon_p=C_\Lambda C_*^{-2}P_*^{-2}+e^{-c/\sqrt\mu}.
\]
At each $Z$ the derivative bound brackets the unique root between
$2/(1+Z^2)\pm C\varepsilon_p$, with $q(1)=1$. Divide by $A_s^2$ and
differentiate the scalar equation through order four. The bounds on the
mixed derivatives of the pressure map give
\[
 A_s\asymp P_*,\qquad
 \left\|q-\frac{2}{1+Z^2}\right\|_{C_Z^4}
 \le C_{\Lambda}C_*^{-2}P_*^{-2}+C e^{-c/\sqrt\mu}.
\]

\paragraph*{Scalar equations during damping.}\label{np:damping}

On $y=\log(R/R_s)\ge0$, prescribe the equatorial slope $d_1(y)$ of the comparison
pair and set $c=\mu(1-d_2)$. Equation \eqref{new:angular-inertia} becomes
\begin{equation}\label{new:positive-Q}
 (\widetilde I_0^\theta)'+(1-c)\widetilde I_0^\theta=\mu(1-d_2)^2\ge0.
\end{equation}

Choose a fixed smooth decreasing $d_1$, initially equal to $d_0$, with
Gaussian decay at infinity and
\begin{equation}\label{new:slow-hazard}
 \sup_{y\ge0}(1+y)|d_1'(y)|\le\eta.
\end{equation}
Choose a fixed smooth nondecreasing \(\vartheta(y)\) equal to zero for
\(y\le1\) and to \(y\) for \(y\ge2\), and set
\(d_1(y)=d_0(1+\vartheta(y))^{-\varepsilon}
\exp[-\vartheta(y)^2/T]\). Thus \(d_1=d_0\) on the entrance collar.
The two contributions to \((1+y)|d_1'|\) are
\(O(\varepsilon)\) and \(O(T^{-\varepsilon/2})\) for
\(0<\varepsilon<1\) and \(T\ge1\). The fixed transition has bounded
\(\vartheta'\).
Choose $\varepsilon$ small, then $T$ large, before choosing $\mu$.

These choices are uniform for $d_0\in[3,5]$. They give
$0\le d_2\le6$, $d_2\le C_{d_1} e^{-y}$ and $-5\mu\le c\le\mu$.
At $Z=1$ the normalized mixed moment stays nonnegative, and the normalized angular
moment stays between two positive constants. Thus \eqref{new:positive-Q}
and the angular moment term in \eqref{new:angular-inertia} give
\begin{equation}\label{new:Q-bound}
 c_{d_1}(e^{-y}+\mu+Z^2)\le \widetilde I_{\rm model}^\theta\le C_{d_1}.
\end{equation}
The homogeneous contribution has a positive lower bound until
$d_2\le1/2$; thereafter the forcing is at least $\mu/4$.
Since $\int_0^\infty d_2<\infty$, for small $\mu$ one has
\[
 A_s e^{-(1/2+\mu)y}\le A(y)
 \le2A_s e^{-(1/2+\mu)y}.
\]

\paragraph*{Inertia and the damping cone.}\label{np:damping-cone}

For the comparison pair, $U^z(R,0)=0$ and $U^\theta(R,0)=2A$. Thus \eqref{pr:moments}
gives $\partial_R M^{z\theta}_{\rm model}(R,0)=-2A^2$. For the damping estimates,
we denote the normalized cumulative angular energy at $Z=0$ by $J_A$ and write
\[
 J_A=-\frac{M^{z\theta}_{\rm model}(R,0)}{2RA^2},\qquad
 M^{z\theta}_{\rm model}(R,Z)=RVZ^2-\frac{2RA^2J_A}{(1+Z^2)^2}.
\]
Since $\partial_y\log A=-1/2-c$, we obtain $\partial_y J_A=1+2cJ_A$.
For the comparison pair we use the same pure growth extrapolation for the angular energy, so
$A\propto R^{1/10}$ gives $J_A=(RA^2)^{-1}\int_0^R A(\rho)^2\,d\rho=5/6$
at the start of the transition. Since $c\le0$ on the fixed transition of length
$L_{\rm tr}$, the damping interval starts with $0<J_A(0)\le5/6+L_{\rm tr}$.
The axial term gives $\partial_y V=d_1^2-V$. In the separated comparison
$P=4p(y)/(1+Z^2)^2$, where $p(y)=-\frac12\int_y^\infty A(s)^2\,ds$.
In particular,
\[
 0\le V\le C_{d_1} e^{-y},\qquad
 J_A\le C_{d_1}(1+y)e^{2\mu y},\qquad |p|\le C A^2.
\]
Substitution gives
\begin{equation}\label{new:axial-inertia}
\begin{aligned}
 \frac{\mathcal I^z}{RF}
 ={}&\frac{Z(1+Z^2)[(d_2-1)d_1-2(1-Z^2)V]}{2A}\\
 &-\frac{4AJ_AZ(1-Z^2)}{(1+Z^2)^2}
   +\frac{4pZ(3-Z^2)}{A(1+Z^2)^2}.
\end{aligned}
\end{equation}
Consequently,
\[
 \frac{|\mathcal I^z|}{RF}\le C_{d_1} A_s(1+y)e^{-(1/2-\mu)y}|Z|,
 \qquad b=\frac{d_1'Z(1+Z^2)}{A}.
\]
For $A_s$ large and $\mu$ small, \eqref{new:slow-hazard} gives for the comparison pair
\begin{equation}\label{new:cone-smallness}
 \frac{|b\mathcal I^z|}{RF}\le C\eta Z^2,\qquad
 b^2\le C_{d_1} A_s^{-2},\qquad
 \mu(\mathcal I^{z\theta})^2
 \le C_{d_1}\mu A_s^2(1+|\log\mu|)^2.
\end{equation}
For the first bound use $|d_1'|(1+J_A)\le C\eta$ and
$|d_1'|(d_1+V)/A^2\le C_{d_1} A_s^{-2}$, since the Gaussian tail absorbs $e^{2\mu y}$.
For the last, use $Z^2/(e^{-y}+\mu+Z^2)^2\le1/[4(e^{-y}+\mu)]$.
We also use the bounds in $y$
\[
 \sup_{y\ge0}\frac{\mu(1+y)^n e^{-(1-2\mu)y}}{(e^{-y}+\mu)^2}
       \le C(1+|\log\mu|)^n,\qquad n=0,2.
\]
Split at $y=|\log\mu|$, using $e^{-y}+\mu\ge e^{-y}$ before the split
and $e^{-y}+\mu\ge\mu$ after it. The same split gives
$e^{-(1-2\mu)y}/(e^{-y}+\mu)\le C$ for the odd inertia.

\paragraph*{The actual incoming data.}\label{np:parity}

Let $U^z=\nu_oU^z_{1o}+\nu_eU^z_{1e}$, with $\nu_o=d_1/d_0$ and fixed Gaussian decay
for $\nu_e$. Since $\partial_Z U^z_{1e}(0)=0$, the actual and comparison $d_2$ agree
throughout damping. Oddness gives $U^z_{1o}-d_0Z=O(\epsilon_{\rm odd}Z^3)$.
Since $U^\theta=A(y)q(Z)$ during damping, define
\[
 \begin{aligned}
 \frac{M^z}{R}&=m_o^z U^z_{1o}+m_e^z U^z_{1e}+r_z,&
 \frac{M^\theta}{R\sqrt{2R}Aq}&=\frac{M^\theta_{\rm model}(R,1)}{R\sqrt{2R}\,A}+r_\theta,\\
 \frac{M^{\theta z}}{R\sqrt{2R}Aq}&=m_o^{\theta z}U^z_{1o}+m_e^{\theta z}U^z_{1e}+r_{\theta z}.
 \end{aligned}
\]
As $q$ is independent of $y$, one has
\[
 \partial_y m_i^z=\nu_i-m_i^z,\qquad
 \partial_y m_i^{\theta z}=\nu_i-(1-c)m_i^{\theta z}\quad(i=o,e).
\]
Thus the three remainders satisfy
\[
 \|r_z\|_{C_Z^4}\le\xi e^{-y},\qquad
 \|r_\theta\|_{C_Z^4}+\|r_{\theta z}\|_{C_Z^4}
       \le C\xi e^{-(1-\mu)y}.
\]
The contribution of $U^z_{1i}$ to $\widetilde I^\theta$ is
\[
 m_i^z(ZU^z_{1i}+L\partial_Z U^z_{1i})
       -m_i^{\theta z}\left[L\partial_Z U^z_{1i}+\left(Z+L\partial_Z\log q\right)U^z_{1i}\right].
\]
The angular linear moment contributes
$M^\theta_{\rm model}(R,1)(1-\tfrac12Z\partial_Z\log q)/(R\sqrt{2R}\,A)$.
Choose the prefix and pressure-shape errors to be at most
$\xi\le\mu^4$ in $C_Z^4$. Then
\[
 \left\|\partial_Z\log q+\frac{2Z}{1+Z^2}\right\|_{C_Z^3}
       \le C\xi.
\]
Taylor's formula and the parity of $U^z_{1o},U^z_{1e}$ give
\[
 |\widetilde I^\theta-\widetilde I_{\rm model}^\theta|
 \le C\epsilon_{\rm odd} Z^2
      +Cj_0 e^{-(1-\mu)y}|Z|
      +C\xi e^{-(1-\mu)y}+C\xi.
\]
Since $e^{-(1-\mu)y}\le C(e^{-y}+\mu)$ and $\xi/\mu\le\mu^3$, division by
$e^{-y}+\mu+Z^2$ bounds this error by $C(j_0+\epsilon_{\rm odd}+\mu^3)$.
Take this below $c_{d_1}/2$ in \eqref{new:Q-bound} to obtain
$\widetilde I^\theta\ge(c_{d_1}/2)(e^{-y}+\mu+Z^2)$.
Taking the odd and even parts of $\mathcal I^z/(RF)$ in $Z$ gives
\[
 \begin{aligned}
 \left|\left(\frac{\mathcal I^z}{RF}\right)_o\right|&\le C_{d_1} A_s(1+y)e^{-(1/2-\mu)y}|Z|,\\
 \left|\left(\frac{\mathcal I^z}{RF}\right)_e\right|&\le C_{d_1}\bigl[j_0/A_s+\xi A_s(1+y)\bigr]e^{-(1/2-\mu)y}.
 \end{aligned}
\]
For the even reference shape, the pressure and angular-energy terms
are odd; the even axial terms and odd axial energy are $O(j_0R_s)$
before division by $RU^\theta$.
A $C_Z^4$ shape error of size $\xi$ contributes
$O(\xi\mathcal E_s(1+y))$ to angular energy and $O(\xi A^2)$ to
pressure, hence $O(\xi A_s(1+y)e^{-(1/2-\mu)y})$ after division,
with the required $Z$ derivatives.

The mixed products retain the Gaussian decay. Indeed,
$b=2(\nu_o'U^z_{1o}+\nu_e'U^z_{1e})/(Aq)$, and the fixed damping functions give
\[
 (1+y)(|\nu_o'|+|\nu_e'|)e^{2\mu y}\le C_{d_1} e^{-y/2}
\]
for small $\mu$. Thus the even shear against the odd inertia satisfies
\[
 \left|\frac{2\nu_e'U^z_{1e}}{Aq}\left(\frac{\mathcal I^z}{RF}\right)_o\right|
 \le C_{d_1} j_0 e^{-y/2}|Z|
 \le \frac{C_{d_1}j_0}{2}(e^{-y}+Z^2).
\]
The same Gaussian bound gives $C_{d_1}(j_0A_s^{-2}+\xi)(e^{-y}+Z^2)$ for
the other mixed product and $C_{d_1}(j_0^2A_s^{-2}+j_0\xi)e^{-y}$ for the
even-even product. The odd-odd product satisfies \eqref{new:cone-smallness}.
Thus
\[
 \frac{|b\mathcal I^z|}{RF}\le C\eta Z^2+\bigl[C_{d_1}j_0+C_{d_1}(A_s^{-2}+\xi)\bigr]\widetilde I^\theta,
\]
\[
 \mu(\mathcal I^{z\theta})^2
 \le C_{d_1}\left[\mu A_s^2(1+|\log\mu|)^2
       +j_0^2/A_s^2+\xi^2 A_s^2(1+|\log\mu|)^2\right].
\]
The odd term uses \eqref{new:cone-smallness} and the even terms use the
bounds in $y$ following it with $n=0,2$. The formula for $b$ above gives $b^2\le C_{d_1} A_s^{-2}$.
Since $-5\mu\le c\le\mu$, we have $1\le a=2+2c\le4$.
Also $\mathcal I^\theta/F=R\widetilde I^\theta\ge c_{d_1} R_s/2\ge16$. With the \hyperref[np:closure]{parameter choices below},
$b^2/a+|b\mathcal I^{z\theta}|+2\mu(\mathcal I^{z\theta})^2\le1/2$ uniformly in $Z$.
Thus Lemma~\ref{oc:cone-test} applies, including at $Z=0$.

\paragraph*{A separated interval.}\label{np:separated}

At the end of damping, $y=\mu^{-1/2}$, we have
\[
 \frac{|\mathcal I^z|}{RF}\le C_{d_1} A_s(1+\mu^{-1/2})e^{-(1/2-\mu)/\sqrt\mu},\qquad \widetilde I^\theta\ge c\mu.
\]
Since $A_s\asymp\mu^{-1/4}$, $|\mathcal I^{z\theta}|=O_N(\mu^N)$ for every fixed $N$.
Multiply $U^z$ by $1-\sigma_{\rm norm}(y-\mu^{-1/2})$; Gaussian decay gives
$O_N(\mu^N)$ velocity, shear and normalized-moment errors.
Then replace the angular shape by
$q(Z)^{1-\sigma_{\rm norm}((y-\mu^{-1/2}-1)/L_f)}$, with base shear $a=3$ and fixed $L_f$ such that
$2\|\sigma_{\rm norm}'\|_\infty\|\log q\|_\infty/L_f<1/2$.
Keep the short shear transitions above $2+\mu$.
During this angular shape change $U^z=0$, so the incoming axial
moment $M_-^z$ is constant in $R$, and
\[
 (\widetilde I^\theta)'+(1-c)\widetilde I^\theta
 =c\left(1-\frac{ZM_-^z+(1-Z^2)\partial_ZM_-^z}{R}\right)
       -\frac12Z\partial_Z\log U^\theta.
\]
Here $c=(a-2)/2$. The fraction in the preceding display is exponentially small.
Since $c\ge c_0>0$ during the shape change and
$-Z\partial_Z\log q\ge-O(\mu^4)$, the forcing is positive.
Thus $\widetilde I^\theta\ge c\mu$. Direct computation gives
\[
 \frac{|\mathcal I^z|}{RF}\le C A_s(1+\mu^{-1/2})e^{-c/\sqrt\mu},\qquad |\mathcal I^{z\theta}|=O_N(\mu^N).
\]
The energy moment from the preceding intervals and its $Z$ derivatives
are $O(\mathcal E_s(1+\mu^{-1/2}))$,
$RU^\theta\asymp R_sA_s e^{1/(2\sqrt\mu)}e^{O(\sqrt\mu)}$, and the pressure
from the outer tail integral is $O((U^\theta)^2)$. Thus both cone conditions hold.

Let $R_{\mathrm E,-}$ be the left endpoint of the next interval, where the velocities are independent of $Z$.
Set $A_{\mathrm E,-}=U^\theta(R_{\mathrm E,-},Z)$ and $y_{\mathrm E}=\log(R/R_{\mathrm E,-})$. Then
\[
 U^\theta=A_{\mathrm E,-}e^{-(1/2+\mu)y_{\mathrm E}},\qquad U_*^z=0,
 \qquad 0\le y_{\mathrm E}\le1/\mu.
\]
Then $\widetilde I^\theta\ge c\mu$. Put $\mathcal E_0=R_{\mathrm E,-}A_{\mathrm E,-}^2$.
We preserve $M^z$ and $M^{\theta z}$ by cancelling their two distinct weights.
Since $U_*^z=0$, the energy increment is quadratic.
We take $U^z=\mu^{1/4}z(Z)H_\theta(y_{\mathrm E})$, where $H_\theta$ is the radial factor
of the axial correction. For a fixed positive flat bump $\varphi$ on $(0,1)$, set
\[
 \Phi(y_{\mathrm E})=\varphi(\mu y_{\mathrm E}),\qquad
 H_\theta(y_{\mathrm E})=U^\theta(\partial_{y_{\mathrm E}}+1/2-\mu)(\partial_{y_{\mathrm E}}+1/2-2\mu)\Phi(y_{\mathrm E}).
\]
Integration by parts gives
\[
 \begin{aligned}
 \int H_\theta\,dR
 &=R_{\mathrm E,-}A_{\mathrm E,-}\left[e^{(1/2-\mu)y_{\mathrm E}}
                (\partial_{y_{\mathrm E}}+1/2-2\mu)\Phi\right]_0^{1/\mu}=0,\\
 \int\sqrt{2R}U^\theta H_\theta\,dR
 &=R_{\mathrm E,-}\sqrt{2R_{\mathrm E,-}}A_{\mathrm E,-}^2\left[e^{(1/2-2\mu)y_{\mathrm E}}
                 (\partial_{y_{\mathrm E}}+1/2-\mu)\Phi\right]_0^{1/\mu}=0.
 \end{aligned}
\]
All endpoint derivatives are zero. The change in the energy moment is
\begin{equation}\label{new:pure-quadratic-energy}
 \Delta M^{z\theta}
 =\frac{\mathcal E_0}{\sqrt\mu}\,b_\mu z(Z)^2,
 \qquad
 b_\mu=\int_0^1e^{-2s}
 \bigl[(1/2-\mu)(1/2-2\mu)\varphi
       +(1-3\mu)\mu\varphi'+\mu^2\varphi''\bigr]^2ds.
\end{equation}
In particular, $b_\mu\to\frac1{16}\int_0^1 e^{-2s}\varphi(s)^2\,ds>0$.
\paragraph*{Exact energy correction.}\label{np:energy}

Figure~\ref{fig:one-axial} shows the axial correction. Let $\mathcal E_{\rm uncor}(Z)$
denote the uncorrected total energy moment ($z\equiv0$). We impose
$\mathcal E_{\rm uncor}(1)=-K\mathcal E_0/\sqrt\mu$ \hyperref[np:closure]{below}.
By \eqref{new:pure-quadratic-energy}, cancelling this moment requires the positive amplitude
\begin{equation}\label{new:positive-root}
 z(Z)=\sqrt{-\frac{\sqrt\mu\,\mathcal E_{\rm uncor}(Z)}{\mathcal E_0 b_\mu}}.
\end{equation}
Put $q_0(Z)=[q(Z)+q(-Z)]/[q(1)+q(-1)]$. Evenness and the $C_Z^2$
bound relative to $2/(1+Z^2)$ give $q_0'(Z)\le-cZ$ for $Z\ge0$.
Thus $q_0(1)=1$ and $q_0\ge1$ on $[-1,1]$.
After the initial shape change, write the comparison angular profile
as $A(R)q_0(Z)^{\vartheta(R)}$. Take $\vartheta=1$ during growth and
damping, then decrease it smoothly to zero. Define the positive even comparison amplitude $z_0$ by
\[
 b_\mu z_0(Z)^2=K+\frac{\sqrt\mu}{2\mathcal E_0}
    \int A(R)^2\bigl[q_0(Z)^{2\vartheta(R)}-1\bigr]\,dR.
\]
The integral covers all intervals after the initial prefix with nonnegative weights.
On $R_s\le R\le R_s e^{1/\sqrt\mu}$ we have $\vartheta=1$ and
$RA(R)^2\asymp R_sA_s^2\asymp\mathcal E_0$.
The remaining growth and transition intervals contribute at most
$C\mathcal E_0$, and the shape interval has fixed logarithmic length. Thus
\[
 b_\mu z_0(1)^2=K,\qquad \partial_Z(b_\mu z_0^2)\le-cZ\quad(0\le Z\le1),\qquad
 \|b_\mu z_0^2-K\|_{C_Z^3}\le C.
\]
This includes the normalized transition energies of size $\sqrt\mu$.

The reserve energy is independent of $Z$ and cancels in the pole difference.
The preceding axial energy is $O(R_s)$ and in the units above contributes
$O(\sqrt\mu/P_*^2)=O(\mu)$ by \eqref{new:parameter-balance}.
The initial prefix is at most $\mu$ in these units for large amplification
radius, and replacing $q$ by $q_0$ costs $O(\xi)$, with $\xi\le\mu^4$.
At $\delta=0$ the matched terminal and central velocities are independent
of $Z$, as is the exterior energy $\epsilon_{\rm em}^2\mathcal J^\theta$ from
\eqref{ue:right-moments} in the finite-energy coordinate. Their normalized pole differences remain $O(\mu)$
for small positive $\delta$ by the \hyperref[np:closure]{parameter continuation below}.
Thus $\|-\sqrt\mu\,\mathcal E_{\rm uncor}/\mathcal E_0-b_\mu z_0^2\|_{C_Z^3}\le C\mu$.
For fixed large $K$ and small $\mu$, the normalized energy
$-\sqrt\mu\,\mathcal E_{\rm uncor}/\mathcal E_0$ is at least $K/2$. Taking square roots gives
\[
 z_0\asymp\sqrt K,\quad z_0'\le0\ (Z\ge0),\quad
 |Zz_0'|/z_0\le C/K,\quad
 \|z-z_0\|_{C_Z^3}\le C\mu/\sqrt K.
\]

\paragraph*{Cone estimates for the energy correction.}\label{np:energy-cone}

Substitution gives the change in $\widetilde I^\theta$,
\[
 \Delta \widetilde I^\theta=-\mu^{5/4}U^\theta(Zz+(1-Z^2)z')\Phi,
\]
and the part of the change in $\mathcal I^z/(RF)$ linear in $z$,
\[
 \mu^{1/4}\left[z(\mu-\partial_{y_{\mathrm E}})-\tfrac12 Z\partial_Z z\right]
             (\partial_{y_{\mathrm E}}+\tfrac12-2\mu)\Phi
 +\mu^{1/4}zH_\theta\frac{ZM_-^z+(1-Z^2)\partial_ZM_-^z}{RU^\theta},
\]
where $M_-^z$ is the constant axial moment before this correction.
The last term in the preceding display is exponentially small, since
$\|M_-^z\|_{C_Z^3}\le C R_s$ and $R/R_s\ge c e^{\mu^{-1/2}}$.
The quadratic part of the change in $\mathcal I^z/(RF)$ is
\[
 \mu^{1/2}\bigl[Zz^2+(1-Z^2)zz'\bigr]H_\theta
              (\partial_{y_{\mathrm E}}+\tfrac12-2\mu)\Phi
 -\frac{2\mu^{1/2}(1-Z^2)zz'}{RU^\theta}
             \int_{R_{\mathrm E,-}}^{R}H_\theta^2\,d\rho.
\]
Both quadratic terms are bounded pointwise by $C_K\mu^{1/2}A_{\mathrm E,-}$, since
$(RU^\theta)^{-1}\int_{R_{\mathrm E,-}}^R H_\theta^2\,d\rho
\le CA_{\mathrm E,-} y_{\mathrm E} e^{-(1/2-\mu)y_{\mathrm E}}\le CA_{\mathrm E,-}$.
Thus the quadratic terms, the preceding exponentially small term,
the value of \(\mathcal I^z/(RF)\) for \(z(\cdot)\equiv0\) and $\Delta \widetilde I^\theta$ are smaller than every fixed power of $\mu$.

Let $B=\mu^{1/4}z_0\Phi/4\ge0$. Differentiation gives $b=-B+O_K(\mu^{5/4})$ and $a=2+2\mu$.
In these bounds the factor $(\partial_{y_{\mathrm E}}+\tfrac12-2\mu)\Phi$ equals
$(\tfrac12-2\mu)\Phi+O(\mu)$, since $\partial_{y_{\mathrm E}}\Phi=\mu\varphi'(\mu y_{\mathrm E})$.
Its $y_{\mathrm E}$ derivative is $O(\mu)$ since $\partial_{y_{\mathrm E}}^2\Phi=\mu^2\varphi''(\mu y_{\mathrm E})$.
Since $\widetilde I^\theta\ge c\mu$, the preceding inertia estimates give
\[
 \mathcal{I}^{z\theta}=\mathcal{I}_{\rm main}^{z\theta}+O_K(\mu^{1/4}),\qquad
 \mathcal{I}_{\rm main}^{z\theta}=-\frac{\mu^{1/4}Zz_0'\Phi}{4\widetilde I^\theta}\ge0,
\]
\[
 \mathcal{I}_{\rm main}^{z\theta}=O_K(\mu^{-3/4}),\qquad
 \mu \mathcal{I}_{\rm main}^{z\theta}\le(C/K)B.
\]
For $a=2+2\mu$, these bounds give
\[
 b^2/a+b\mathcal{I}^{z\theta}+2\mu (\mathcal{I}^{z\theta})^2
 \le -(1-2C/K)B\mathcal{I}_{\rm main}^{z\theta}+O_K(\mu^{1/2})\le\tfrac12
\]
for fixed large $K$ and then small $\mu$. Since $\mathcal I^\theta/F=R\widetilde I^\theta\ge c\mu R_{\mathrm E,-}\ge16$
and $\kappa\ge a>2$, Lemma~\ref{oc:cone-test} gives the strict cone.

\paragraph*{The central angular reserve.}\label{np:reserve}

On the reserve interval the normalized angular moment approaches the
central entrance value while the shear remains strictly admissible.
Set \(R_{\mathrm E,+}=R_{\mathrm E,-}e^{1/\mu}\), \(A_{\mathrm E,+}=A_{\mathrm E,-}e^{-1/(2\mu)-1}\) and \(y_{\mathrm{res}}=\log(R/R_{\mathrm E,+})\).
Let $L_{\rm res}>1$ be the logarithmic length of the reserve, so its
outer radius is $R_{\mathrm{res},+}=R_{\mathrm E,+}e^{L_{\rm res}}$. On $0\le y_{\mathrm{res}}\le L_{\rm res}$ use the rate $r_0$ of
\eqref{nc:paper-rates} and the cutoff $\sigma_{\rm norm}$ of \eqref{pr:cutoffs}. Prescribe
\[
 c(y_{\mathrm{res}})=\mu+(r_0-\mu)\sigma_{\rm norm}(y_{\mathrm{res}}),\qquad U^z=0,\qquad
 U^\theta=A(y_{\mathrm{res}})=A_{\mathrm E,+}\exp\!\left[-\tfrac12y_{\mathrm{res}}-\int_0^{y_{\mathrm{res}}} c(v)\,dv\right].
\]
All derivatives match at $y_{\mathrm{res}}=0$, by the flatness of $\sigma_{\rm norm}$ there and of
$H_\theta$ at the end of the energy window. For $y_{\mathrm{res}}\ge1$ one has
$A(y_{\mathrm{res}})=A(1)e^{-(1/2+r_0)(y_{\mathrm{res}}-1)}$. The two vanishing integrals of $H_\theta$
give the original values of $M^z,M^{\theta z}$ at $R_{\mathrm E,+}$. Both are constant
on the reserve.

At $\delta=0$ we have $\mathcal I^z/(RF)=A\widetilde I^z$ on the reserve.
At $\Phi=0$ the energy correction gives $\widetilde I^\theta(0)\gtrsim\mu$ and
$|\mathcal{I}^{z\theta}(0)|\lesssim_K\mu^{1/4}$. Since $A$ is independent of $Z$,
\eqref{oc:radial-system} gives
\[
 \partial_{y_{\mathrm{res}}}\widetilde I^\theta+(1-c)\widetilde I^\theta=c\left(1-\frac{ZM^z+L\partial_ZM^z}{R}\right),\qquad
 \partial_{y_{\mathrm{res}}}(A\widetilde I^z)=(c-\tfrac12)A\widetilde I^z+A(Z+\Pi).
\]
The fraction in the first equation is $O(R_s/R_{\mathrm E,+})$, since $M^z$ stays constant on the reserve. At the value selected by the pressure
equation, $P=-\frac12\int_R^\infty(U^\theta)^2\,d\rho/\rho$. The part up
to the end of the reserve is independent of $Z$ and bounded by $A(y_{\mathrm{res}})^2$,
since $\partial_{y_{\mathrm{res}}}\log A\le-1/2$. The remaining tail has $C_Z^1$ norm
$O_{N_{\rm mid}}((U_b^\theta)^2)$ by the \hyperref[np:closure]{pressure bounds below},
where $N_{\rm mid}=\log(R_{\rm mid}/R_b)$ is the logarithmic length of the central auxiliary connection.
Consequently $|\Pi|\le C$ throughout the reserve. The first equation gives
$\widetilde I^\theta(y_{\mathrm{res}})\ge e^{-y_{\mathrm{res}}}\widetilde I^\theta(0)$, hence $R\widetilde I^\theta\ge R_{\mathrm E,+}\widetilde I^\theta(0)\gg1$. Integrating the second
equation over the fixed length $L_{\rm res}$ gives
\[
 \sup_{0\le y_{\mathrm{res}}\le L_{\rm res}} |\mathcal{I}^{z\theta}(y_{\mathrm{res}})|\le C_{L_{\rm res},N_{\rm mid}}
       \left(|\mathcal{I}^{z\theta}(0)|+A_{\mathrm E,+}/\widetilde I^\theta(0)\right)=o(1)\qquad(\mu\downarrow0).
\]
Here $b=0$, $\kappa=a=2+2c\in(2,4)$, and $(a-2)(\mathcal{I}^{z\theta})^2=o(1)$.
Since $R\widetilde I^\theta\ge16$, Lemma~\ref{oc:cone-test} gives the strict cone throughout the reserve.

At the reserve exit $R_{\mathrm{res},+}$ put $A_{\mathrm{res},+}=A(L_{\rm res})$.
By \eqref{pr:moments}, $x=M^\theta/(R\sqrt{2R}U^\theta)$ solves
$\partial_{y_{\mathrm{res}}}x=1-(1-c)x$. On the preceding power interval its deviation
from $(1-\mu)^{-1}$ contracts by $e^{-(1-\mu)/\mu}$. Integration over \(0\le y_{\mathrm{res}}\le1\) bounds \(x(1)\) in \(C_Z^3\). On the constant part,
\[
 x(L_{\rm res})-\frac1{1-r_0}
 =e^{-(1-r_0)(L_{\rm res}-1)}
      \left(x(1)-\frac1{1-r_0}\right).
\]
Choose $L_{\rm res}$ so the right side is smaller than $\tau/10$ in $C_Z^3$
before the short transition to the central inlet. The same contraction
brings $\widetilde I^\theta$ close to $r_0/(1-r_0)$ with an additional $O(R_s/R_{\mathrm E,+})$ error
from the axial moment in the first equation. Moreover,
\[
 \left\|\frac{M^z(R_{\mathrm{res},+})}{R_{\mathrm{res},+}A_{\mathrm{res},+}}\right\|_{C_Z^3}
 +\left\|\frac{M^{\theta z}(R_{\mathrm{res},+})}{R_{\mathrm{res},+}\sqrt{2R_{\mathrm{res},+}}A_{\mathrm{res},+}^2}\right\|_{C_Z^3}
 \le C A_s^{-1}e^{-c/\mu}+o_{C_*\to\infty}(1).
\]
Indeed, $M^z=O(R_s)$ and $M^{\theta z}=O(R_s\sqrt{R_s}A_s)$. Their
denominators grow on the separated interval as $e^{(1/2-\mu)y_{\mathrm E}}$
and $e^{(1/2-2\mu)y_{\mathrm E}}$, respectively. The correction $H_\theta$
preserves both moments and the fixed reserve changes these denominators
by bounded factors. The same angular decay gives $A_{\mathrm{res},+}\le CP_*e^{-c/\mu}$.

\paragraph*{Matching the inlet velocities and choosing the parameters.}\label{np:closure}
For this transition choose \(y_{b,-}<0\) close to zero and set \(R_b=R_{\mathrm{res},+}e^{-y_{b,-}}\).
Then \(y_b=\log(R/R_b)\) runs from \(y_{b,-}=\log(R_{\mathrm{res},+}/R_b)\) to \(0\).
At \(R_{\mathrm{res},+}\) the velocities are \(A_{\mathrm{res},+}>0\) and zero, both independent of \(Z\).
In the backward central model \eqref{nc:paper-incoming-control}, set
\(D_{\rm in}=K_{\rm in}=0\) and determine the inlet ratio
\(B_{\rm in}=B(0)=U^z(R_b)/U^\theta(R_b)\) at fixed \(\epsilon_{\rm em}\).
The two corresponding forward moments retain their actual, generally nonzero, values.

Continue the model shears \(a_c(y_b),b_c(y_b;B_{\rm in},\epsilon_{\rm em})\) across \(y_b=0\)
to \([y_{b,-},0]\) for sufficiently small \(|y_{b,-}|\). This uses the analytic germ of \(Q\) at
zero and the ordinary differential equations \eqref{nc:paper-state-odes}.
The angular shear is independent of both parameters.
The terminal control equations \eqref{nc:paper-functionals}
are linear with the fixed invertible two-column matrix. Hence
\(b_c\) is affine in \(B_{\rm in}\) and \(\sqrt{1+\epsilon_{\rm em}^2}\), with bounded derivatives for
fixed \(N_{\rm mid}\). Take a smooth cutoff \(\chi(y_b)\) equal to zero near
\(y_{b,-}\) and one near zero, with \(0\le\chi\le1\), and prescribe
\[
 a=(1-\chi)(2+2r_0)+\chi a_c,\qquad b=\chi b_c .
\]
Define \(g_a(y_b)=\exp(\tfrac12\int_{y_{b,-}}^{y_b}(1-a(v))\,dv)\).
Since \(a=1-2\partial_{y_b}\log U^\theta\) and
\(b=2\partial_{y_b}U^z/U^\theta\), integration gives
\[
 A_b=A_{\mathrm{res},+}g_a(0),\qquad
 U^\theta(y_b)=A_{\mathrm{res},+}g_a(y_b),\qquad
 U^z(y_b)=\frac{A_{\mathrm{res},+}}{2}\int_{y_{b,-}}^{y_b} g_a(v)\chi(v)b_c(v;B_{\rm in},\epsilon_{\rm em})\,dv .
\]
Both endpoint velocities match precisely when
\[
 B_{\rm in}=\frac{1}{2g_a(0)}\int_{y_{b,-}}^{0}
       g_a(v)\chi(v)b_c(v;B_{\rm in},\epsilon_{\rm em})\,dv
       =:\mathcal B_{y_{b,-}}(B_{\rm in},\epsilon_{\rm em}).
\]
The ratio \(g_a(v)/g_a(0)\) stays bounded independently of \(A_{\mathrm{res},+}\).
Consequently \(|\partial_{B_{\rm in}}\mathcal B_{y_{b,-}}|\le C_{N_{\rm mid}}|y_{b,-}|\)
and \(|\mathcal B_{y_{b,-}}(0,\epsilon_{\rm em})|\le C_{N_{\rm mid}}|y_{b,-}|\) for small \(\epsilon_{\rm em}\).
Choosing \(C_{N_{\rm mid}}|y_{b,-}|<1/2\) gives a unique
\(B_{\rm in}=O_{N_{\rm mid}}(|y_{b,-}|)\) depending smoothly on \(\sqrt{1+\epsilon_{\rm em}^2}\).
The solution \(B_{\rm in}\) is independent of \(Z\). The derivative acts by scalar
multiplication, giving the same bound in \(C_Z^3\). Matching the velocities at the endpoints
and the shears on their collars matches all jets.

The multiplier \(\eta(y_b)=W_1/W(y_b)\) in \eqref{nc:paper-Gamma-profile} depends only on \(y_b=\log(R/R_b)\).
Since \(W(0)=1\) by \eqref{nc:paper-angular}, \(\eta(0)=W_1\). At \((R_b,0)\) the equation
\(A_b=h_\infty W_1 R_b^{-(1+\delta)/2}H_e(2/R_b)\) determines \(h_\infty>0\) by \eqref{ue:heat-factor}.
The entrance adjustment \eqref{nc:paper-entrance-lift} matches the model velocity jets for every \(Z\).
The angular profile is independent of the exterior axial coefficient.

The shear on the short interval remains in the common strict
neighborhood of the reserve and central entrance shears:
\(a_c(0)=2+2r_0+O(N_{\rm mid}^{-1})\), and \(b_c(0;0,\epsilon_{\rm em})=O(W_1^{-1})\) by
\eqref{nc:paper-control}--\eqref{nc:paper-model}. Continuity for fixed
\(N_{\rm mid}\) and the bound on \(B_{\rm in}\) preserve these estimates up to
\(O_{N_{\rm mid}}(|y_{b,-}|)\). In the units \eqref{nc:paper-moment-units} every
density integral on this interval is \(O_{N_{\rm mid}}(|y_{b,-}|)\) in \(C_Z^3\):
\(dR=R\,dy_b\), \(R/R_b\asymp1\), and the two velocities divided by
\(A_b\) have bounded \(C_Z^3\) norms. Thus the pre-existing normalized
angular moment error and the two small forward axial and mixed moments
change by \(O_{N_{\rm mid}}(|y_{b,-}|)\), independently of \(C_*\).
Take this error below \(\tau/10\). By \eqref{nc:paper-entrance-targets},
the backward axial and mixed targets are small in precisely these units.
With pressure and energy matched below, substitution in
\eqref{pr:moment-inertia} changes the normalized inertia
\(\mathcal I/(RF)\) over the short interval by
\(O_{N_{\rm mid}}(\tau+|y_{b,-}|+R_b^{-1}+(R_b A_b)^{-1})+o_{\delta\to0}(1)\) relative to
the central entrance. Indeed the terms are linear in the normalized
moments and their first \(Z\)-derivatives, with bounded multipliers
when \(R/R_b\) and \(U^\theta/A_b\) stay close to one and \(A_b\le1\).
For the potentially large constant part of \(M^{z\theta}\), the
coefficient in \eqref{pr:moment-inertia} contains \(\delta\);
its contribution tends to zero at fixed reference parameters as
\(\delta\downarrow0\).
The positive central margins in \eqref{nc:paper-model-margins} and
\eqref{nc:paper-relative-stress} therefore give strict admissibility.
Lemma~\ref{nc:paper-five-window} subsequently removes the three
remaining forward moment errors.

\paragraph{Pressure and exterior energy.}
The transition and central angular profiles are independent of the
exterior axial coefficient and of \(z(Z)\). After dividing the angular
velocities by \(A_b\), they and their fixed \(Z\)-derivatives are
bounded on the fixed logarithmic interval by \eqref{nc:paper-heat-estimates}
and \eqref{nc:paper-entrance-lift}. The same bounds hold after a
\(y_{\mathrm{g},+}\)-derivative: \(R_b\) is a fixed multiple of \(R_s\),
\(\partial_{y_{\mathrm{g},+}}\log R_s=1\), and \(\partial_{y_{\mathrm{g},+}}\log A_b\) is bounded
by the growth formula for \(A_s\) and the fixed length transitions.
In the exterior use \(R=R_{\rm mid}\,v\), \(v\ge1\): the normalized pressure
integrand and its differentiated versions are bounded by \(C v^{-2}\).
Thus the transition, central profile and exterior contribute \(O_{N_{\rm mid}}(A_b^2)\)
to \(\mathcal R_p\), also after these derivatives, as required by
the \hyperref[np:amplification]{pressure inverse above}.

Put $e=\epsilon_{\rm em}^2/\delta$ and
$\mathcal E_h=R_{\rm mid}[U_e^\theta(R_{\rm mid},1)]^2>0$.
By Lemma~\ref{ue:regularized-tail}, the exterior energy at $Z=1$ is $e\mathcal E_h/2$.
At $\delta=0$ the finite velocities have $\epsilon_{\rm em}=0$ and are
independent of $e$. With $z=0$, the separated interval contributes
\[
 -\frac12\int_0^{1/\mu}\mathcal E_0e^{-2\mu y_{\mathrm E}}\,dy_{\mathrm E}
       =-\frac{1-e^{-2}}{4\mu}\mathcal E_0.
\]
The reserve has $S=RA^2$, $S(0)=e^{-2}\mathcal E_0$ and $\partial_{y_{\mathrm{res}}}S=-2cS$,
hence contributes $O(\mathcal E_0)$ over its fixed length. The preceding
axial energy is $O(R_s)$ and the remaining finite terms are
$O(\mathcal E_0(1+\mu^{-1/2}))$. Thus the uncorrected finite energy at
$Z=1$ lies between $-C\mathcal E_0/\mu$ and $-c\mathcal E_0/\mu$.
Adding $e\mathcal E_h/2$ reaches $-K\mathcal E_0/\sqrt\mu$ at a unique
$e\asymp\mathcal E_0/(\mu\mathcal E_h)>0$. Equation \eqref{new:positive-root}
then supplies the pointwise correction.

\paragraph{Parameter order and continuation.}
The $Z^2$ coefficient in \eqref{new:angular-inertia} has a positive
lower bound independent of the damping rate. Use it to fix $\eta$,
then the damping function, uniformly for $d_0\in[3,5]$.
Fix the central length, the reserve length and the short transition
for the tolerance $\tau$, and take $K$ as in the \hyperref[np:energy]{energy correction}.
Choose $j_0,\epsilon_{\rm odd}$ below the transition and parity tolerances.
These tolerances are independent of $\mu$.
Then take $\mu$ small with
\begin{equation}\label{new:parameter-balance}
 \mu P_*^4\text{ fixed and sufficiently small},\qquad A_s\asymp P_*.
\end{equation}
Then $A_s\to\infty$ and $\mu A_s^2(1+|\log\mu|)^2\to0$, as required.
With the other core parameters
fixed, choose $C_*$ large enough that
$\xi\le\mu^4$, the normalized errors from the preceding intervals
are below the fixed tolerance, and $R_b^{-1}+(R_bU_b^\theta)^{-1}$ is sufficiently small.
Indeed $U_b^\theta\asymp A_{\mathrm{res},+}\le CP_*e^{-c/\mu}$, and for fixed $\mu,P_*$
the growth radius from the pressure equation tends to infinity with $C_*$.
Keep this $C_*$ fixed in the continuation from $\delta=0$.

The incoming data, finite integrals and reconstructed moments are
smooth in $\delta$. Lemma~\ref{ue:regularized-tail} gives the same
extension for all five exterior primitives. In particular,
$M_e^\theta|_{\delta=0}=\sqrt2h_\infty(R+L)$, whose $L$ term gives
the included $O(R_b^{-1})$ normalized entrance error.
The \hyperref[np:amplification]{pressure inverse} gives an invertible block
\(\partial_{(y_{\mathrm{g},+},q)}\mathcal P\) on variations of \(y_{\mathrm{g},+}\) and on \(C_Z^4\) variations of \(q\) vanishing at the pole,
independent of \(e,z\). Use the uncorrected
pole energy and the pointwise energy divided by \(\mathcal E_0/\sqrt{\mu}\)
for the last two rows.
At \(\delta=0\) their joint derivative in \(((y_{\mathrm{g},+},q),e,z)\) is
\[
 \begin{pmatrix}
 \partial_{(y_{\mathrm{g},+},q)}\mathcal P&0&0\\
 *&\mathcal E_h/2&0\\
 *&*&2b_\mu z(\cdot)
 \end{pmatrix}.
\]
By \eqref{new:positive-root}, multiplication by the last diagonal entry has a bounded inverse in \(C_Z^3\). The two vanishing integrals preceding \eqref{new:pure-quadratic-energy} show that varying \(z\) preserves the central inlet and the first two equations.
For \(\delta>0\), the exterior factor \(\sqrt{1+\epsilon_{\rm em}^2}=\sqrt{1+\delta e}\) has \(e\)-derivative
\(O(\delta)\); the differentiated short-interval equation has the
bounded inverse \((1-\partial_{B_{\rm in}}\mathcal B_{y_{b,-}})^{-1}\).
Thus the finite central and transition energy has \(e\)-derivative
\(O_{N_{\rm mid}}(\delta\mathcal E_b)\),
where $\mathcal E_b=R_b(U_b^\theta)^2$.
For small $\delta$, the energy derivative is at least $\mathcal E_h/4$.
We solve all three equations using the scalar and pointwise inverses. Thus
\[
 \epsilon_{\mathrm{em}}^2\asymp
      \frac{\delta\mathcal E_0}{\mu\mathcal E_h}>0.
\]
Reduce $\delta$ so that the exterior coefficient remains in the
permitted interval and the additional normalized $C_Z^3$ energy
error is at most $\mu$. Continuity on the fixed finite intervals
preserves the strict inequalities defining relaxed admissibility,
the strict cone margins on the reserve and terminal transition, and the incoming
hypotheses. The two zero moments of $H_\theta$ persist,
since their radial weights and the separated power interval are
unchanged. Hence the lemma follows.
\end{proof}
\subsubsection{The auxiliary connection from \texorpdfstring{\(R_b\) to \(R_{\rm mid}\)}{Rb to Rmid}}
\label{nc:paper-section}

Here we match the exact exterior at \(R_{\rm mid}=e^{N_{\rm mid}}R_b\) to
construct an auxiliary profile. Two scalar controls set the axial and
mixed terminal moments. The final profile meets the mixed
power at this radius, with the moments obtained by integrating back from \(R_e\).

For the auxiliary model, choose a nonincreasing weight $Q\in C^\infty([0,1])$, analytic on $[0,1)$,
such that
\begin{equation}\label{nc:paper-shape}
 Q>0\quad(0\le s<1),\qquad
 I:=\int_0^1\frac{Q(s)}{1+Q(s)}\,ds>\frac12.
\end{equation}
Assume $e^{\gamma/(1-s)^2}Q(s)$ extends smoothly and positively to
$s=1$ for some fixed $\gamma>0$.  For example,
$Q(s)=4\exp\{-\gamma[(1-s)^{-2}-1]\}$ has these properties for
sufficiently small $\gamma>0$, since $I\to4/5>1/2$.

Let $y_b=\log(R/R_b)=N_{\rm mid}s$.
Fix $N_{\rm mid}$ sufficiently large before the small parameters below, and define
\begin{equation}\label{nc:paper-angular}
 \begin{alignedat}{2}
 W(y_b)&=\frac{1+Q(s)}{1+Q(0)}
       \exp\left(N_{\rm mid}\int_0^s\frac{Q(t)}{1+Q(t)}\,dt\right),
 &\qquad W_1&=W(N_{\rm mid})=\frac{e^{N_{\rm mid}I}}{1+Q(0)},\\
 c&=\frac{Q+\partial_sQ/N_{\rm mid}}{1+Q},
 &\qquad a&=2+2c.
 \end{alignedat}
\end{equation}
Here and below, an unmarked \(Q\) is evaluated at \(s=y_b/N_{\rm mid}\).
Then $\partial_{y_b}W=cW$ and
\begin{equation}\label{nc:paper-angular-ode}
 \partial_{y_b}Q=c-(1-c)Q.
\end{equation}
Define
\begin{equation}\label{nc:paper-rates}
 r_0=\frac{Q(0)}{1+Q(0)}>\frac12 .
\end{equation}

For two coefficients $A_{N_{\rm mid}},C_{N_{\rm mid}}$, set
\begin{equation}\label{nc:paper-control}
 e(y_b)=\frac{(1-y_b/N_{\rm mid})^2y_b e^{-y_b}}{W_1}
                  (A_{N_{\rm mid}}+C_{N_{\rm mid}}y_b),\qquad E=Qe.
\end{equation}
Solve
\begin{equation}\label{nc:paper-state-odes}
 \begin{gathered}
 \partial_{y_b}G=cG-E/2,\qquad B=G-E,
 \\
 \partial_{y_b}K+(1/2-2c)K=B,\qquad G(0)=K(0)=0.
 \end{gathered}
\end{equation}
For $U_b^\theta>0$, define
\begin{equation}\label{nc:paper-model}
 U^\theta=U_b^\theta e^{-y_b/2}/W(y_b),\qquad
 U^z=BU^\theta,\qquad
 b=-G+(a-2)E-2\partial_{y_b}E.
\end{equation}
Thus $a=1-2\partial_{y_b}\log U^\theta$ and
$b=2\partial_{y_b}U^z/U^\theta$. The three moment equations
\eqref{pr:moments} hold with
\begin{equation}\label{nc:paper-linear-moments}
 \frac{M^\theta}{R\sqrt{2R}U^\theta}=1+Q,\qquad
 \frac{M^z}{RU^\theta}=2G,\qquad
 \frac{M^{\theta z}}{R\sqrt{2R}(U^\theta)^2}=K,
\end{equation}
by $2\partial_{y_b}G+(1-2c)G=B$ and \eqref{nc:paper-angular-ode}.

\begin{lemma}[Two endpoint controls]\label{nc:paper-control-lemma}
Let $\epsilon_{\rm em}^2$ range over a fixed sufficiently small nonnegative
interval. For sufficiently large $N_{\rm mid}$, there are unique coefficients
$A_{N_{\rm mid}},C_{N_{\rm mid}}$ in \eqref{nc:paper-control}, uniformly bounded in
$N_{\rm mid},\epsilon_{\rm em}^2$, such that
\begin{equation}\label{nc:paper-terminal}
 G(N_{\rm mid})=B(N_{\rm mid})=\frac{\sqrt{1+\epsilon_{\rm em}^2}}{\sqrt2},
 \qquad K(N_{\rm mid})=\sqrt{2(1+\epsilon_{\rm em}^2)}.
\end{equation}
Moreover, $A_{N_{\rm mid}}>0$ and $C_{N_{\rm mid}}<0$ throughout this interval.
\end{lemma}
\begin{proof}
The two terminal conditions on \(G,K\) are linear in
\(A_{N_{\rm mid}},C_{N_{\rm mid}}\). We compute their matrix.
Integrating \eqref{nc:paper-state-odes} gives
\begin{equation}\label{nc:paper-state-integrals}
 G(y_b)=-\frac12W(y_b)\int_0^{y_b} E(t)/W(t)\,dt,\qquad
 K(y_b)=W(y_b)^2e^{-y_b/2}\int_0^{y_b} e^{t/2}B(t)/W(t)^2\,dt.
\end{equation}
Define
\begin{equation}\label{nc:paper-kernels}
 \mathcal K_1(y_b)=W(y_b)^{-1},\qquad
 \mathcal K_2(y_b)=-\frac1{2W(y_b)}\int_{y_b}^{N_{\rm mid}}\frac{e^{t/2}}{W(t)}\,dt
                         -\frac{e^{y_b/2}}{W(y_b)^2}.
\end{equation}
Substitution and Fubini give
\begin{equation}\label{nc:paper-functionals}
 G(N_{\rm mid})=-\frac{W_1}{2}\int_0^{N_{\rm mid}}\mathcal K_1E\,dy_b,
 \qquad \frac{e^{N_{\rm mid}/2}K(N_{\rm mid})}{W_1}
                       =W_1\int_0^{N_{\rm mid}}\mathcal K_2E\,dy_b.
\end{equation}
Hence
\begin{equation}\label{nc:paper-kernel-ratio}
 \partial_{y_b}(\mathcal K_2/\mathcal K_1)=c(y_b)e^{y_b/2}/W(y_b).
\end{equation}

Monotonicity of $Q/(1+Q)$ gives
\begin{equation}\label{nc:paper-weight-bounds}
 W(y_b)\ge\frac{e^{Iy_b}}{1+Q(0)},\qquad
 |\mathcal K_1(y_b)|\le Ce^{-Iy_b},\qquad
 |\mathcal K_2(y_b)|\le Ce^{-(2I-1/2)y_b}.
\end{equation}
For fixed $y_b$, $Q(y_b/N_{\rm mid})\to Q(0)$ and $W(y_b)\to e^{r_0y_b}$.
For the columns, substitute
$W_1E=Q(y_b/N_{\rm mid})(1-y_b/N_{\rm mid})^2y_b e^{-y_b}$ and its product with $y_b$.
By \eqref{nc:paper-weight-bounds}, an integrable exponential times a
polynomial bounds both integrands. Dominated convergence gives
\begin{equation}\label{nc:paper-limit-matrix}
 \begin{pmatrix}
 -\dfrac{Q(0)}{2(r_0+1)^2}&-\dfrac{Q(0)}{(r_0+1)^3}\\[8pt]
 -\dfrac{Q(0)r_0}{(r_0-1/2)(2r_0+1/2)^2}
 &-\dfrac{2Q(0)r_0}{(r_0-1/2)(2r_0+1/2)^3}
 \end{pmatrix}.
\end{equation}
Its two ratios of the second entry to the first are the weighted means
\[
 \frac{\int_0^\infty y_b^2e^{-(r_0+1)y_b}\,dy_b}
      {\int_0^\infty y_b e^{-(r_0+1)y_b}\,dy_b}
 =\frac2{r_0+1}
 >\frac2{2r_0+1/2}
 =\frac{\int_0^\infty y_b^2e^{-(2r_0+1/2)y_b}\,dy_b}
      {\int_0^\infty y_b e^{-(2r_0+1/2)y_b}\,dy_b}.
\]
The limiting matrix is invertible. Its finite-\(N_{\rm mid}\) inverses are
uniformly bounded for all sufficiently large \(N_{\rm mid}\).
Since $e^{N_{\rm mid}/2}/W_1\to0$, the second target in
\eqref{nc:paper-functionals} tends to zero. Hence
\begin{equation}\label{nc:paper-coefficient-limits}
 A_{N_{\rm mid}}\longrightarrow
 \frac{\sqrt{2(1+\epsilon_{\rm em}^2)}(r_0+1)^3}
                       {Q(0)(r_0-1/2)}>0,
 \qquad
 C_{N_{\rm mid}}\longrightarrow-\left(r_0+\frac14\right)\lim A_{N_{\rm mid}}<0.
\end{equation}
The convergence is uniform in the stated parameter interval.
\end{proof}

\begin{lemma}[Instability and cone margins]\label{nc:paper-model-cone}
Fix $N_{\rm mid}$ sufficiently large. Then
\begin{equation}\label{nc:paper-model-margins}
 \begin{gathered}
 a\ge3/2,\qquad a+b^2/a-2\ge\kappa_*>0,\\
 \frac{a-be}{\sqrt{(1+e^2)(a^2+b^2)}}\ge\frac13,\qquad
 \frac{2-b^2/a-2be-(a-2)e^2}{1+e^2}\ge1.
 \end{gathered}
\end{equation}
The terminal velocity and moment differences from the mixed critical
power profile are flat.
\end{lemma}
\begin{proof}
By the coefficient bounds,
\begin{equation}\label{nc:paper-small-control}
 \|e\|_{C_{y_b}^m}+\|E\|_{C_{y_b}^m}\le C_m/W_1.
\end{equation}
Since $W(y_b)/W(t)\le1+Q(0)$ for $y_b\le t$, backward integration gives
\begin{equation}\label{nc:paper-shear-limit}
 \left\|G-\frac{\sqrt{1+\epsilon_{\rm em}^2}W}{\sqrt2W_1}\right\|_\infty
 +\left\|b+\frac{\sqrt{1+\epsilon_{\rm em}^2}W}{\sqrt2W_1}\right\|_\infty
                \le CN_{\rm mid}/W_1\longrightarrow0.
\end{equation}
The definition of \(c\) gives the first bound in
\eqref{nc:paper-aW}.
\begin{equation}\label{nc:paper-aW}
\begin{gathered}
 2-2\|\partial_sQ\|_\infty/N_{\rm mid}\le a<4,\\
 \frac{W(y_b)}{W_1}=(1+Q(y_b/N_{\rm mid}))
       \exp\!\left(-\int_{y_b}^{N_{\rm mid}}\frac{Q(v/N_{\rm mid})}{1+Q(v/N_{\rm mid})}\,dv\right)
       \le1+\|\partial_sQ\|_\infty/N_{\rm mid} .
\end{gathered}
\end{equation}
For the second bound, an interior maximum of $W$ has $c=0$ and
$Q\le\|\partial_sQ\|_\infty/N_{\rm mid}$. At the two endpoints the ratio
\(W/W_1\) equals \(1/W_1<1\) and \(1\), respectively, for large \(N_{\rm mid}\).
Thus $a\ge19/10$, $|b|<4/5$, and $|e|<1/100$ for large $N_{\rm mid}$
and small $\epsilon_{\rm em}^2$. Hence
$a-be>3/2$ and $2-b^2/a-2be-(a-2)e^2>3/2$, giving the cone bounds.

If $c>0$, then $a>2$. If $c\le0$, monotonicity gives
$Q\le\|\partial_sQ\|_\infty/N_{\rm mid}$, so the integral in
\eqref{nc:paper-aW} is at most \(\|\partial_sQ\|_\infty\).
Thus \(W/W_1\ge\exp(-\|\partial_sQ\|_\infty)\).
By \eqref{nc:paper-shear-limit}, \(|b|\) has a positive lower bound
independent of large \(N_{\rm mid}\) in this region. Since \(a<4\) and
\(a-2\ge-2\|\partial_sQ\|_\infty/N_{\rm mid}\), one has \(a+b^2/a-2>0\)
for all sufficiently large \(N_{\rm mid}\).
Compactness gives $\kappa_*>0$ for fixed $N_{\rm mid}$. Since $c,E$ are
flat at $N_{\rm mid}$, backward integration
of \eqref{nc:paper-state-odes} gives flatness of
$W/W_1-1$, $B-\sqrt{1+\epsilon_{\rm em}^2}/\sqrt2$,
$2G-\sqrt{2(1+\epsilon_{\rm em}^2)}$, and
$K-\sqrt{2(1+\epsilon_{\rm em}^2)}$.
\end{proof}

\paragraph{Small incoming data.}
For small $B_{\rm in},D_{\rm in},K_{\rm in}$, replace
\eqref{nc:paper-control} by
\begin{equation}\label{nc:paper-incoming-control}
 e=(1-y_b/N_{\rm mid})^2e^{-y_b}
 \left\{\frac{D_{\rm in}-2B_{\rm in}}{2Q(0)}
           +\frac{y_b(A_{N_{\rm mid}}+C_{N_{\rm mid}}y_b)}{W_1}\right\},
 \qquad E=Qe,
\end{equation}
and take $G(0)=D_{\rm in}/2$, $K(0)=K_{\rm in}$.
Solve the same matrix equation for \eqref{nc:paper-terminal}.
The coefficients are affine in $B_{\rm in},D_{\rm in},K_{\rm in}$
and $\sqrt{1+\epsilon_{\rm em}^2}$. Then
\begin{equation}\label{nc:paper-incoming}
 B(0)=B_{\rm in},\qquad 2G(0)=D_{\rm in},\qquad K(0)=K_{\rm in}.
\end{equation}
For small incoming data the margins persist. Backward integration bounds
changes in $B,2G,K$, after subtracting their terminal values, by
$C(|B_{\rm in}|+|D_{\rm in}|+|K_{\rm in}|)Q$.
One $y_b$ derivative replaces $Q$ by $Q+|\partial_{y_b}Q|$.
Differentiating the coefficients gives the corresponding $Z$ bounds;
small $C_Z^2$ norms preserve the cone. Constants may depend on $N_{\rm mid}$.
The angular profile is independent of these data,
$\epsilon_{\rm em}^2$, and the incoming energy moment.

\paragraph*{The exact exterior and its flat stress edge.}
\label{nc:paper-Gamma-section}

For the exterior
\eqref{ue:heat-factor}--\eqref{ue:right-moments}, define
\begin{equation}\label{nc:paper-Gamma-profile}
 \eta(y_b)=\frac{W_1}{W(y_b)}
       =\frac{e^{\int_{y_b}^{N_{\rm mid}}\frac{Q(v/N_{\rm mid})}{1+Q(v/N_{\rm mid})}\,dv}}{1+Q(y_b/N_{\rm mid})},\qquad
 U_\Gamma^\theta=\eta U_e^\theta,\qquad
 U_\Gamma^z=\eta B U_e^\theta.
\end{equation}
At \(R_{\rm mid}\) fix the five right moments by \eqref{ue:right-moments} and integrate
\eqref{pr:moments} backward.

\paragraph{Gaussian estimates at the outer endpoint.}
Let $d_{\mathrm{aux}}=N_{\rm mid}-y_b=\log(R_{\rm mid}/R)$ and $c_{\rm aux}=\gamma N_{\rm mid}^2$. In this paragraph
write $Q$ as a function of $d_{\mathrm{aux}}$. Its endpoint condition gives
$Q(d_{\mathrm{aux}})=\rho_0(d_{\mathrm{aux}})e^{-c_{\rm aux}/d_{\mathrm{aux}}^2}$ with $\rho_0$ smooth and $\rho_0(0)>0$.
For any $\nu\in\mathbb R$ and smooth $f=f(d_{\mathrm{aux}},Z)$,
there are smooth $f_1,f_2$ up to $d_{\mathrm{aux}}=0$ such that
\begin{equation}\label{nc:paper-gaussian}
\begin{aligned}
 \partial_{d_{\mathrm{aux}}}(Qd_{\mathrm{aux}}^\nu f)&=Qd_{\mathrm{aux}}^{\nu-3}f_1,\qquad
 &\int_0^{d_{\mathrm{aux}}} Q(t)t^\nu f(t,Z)\,dt&=Q(d_{\mathrm{aux}})d_{\mathrm{aux}}^{\nu+3}f_2(d_{\mathrm{aux}},Z),\\
 |\partial_{y_b}Q|&\asymp Qd_{\mathrm{aux}}^{-3},&
 |\partial_{y_b}^i Q|&\le C_iQd_{\mathrm{aux}}^{-3i}.
\end{aligned}
\end{equation}
The derivative formula follows from
$d_{\mathrm{aux}}^3\partial_{d_{\mathrm{aux}}}\log Q=2c_{\rm aux}+d_{\mathrm{aux}}^3\partial_{d_{\mathrm{aux}}}\log\rho_0$.
The primitive formula is \eqref{fin:gaussian-smooth-quotient} with
$j=-\nu$ and amplitude $\rho_0 f$. Both formulas allow smooth auxiliary
parameters. If each mixed derivative of $h$ is bounded by
$C d_{\mathrm{aux}}^{-M}e^{-c_4/d_{\mathrm{aux}}^4}$ near zero, with $c_4>0$ and constants allowed
to depend on the derivative, then $h/(Q^\alpha d_{\mathrm{aux}}^\beta)$ extends
smoothly by zero for every fixed $\alpha,\beta\in\mathbb R$. Indeed,
each differentiated quotient is bounded by a power of $d_{\mathrm{aux}}^{-1}$ times
$\exp(-c_4/d_{\mathrm{aux}}^4+\alpha c_{\rm aux}/d_{\mathrm{aux}}^2)$, which vanishes to every order.

\begin{proposition}[Exact exterior lift]\label{nc:paper-Gamma}
Suppose $\delta$, $R_b^{-1}$, $\sup U_e^\theta$ and the incoming
$C_Z^2$ norms are sufficiently small, with thresholds depending on the
fixed curve. Then \eqref{nc:paper-Gamma-profile} defines a pair with
positive swirl, negative angular shear, and strictly admissible nonzero
stress on $[R_b,R_{\rm mid})$. It has all exterior terminal derivatives and
Gaussian flat stress at $R_{\rm mid}$.
\end{proposition}
\begin{proof}
First hold the incoming data fixed under $Z$ differentiation. The state equations give
\begin{equation}\label{nc:paper-Gamma-differences}
 \partial_{y_b}\eta=-c\eta,\qquad
 \eta B-\frac{\sqrt{1+\epsilon_{\rm em}^2}}{\sqrt2}
       =-\eta E+\frac12\int_{y_b}^{N_{\rm mid}}\eta E\,dt,
 \qquad \partial_{y_b}(\eta B)=-\eta[(1/2-c)E+\partial_{y_b}E].
\end{equation}
By monotonicity, $\eta-1$, $\eta B-\sqrt{1+\epsilon_{\rm em}^2}/\sqrt2$
and their integrals from $y_b$ to $N_{\rm mid}$ are bounded by $CQ$;
their first derivatives are bounded by $CQ(1+k)$, where
\begin{equation}\label{nc:paper-edge-weight}
 k=-2\partial_{y_b}Q/Q\ge0,\qquad
 \mathcal W=Q\left(R+\frac{k}{1+Q}\right).
\end{equation}
The product $kQ$ is smooth and bounded.

In the positive integral for the heat factor $H_e(\xi)$ in Lemma~\ref{ue:regularized-tail},
$\xi H_e'/H_e$ is the weighted mean of
$-(\delta/2)\xi v/(1+\xi v)$, so $-\delta/2<\xi H_e'/H_e\le0$ for $\delta>0$.
Since $a_e=2+\delta+2\xi H_e'/H_e$, the slope bound follows.
The differentiated expansion of $H_e(2L/R)$ in the same lemma and
$(\rho/R)^{-\delta/2}=1+O_{N_{\rm mid}}(\delta)$ give
\begin{equation}\label{nc:paper-heat-estimates}
 \begin{gathered}
 2<a_e:=1-2\partial_{y_b}\log U_e^\theta\le2+\delta,\qquad
 \|\partial_Z^m\log U_e^\theta\|_\infty\le C_m\delta/R\quad(m\ge1),\\
 U_e^\theta(\rho,Z)/U_e^\theta(R,Z)
       =(\rho/R)^{-1/2}[1+O_{C_Z^m}(\delta)]
                 \quad(R\le\rho\le e^{N_{\rm mid}}R).
 \end{gathered}
\end{equation}
Here $R\ge1$, $|Z|\le1$, and constants may depend on $N_{\rm mid}$.

Denote differences from the exterior by $\Delta$. Substituting
\eqref{nc:paper-Gamma-differences} and \eqref{nc:paper-heat-estimates}
into the backward integrals of \eqref{pr:moments} gives
\begin{equation}\label{nc:paper-moment-estimates}
 \begin{alignedat}{2}
 \|\Delta M^\theta\|_{C_Z^m}&\le C_mR\sqrt{2R}U_e^\theta Q,&
 \qquad\|\Delta M^z\|_{C_Z^m}&\le C_mRU_e^\theta Q,\\
 \|\Delta M^{\theta z}\|_{C_Z^m}&\le C_mR\sqrt{2R}(U_e^\theta)^2Q,&
 \qquad\|\Delta M^{z\theta}\|_{C_Z^m}&\le C_mR(U_e^\theta)^2Q,\\
 \|\Delta P\|_{C_Z^m}&\le C_m(U_e^\theta)^2Q.
 \end{alignedat}
\end{equation}
Write $\xi=(B_{\rm in},D_{\rm in},K_{\rm in})$. At fixed $\xi$ the two
linear moment differences gain a factor $C_m\delta/R$ under positive-order
$Z$ differentiation. When $\xi=\xi(Z)$ this gain still holds for
$\Delta M^\theta$ since $\eta$ is independent of $\xi$. For
$\partial_Z\Delta M^z$ the chain rule adds
$-\int_R^{R_{\rm mid}}\eta U_e^\theta\,\partial_\xi B\cdot\partial_Z\xi\,d\rho$.
Since $e$ depends affinely on $\xi$, differentiating
\eqref{nc:paper-Gamma-differences} with $E=Qe$ gives $|\partial_\xi B|\le C Q$.
Hence $|\partial_Z\Delta M^z|\le CRU_e^\theta Q(\delta/R+|\partial_Z\xi|)$.

Subtract the exterior identity \(\mathcal I_e+\mathcal S_e=0\) in
\eqref{pr:moment-inertia}. Since \([(1-\delta)Z+L\partial_Z]Z=L_\delta\),
the additive \(Z\) in \(M_e^z\) contributes \(1\).
Thus the part of the linear inertial coefficient independent of the
exterior velocities is $-R/L_\delta+1$. Its value at $\delta=0$ is
$-R+1$.
In the linear velocity terms, set the heat ratios and explicit
$\delta$ factors to their $\delta=0$ values. Call the resulting stress
$\mathcal L_0$. The exterior energy integral is evaluated with $\delta>0$.
For $U_{e,0}^\theta$ proportional to $R^{-1/2}$, the moment equations
and their zero terminal increments give
\begin{equation}\label{nc:paper-linear-identities}
 \begin{aligned}
 -R(\eta-1)+\frac{\Delta M_0^\theta}{\sqrt{2R}U_{e,0}^\theta}
       &=R\eta Q,\\
 -R\left(\eta B-\frac{\sqrt{1+\epsilon_{\rm em}^2}}{\sqrt2}\right)
       +\frac{\Delta M_0^z}{2U_{e,0}^\theta}&=R\eta E.
 \end{aligned}
\end{equation}
These identities follow from $\partial_{y_b}[R\eta(1+Q)]=R\eta$ and
$2\partial_{y_b}G+(1-2c)G=B$.
For $F=U_\Gamma^\theta/\sqrt{2R}$, the exact shear difference is
\[
 \begin{aligned}
 \frac{\mathcal S^\theta-\mathcal S_e^\theta}{F}
   &=-2c+a_e(\eta^{-1}-1),\\
 \frac{\mathcal S^z-\mathcal S_e^z}{F}
   &=b+(2-a_e)B
       +(a_e-1)\frac{\sqrt{1+\epsilon_{\rm em}^2}}{\sqrt2\eta}.
 \end{aligned}
\]
Here $b=2\partial_{y_b}B-(1+2c)B$ is the model shear.
Before adding this shear, the frozen linear inertial difference,
divided by $F$, is
\[
 \left(RQ+1-\eta^{-1},\quad
 RE+B-\frac{\sqrt{1+\epsilon_{\rm em}^2}}{\sqrt2\eta}\right).
\]
The state equations also give $B+b=(a-3)E-2\partial_{y_b}E$.
Adding the frozen shear therefore gives
\begin{equation}\label{nc:paper-linear-stress}
 \frac{\mathcal L_0}{F}
       =(RQ+\eta^{-1}-1-2c,\;RE+(a-3)E-2\partial_{y_b}E).
\end{equation}
In the differences between the actual and frozen shears, \(-2c\)
cancels in the angular component and \(2\partial_{y_b} B\) cancels in the
axial component. The remaining differences are
\[
 (a_e-2)(\eta^{-1}-1),\qquad
 (a_e-2)\left(
 \frac{\sqrt{1+\epsilon_{\rm em}^2}}{\sqrt2\eta}-B\right).
\]
Both are $O_{C_Z^m}(\delta Q)$ by
\eqref{nc:paper-Gamma-differences} and
\eqref{nc:paper-heat-estimates}. No derivative of $Q$ occurs here.
At fixed $\xi$, differentiated linear moments have the gain just proved.
Substituting the heat ratios of \eqref{nc:paper-heat-estimates} into the
remaining backward integrals costs $O_{C_Z^m}(\delta RQ)$.
Each remaining term has another velocity factor.
By \eqref{nc:paper-moment-estimates},
$\Delta M^z$, $\Delta M^{\theta z}/(\sqrt R U_e^\theta)$,
$\Delta M^{z\theta}/U_e^\theta$, and $R\Delta P/U_e^\theta$ are
$O_{C_Z^{m+1}}(RU_e^\theta Q)$.
For variable $\xi$ evaluate \eqref{nc:paper-linear-stress} at $\xi(Z)$.
The additional terms in \eqref{pr:moment-inertia} arise from
$\partial_\xi\Delta M^j\cdot\partial_Z\xi$ for $j=z,\theta z,z\theta$.
In the density derivative one axial factor becomes
$\eta U_e^\theta\partial_\xi B$, and the right moment is independent of $\xi$.
Thus each differentiated moment has the corresponding bound in
\eqref{nc:paper-moment-estimates}. Dividing by $F$ gives at most
$C(1+\sup U_e^\theta)RQ|\partial_Z\xi|$.
The densities are at most quadratic in the bounded inlet parameters. The
product rule therefore gives
\begin{equation}\label{nc:paper-stress-remainder}
 \| (\mathcal T-\mathcal L_0)/F\|_{C_Z^m}
                        \le C_m(\delta+\sup U_e^\theta+\|\partial_Z\xi\|_{C_Z^m})RQ.
\end{equation}

Since $\int_{y_b}^{N_{\rm mid}}\frac{Q(v/N_{\rm mid})}{1+Q(v/N_{\rm mid})}\,dv\le N_{\rm mid}Q(y_b/N_{\rm mid})$,
$|(\eta^{-1}-1)/Q|\le C$.  Dividing \eqref{nc:paper-linear-stress}
by $Q$ gives
\begin{equation}\label{nc:paper-edge-rescaling}
 \begin{aligned}
 \frac{\mathcal L_0^\theta}{FQ}
   &=R+\frac{k}{1+Q}+\frac{\eta^{-1}-1}{Q}-\frac2{1+Q},\\
 \frac{\mathcal L_0^z}{FQ}
   &=\left(R+\frac{k}{1+Q}\right)e
                      +\frac{kQ}{1+Q}e+(a-3)e-2\partial_{y_b}e.
 \end{aligned}
\end{equation}
The remaining terms are bounded. Hence
\begin{equation}\label{nc:paper-relative-stress}
 \begin{gathered}
 \frac{\mathcal T}{F\mathcal W}=(1,e)+\mathcal R,\qquad
 \|\mathcal R\|_{C_Z^m}\le C_m(\delta+\sup U_e^\theta+R_b^{-1}+\|\partial_Z\xi\|_{C_Z^m}),\\
 a_\delta=a+(a_e-2),\qquad b_\delta=b-(a_e-2)B.
 \end{gathered}
\end{equation}
For normalized stress $(t_\theta,t_z)$, require positivity of
\[
 a,\quad a^2+b^2-2a,\quad at_\theta-bt_z,\quad
 2a(at_\theta-bt_z)^2-(a^2+b^2-2a)(bt_\theta+at_z)^2.
\]
At $(t_\theta,t_z)=(1,e)$ the last expression equals
$a(a^2+b^2)[2-b^2/a-2be-(a-2)e^2]$.
The margins in Lemma~\ref{nc:paper-model-cone} persist under
\eqref{nc:paper-relative-stress}. Thus
\begin{equation}\label{nc:paper-stress-lower}
 |\mathcal T|\ge\tfrac12F\mathcal W>0\quad(R<R_{\rm mid}).
\end{equation}
The constant at order $m$ may depend on bounds for $\xi$ through order
$m+1$. The strict cone requires only the zero-order estimate.
At the endpoint, \eqref{nc:paper-Gamma-profile},
\eqref{nc:paper-Gamma-differences} and the primitive rule
\eqref{nc:paper-gaussian} give smooth quotients by $Q$ for
$\eta-1$, $\eta B-\sqrt{1+\epsilon_{\rm em}^2}/\sqrt2$ and
$\partial_\xi B$. Each moment density difference is therefore $Q$
times a smooth coefficient. Its zero-terminal primitive, and likewise
$\partial_\xi\Delta M^j$ for $j=z,\theta z,z\theta$, have the factor
$Qd_{\mathrm{aux}}^3$. In particular all exterior terminal derivatives agree.
The lower bound in \eqref{nc:paper-stress-lower} and
\eqref{nc:paper-relative-stress} give uniform relative cone margins
up to this auxiliary endpoint.
\end{proof}

\paragraph{Exact entrance derivatives.}
Set $U_b^\theta=U_\Gamma^\theta(R_b,0)$ and
$U_\ell^\theta=U_b^\theta e^{-y_b/2}/W(y_b)$,
$U_\ell^z=BU_\ell^\theta$.  Define
\begin{equation}\label{nc:paper-entrance-weight}
 \Psi(s)=\frac{e^{-(1+\gamma)/(1-s)^2}}
       {e^{-(1+\gamma)/s^2}+e^{-(1+\gamma)/(1-s)^2}},\qquad
 \Psi(0)=1,\quad\Psi(1)=0,
\end{equation}
and define
\begin{equation}\label{nc:paper-entrance-lift}
 \widetilde U^\theta=U_\Gamma^\theta
       +\Psi^2(U_\ell^\theta-U_\Gamma^\theta),\qquad
 \widetilde U^z=U_\Gamma^z+\Psi^2(U_\ell^z-U_\Gamma^z).
\end{equation}
Since \(1-\Psi^2\) is flat at \(s=0\) and \(\Psi^2\) is flat at
\(s=1\), the pair has the model entrance jets and the exterior
terminal jets. It is analytic in \(R\) on the open interval. By \eqref{nc:paper-heat-estimates},
$\|U_\ell^{\theta,z}-U_\Gamma^{\theta,z}\|_{C^m_{s,Z}}
\le C_m\delta U_b^\theta$. Since $2(1+\gamma)>\gamma$, the derivatives
of $\Psi^2$ are $O(Q)$ at the outer edge. Thus
\eqref{nc:paper-relative-stress} changes by $O_{C_Z^m}(\delta)$,
preserving the relative cone margins and all exterior terminal jets.

\paragraph{The five entrance targets.}
For \eqref{nc:paper-entrance-lift}, integrate backward from the fixed
right moments. Equations \eqref{nc:paper-heat-estimates} and
\eqref{ue:right-moments} give
\begin{equation}\label{nc:paper-entrance-targets}
 \begin{aligned}
 \frac{M^\theta(R_b)}{R_b\sqrt{2R_b}U_b^\theta}
     &=1+Q(0)+O(\delta+R_b^{-1}),\\
 \frac{M^z(R_b)}{R_bU_b^\theta}
     &=D_{\rm in}+\frac{Z}{R_bU_b^\theta}+O(\delta),\\
 \frac{M^{\theta z}(R_b)}{R_b\sqrt{2R_b}(U_b^\theta)^2}
     &=K_{\rm in}+O(\delta).
 \end{aligned}
\end{equation}
The contribution of the additive $Z$ in $M_e^z$ equals the second term.
For moment correction,
take $\delta+R_b^{-1}+(R_bU_b^\theta)^{-1}$ sufficiently small.
The exact pressure and energy targets are
\begin{equation}\label{nc:paper-quadratic-targets}
 \begin{aligned}
 P(R_b)&=P_e(R_{\rm mid})-\int_{R_b}^{R_{\rm mid}}\frac{(\widetilde U^\theta)^2}{2R}\,dR,
       &M^p(R_b)&=P(R_b)-P_0,\\
 M^{z\theta}(R_b)&=-\epsilon_{\rm em}^2\mathcal J^\theta(R_{\rm mid},Z)
       -\int_{R_b}^{R_{\rm mid}}[(\widetilde U^z)^2-(\widetilde U^\theta)^2/2]\,dR.
 \end{aligned}
\end{equation}
The angular velocity and pressure target are independent of
$\epsilon_{\rm em}^2$ and the incoming energy moment.  The axial
velocity depends only on the three incoming data and
$\sqrt{1+\epsilon_{\rm em}^2}$.

\paragraph*{A uniform five-moment window.}
\label{nc:paper-window-section}

\begin{lemma}[Five-moment correction on a fixed window]\label{nc:paper-five-window}
Let $J\Subset(0,\infty)$ be a fixed interval with nonempty interior. For all sufficiently large $N_{\rm mid}>\sup J$, five smooth variations supported in $J$ have a normalized moment matrix with uniformly bounded inverse. This bound persists for small incoming data, small $\delta$, and \eqref{nc:paper-entrance-lift}. With the right moments fixed we can correct small normalized left moment errors. With nearby left moments prescribed we can recover the original right moments. Both corrections preserve the endpoint jets of $U^\theta,U^z$ and the strict cone.
\end{lemma}
\begin{proof}
Use the moment units
\begin{equation}\label{nc:paper-moment-units}
 R_b\sqrt{2R_b}U_b^\theta,\quad R_bU_b^\theta,\quad
 R_b\sqrt{2R_b}(U_b^\theta)^2,\quad
 R_b(U_b^\theta)^2,\quad (U_b^\theta)^2
\end{equation}
for $M^\theta,M^z,M^{\theta z},M^{z\theta},M^p$.
For zero incoming data, $U^\theta/U_b^\theta\to e^{-(1/2+r_0)y_b}$ and
$U^z/U_b^\theta\to0$ on $J$, with convergence of all derivatives.
For variations of size $U_b^\theta$, the three limiting angular weights are
\begin{equation}\label{nc:paper-window-weights}
 e^{3y_b/2},\quad -e^{(1/2-r_0)y_b},\quad e^{-(1/2+r_0)y_b},
\end{equation}
and the two axial rows have weights $e^{y_b},e^{(1-r_0)y_b}$.
The other blocks vanish, by \eqref{pr:moments} and $dR=R_be^{y_b}\,dy_b$.
Use a common smooth cutoff supported in $J$ and positive on a subinterval.
The exponents within each block are distinct, so its rows are independent.
Lemma~\ref{pr:gram} gives a uniform inverse for large $N_{\rm mid}$ and then for
small incoming data and the $O(\delta)$ heat and entrance corrections.

In the units \eqref{nc:paper-moment-units}, the five increments satisfy
$A\alpha+\mathcal Q(\alpha,\alpha)=m_{\rm def}$ as in \eqref{rp:quadratic-moment-map}.
With $R_b$ fixed as the initial endpoint, $m_{\rm def}$ is the prescribed
minus the computed moment vector at $R_{\rm mid}$. With the moments at
$R_{\rm mid}$ fixed, it is the computed minus the prescribed vector at $R_b$.
Lemma~\ref{pr:gram} and the product rule bound the inverse and $\mathcal Q$
in $C_Z^m$. Lemma~\ref{pr:quadratic} then gives
$\|\alpha\|_{C_Z^m}\le C_m\|m_{\rm def}\|_{C_Z^m}$.

In these units, each partial moment integral and the velocity change with one $y_b$ derivative are bounded by $C\|m_{\rm def}\|_{C_Z^m}$. By \eqref{pr:moment-inertia}, the changes in $(a,b)$ and $\mathcal T/(RF)$ on compact $y_b$ intervals are $O(\|m_{\rm def}\|_{C_Z^2})$, uniformly for $R_b\ge1$ and $0<U_b^\theta\le1$ under the stated smallness conditions. Hence the strict margins persist.

Both corrections preserve the tangential endpoint jets and the full reconstruction to the right of $J$, including the flat edge. Backward integration changes the left moments by $-m_{\rm def}$ and may change $U^r,P,\mathcal I,\mathcal T$ there. Forward integration preserves the prescribed left reconstruction.
\end{proof}
\subsubsection{Strict admissibility and analytic matching of the five moments}
\label{rr:section}

\paragraph*{A strict shear loop with the prescribed mean.}
\label{rr:loop-section}

We vary the shear along a strictly admissible periodic curve. Preserving
its mean permits a velocity correction of order \(N^{-1}\) at radial frequency \(N\).

Let $\mathcal K=[R_1,R_2]\times[-1,1]$ be a compact interior
rectangle, with moments and inertia given by
\eqref{pr:moments}--\eqref{pr:moment-inertia}.

\begin{lemma}\label{rr:loop}
Suppose $F>0$, $a>0$, and the pair is relaxed admissible on
$\mathcal K$ and strictly admissible near its radial boundaries.
Then there is a smooth one-periodic loop $(a_L,b_L)$ with
\begin{equation}\label{rr:loop-mean}
 \int_0^1(a_L,b_L)\,d\varphi=(a,b).
\end{equation}
With $\mathcal I$ fixed in phase, $F(-a_L,b_L)$ is strictly admissible
at every phase and equals the original shear near the radial boundaries.
\end{lemma}
\begin{proof}
Put
\[
 t_0=-b/a,\qquad
 H(t)=\frac{\mathcal I^\theta+t\mathcal I^z}{F},\qquad
 J(t)=\frac{\mathcal I^z-t\mathcal I^\theta}{F}.
\]
For $\mathcal S'=-Fa'(1,t)$ and $v=a'(1+t^2)$, one has
\begin{equation}\label{rr:frozen-cone}
\begin{aligned}
 - (\mathcal I+\mathcal S')\cdot\mathcal S'
       &=F^2a'[H(t)-v],\\
 2\bigl((\mathcal I+\mathcal S')\cdot\mathcal S'\bigr)^2
 -(v-2)\bigl((\mathcal I+\mathcal S')\cdot
                        (\mathcal S')^\perp\bigr)^2
       &=F^4(a')^2\{2[H(t)-v]^2-(v-2)J(t)^2\}.
\end{aligned}
\end{equation}
If $\kappa\le2$, \eqref{pr:relaxed} gives $H(t_0)>2$; otherwise
the direction condition gives $H(t_0)>\kappa>2$. Thus
$m:=\min_{\mathcal K}(H(t_0)-2)>0$.

Use the smooth step \(\sigma_{\rm norm}(s)\) from \eqref{pr:cutoffs}, for which
\(\sigma_{\rm norm}(1/2)=1/2\).
Take $0<\eta\le\min(1/2,m/8)$, which will be reduced later, and set
\begin{equation}\label{rr:loop-amplitude}
 q=\begin{cases}
 \displaystyle\sigma_{\rm norm}\!\left(\frac{2+\eta-\kappa}{\eta}\right)
       \sqrt{\frac{2+2\eta-\kappa}{2a}},&\kappa<2+\eta,\\
 0,&\kappa\ge2+\eta,
 \end{cases}
 \qquad v=\kappa+2aq^2.
\end{equation}
Take $\eta$ smaller than half the minimum of $\kappa-2$ on the
radial boundaries. Then $q$ is smooth and vanishes near them.
Moreover
\[
 q\le\sqrt{\frac3{2\min_{\mathcal K}a}},\qquad
 v>2+\eta/4,
 \qquad v\le2+2\eta\quad\text{where }q\ne0.
\]
If $\kappa\le2$, then $v=2+2\eta$. If
$2\le\kappa\le2+\eta/2$, use $\sigma_{\rm norm}\ge1/2$; otherwise use $v\ge\kappa$.

For $0\le\psi\le2\pi$, define
\begin{equation}\label{rr:poisson-loop}
 \varrho=\frac{q\mathcal I^z/F}
             {\sqrt{(m/4)^2+(q\mathcal I^z/F)^2}},\qquad
 t(\psi)=t_0+2q\sqrt{1-\varrho^2}
       \frac{\cos\psi-\varrho}{1-2\varrho\cos\psi+\varrho^2}.
\end{equation}
The formulas are smooth at $q\mathcal I^z=0$.
The Poisson kernel satisfies
\[
 P_\varrho(\psi)
 =\frac{1-\varrho^2}{1-2\varrho\cos\psi+\varrho^2}
 =1+2\sum_{n\ge1}\varrho^n\cos(n\psi),\qquad
 \left\langle(P_\varrho-1)^2\right\rangle
       =\frac{2\varrho^2}{1-\varrho^2},
\]
where brackets denote the mean over $[0,2\pi]$.
The Poisson kernel identity and
\((\mathcal I^z/F)(t-t_0)=(m/4)(P_\varrho-1)\) give
\begin{equation}\label{rr:poisson-averages}
 \langle t\rangle=t_0,\qquad
 \langle(t-t_0)^2\rangle=2q^2,\qquad
 H(t)\ge2+3m/4.
\end{equation}
At \(\mathcal I^z=0\), \(t=t_0+2q\cos\psi\) gives the same identities.
Since \(|t-t_0|\le 2q+16|\mathcal I^z|q^2/(mF)\), the bounds for \(t,J(t)\)
are independent of \(\eta\). Reduce \(\eta\) until \(2\eta\sup|J(t)|^2\le m^2/4\).

Define the increasing phase and the loop by
\begin{equation}\label{rr:phase}
 \varphi(\psi)=\frac{a}{2\pi v}
        \int_0^\psi(1+t(\vartheta)^2)\,d\vartheta,
 \qquad
 (a_L,-b_L)(\varphi)
   =\frac{v}{1+t(\psi(\varphi))^2}(1,t(\psi(\varphi))).
\end{equation}
By \eqref{rr:poisson-averages}, $a\langle1+t^2\rangle=v$.
Thus the phase increases by one per period, and changing variables
gives \eqref{rr:loop-mean}. The instability parameter is $v>2$.
For $q\ne0$,
\[
 H(t)-v\ge m/2,\qquad
 2[H(t)-v]^2-(v-2)J(t)^2\ge m^2/4.
\]
Hence \eqref{rr:frozen-cone} gives strict admissibility for $q\ne0$;
for $q=0$ the shear is unchanged. Compactness gives positive margins.
\end{proof}

\paragraph*{A finite-frequency velocity correction.}
\label{rr:finite-frequency-section}
Integrating the oscillating shear gives the following correction to the velocities.

\begin{lemma}\label{rr:finite-frequency}
Let $(a_L,b_L)$ be the loop in Lemma~\ref{rr:loop}.
Then there is a velocity pair equal to the given pair near
the radial boundaries, with $O(N^{-1})$ velocity and cumulative
moment errors in every fixed $C_Z^m$ norm, uniformly in $R$.
Its normalized shear differs by $O(N^{-1})$ from the loop at
phase $N\log(R/R_1)$.
\end{lemma}
\begin{proof}
Let $A_L,B_L$ be the zero-mean periodic primitives
\[
 \partial_\varphi A_L=-\tfrac12(a_L-a),\qquad
 \partial_\varphi B_L=\tfrac12U^\theta(b_L-b).
\]
By \eqref{rr:loop-mean}, they exist and vanish near the radial
boundaries. Extend the primitives by zero and set
\begin{equation}\label{rr:finite-velocity}
 U_N^\theta=U^\theta e^{A_L(R,Z,N\log(R/R_1))/N},\qquad
 U_N^z=U^z+B_L(R,Z,N\log(R/R_1))/N.
\end{equation}
Holding the phase fixed in \(\partial_R A_L\) and \(\partial_R B_L\) gives
\begin{equation}\label{rp:finite-shear}
 a_N=a_L-2R\partial_R A_L/N,\qquad
 b_N=e^{-A_L/N}
           \left(b_L+\frac{2R\partial_R B_L}{NU^\theta}\right).
\end{equation}
Since the phase is independent of $Z$, these identities give the
velocity and shear bounds. The first radial derivatives of the velocities
are bounded in $C_Z^m$ independently of $N$. Integrating the density differences
in \eqref{pr:moments} from the same incoming values gives
\begin{equation}\label{rr:partial-moment-bound}
 \sup_{R_1\le R\le R_2}
   \|M_N^j(R)-M^j(R)\|_{C_Z^m}\le C_m/N,
 \qquad j=\theta,z,\theta z,z\theta,p.
\end{equation}
\end{proof}

\paragraph{Matching the power at \(R_{\rm mid}\).}
The auxiliary profile ends at the heat exterior. We now change its endpoint
derivatives to those of the power and correct the resulting five moment errors
on the earlier window \(J\).
\begin{lemma}[Matching a left profile to the exact power]\label{pow:insertion}
Suppose a smooth auxiliary pair is strictly admissible on
\([R_{ac},R_{\rm mid})\), has the exterior moments at \(R_{\rm mid}\),
and continues there as the exact exterior. On the final left collar
assume the pair and its backward moments are those of
\eqref{nc:paper-entrance-lift} and Proposition~\ref{nc:paper-Gamma}.
Assume an earlier strict window \(J\) has the five variations of
Lemma~\ref{nc:paper-five-window}. For some fixed small \(\ell,\epsilon_0>0\)
and all sufficiently small \(\delta>0\) in the parameter family of
Lemma~\ref{np:preparation}, the left profile can be changed to match the
power jets and the five targets \eqref{pow:inlet-moments} at \(R_{\rm mid}\).
The inner data and axis pressure are preserved, and the stress is strictly
admissible throughout \([R_{ac},R_e)\), with a positive margin on
\([R_{ac},R_{\rm mid}]\).
\end{lemma}
\begin{proof}
Let \(\chi\) rise smoothly from zero to one on
\([e^{-2\ell}R_{\rm mid},e^{-\ell}R_{\rm mid}]\) and remain constant
on either side. Subtract \(\chi(U_e^\theta-V_p)(1,U_e^z/U_e^\theta)\) from the auxiliary velocities.
Since the auxiliary pair has the exterior jets at \(R_{\rm mid}\), the new
pair has the jets of the middle power in \eqref{pow:profiles}. Integrate the new moments backward from
\eqref{pow:inlet-moments}. Between \(J\) and the start of the switch the
moment change is constant in \(R\). By \eqref{pr:moments} it is
\(O_\ell(\epsilon_0)\) in the fixed moment units, including two \(Z\)
derivatives. The forward correction in Lemma~\ref{nc:paper-five-window}
realizes it on \(J\) and preserves the data to the left. From \(J\) to the
end of the switch the changes of velocity in \(C^1\) and of moments in
\(C_Z^1\) are \(O_\ell(\epsilon_0)\). For fixed \(\ell\), compactness and
\eqref{nc:paper-relative-stress} give a uniform strict margin there.
The moment units and inverse on \(J\) stay bounded at \(\delta=0\)
by continuation in \(\epsilon_{\rm em}^2/\delta\)
(Lemmas~\ref{np:preparation} and \ref{ue:regularized-tail}).

On \([e^{-\ell}R_{\rm mid},R_{\rm mid}]\) we have \(\chi=1\).
Apply the calculation of \eqref{pow:stress-comparison} to the difference from the
auxiliary pair, replacing \(M_e^z\) in \(b_0\) by its auxiliary value.
The auxiliary velocities and their first \(Z\) derivatives are
\(O(U_e^\theta)\) by smooth dependence and matching of the exterior jets.
Thus the same five moment bounds apply. By continuity from \(R_{\rm mid}\),
\(b_0>1/2\) on this collar for small \(\ell\). The stress increment again
lies strictly inside the cone of the new shear. The old nonzero stress
keeps its relative margin by \eqref{nc:paper-relative-stress} and
closeness of the two shears. For the fixed new shear the two inequalities
\[
 -\mathcal T\cdot\mathcal S
 \ \pm\sqrt{(\kappa-2)/2}\,\mathcal T\cdot\mathcal S^\perp>0
\]
are linear in \(\mathcal T\). Thus the sum of the old stress and the increment
is strictly admissible, also at \(R_{\rm mid}\), where the old stress vanishes.
\end{proof}

\paragraph*{Radial analytic approximation with all endpoint derivatives.}
\label{rr:analytic-section}

Define on the open interval
\begin{equation}\label{rr:rho}
 \rho(R)=\exp\left(-\frac1{(R-R_-)^4}
                         -\frac1{(R_+-R)^4}\right),
\end{equation}
and set $\rho(R_-)=\rho(R_+)=0$.

\begin{lemma}\label{rr:analytic}
Let $U$ be smooth on $[R_-,R_+]\times[-1,1]$. For every sufficiently small
$\varepsilon>0$ there is a smooth $\widetilde U$, real analytic in $R$ on
$(R_-,R_+)$, with the same endpoint derivatives and, for each $i,k\ge0$,
\begin{equation}\label{rr:analytic-error}
 |\partial_R^i\partial_Z^k(\widetilde U-U)|\le C_{ik}\varepsilon\sqrt\rho.
\end{equation}
\end{lemma}
For the variable-width Gaussian formula and its proof, see Appendix~\ref{app:tools}.

\paragraph*{Analytic preparation and the moment matrices.}
\label{rr:gram-section}\label{rr:realization-section}
Fix the axial primitive \(M^z\), so \(U^z=\partial_R M^z\) is known.
On a fixed radial interval, write \([M]\) for the terminal value of a moment
minus its incoming value. With \(\Gamma=\sqrt{2R}\,U^\theta\),
\eqref{pr:moments} gives two linear and two quadratic moments,
\[
 \int\!\left(\Gamma,U^z\Gamma,\frac{\Gamma^2}{R},\frac{\Gamma^2}{R^2}\right)dR
 =\left([M^\theta],[M^{\theta z}],
        4\int(U^z)^2\,dR-4[M^{z\theta}],4[M^p]\right).
\]
The power jets give a uniform inverse for these four angular moments.

\begin{lemma}[Rank on interior windows]\label{rp:analytic-local-rank}
Let a smooth pair be radially analytic in the open interval. At one endpoint,
suppose each component has the jets of a positive constant multiple of \(R^{-p}\), with the same \(p>0\).
The four nonzero angular entries of
\eqref{rr:rows}, in the order \(\theta,\theta z,z\theta,p\), are
\begin{equation}\label{rp:four-angular-rows}
 \sqrt{2R},\qquad \sqrt{2R}\,U^z,\qquad -U^\theta,\qquad U^\theta/R .
\end{equation}
On each fixed interior window, these four functions and the five rows
\eqref{rr:rows} each have a Gram matrix with uniformly positive least
eigenvalue in \(Z\), for any smooth nonnegative weight positive on an open subinterval.
\end{lemma}
\begin{proof}
A constant relation on that subinterval extends by radial analyticity to the
whole piece. At the power endpoint the four angular entries have the jets
of nonzero multiples of \(R^{1/2},R^{1/2-p},R^{-p},R^{-p-1}\).
Their derivatives of orders zero through three in \(\log R\) give an
invertible Vandermonde matrix. In a relation among the five rows, the
angular component therefore leaves only the axial row \((0,1)\), whose
coefficient must also vanish. Compactness in \(Z\) gives the uniform bounds.
\end{proof}

\begin{lemma}[Analytic preparation]\label{rr:realization}
Let a smooth pair on \([R_{ac},R_{\rm mid}]\times[-1,1]\) have positive \(F\),
negative angular shear, and nonzero strictly admissible stress throughout
the closed rectangle. At \(R=R_{\rm mid}\), suppose each component has the jets of
a positive constant multiple of \(R^{-p}\), with the same \(p>0\).
There is a pair, radially analytic on \((R_{ac},R_{\rm mid})\), with the same endpoint jets and
the same five actual moments at \(R_{\rm mid}\). Its velocities are arbitrarily
close in every prescribed finite \(C^m\) norm and its stress remains
strictly admissible on the closed rectangle.
\end{lemma}
\begin{proof}
Apply Lemma~\ref{rr:analytic} to \(U^\theta\) and \(M^z\). The endpoint values
of \(\widetilde M^z\) fix the axial integral exactly. With one extra radial
derivative, \eqref{rr:analytic-error} gives
\(\|\partial_R\widetilde M^z-U^z\|_{C^m}=O_m(\varepsilon)\) and equality of all endpoint jets.

Evaluate the four entries \(\omega_j^\theta\) in \eqref{rp:four-angular-rows}
at \((\widetilde U^\theta,\partial_R\widetilde M^z)\) and set
\begin{equation}\label{rr:single-formula}
 U_*^z=\partial_R\widetilde M^z,\qquad
 U_*^\theta=\widetilde U^\theta+\rho\sum_j c_j(Z)\omega_j^\theta,\qquad
 G_{ij}=\int_{R_{ac}}^{R_{\rm mid}}\rho\,\omega_i^\theta \omega_j^\theta\,dR .
\end{equation}
Before smoothing, \(\rho>0\) throughout the interior, so a Gram null vector
would give a relation up to the power endpoint. The Vandermonde calculation in
Lemma~\ref{rp:analytic-local-rank} excludes it. Compactness in \(Z\) gives a
uniform inverse, which persists for small \(\varepsilon\).
With the incoming moments fixed, \eqref{pr:moments} gives the change
\(Gc+\mathcal Q(c,c)\), quadratic only in the energy and pressure rows.
Lemma~\ref{pr:quadratic} removes the \(O_m(\varepsilon)\) errors in
\(C_Z^m\) with \(c=O_m(\varepsilon)\). The factor \(\rho\) preserves the jets
and radial analyticity. The \(C^1\) bound in \eqref{pr:moment-inertia}
preserves the strict margins, and \eqref{pr:moments} gives the reconstructed jets.
\end{proof}

\begin{theorem}[Connection through an exact mixed power]
\label{thm:connection}\label{cl:leading}
\label{sc:base-theorem}\label{sc:base-leading-theorem}
Choose \(Q,\gamma,N_{\rm mid}\) as in Section~\ref{nc:paper-section} and fix
\(\mu P_*^4>0\) sufficiently small. For sufficiently large \(P_*\) one can
choose the inner parameters, \(0<\delta<1/200\), \(0<\ell,\epsilon_0<1\),
\(h_\infty,\epsilon_{\rm em}>0\) and \(R_b>5\) so that
\eqref{rr:single-formula} and \eqref{pow:profiles} give the leading background
of Theorem~\ref{intro:background}\textup{(i)}, with
\(R_{\rm mid}=e^{N_{\rm mid}}R_b\) and \(R_e=e^{2\ell}R_{\rm mid}\).
Its five moments at \(R_e\) are \eqref{ue:right-moments}.
The leading stress satisfies \eqref{intro:stress-bounds} and
\eqref{intro:uniform-cone} with the subscript \(B\) omitted, and \(F,-\mathcal S^\theta\)
have positive uniform lower and upper bounds on \([R_{\rm in},R_e]\times[-1,1]\).
For \(d=\log(R_e/R)\) near the outer edge,
\begin{equation}\label{conn:outer-bounds}
 |\mathcal T|\ge c\,d^{-3}e^{-c_h/d^2},\qquad
 |(R\partial_R)^i\partial_Z^k\mathcal T|
     \le C_{ik}d^{-3i-3}e^{-c_h/d^2},
 \qquad c_h=\ell^2>0.
\end{equation}

On each fixed nonempty open window in $(R_{ac},R_e)$, the
five-moment and four-angular-moment Gram matrices are positive
definite, with inverse bounds uniform in \(Z\in[-1,1]\) and independent
of the physical scale.
\end{theorem}

\begin{proof}
Fix the tolerance in Lemma~\ref{nc:paper-five-window}, then choose \(\Lambda,C_*\)
in the order of Lemma~\ref{np:preparation}. Apply Lemma~\ref{en:passage} with
\((R_1,R_2,a_*)=(R_{ac},5,4/5)\) and a sufficiently small \(C_Z^4\) perturbation.
Lemma~\ref{np:preparation} gives the prescribed pressure and energy
and makes the three remaining errors arbitrarily small in \(C_Z^3\).

Choose \(\ell,\epsilon_0,\delta\) in that order by Lemma~\ref{pow:insertion}.
Lemmas~\ref{rr:loop} and \ref{rr:finite-frequency} preserve the boundary shears.
Take \(N\) large enough that the \(O(N^{-1})\) errors in shear, inertia and
moments fit the strict cone margin and the units \eqref{nc:paper-moment-units}.
Lemma~\ref{nc:paper-five-window} restores the five moments on \(J\), so the
data to its right agree with Proposition~\ref{nc:paper-Gamma}.
Lemma~\ref{pow:insertion} now matches the power jets and \eqref{pow:inlet-moments}
at \(R_{\rm mid}\) with a strict margin up to that endpoint.
Lemma~\ref{rr:realization} preserves these data when making the left profile
analytic. Moment matching (Lemma~\ref{pr:matching}) then gives the exterior reconstruction.

The edge estimates \eqref{ic:actual-angle}, \eqref{ic:inner-flat-size} and
\eqref{pow:edge-data} give the relative margins at both edges and
\eqref{conn:outer-bounds} at \(R_e\). Compactness gives the interior margins
and the bounds for \(F,-\mathcal S^\theta\). Every open window contains a
subinterval of one analytic piece. Each piece has the power jets at one
endpoint, so Lemma~\ref{rp:analytic-local-rank} gives both rank assertions.
\end{proof}

\subsection{Exact physical fields from finite corrections}
\label{fin:section}

We start with the five matched leading moments of Theorem~\ref{cl:leading}
and insert the exact core by a streamfunction. Four angular variations
correct the remaining physical moments, and radial integration gives
the pressure and stress. Choose two open intervals
\begin{equation}\label{fin:reserved-windows}
 \overline I_2\cup\overline I_3\subset(R_{ac},R_{\rm mid}),
 \qquad \sup I_2<\inf I_3.
\end{equation}
We correct the physical moments on \(I_2\) and reserve \(I_3\) for the pulses.
On their closures the Gaussian weight is bounded below. Fix the inverse bounds
from Lemma~\ref{rp:analytic-local-rank} and the cone margins before choosing $\lambda_*$.

\paragraph{Inserting the exact core}
\label{fin:seed-section}

Let $(u_c,p_c)$ be the exact core of Section~\ref{ex:appendix}, and
let $\bar{\lambda}=\lambda^{2\delta}$. Treat $(\lambda,R,Z)$ as independent
variables. For $a\in\mathbb R$, define
\begin{equation}\label{fin:operators}
\begin{aligned}
 \mathcal T_aG&=\frac1{1-\delta Z^2}
 \left[-\frac12(a+\lambda\partial_\lambda)G
       +\frac{1-\delta}{2}Z\partial_ZG+R\partial_RG\right],\\
 \mathcal Z_aG&=\frac1{1-\delta Z^2}
 \left[Z(a+\lambda\partial_\lambda)G
       +(1-Z^2)\partial_ZG-2ZR\partial_RG\right].
\end{aligned}
\end{equation}
By Section~\ref{pr:section},
\begin{equation}\label{fin:chain-rules}
\begin{aligned}
 \partial_t(\lambda^aG)&=\lambda^{a-2}\mathcal T_aG, &
 \qquad\partial_z(\lambda^aG)&=\lambda^{a-1+\delta}\mathcal Z_aG,\\
 \partial_z^2(\lambda^aG)
   &=\lambda^{a-2+2\delta}\mathcal Z_{a-1+\delta}\mathcal Z_aG, &
 \qquad\partial_r(\lambda^aG)&=\lambda^{a-1}\sqrt{2R}\,\partial_RG.
\end{aligned}
\end{equation}
On functions of $\bar{\lambda}$, $\lambda\partial_\lambda=2\delta\bar{\lambda}\partial_{\bar{\lambda}}$.
On a fixed closed annulus in $R>0$, $\partial_R$ and $R\partial_R$ give
equivalent fixed-order bounds, with constants depending on the annulus.
Write $O_m(\bar{\lambda})$ for a bound $C_m\bar{\lambda}$ after at most $m$ compositions of
$\lambda\partial_\lambda,R\partial_R,\partial_Z$.

By \eqref{ex:scale-remainder}, each normalized core increment in
\eqref{ex:normalized}, divided by $\bar{\lambda}$, extends holomorphically
to $\bar{\lambda}=0$ with bounded derivatives on smaller neighborhoods.
Write $(u_{c,(0)},p_{c,(0)})$
for the corresponding leading core on $R\le4.1/\Lambda$.

Define $\Psi=r^{-1}\int_0^r s u^z(t,s,z)\,ds$. Then
$\partial_r(r\Psi)=ru^z$, and incompressibility gives $u^r=-\partial_z\Psi$.
Choose $R_{\rm in}<R_c<R_{ac}$ and a smooth $\chi_c(R)$ equal to one near
$[0,R_{\rm in}]$ and supported in $[0,R_c)$. Choose smooth $\chi_h(R)$ equal to zero near
the core and one on the terminal collar and exterior, with transition
disjoint from $I_2\cup I_3$. Define
\begin{equation}\label{fin:seed}
\begin{aligned}
 \Psi_*&=\Psi_{(0)}+\chi_c(\Psi_c-\Psi_{c,(0)}),\\
 u_*^\theta&=u_{(0)}^\theta
                 +\chi_c(u_c^\theta-u_{c,(0)}^\theta),\qquad
 p_*=p_{(0)}+\chi_c(p_c-p_{c,(0)})-\frac{\chi_h(R)}{2r^2},\\
 u_*^z&=(\partial_r+r^{-1})\Psi_*,\qquad
 u_*^r=-\partial_z\Psi_*.
\end{aligned}
\end{equation}
Since $r\Psi=0$ at the axis, integrating \eqref{fin:seed} gives
\begin{equation}\label{fin:axial-primitive}
 \int_0^r s(u_*^z-u_{(0)}^z)\,ds
     =r\chi_c(\Psi_c-\Psi_{c,(0)}).
\end{equation}
Thus the axial primitive is unchanged for $R\ge R_c$. The product rule gives
\begin{equation}\label{fin:seed-product-rule}
 u_*^z-u_{(0)}^z=\chi_c(u_c^z-u_{c,(0)}^z)
                    +(\partial_r\chi_c)(\Psi_c-\Psi_{c,(0)}),
 \qquad
 u_*^r-u_{(0)}^r=\chi_c(u_c^r-u_{c,(0)}^r)
                    -(\partial_z\chi_c)(\Psi_c-\Psi_{c,(0)}).
\end{equation}
The streamfunction gives $\operatorname{div}u_*=0$. The fields agree with
the exact core on $R\le R_{\rm in}$ and with \eqref{ue:physical} on $R\ge R_e$.
By \eqref{ex:remainder-factor}, the normalized velocity increments are
$O_m(\bar{\lambda})$. The added normalized pressure $-\bar{\lambda}\chi_h/(4R)$
is also $O_m(\bar{\lambda})$ on fixed profile ranges and zero near the axis.

Let \(\mathsf T_{(0)}=\lambda^{-2-\delta}\mathcal T_{(0)}\) be the
physical leading stress and put
\(E_* =\NS(u_*,p_*)+\mathcal D\mathsf T_{(0)}\).
The leading identities \eqref{pr:leading-residual}, the exact core, and
\eqref{fin:chain-rules} give, on the fixed stress annulus,
\begin{equation}\label{fin:seed-residual}
 \lambda^{3+\delta}E_*^\theta=O_m(\bar{\lambda}),\qquad
 \lambda^{3+\delta}E_*^z=O_m(\bar{\lambda}),\qquad
 \lambda^{3+2\delta}E_*^r=O_m(\bar{\lambda}).
\end{equation}
Indeed, the normalized axial viscosity is $\lambda^{2\delta}$;
the radial inertial and viscous terms are $O(\lambda^{-3})$, hence
$O(\bar{\lambda})$ in the unit $\lambda^{-3-2\delta}$. The polynomial
Navier--Stokes operator preserves these bounds under $O_m(\bar{\lambda})$
increments, since $1-\delta Z^2\ge1-\delta$. The residual is zero on
$R\le R_{\rm in}$ and $R\ge R_e$.

\paragraph{Physical moments.}
\label{fin:moments-section}

Let $(u,p)$ be smooth and axisymmetric, with $\operatorname{div}u=0$, and equal the mixed exterior for large $r$. At fixed $(t,z)$, define
\begin{equation}\label{fin:moments}
\begin{aligned}
 m_1(u,p)&=\int_0^\infty
       (ru^z-2^{\delta/2}\sqrt{1+\epsilon_{\mathrm{em}}^2} h_\infty r^{-\delta})\,dr,\\
 m_2(u,p)&=\int_0^\infty
       (r^2u^\theta-2^{(1+\delta)/2}h_\infty r^{1-\delta})\,dr,\\
 m_3(u,p)&=\int_0^\infty\NS^r(u,p)\,dr,\\
 m_4(u,p)&=\int_0^\infty
       (r^2u^\theta u^z-2^{\delta+1/2}\sqrt{1+\epsilon_{\mathrm{em}}^2} h_\infty^2r^{-2\delta})\,dr,\\
 m_5(u,p)&=\lim_{b\to\infty}\left\{
       \int_0^b r\left((u^z)^2+p+\frac r2\NS^r(u,p)\right)\,dr
                                  +\frac12\log b\right\}.
\end{aligned}
\end{equation}
Since $0<\delta<1/2$, the subtracted powers are integrable at zero.
At infinity, the Gamma error decays two powers faster and
$p=-1/(2r^2)+O(r^{-2-2\delta})$. Thus the finite parts and all fixed
$(t,z)$ derivatives converge uniformly on compact $(t,z)$-sets.

\begin{lemma}[Exact conservative moment identities]
\label{fin:moment-identities}
Suppose $(u,p)$ is as above and $\mathsf T$ is smooth and supported
in an annulus. Let $E=\NS(u,p)+\mathcal D\mathsf T$. Then
\begin{equation}\label{fin:compatibility-identities}
\begin{aligned}
 \int_0^\infty E^r\,dr&=m_3,\\
 \int_0^\infty r^2E^\theta\,dr
   &=(\partial_t-\partial_z^2)m_2+\partial_zm_4,\\
 \int_0^\infty rE^z\,dr
       +\frac12\partial_z\int_0^\infty r^2E^r\,dr
   &=(\partial_t-\partial_z^2)m_1+\partial_zm_5.
\end{aligned}
\end{equation}
\end{lemma}
\begin{proof}
Incompressibility gives the exact conservative forms
\begin{equation}\label{fin:conservative}
\begin{aligned}
 r^2\NS^\theta
  ={}&\partial_t(r^2u^\theta)+\partial_r(r^2u^ru^\theta)
      +\partial_z(r^2u^zu^\theta)\\
     &-\partial_r(r^2\partial_ru^\theta-ru^\theta)
      -\partial_z^2(r^2u^\theta),\\
 r\NS^z
  ={}&\partial_t(ru^z)+\partial_r(ru^ru^z)
      +\partial_z\bigl(r((u^z)^2+p)\bigr)\\
     &-\partial_r(r\partial_ru^z)-\partial_z^2(ru^z).
\end{aligned}
\end{equation}
Integrate over finite radial intervals and pass to the limits. Axis
regularity removes the inner boundary terms. At infinity,
$r^2u^ru^\theta,r^2\partial_ru^\theta-ru^\theta=O(r^{-\delta})$
and the radial fluxes in the axial equation are $O(r^{-1-\delta})$. The subtracted
terms are independent of $t,z$, and the stress boundary terms
$[r^2\mathsf T^\theta]_0^\infty,[r\mathsf T^z]_0^\infty$ vanish.
Adding $\frac12\partial_z\int r^2\NS^r\,dr$ gives the last identity.
\end{proof}

The targets for this exterior are
\begin{equation}\label{fin:moment-targets}
 m_1=z,\qquad m_2=m_3=m_4=m_5=0.
\end{equation}
By \eqref{fin:axial-primitive} and the streamfunction identity preceding
\eqref{pr:radial}, the additive $Z$ in \eqref{ue:right-moments} gives
$m_1(u_*,p_*)=\lambda^{1-\delta}Z=z$. With
$\varepsilon_\lambda=\bar{\lambda}(1+|\log\lambda|)$ for $0<\lambda\le1$, the remaining defects satisfy
\begin{equation}\label{fin:moment-units}
 \left(\lambda^{-2+\delta}m_2,
       \lambda^{2+2\delta}m_3,
       \lambda^{-1+2\delta}m_4,
       \lambda^{2\delta}m_5\right)(u_*,p_*)
       =O_m(\varepsilon_\lambda).
\end{equation}
Here $O_m(\varepsilon_\lambda)$ includes normalized derivatives through order $m$.
For the leading fields $(u_{(0)},p_{(0)})$, matching at $R_e$ and taking
finite parts gives $m_2(u_{(0)},p_{(0)})=m_4(u_{(0)},p_{(0)})=0$ and
\[
 \int_0^\infty r\bigl((u^z_{(0)})^2+p_{(0)}\bigr)\,dr
 =\lambda^{-2\delta}
   \left[M^{z\theta}(\infty,Z)+\lim_{R\to\infty}RP(R,Z)\right]=0.
\]
Indeed, for fixed $\delta>0$, the right energy primitive tends to
zero and $P=O(R^{-1-\delta})$. Inserting the core changes the normalized
$m_2,m_4$ by $O_m(\bar{\lambda})$. Equation \eqref{fin:swirl-equations} cancels these defects.
By \eqref{fin:seed-residual}, the compact radial residual
contributes $O_m(\bar{\lambda})$ to normalized $m_3,m_5$. The sink term is
\[
\begin{aligned}
 \lim_{b\to\infty}\left\{-\frac12\int_0^b
   \chi_h(r^2/(2\lambda^2))\frac{dr}{r}+\frac12\log b\right\}
 =\frac12\log\lambda+
 \frac14\left[\log(2R_e)-\int_0^{R_e}\frac{\chi_h(R)}R\,dR\right].
\end{aligned}
\]
Differentiating under the integrals gives \eqref{fin:moment-units},
since every fixed normalized derivative of $\bar{\lambda}\log\lambda$ is
$O(\varepsilon_\lambda)$.

\paragraph{Four physical moment equations}
\label{fin:repair-section}

Use $x=R/R_b$ and
$A=U_{(0)}^\theta(R_b,0)>0$. Define
$\phi_0(x,Z)=U_{(0)}^\theta(R_bx,Z)/A$ and
$\psi_0(x,Z)=U_{(0)}^z(R_bx,Z)/A$.
Let $J_i=I_i/R_b$. Choose $\rho_2$ smooth and positive on $J_2$, flat
and zero outside it, and $\rho_3\in C_c^\infty(J_3)$ nonnegative
and positive on a nonempty open subinterval. For a normalized pair
$(\phi,\psi)$ put
\begin{equation}\label{fin:angular-basis}
 r_{\phi,\psi}=(\sqrt x,\phi/x,\sqrt x\psi,\phi)^T.
\end{equation}
On each window \(J_i\), \(i=2,3\), use the four fixed test functions
\begin{equation}\label{fin:mixed-gram}
 \xi_j=\rho_i(r_{\phi_0,\psi_0})_j,\quad
 G_i(Z)=\int_{J_i}\rho_i\,r_{\phi_0,\psi_0}r_{\phi_0,\psi_0}^{T}\,dx,\quad
 g_i=\inf_{|Z|\le1}\lambda_{\min}(G_i(Z)).
\end{equation}

\begin{lemma}[Inverse using four angular profiles]
\label{fin:mixed-inverse}
For $i=2,3$, $g_i>0$. For $(\phi,\psi)$ sufficiently close in $C^0(J_i)$
to $(\phi_0,\psi_0)$ uniformly in $Z$, the physical moment matrix
\begin{equation}\label{fin:swirl-jacobian}
 B_{j\ell}=\int_{J_i}(r_{\phi,\psi})_j\xi_\ell\,dx
\end{equation}
satisfies $\|B^{-1}\|\le2/g_i$. For smooth data, each fixed parameter
derivative of $B^{-1}$, including compositions of
$\lambda\partial_\lambda,\partial_Z$, is bounded in terms of $g_i^{-1}$
and the corresponding finite collection of row and test-function seminorms.
\end{lemma}
\begin{proof}
Lemma~\ref{rp:analytic-local-rank} gives $g_i>0$ after the fixed rescalings by $R_b,A$.
Since $J_i$ stays away from zero,
$\|B-G_i\|\le C_i(\|\phi-\phi_0\|_\infty+\|\psi-\psi_0\|_\infty)$.
For $\|B-G_i\|\le g_i/2$, the inverse estimate in Lemma~\ref{pr:gram} applies
uniformly in $Z$. Derivatives of the integral in \eqref{fin:swirl-jacobian}
act on both the rows and the test functions, giving the stated seminorms.
\end{proof}

For the initial pair $\lambda^{1+\delta}(u_*^\theta,u_*^z)/A$ on $I_2$,
$B=G_2$ since $R_c<R_{ac}$. Set
\begin{equation}\label{fin:angular-correction}
 w^\theta=\lambda^{-1-\delta}A h(x,\lambda,Z),\qquad
 h=\sum_{j=1}^4 b_j(\lambda,Z)\xi_j(x,Z),\qquad
 u_B=u_*+w^\theta e_\theta.
\end{equation}
Solve the four exact equations
\begin{equation}\label{fin:swirl-equations}
\begin{aligned}
 \int_0^\infty r^2w^\theta\,dr&=-m_2(u_*,p_*),\\
 \int_0^\infty
   \frac{2u_*^\theta w^\theta+(w^\theta)^2}{r}\,dr
                                      &=m_3(u_*,p_*),\\
 \int_0^\infty r^2u_*^z w^\theta\,dr&=-m_4(u_*,p_*),\\
 \int_0^\infty r[2u_*^\theta w^\theta+(w^\theta)^2]\,dr
                                      &=2m_5(u_*,p_*).
\end{aligned}
\end{equation}
Only $-(u^\theta)^2/r$ depends on swirl in the radial equation.
Hence \eqref{fin:swirl-equations} sets $m_2,m_3,m_4,m_5$ to zero
algebraically, with $u^r,u^z$ fixed.

After dividing by the physical and amplitude factors, the equations
are $B(Z)b+N(b)=d(\lambda,Z)$, where
\begin{equation}\label{fin:swirl-quadratic}
 N(b)=\left(0,\ \frac12\int_{J_2}\frac{h^2}{x}\,dx,
                \ 0,\ \frac12\int_{J_2}h^2\,dx\right)^T
\end{equation}
and
\begin{equation}\label{fin:mixed-target-units}
 d=\left(
 -\frac{\lambda^{-2+\delta}m_2}{\sqrt2 R_b^{3/2}A},\quad
 \frac{\lambda^{2+2\delta}m_3}{A^2},\quad
 -\frac{\lambda^{-1+2\delta}m_4}{\sqrt2 R_b^{3/2}A^2},\quad
 \frac{\lambda^{2\delta}m_5}{R_bA^2}\right)^T.
\end{equation}
Here we used $r^2dr=\lambda^3\sqrt{2R}\,dR$, $rdr=\lambda^2dR$
and divided the last equation in \eqref{fin:swirl-equations} by two.
In fixed Euclidean norms choose finite constants
\begin{equation}\label{fin:ift-gate}
\begin{aligned}
 L_B&=\sup_{|Z|\le1}\|B(Z)^{-1}\|, &
 \qquad|N(b)|&\le C_Q|b|^2, & \qquad\|DN(b)\|&\le2C_Q|b|,\\
 4L_BC_Q\varrho&\le1, & \qquad|d|&\le\frac{\varrho}{2L_B}.
\end{aligned}
\end{equation}
By \eqref{fin:moment-units}, $d=O_m(\varepsilon_\lambda)$, so small $\lambda_*$
gives \eqref{fin:ift-gate} and $8L_B^2C_Q|d|\le1$.
Lemma~\ref{pr:quadratic} gives the unique small solution $|b|\le2L_B|d|\le\varrho$.
Since $\|(B+DN(b))^{-1}\|\le2L_B$, differentiation gives
$b=O_m(\varepsilon_\lambda)$ through $Z=\pm1$.

Consequently,
\begin{equation}\label{fin:five-zero}
\begin{aligned}
 m_1(u_B,p_*)&=z, & \qquad m_j(u_B,p_*)&=0\quad(2\le j\le5),\\
 \lambda^{1+\delta}(u_B^\theta-u_{(0)}^\theta,
                    u_B^z-u_{(0)}^z)&=O_m(\varepsilon_\lambda), &
 \qquad\lambda(u_B^r-u_{(0)}^r)&=O_m(\varepsilon_\lambda).
\end{aligned}
\end{equation}
The corrections are supported in $I_2$, so the core, exterior,
and velocity on $I_3$ are unchanged.

\paragraph{The realization map.}\label{fin:realization-map}

For the three double-mean defects of Section~\ref{pulse:sec:3.1}, write
$P=\int_0^\infty\langle f^r\rangle_{\theta,Y}\,dr$,
$J^{\theta z}$ and $J^{z\theta}$.
The linearized moment equations at the background, for $v=be_\theta$, are
\begin{equation}\label{fin:realization-rows}
\begin{aligned}
 \int r^2b\,dr&=0,
 &\qquad 2\int u_B^\theta b\,\frac{dr}{r}&=P,\\
 \int r^2u_B^z b\,dr&=-J^{\theta z},
 &\qquad \int r u_B^\theta b\,dr&=J^{z\theta}.
\end{aligned}
\end{equation}
Indeed, the linearized radial force is $-2u_B^\theta b/r$.
Thus the radial defect changes by $-P$. After positive row normalization
the matrix is \eqref{fin:swirl-jacobian} on $I_3$.

\paragraph{Exact pressure and stress reconstruction}
\label{fin:reconstruction-section}

Keep the temporary pressure \(p_*\), and define
\begin{equation}\label{fin:raw-error}
 E=\NS(u_B,p_*)+\mathcal D\mathsf T_{(0)}.
\end{equation}
It is smooth, supported in $[R_{\rm in},R_e]$, and satisfies
\eqref{fin:seed-residual} with $O_m(\varepsilon_\lambda)$ by
\eqref{fin:five-zero}. Lemma~\ref{fin:moment-identities} gives
\begin{equation}\label{fin:three-zero}
 \int_0^\infty E^r\,dr=0,\qquad
 \int_0^\infty r^2E^\theta\,dr=0,\qquad
 \int_0^\infty rE^z\,dr
      +\frac12\partial_z\int_0^\infty r^2E^r\,dr=0.
\end{equation}
Define the finite radial integrals below, omitting $(t,s,z)$
in their integrands:
\begin{equation}\label{fin:quadratures}
\begin{aligned}
 q(t,r,z)&=\int_r^\infty E^r\,ds=-\int_0^r E^r\,ds,\\
 \Delta\mathsf T^\theta(t,r,z)
     &=-\frac1{r^2}\int_0^r s^2E^\theta\,ds
       =\frac1{r^2}\int_r^\infty s^2E^\theta\,ds,\\
 \Delta\mathsf T^z(t,r,z)
     &=-\frac1r\int_0^r s(E^z+\partial_z q)\,ds
       =\frac1r\int_r^\infty s(E^z+\partial_z q)\,ds,\\
 p_B&=p_*+q,\qquad
 \mathsf T_B=\mathsf T_{(0)}+\Delta\mathsf T.
\end{aligned}
\end{equation}
By Fubini's theorem,
\begin{equation}\label{fin:pressure-moment}
 \int_0^\infty rq\,dr
       =\frac12\int_0^\infty r^2E^r\,dr.
\end{equation}
By \eqref{fin:three-zero}, $q$ and both stress primitives vanish outside the
annulus and extend smoothly by zero, since their sources are smooth there.
In particular $q$ preserves both endpoint pressures. Differentiation gives $\partial_r q=-E^r$,
$(\partial_r+2/r)\Delta\mathsf T^\theta=-E^\theta$, and
$(\partial_r+1/r)\Delta\mathsf T^z=-(E^z+\partial_z q)$. Hence \eqref{fin:raw-error} gives
\begin{equation}\label{fin:exact-equation}
 \nabla\cdot u_B=0,\qquad
 \NS(u_B,p_B)=-\mathcal D\mathsf T_B.
\end{equation}
Moreover,
\[
 m_5(u_B,p_B)-m_5(u_B,p_*)
 =\int_0^\infty\left(rq+\frac{r^2}{2}\partial_rq\right)dr
 =\left[\frac{r^2q}{2}\right]_0^\infty=0.
\]
Thus the five physical moment conditions hold for $(u_B,p_B)$.

Integrating \eqref{fin:seed-residual} radially gives
\(q=\lambda^{-2-2\delta}O_m(\varepsilon_\lambda)\); hence
\(\partial_z q=\lambda^{-3-\delta}O_m(\varepsilon_\lambda)\).
The stress primitives then give
\begin{equation}\label{fin:uniform-increments}
 \lambda^{2+2\delta}(p_B-p_{(0)})=O_m(\varepsilon_\lambda),\qquad
 \lambda^{2+\delta}\Delta\mathsf T=O_m(\varepsilon_\lambda)
 \quad (R_{\rm in}\le R\le R_e).
\end{equation}
Use the fixed integration variable $R$; the kernels are smooth since
$R_{\rm in}>0$. Physical estimates follow from \eqref{fin:chain-rules}.

\paragraph{Relative estimates at the inner stress edge.}\label{fin:edges-section}
Let $K$ be given by \eqref{ic:detail-data-sizes}.
Choose a fixed collar \(0\le y\le y_0<1/(2C_*K^{100})\),
\(y=\log(R/R_{\rm in})\), on which \(\chi_c=1\) and $R<R_{ac}$.
The final leading profile is unchanged there. Write
\(\omega_a(y)=1-\chi_b(y)\), with \(\chi_b\) from \eqref{ic:exit-cutoff}.
The weight vanishes at $R=R_{\rm in}$ and satisfies $\omega_a'\ge0$ and
\begin{equation}\label{fin:inner-weight}
 \omega_a(y)\asymp e^{-1/(C_*^2K^{200}y^2)},\qquad
 |\partial_y^j\omega_a(y)|\le C_j\omega_a(y)y^{-q_j},\qquad
 \int_0^y(\omega_a+\omega_a')(s)\,ds\le(1+y_0)\omega_a(y).
\end{equation}
Equation~\eqref{ic:inner-flat-size} also gives
\(|\mathcal T_{(0)}|\asymp\omega_a\) on this collar,
uniformly in \(Z\).

Write $d=u_{(0)}-u_{c,(0)}$ and let $d_{\rm n}$ denote its normalized
components. In the zero-data difference equations
\eqref{ic:actual-bridge} used for \eqref{ic:flat-error},
$\partial_y\log(F_{(0)}/F_c)=-\omega_a\partial_y\log F_c$ and $\omega_a$ is
independent of $Z$. One radial integration and two differentiations
give, for every fixed $k$,
\begin{equation}\label{fin:edge-difference-bounds}
 \begin{aligned}
 |\partial_Z^k d_{\rm n}|&\le C_k\int_0^y\omega_a(s)\,ds,\\
 |\partial_y\partial_Z^k d_{\rm n}|&\le C_k\omega_a(y),\qquad
 |\partial_y^2\partial_Z^k d_{\rm n}|
       \le C_k(\omega_a(y)+\omega_a'(y)).
 \end{aligned}
\end{equation}
The radial component satisfies the same bounds by \eqref{pr:radial},
with zero axial moment difference at $R_{\rm in}$. Equation~\eqref{fin:gaussian-primitive}
also gives $\int_0^y\omega_a\le Cy^3\omega_a(y)$.

Let $h=u_c-u_{c,(0)}$. Here $u_*=u_c+d$ and the moment corrections
vanish. Direct computation gives
\begin{equation}\label{fin:inner-residual-identity}
 E=\NS(u_{(0)},p_{(0)})+\mathcal D\mathsf T_{(0)}
       -\NS(u_{c,(0)},p_{c,(0)})
       +(h\cdot\nabla)d+(d\cdot\nabla)h.
\end{equation}
The first three tangential terms give $-\partial_z^2d^{\theta,z}$.
The radial pressure cancels the centrifugal term, leaving
\[
 \partial_tu^r+u^r\partial_ru^r+u^z\partial_zu^r
 -\left(\partial_r^2+r^{-1}\partial_r-r^{-2}+\partial_z^2\right)u^r.
\]
Its difference uses at most two $y$ derivatives, with unit $\lambda^{-3}$.
For the swirl, \eqref{ex:scale-remainder} bounds
$\lambda^{2+\delta}h^\theta/(r\bar{\lambda})=\lambda^{2-\delta}h^\theta/r$. Thus
\[
 -\frac{2h^\theta d^\theta}{r}
     =-2\lambda^{-3}\Bigl(\lambda^{2-\delta}\frac{h^\theta}{r}\Bigr)
                      (\lambda^{1+\delta}d^\theta).
\]
By \eqref{ex:remainder-factor}, \eqref{fin:chain-rules}, and
\eqref{fin:edge-difference-bounds}, the other cross terms give
\begin{equation}\label{fin:inner-source}
 \lambda^3|E^r|+\lambda^{3-\delta}\bigl(|E^\theta|+|E^z|\bigr)
    \le C(\omega_a+\omega_a').
\end{equation}
The same estimate holds after fixed compositions of $\lambda\partial_\lambda,\partial_Z$
by \eqref{ex:scale-remainder} and the edge bounds \eqref{fin:inner-weight}.

In integrals with respect to $\rho$, we evaluate $E$ at $(t,\lambda\sqrt{2\rho},z)$.
Let $P_c=\lambda^2q$. Integrating \eqref{fin:three-zero} from $R=R_{\rm in}$ gives
\begin{equation}\label{fin:inner-pressure}
 P_c(\lambda,R,Z)
   =-\lambda^3\int_{R_{\rm in}}^R\frac{E^r}{\sqrt{2\rho}}\,d\rho,
 \qquad
 |P_c|+|\lambda \partial_\lambda P_c|+|\partial_Z P_c|\le C\omega_a,\quad
 |\partial_R P_c|\le C(\omega_a+\omega_a').
\end{equation}
Hence
\begin{equation}\label{fin:inner-pressure-z}
 \partial_zq=\frac{\lambda^{-3+\delta}}{1-\delta Z^2}
       \bigl(-2ZP_c+Z\lambda \partial_\lambda P_c+(1-Z^2)\partial_Z P_c-2ZR\partial_R P_c\bigr),
 \qquad
 |\partial_zq|\le C\lambda^{-3+\delta}(\omega_a+\omega_a').
\end{equation}
Integrating \eqref{fin:inner-source} and \eqref{fin:inner-pressure-z}
in the stress primitives and using \eqref{fin:inner-weight} gives
\begin{equation}\label{fin:inner-relative}
 |\lambda^{2+\delta}\Delta\mathsf T|
       \le C\bar{\lambda}\omega_a
       \le C\bar{\lambda}|\mathcal T_{(0)}|.
\end{equation}
Differentiating \eqref{fin:quadratures} and using \eqref{fin:inner-weight}
gives the bound $C_m\bar{\lambda}\omega_a y^{-q_m}$ at every fixed normalized order.

To extend the unit stress direction to $R_{\rm in}$, divide out $\omega_a$.
The zero-data equations \eqref{ic:actual-bridge} and the Gaussian primitive
\eqref{fin:gaussian-smooth-quotient} give $\partial_y d_{\rm n}=\omega_a$
times a smooth coefficient. Thus each scaled residual in \eqref{fin:inner-source}
is a sum of $\omega_a$ and $\omega_a'$ times smooth coefficients.
For a smooth $g$, integration by parts gives
\[
 \int_0^y\omega_a'(s)g(s)\,ds
   =\omega_a(y)g(y)-\int_0^y\omega_a(s)g'(s)\,ds.
\]
The last integral is $y^3\omega_a$ times a smooth function by
\eqref{fin:gaussian-smooth-quotient}. Applying these identities to
\eqref{fin:inner-pressure} and then \eqref{fin:quadratures} shows that
$\lambda^{2+\delta}\Delta\mathsf T/(\bar{\lambda}\omega_a)$ is smooth.
Likewise \eqref{ic:detail-moment-differences} gives a smooth $e/\omega_a$ in
$\mathcal T_{(0)}/\omega_a=F(\bar{\mathcal I}/\bar F+e/\omega_a)$.
By \eqref{fin:inner-relative} and the leading lower bound, the sum
$\lambda^{2+\delta}\mathsf T_B/\omega_a$ has a nonzero limit for small $\bar{\lambda}$.
Its normalization therefore extends smoothly to $R_{\rm in}$.

\paragraph{The outer edge}
\label{fin:outer-section}

Fix the final profile \eqref{pow:profiles}, its radii, its exterior
amplitude, and $\delta>0$ before the physical scale.
On an outer collar beyond the compact physical corrections and pressure
cutoff transitions,
\begin{equation}\label{fin:outer-seed}
\begin{aligned}
 u^\theta&=\lambda^{-1-\delta}U^\theta, &
 \qquad u^z&=\lambda^{-1-\delta}U^z,\\
 u^r&=\lambda^{-1}U^r, &
 \qquad p_*&=\lambda^{-2-2\delta}P-\frac1{2r^2},\\
 \mathsf T_{(0)}&=\lambda^{-2-\delta}\mathcal T, &
 \qquad R&=\frac{r^2}{2\lambda^2}, &
 y&=\log(R/R_e).
\end{aligned}
\end{equation}
The profile pressure and moments match the exterior at $R_e$,
and \eqref{fin:three-zero} holds.

\paragraph{The radial residual.}

Let $V=U^r-U_e^r$, where $U_e^r=-1/\sqrt{2R}$. On this collar
$d=\log(R_e/R)=-y$. By \eqref{pow:edge-data}, the cutoff $\sigma_{\rm norm}(d/\ell)$ divided by
$e^{-\ell^2/d^2}$ extends smoothly to a positive value at $R=R_e$.
The exterior fields \eqref{ue:physical} are independent of physical $z$.
By incompressibility and the right axial primitive,
\begin{equation}\label{fin:outer-radial-difference}
 V=\frac{2ZR(U^z-U_e^z)-\big((1-\delta)Z+(1-Z^2)\partial_Z\big)(M^z-M_e^z)}
                {(1-\delta Z^2)\sqrt{2R}}.
\end{equation}
For a profile \(G(y,Z)\), the derivatives in \eqref{fin:operators} reduce to
\begin{equation}\label{fin:outer-operators}
 \lambda\partial_\lambda G=0,\qquad R\partial_R G=\partial_yG .
\end{equation}
By \eqref{fin:chain-rules}, with $\bar{\lambda}=\lambda^{2\delta}$,
\begin{equation}\label{fin:outer-radial-residual}
\begin{aligned}
 \lambda^3 E^r={}&\mathcal T_{-1}V
 +\sqrt{2/R}\{U_e^r\partial_y V+(\partial_y U_e^r)V+V\partial_y V\}
 +U^z\mathcal Z_{-1}V\\
 &-\frac2R \partial_y^2 V+\frac{V}{2R}
       -\bar{\lambda}\mathcal Z_{-2+\delta}\mathcal Z_{-1}V.
\end{aligned}
\end{equation}
Indeed, $\partial_R P=(U^\theta)^2/(2R)$ cancels the centrifugal term.
Subtracting the exact radial equation for $u_e^r=-1/r$ removes the sink
pressure. The radial vector Laplacian of $\lambda^{-1}V$ is
$\lambda^{-3}(2\partial_y^2V/R-V/(2R))$, and \eqref{fin:chain-rules}
gives the axial Laplacian. Thus $\lambda^3 E^r$ uses at most two $y$ derivatives and
depends on $\lambda$ only through $\bar{\lambda}$.

\paragraph{Tangential sources.}

By \eqref{pr:leading-residual}, $E^{\theta,z}=-\partial_z^2u^{\theta,z}$.
Factoring out the $z$-independent $u_e^\theta$ before differentiating gives
\begin{equation}\label{fin:outer-tangential-general}
 \lambda^{3-\delta}\begin{pmatrix}E^\theta\\E^z\end{pmatrix}
 =-U_e^\theta\,
       \mathcal Z_{-1+\delta}\mathcal Z_0
                 \left[\frac{1}{U_e^\theta}
                 \begin{pmatrix}U^\theta-U_e^\theta\\ U^z-U_e^z\end{pmatrix}\right].
\end{equation}
The operators act componentwise with bounded smooth coefficients on $|Z|\le1$.

\paragraph{Gaussian bounds and the unit stress direction.}
By \eqref{pow:edge-data}, all three velocity differences are $\sigma_{\rm norm}(d/\ell)$ times smooth coefficients. Since $U_e^\theta>0$ on the collar, division by $U_e^\theta$ preserves this property.
Each $y$ derivative costs at most $d^{-3}$ by \eqref{nc:paper-gaussian}.
The two derivatives in the physical residual formulas give
\begin{equation}\label{fin:outer-source-bound}
 (\lambda^3 E^r,\lambda^{3-\delta}E^\theta,\lambda^{3-\delta}E^z)
 =\sigma_{\rm norm}(d/\ell)d^{-6} B(d,Z,\bar{\lambda}),
 \qquad B\in C^\infty,\quad 0\le\bar{\lambda}\le1 .
\end{equation}
By \eqref{fin:three-zero}, integrating from the outer edge gives
\begin{equation}\label{fin:outer-pressure}
\begin{aligned}
 q(t,r,z)&=\int_r^\infty E^r(t,s,z)\,ds
           =\lambda^{-2}P_c(\lambda,R,Z),\\
 P_c&=\lambda^3\int_R^{R_e}\frac{E^r}{\sqrt{2\rho}}\,d\rho, &
 \partial_y P_c&=-\sqrt{R/2}\,\lambda^3 E^r.
\end{aligned}
\end{equation}
The same primitive rule and \eqref{fin:outer-pressure} give
\begin{equation}\label{fin:outer-pressure-bound}
 P_c=\sigma_{\rm norm}(d/\ell)d^{-3} B_p(d,Z,\bar{\lambda}),\qquad
 \partial_yP_c=\sigma_{\rm norm}(d/\ell)d^{-6}\widetilde B_p(d,Z,\bar{\lambda}),
 \qquad B_p,\widetilde B_p\in C^\infty .
\end{equation}
By the chain rule,
\begin{equation}\label{fin:outer-pressure-z}
 \partial_z q=\lambda^{-3+\delta}\mathcal Z_{-2}P_c,\qquad
 \mathcal Z_{-2}P_c=\frac{-2ZP_c+Z\lambda\partial_\lambda P_c+(1-Z^2)\partial_Z P_c-2Z\partial_y P_c}
                         {1-\delta Z^2}.
\end{equation}
Hence $(\mathcal Z_{-2}P_c)/(\sigma_{\rm norm}(d/\ell)d^{-6})$ extends smoothly to $R=R_e$, including at $\bar{\lambda}=0$.

Use the backward stress formulas in \eqref{fin:quadratures} and
change variables $s=\lambda\sqrt{2\rho}$:
\begin{equation}\label{fin:outer-source}
\begin{aligned}
 \lambda^{2+\delta}\Delta\mathsf T^\theta
   &=\frac{\bar{\lambda}}{2R}
            \int_R^{R_e}\sqrt{2\rho}\,\lambda^{3-\delta}E^\theta\,d\rho,\\
 \lambda^{2+\delta}\Delta\mathsf T^z
   &=\frac{\bar{\lambda}}{\sqrt{2R}}
            \int_R^{R_e}(\lambda^{3-\delta}E^z+\mathcal Z_{-2}P_c)\,d\rho.
\end{aligned}
\end{equation}
Applying \eqref{nc:paper-gaussian} once more and using the leading
stress in \eqref{pow:edge-data}, we obtain
\begin{equation}\label{fin:outer-relative}
 \mathcal T=\sigma_{\rm norm}(d/\ell)d^{-3}\mathbf A(d,Z),\qquad
 \lambda^{2+\delta}\Delta\mathsf T
       =\bar{\lambda} \sigma_{\rm norm}(d/\ell)d^{-3}\mathbf A_1(d,Z,\bar{\lambda}).
\end{equation}
Both coefficients are smooth at $R_e$, where $A^\theta(0,Z)$ has a positive minimum by \eqref{pow:edge-data}. Thus $|\lambda^{2+\delta}\Delta\mathsf T|\le C_{\rm rel}\bar{\lambda}|\mathcal T|$ near the edge and, by compactness and $|\mathcal T|>0$, on the rest of the collar.
For small $\bar{\lambda}$, $(\mathbf A+\bar{\lambda}\mathbf A_1)^\theta>0$ near $R=R_e$. Normalization extends both unit stress directions smoothly, and
$|\lambda^{2+\delta}\mathsf T_B|\ge c\sigma_{\rm norm}(d/\ell)d^{-3}$ there.
For every fixed normalized derivative order $m$,
\begin{equation}\label{fin:outer-derivatives}
 |\mathscr D^{\le m}(\lambda^{2+\delta}\Delta\mathsf T)|
       \le C_m\bar{\lambda} e^{-\ell^2/d^2}d^{-N_m}
\end{equation}
for some finite $N_m$. The normalized pressure correction
$\lambda^{2+2\delta}q=\bar{\lambda} P_c$ satisfies this bound as well. The
three normalized velocity differences are $\sigma_{\rm norm}(d/\ell)$ times smooth functions.
Consequently,
\[
 \lim_{R\uparrow R_e}
 \partial_\lambda^\ell\partial_R^j\partial_Z^k
 \bigl(q,u_B-u_e,\Delta\mathsf T\bigr)=0
 \qquad(\ell,j,k\ge0),
\]
for every fixed $0<\lambda\le\lambda_*$ and $|Z|\le1$.

\paragraph{Admissibility and the scale.}
\label{fin:completion-section}

Between the endpoint collars, $|\mathcal T_{(0)}|$ has a positive minimum.
The compact estimates, \eqref{fin:inner-relative}, and \eqref{fin:outer-relative}
give fixed constants $C_{\rm rel},C_{\rm vel}$ such that
\begin{equation}\label{fin:global-relative}
 |\mathcal T_B-\mathcal T_{(0)}|
       \le C_{\rm rel}\varepsilon_\lambda|\mathcal T_{(0)}|,
 \qquad
 |F_B-F_{(0)}|+|\mathcal S_B-\mathcal S_{(0)}|
       \le C_{\rm vel}\varepsilon_\lambda .
\end{equation}
Here $\mathcal T_B=\lambda^{2+\delta}\mathsf T_B$,
$F_B=\lambda^{2+\delta}u_B^\theta/r$, and
\[
 \mathcal S_B=\lambda^{2+\delta}
   \bigl(r\partial_r(u_B^\theta/r),\ \partial_ru_B^z\bigr).
\]
Define the leading versions similarly. The edge and interior bounds give
\begin{equation}\label{fin:edge-bounds}
 |\mathcal T_B|\ge c\omega,\qquad
 |\mathscr D^{\le m}\mathcal T_B|
       \le C_m\omega d_e^{-q_m},
 \qquad
 \omega=\exp\!\left(-\frac1{C_*^2K^{200}\log^2(R/R_{\rm in})}
                    -\frac{c_h}{\log^2(R_e/R)}\right),
\end{equation}
with \(c_a=1/(C_*^2K^{200})\), \(c_h=\ell^2\) and the edge distance
\(d_e\) of \eqref{intro:edge-weight}.
For $C_{\rm rel}\varepsilon_\lambda\le1/2$, use the relative lower bound
and the comparability of $\omega$ with each endpoint weight.
Gaussian decay gives smooth zero extensions of all derivatives.

The edge estimates and compactness give a uniform leading admissibility
margin. Since the regular swirl and negative angular shear have positive
minima on the closed annulus, the functions
\[
 \kappa-2,\qquad
 -\frac{\mathcal T}{|\mathcal T|}\cdot
       \frac{\mathcal S}{|\mathcal S|},\qquad
 1-\sqrt{\frac{\kappa-2}{2}}
       \left|\frac{\mathcal T\cdot\mathcal S^\perp}
                   {\mathcal T\cdot\mathcal S}\right|
\]
are Lipschitz in $(F,\mathcal S,\mathcal T/|\mathcal T|)$ on a fixed
neighborhood of the leading data. Let $c_0>0$ bound their leading values below,
$\eta_{\rm cone}>0$ be the neighborhood size, and $L_{\rm cone}$ a Lipschitz bound.
For $C_{\rm rel}\varepsilon_\lambda\le1/2$, the unit stress changes by at most
$4C_{\rm rel}\varepsilon_\lambda$. Hence \eqref{fin:global-relative} preserves
uniform admissibility for small $\varepsilon_\lambda$.

Let $C_{\rm mom}$ bound $|d|/\varepsilon_\lambda$ in
\eqref{fin:ift-gate}.  For $0<\lambda\le1$,
\[
 \varepsilon_\lambda=\lambda^{2\delta}(1+|\log\lambda|)
       \le\left(1+\frac1{e\delta}\right)\lambda^\delta.
\]
It is enough to choose
\begin{equation}\label{fin:finite-gate}
\begin{aligned}
 \lambda_*&\le1,\qquad \lambda_*^{2\delta}\le\bar{\lambda}_0/2,\\
 \left(1+\frac1{e\delta}\right)\lambda_*^\delta&\le
 \min\left\{
 \frac{\varrho}{2L_BC_{\rm mom}},\quad
 \frac1{2C_{\rm rel}},\quad
 \frac{\min\{\eta_{\rm cone},c_0/(2L_{\rm cone})\}}
      {2(C_{\rm vel}+4C_{\rm rel})}
 \right\}.
\end{aligned}
\end{equation}
These constants depend on the fixed profile derivatives, cutoffs,
core quotient, final leading Gram matrix, and radial integrals. Higher fixed
orders require only additional seminorms, with the same $\lambda_*$.

\begin{proof}[Completion of Theorem~\ref{intro:background}]
Theorem~\ref{cl:leading} gives the leading profile.
Equations \eqref{fin:seed}, \eqref{fin:five-zero} and \eqref{fin:quadratures}
give the five physical moments, the exact equation \eqref{fin:exact-equation}
and the prescribed fields at both ends. The velocity corrections are supported
in $R<R_c$ and $I_2$, leaving $I_3$ and the power and transition profiles unchanged.
Estimates \eqref{fin:global-relative}--\eqref{fin:finite-gate} give the
weighted stress bounds and uniform admissibility.

In the core the regular quotients \eqref{ex:normalized} are analytic.
The radial and angular velocities in normalized units are $\sqrt{2R}$ times
their regular quotients, so $R\partial_R$ and $\lambda\partial_\lambda
=2\delta\bar{\lambda}\partial_{\bar{\lambda}}$ preserve the bounds. In the exterior the
normalized velocities are independent of $\lambda$. In \eqref{ue:heat-factor},
each fixed composition of $R\partial_R$ and $\partial_Z$ applied to the
integrand defining $H_e(2L/R)$ is bounded by $Ce^{-v}v^{\delta/2-1}(1+v)^N$
for some $N$, uniformly on $R\ge R_e$, $|Z|\le1$.
Thus $\mathscr D$ derivatives of order at most $m$ of the tangential
velocities are $O_m(R^{-(1+\delta)/2})$ and those of $U_e^r=-1/\sqrt{2R}$
are $O_m(R^{-1/2})$. The velocity error is $O_m(\bar{\lambda})$ in the core
by Proposition~\ref{ex:remainder}, $O_m(\varepsilon_\lambda)$ on the
annulus by \eqref{fin:five-zero}, and zero in the exterior. Together these
bounds give the stated estimate on $0\le R\le R_0$ for every fixed $R_0>0$.

At $Z=\pm1$ all fixed derivatives extend in physical coordinates since
$\det\partial(t,z)/\partial(\lambda,Z)=-2\lambda^{2-\delta}(1-\delta Z^2)\ne0$.
Near $t=1,z=0,r>0$, the one-sided extension follows from the exact exterior
Gamma integral with $r$ bounded away from zero. Axis regularity follows
from the exact core, and \eqref{pr:tensor} gives the symmetric stress tensor.
\end{proof}

\subsubsection{Data for oscillatory correction}
\label{fin:partII-interface}

\paragraph{Estimates on dyadic bands.}
On a band $1/2<2^k\lambda<2$, introduce coordinates
\[
 \widehat R=2^{2k-1}r^2,\qquad
 \widehat Z=2^{(1-\delta)k}z,\qquad
 \widehat t=2^{2k}(1-t).
\]
In these coordinates, define
\begin{equation}\label{fin:dyadic-normalization}
\begin{aligned}
 \widehat u_B&=2^{-(1+\delta)k}u_B,
 &\qquad\widehat F_B&=2^{-(2+\delta)k}u_B^\theta/r,\\
 \widehat S_B&=2^{-(2+\delta)k}
       \bigl(r\partial_r(u_B^\theta/r),\partial_ru_B^z\bigr),
 &\qquad\widehat T_B&=2^{-(2+\delta)k}\mathsf T_B.
\end{aligned}
\end{equation}
Define the leading fields similarly. On this band and the fixed stress
annulus, hatted derivatives are bounded combinations of
$\lambda\partial_\lambda,R\partial_R,\partial_Z$. Moreover,
\[
 \varepsilon_\lambda
       =\lambda^{2\delta}(1+|\log\lambda|)
       \le C(1+k)2^{-2\delta k}.
\]
Hence
\begin{equation}\label{fin:dyadic-input}
\begin{aligned}
 (\widehat F_B,\widehat S_B,\widehat T_B)
  &=(\widehat F_{(0)},\widehat S_{(0)},\widehat T_{(0)})
       +O_{C^m}((1+k)2^{-2\delta k}),\\
 |\widehat T_B-\widehat T_{(0)}|
       &\le C(1+k)2^{-2\delta k}|\widehat T_{(0)}|,\\
 \widehat u_B^r&=O_{C^m}(2^{-\delta k}).
\end{aligned}
\end{equation}
The stress difference also satisfies the $\omega d_e^{-q_m}$
derivative bound. Equivalently, in units of squared tangential velocity,
\[
 2^{-2(1+\delta)k}(\mathsf T_B-\mathsf T_{(0)})
   =O_{C^m}\bigl((1+k)2^{-3\delta k}
                   \omega d_e^{-q_m}\bigr).
\]
The coefficient classes allow polynomial factors in $k$. Thus a sufficiently
large initial band index $k$ preserves strict admissibility, with zero equation residual.

\begin{proposition}[Quantitative data on the stress annulus]\label{fin:mixed-interface}
The constructed background has the following properties on $0<\lambda\le\lambda_*$
after $\lambda_*$ is decreased once.
\begin{enumerate}[label=\textup{(B\arabic*)},leftmargin=3.1em]
\item The Cartesian fields are smooth away from the vertex, including the
regular axis and the one-sided terminal time, and
$\operatorname{div}u_B=0$, $\NS(u_B,p_B)=-\mathcal D\mathsf T_B$.
Thus $E_B=0$.
\item On each bounded $R$ interval the four components
$\lambda u_B^r,\lambda^{1+\delta}u_B^\theta,\lambda^{1+\delta}u_B^z,
\lambda^{2+2\delta}p_B$ are bounded after any fixed composition of
$\lambda\partial_\lambda,\sqrt{2R}\partial_R,\partial_Z$.
\item The stress is supported in $R_{\rm in}\le R\le R_e$, extends smoothly by
zero there, and satisfies \eqref{intro:edge-weight}--\eqref{intro:stress-bounds}.
\item The normalized swirl and negative angular shear have uniform positive
upper and lower bounds on the closed annulus. The cone margin
\eqref{intro:uniform-cone} holds with a single $c_{\rm cone}>0$.
\item On $I_3$ the normalized tangential velocity pair is exactly
$(\phi_0,\psi_0)$ in units $A\lambda^{-1-\delta}$. For $x=R/R_b$ the
four rows $(\sqrt{x},\phi_0/x,\sqrt{x}\psi_0,\phi_0)$ have a Gram
matrix with a positive lower eigenvalue uniform in $Z$.
\end{enumerate}
The core supplies the axial trace \eqref{intro:blowup-trace}.
For $R\ge R_e$ the exterior is exactly \eqref{ue:physical}.
\end{proposition}

\begin{proof}
Theorem~\ref{intro:background} gives \textup{(B1)}, \textup{(B3)},
\textup{(B4)}, the axis trace and the exterior. For \textup{(B2)}, put
$r_{\mathrm{sim}}=\sqrt{2R}$. Then $\sqrt{2R}\partial_R=\partial_{r_{\mathrm{sim}}}$, and in the core the
normalized radial and angular velocities are $r_{\mathrm{sim}}$ times smooth functions of
$r_{\mathrm{sim}}^2$, while the axial velocity and pressure are smooth in $r_{\mathrm{sim}}^2$.
Away from the axis use \eqref{fin:five-zero} and \eqref{fin:uniform-increments}
on the annulus and \eqref{ue:physical} in the exterior.
On the unchanged window $I_3$, Lemma~\ref{fin:mixed-inverse} gives \textup{(B5)}.
\end{proof}

Section~\ref{pulse:sec:3.3} uses \(\xi_j=\rho_3(r_{\phi_0,\psi_0})_j\) on \(I_3\) for the
nonlinear moment correction at initialization. Lemma~\ref{fin:mixed-inverse}
also gives the inverse for the subsequent linear correction at the initialized mean.
\section{Geometry and primary pulses}\label{pulse:sec:2}
Proposition~\ref{fin:mixed-interface} gives \(\NS(u_B,p_B)=-\mathcal D\mathsf T_B\), so \(E_B\) vanishes in the force identities. On a band \(\lambda\asymp2^{-k}\), use the units
\[
 \varepsilon_k=2^{-\delta k},\qquad U_k=2^{(1+\delta)k},\qquad F_k=2^{(3+\delta)k}
\]
for the small parameter, tangential velocity and force, respectively. The pressure and radial stress have units \(U_k^2\) and \(U_k^2\varepsilon_k\), respectively. The desired velocity correction has size \(U_k\sqrt{\varepsilon_k}\) in the annulus \(R_{\rm in}<R<R_e\).
\pulsebridge{2}
We choose compatible labels and finite Fourier windows on each band,
then construct the two divergence-free pulse families.
Corollary~\ref{pulse:thm:2.16} gives their size and mean flux.

\subsection{Coordinates, covers and coefficient domains}\label{pulse:sec:2.1}
Fix a scale \(2^{-k}\), where \(k\ge k_0\) and \(k_0>1\) is a fixed sufficiently large integer. Set \(\delta_0=10^{-5}\) for the auxiliary radial derivative estimates. At order \(m\), \(P_m(k)\) denotes a polynomial bound in \(k\) whose coefficients depend only on the fixed geometry and finitely many background seminorms from Proposition~\ref{fin:mixed-interface}.

\subsubsection{Similarity and band coordinates}\label{pulse:sec:2.1.1}
\pulsebridge{2.1.1}

For a function of \((\lambda,R,Z)\), the chain rule gives
\(\partial_j=(\partial_j\lambda)\partial_\lambda+
(\partial_jR)\partial_R+(\partial_j Z)\partial_Z\), with the physical derivatives in \eqref{pr:derivatives}.

On the band \(\mathcal B_k=\{1/2<2^k\lambda<2\}\) introduce the coordinates
\begin{equation}
\widehat R=2^{2k}r^2/2,\quad \widehat Z=2^{(1-\delta)k}z,\quad \widehat T=2^{2k}(1-t),\quad \sigma=2^k\lambda.
\label{pulse:eq:2.5}
\end{equation}
Here \(\sigma\) is the similarity scale \(\lambda\) in units of \(2^{-k}\), so \(1/2<\sigma<2\) on \(\mathcal B_k\). We have
\begin{subequations}
\begin{equation}
\widehat R=\sigma^2 R,\qquad \widehat Z=\sigma^{1-\delta}Z,\qquad \widehat T=\sigma^2(1-Z^2).
\label{pulse:eq:2.6}
\end{equation}

We use the enlarged annular chart
\begin{equation}\label{pulse:band-chart}
  \tfrac18\le \sigma=2^k\lambda\le8,\qquad
  R_{\rm in}/4\le R\le4R_e,\qquad |Z|\le11/10.
\end{equation}
\end{subequations}
On \eqref{pulse:band-chart} the absolute Jacobian of
\((\sigma,Z)\mapsto(\widehat T,\widehat Z)\) is \(2\sigma^{2-\delta}(1-\delta Z^2)\),
and \(\partial_R\widehat R=\sigma^2\). Both are bounded above and away from
zero uniformly in \(k\). Thus \eqref{pulse:eq:2.6} and its inverse have bounded
derivatives on this compact chart and give equivalent Sobolev norms. At fixed \(k\),
\begin{equation}
\begin{alignedat}{2}2^{-k}\partial_r&=\sqrt{2\widehat R}\,\partial_{\widehat R},\qquad&2^{-(1-\delta)k}\partial_z&=\partial_{\widehat Z},\\2^{-2k}\partial_t&=-\partial_{\widehat T},&2^{-k}r^{-1}\partial_\theta&=(2\widehat R)^{-1/2}\partial_\theta.\end{alignedat}
\label{pulse:eq:2.7}
\end{equation}

In these coordinates set \(\widehat u_B=2^{-(1+\delta)k}u_B\), \(\widehat F_B=\widehat u_B^\theta/\sqrt{2\widehat R}\), and \(\widehat T_B=2^{-(2+\delta)k}\mathsf T_B\). Define \(\widehat S_B=(2\widehat R\partial_{\widehat R}\widehat F_B,\sqrt{2\widehat R}\partial_{\widehat R}\widehat u_B^z)\). The normalization in \eqref{intro:physical-shear} and the band estimates \eqref{fin:dyadic-input} give
\begin{equation}
\begin{gathered}
(\widehat F_B,\widehat S_B,\widehat T_B)=\sigma^{-2-\delta}(F_B,\mathcal S_B,\mathcal T_B),\\
\|\widehat u_B^{\theta,z}\|_{C^m}\le C_m,\qquad \|\widehat u_B^r\|_{C^m}\le C_m\varepsilon_k .
\end{gathered}\label{pulse:eq:2.8}
\end{equation}
These bounds hold on the part of \eqref{pulse:band-chart} with \(|Z|\le1\), the physical range of the background. On the stress annulus the common positive factor preserves \(\kappa_B\) and the cone margins.

\subsubsection{Auxiliary covers and directional inverses}\label{pulse:sec:2.1.2}
\pulsebridge{2.1.2}

Equip \(\mathbb T_Y^2=\mathbb R^2/\mathbb Z^2\) with normalized Haar measure. For the pulse construction we use the covering matrix
\begin{equation}
E=\begin{pmatrix}3&1\\1&5\end{pmatrix},\quad \Lambda_r=4-\sqrt2,\quad\Lambda_t=4+\sqrt2,\quad v_r=(1,1-\sqrt2)^{\mathsf T},\quad v_t=(\sqrt2-1,1)^{\mathsf T}.
\label{pulse:eq:2.9}
\end{equation}
Thus \(Ev_j=\Lambda_jv_j\), \(v_r\cdot v_t=0\) and \(\det E=14\). Its powers amplify the radial and time frequencies at different rates. To compare them on a physical band set
\begin{equation}
\rho_0=\frac{\log\Lambda_r}{\log\Lambda_t}>\frac12,\quad d_0=(2+\delta)\rho_0-\delta_0\delta>0,\quad Y=t v_t+r^{d_0}v_r\pmod{\mathbb Z^2},\quad Y_i=E^iY.
\label{pulse:eq:2.10}
\end{equation}
Regard \(Y_i\) as an independent torus variable when taking averages. After evaluation at \(Y_i=E^i(tv_t+r^{d_0}v_r)\), the physical time and radial derivatives act on a coefficient through
\begin{equation}
D_t^{(i)}=\partial_t+\Lambda_t^i v_t\cdot\nabla_{Y_i},\qquad D_r^{(i)}=\partial_r+d_0r^{d_0-1}\Lambda_r^i v_r\cdot\nabla_{Y_i}.
\label{pulse:eq:2.11}
\end{equation}
On the right, \(\partial_t,\partial_r\) hold \(Y_i\) fixed. The only variable coefficient depends on \(r\) alone, so \(D_t^{(i)}D_r^{(i)}=D_r^{(i)}D_t^{(i)}\). Axial and angular derivatives hold the auxiliary variable fixed even after physical evaluation. Radial auxiliary oscillations are supported away from the axis.

Choose
\begin{equation}
i(k)=\left\lfloor\log_{\Lambda_t}\frac{2^{(2+\delta)k}}{k^2}\right\rfloor,\qquad A_k=2^{(2+\delta)k}\Lambda_t^{-i(k)},\qquad M_k=2^{-kd_0}\Lambda_r^{i(k)}.
\label{pulse:eq:2.12}
\end{equation}
For large \(k\), \(i(k)\ge 1\). The floor and the definition of \(d_{0}\) give
\begin{equation}
k^2\le A_k<\Lambda_t k^2,\quad M_k=\varepsilon_k^{-\delta_0}A_k^{-\rho_0},\quad \Lambda_t^{-i(k)}2^{(3+\delta)k}/2^{(1+\delta)k}=A_k\varepsilon_k.
\label{pulse:eq:2.13}
\end{equation}

In band coordinates,
\begin{equation}
\begin{aligned}
D_t^{(i)}&=-2^{2k}\partial_{\widehat T}+\Lambda_t^i v_t\cdot\nabla_{Y_i},\\
D_r^{(i)}&=2^k\sqrt{2\widehat R}\,\partial_{\widehat R}
 +2^{(1-d_0)k}\Lambda_r^i d_0(2\widehat R)^{(d_0-1)/2}v_r\cdot\nabla_{Y_i}.
\end{aligned}\label{pulse:eq:2.14}
\end{equation}
Every derivative of the normalized coefficient retains just one factor \(M_{k}\).

\noindent\textbf{Lemma \PulseResult{2.1} (Directional divisor).}
\textit{For \(v=v_r,v_t\) and \(n\in\mathbb Z^2\setminus\{0\}\), one has
\(|v\cdot n|\ge [(1+2\sqrt2)(1+|n|)]^{-1}\).}

\noindent\textbf{Proposition \PulseResult{2.2} (Zero-mean inverse).}
\textit{For smooth \(f\) of zero \(Y\)-mean, the unique zero-mean function with directional derivative \(f\) along \(v\) is}
\begin{equation}
\bigl[(v\cdot\nabla_Y)^{-1}f\bigr](Y)=\sum_{n\ne0}\frac{\widehat f(n)}{2\pi i\,v\cdot n}\,e^{2\pi i n\cdot Y}.
\label{pulse:eq:2.15}
\end{equation}
\textit{The inverse satisfies \(\|(v\cdot\nabla_Y)^{-1}f\|_{H^m}\le C_m\|f\|_{H^{m+1}}\) in the joint slow, auxiliary and carrier variables. On \(\mathbb T^2\) one has \(\|(v\cdot\nabla_Y)^{-1}f\|_{C^m}\le C_m\|f\|_{C^{m+3}}\). It preserves reality and support in the passive variables, and commutes with their derivatives.}

For cover levels \(k_1\le k_2\), let \(f\) be any smooth function of \(Y_{k_2}\) on the level-\(k_2\) torus. With the fixed covering matrix \(E\) in \eqref{pulse:eq:2.9}, its pullback to the level-\(k_1\) torus is \(f(E^{k_2-k_1}Y_{k_1})\). The map \(n\mapsto(E^{k_2-k_1})^{\mathsf T}n\) on Fourier indices is injective, so
\begin{equation}
\begin{gathered}
\int f(E^{k_2-k_1}Y_{k_1})\,dY_{k_1}=\int f(Y_{k_2})\,dY_{k_2},\qquad
\bigl\|f(E^{k_2-k_1}Y_{k_1})\bigr\|_{L^2}=\|f\|_{L^2},\\
\bigl\|f(E^{k_2-k_1}Y_{k_1})\bigr\|_{H^m}\le C_m\|E\|^{m(k_2-k_1)}\|f\|_{H^m}.
\end{gathered}\label{pulse:eq:2.16}
\end{equation}

\noindent\textbf{Lemma \PulseResult{2.3}} (Cover compatibility)\textbf{.} \textit{Let \(\Lambda_v\) be the eigenvalue of \(E\) corresponding to \(v\). If \(f\) has zero mean, then}
\begin{equation}
\Lambda_v^{-k_1}(v\cdot\nabla_{Y_{k_1}})^{-1}\!\bigl[f(E^{k_2-k_1}Y_{k_1})\bigr]
=\Lambda_v^{-k_2}\bigl[(v\cdot\nabla_{Y_{k_2}})^{-1}f\bigr](E^{k_2-k_1}Y_{k_1}).\label{pulse:eq:2.17}
\end{equation}
\textit{Products, physical derivatives, and angular and auxiliary averages commute with compatible pullback.}

\begin{proof}[Proof of Lemma~\ref{pulse:thm:2.1}, Proposition~\ref{pulse:thm:2.2}, and Lemma~\ref{pulse:thm:2.3}]
The product of \(v\cdot n\) with its algebraic conjugate is
\((n_1+n_2)^2-2n_2^2\) or \((n_2-n_1)^2-2n_1^2\), a nonzero integer.
The conjugate is bounded by \((2+\sqrt2)|n|\), giving the divisor bound.
Thus \eqref{pulse:eq:2.15} converges for smooth \(f\) and gives the unique
zero-mean inverse. Parseval's identity gives the joint \(H^m\) estimate,
and \(H^{m+2}(\mathbb T^2)\hookrightarrow C^m\) gives the \(C^m\) estimate.
The multiplier respects conjugation and acts pointwise in the passive variables.

The chain rule for \(f(E^{k_2-k_1}Y_{k_1})\) multiplies the derivative along
\(v\) by \(\Lambda_v^{k_2-k_1}\), giving \eqref{pulse:eq:2.17} by uniqueness.
Use \eqref{pulse:eq:2.11} for physical derivatives and
\eqref{pulse:eq:2.16} for Haar averages. Products commute with pullback.
\end{proof}

\subsubsection{Partitions, buffered labels and pulse clocks}\label{pulse:sec:2.1.3}
\pulsebridge{2.1.3}

Choose a nonnegative \(\chi\in C_c^\infty((-1,1))\), positive on \([-1/2,1/2]\). Divide \(\chi(s)\) by \(\bigl(\sum_{n\in\mathbb Z}\chi(s-n)^2\bigr)^{1/2}\). The denominator is positive and \(1\)-periodic, so the normalized cutoff satisfies \(\sum_{n\in\mathbb Z}\chi(s-n)^2=1\). On band \(k\) use the grid of spacing \(k^{-6}\) and set
\begin{equation}
\zeta_k(\lambda)=\chi(k+\log_2\lambda),\qquad
\zeta_{k,\alpha}(X)=\prod_{\ell=1}^{3}\chi(k^6X_\ell-\alpha_\ell),\qquad
X=(\sqrt{2\widehat R},\widehat Z,\widehat T).\label{pulse:eq:2.18}
\end{equation}
Their squares form partitions of unity. At most two bands and eight boxes per band meet at a point. For each \((k,\alpha)\) use the single outer box \(\lvert k^6X-\alpha\rvert_\infty<1\). The distance from \(\operatorname{supp}\zeta_{k,\alpha}\) to the box boundary is at least \(k^{-6}\operatorname{dist}(\operatorname{supp}\chi,\{-1,1\})>0\). Its \(X\)-derivatives of order \(\beta\) are \(O(k^{6|\beta|})\). A fixed compact chart meets \(O(k^{18})\) boxes. In \(r,z,t\) their widths are
\begin{equation}
O(2^{-k}k^{-6}),\qquad O(2^{-(1-\delta)k}k^{-6}),\qquad O(2^{-2k}k^{-6}).\label{pulse:eq:2.19}
\end{equation}

The uniform \(C^2\) bound in \eqref{pulse:eq:2.8} and the lower bound for \(\sqrt{2\widehat R}\) give uniform \(C^1\) bounds for \(\widehat F_B,\widehat S_B\) and \(1/\sqrt{2\widehat R}\). At \(Z=\pm1\) a fixed-width smooth extension preserves these bounds. Thus their oscillation on each outer box is \(O(k^{-6})\). For each box meeting the physical stress annulus we choose its reference point in the intersection.

The two pulse labels on each box are \(\gamma=(k,\alpha,\pm)\). In paired formulas the upper signs refer to the \(+\) pulse and the lower signs to the \(-\) pulse. We choose all slow cutoffs for a label inside its outer box. The gap to the boundary accommodates any fixed number of nested cutoffs, with gaps comparable to \(k^{-6}\) for that number. Labels are adjacent if their outer boxes and broad band supports meet in physical variables.

\noindent\textbf{Lemma \PulseResult{2.4}} (Bounded degree)\textbf{.} \textit{The graph has uniformly bounded degree. Native neighbors satisfy} \(|k-k'|\le2\) \textit{and} \(|i(k)-i(k')|\le2\)\textit{. Fixed broader bands replace these by fixed constants.}

\textit{Proof.} Band intersection bounds \(|k-k'|\). At a common physical point,
\begin{equation}
X'=\operatorname{diag}(2^{k'-k},2^{(1-\delta)(k'-k)},2^{2(k'-k)})X.\label{pulse:eq:2.20}
\end{equation}
These ratios and \((k^\prime/k)^6\) are bounded, so only boundedly many grid
centers meet a given box. The derivative of the expression inside the floor
in \eqref{pulse:eq:2.12} is \(((2+\delta)\log 2-2/k)/\log\Lambda_t\in(0,1)\)
for large \(k\). Taking floors gives the stated level bound. \hfill\(\square\)

\noindent\textbf{Proposition \PulseResult{2.5}} (Separated outer labels)\textbf{.} \textit{There are finitely many centers} \(O_1,\ldots,O_Q\in\mathbb T^2\)\textit{, a graph coloring} \(c(\gamma )\)\textit{, and} \(r_{0}>0\)\textit{. The rectangles} \(\mathcal R_\gamma=\{O_{c(\gamma)}+l_r v_r+l_t v_t:|l_r|,|l_t|<r_0\}\) \textit{are injective charts. Adjacent labels satisfy}
\begin{equation}
E^{-i(k)}\mathcal R_\gamma\cap E^{-i(k')}\mathcal R_{\gamma'}=\varnothing.\label{pulse:eq:2.21}
\end{equation}
\textit{Nested inner rectangles have fixed native buffers. Across} \(O(\log k)\) \textit{lower-cover levels, branch counts and inverse buffer widths are polynomial in} \(k\)\textit{.}

\textit{Proof.} Greedy coloring uses at most one more color than the degree bound.
Choose each center outside the images and preimages of earlier centers under
\(E^d\) for \(0\le d\le d_*\), where \(d_*\) bounds the level difference on an edge.
This excludes only finitely many points. Choose \(r_0\) within the injectivity
radius with \(C(1+\Lambda_t^{d_*})r_0\) smaller than the resulting positive
separation. This gives \eqref{pulse:eq:2.21}.
Over \(d\) levels there are \(14^d\) inverse branches with contraction factors
\(\Lambda_r^{-d},\Lambda_t^{-d}\), so \(d=O(\log k)\) costs only powers of \(k\).
\hfill\(\square\)

When the slow supports of two distinct labels meet, \eqref{pulse:eq:2.21}
separates their auxiliary rectangles on every common cover. Their products
therefore vanish even after differentiation and physical evaluation.
On a native rectangle set
\begin{equation}
\tau_\gamma=A_k(l_t+r_0),\qquad T_k=2r_0 A_k\asymp k^2.\label{pulse:eq:2.22}
\end{equation}
Its base-cover branches are \(E^{-i}(O_\gamma+n)+\Lambda_r^{-i}l_r v_r+\Lambda_t^{-i}l_t v_t\), \(n\in\mathbb Z^2/E^i\mathbb Z^2\). Orthogonality gives \(l_t(Y_i(t,r))=\Lambda_t^i t-v_t\cdot(O_\gamma+n)/|v_t|^2\). The corresponding \(l_{r}\) depends on \(r\) alone. Hence
\begin{equation}
2^{-(2+\delta)k}D_t\tau_\gamma=1,\quad D_r\tau_\gamma=\partial_z\tau_\gamma=0,\quad D_tl_r=\partial_zl_r=0.\label{pulse:eq:2.23}
\end{equation}

Choose nonnegative cutoffs \(\chi_{\rm aux}\in C_c^\infty((-r_0/3,r_0/3))\) and \(\chi_\tau\in C_c^\infty((1/6,5/6))\), with \(\chi_{\rm aux}\not\equiv0\) and \(\chi_\tau=1\) on \([3/10,7/10]\). Set
\begin{equation}
\chi_\gamma=\zeta_k\zeta_{k,\alpha}\,\chi_{\rm aux}(l_r)\chi_\tau(\tau_\gamma/T_k).\label{pulse:eq:2.24}
\end{equation}
The larger cutoffs have polynomial derivative bounds and extend smoothly by zero from the rectangles. The Haar Jacobian is
\begin{equation}
J_E=|\det[v_r,v_t]|=|v_r||v_t|=4-2\sqrt2,\qquad dY_i=J_E\,dl_r\,dl_t=\frac{2r_0J_E}{T_k}\,dl_r\,d\tau.\label{pulse:eq:2.25}
\end{equation}

For a matrix \(H(\tau )\) independent of \(l_{r}\),
\begin{equation}
\int_{\mathbb T^2}\chi_{\rm aux}(l_r)^2\chi_\tau(\tau_\gamma/T_k)^2H(\tau_\gamma)\,dY_i
=\frac{2r_0J_E}{T_k}\left(\int\chi_{\rm aux}(l_r)^2\,dl_r\right)\int_0^{T_k}\chi_\tau(\tau/T_k)^2H(\tau)\,d\tau.\label{pulse:eq:2.26}
\end{equation}
Haar invariance preserves \eqref{pulse:eq:2.26} on lower covers.

By \eqref{pulse:eq:2.23} the native clock \(\tau_\gamma\) increases with physical
time at rate \(2^{(2+\delta)k}\), while \(\widehat T=2^{2k}(1-t)\) decreases.
Thus \(T_k\asymp k^2\) represents a physical interval of length
\(O(2^{-(2+\delta)k}k^2)\). By \eqref{pulse:eq:2.8} and \eqref{pulse:eq:2.13},
the slow and transverse auxiliary displacements are \(O(k^2\varepsilon_k)\) and
\(O(P(k)\varepsilon_k^{1-\delta_0})\). Both fit inside the polynomial buffers.

\subsubsection{A short smooth extension}\label{pulse:sec:2.1.4}
\pulsebridge{2.1.4}

\noindent\textbf{Lemma \PulseResult{2.6} (Half-space extension).}
\textit{Seeley's real linear operator extends smooth data smoothly across a
half-space boundary. Its \(H^m\to H^m\) norm on a strip of width \(0<h<1\)
is at most \(C_mh^{-m}\).
It acts pointwise in the passive variables, preserves their supports, and
commutes with their derivatives.}

Appendix~\ref{app:pulse-tools} proves the lemma. For \(b_j=2^j\), the
coefficients of the extension satisfy, for every integer \(m\ge0\),
\begin{equation}
 |a_j|\le C2^{-j(j+1)/2},\qquad \sum_j |a_j|b_j^m<\infty,\qquad
 \sum_j a_j(-b_j)^m=1 .
\label{pulse:eq:2.27}
\end{equation}
The cutoff \(\chi\) in \eqref{pulse:eq:2.18} equals one near zero.
For \(Z=1+x>1\), set
\begin{equation}
(\operatorname{Ext}_h f)(1+x)=\sum_j a_j\chi(b_j x/h)f(1-b_j x).\label{pulse:eq:2.28}
\end{equation}
At \(Z=-1-x\), use \(-1+b_jx\) as the last argument. Keep \(f\) on \([-1,1]\). The extension samples within \(2h\) of \(Z\) and vanishes for \(|Z|\ge1+h\). For its width and the edge cutoffs, set
\begin{equation}
\kappa=1+\log(1/\lambda),\qquad h(\kappa)=c_{\rm ext}(1+\kappa)^{-7},\qquad 0<c_{\rm ext}<1/100.\label{pulse:eq:2.29}
\end{equation}
We use this extension at fixed \((\lambda, R, Y, \eta )\) on every window.

\noindent\textbf{Proposition \PulseResult{2.7}} (Support and norms of short extension)\textbf{.} \textit{On band} \(k\)\textit{,} Ext \(:H^{m}\to H^{m}\) \textit{has norm} \(P_{m}(k)\)\textit{. It preserves} \(\lambda, R\) \textit{supports and radial independence. It commutes with cover pullbacks, averages and conjugation. Its slow support propagation is} \(O(k^{-7})\)\textit{. This fits inside every fixed number of} \(k^{-6}\) \textit{label collars.}

\textit{Proof.} On a band \(\kappa\asymp k\), so differentiating the variable
width in \eqref{pulse:eq:2.28} costs only fixed powers of \(k\) and \(b_j\).
Equation \eqref{pulse:eq:2.27} sums these bounds. The passive variables stay fixed,
and \eqref{pulse:eq:2.6} carries the sampling distance \(2h\) to
\(O(k^{-7})=o(k^{-6})\) in slow coordinates. \hfill\(\square\)

With time as the normal variable and fixed width, the same extension preserves terminal jets and compact spatial support.

\subsubsection{The full coefficient torus}\label{pulse:sec:2.1.5}
\pulsebridge{2.1.5}

Choose a fixed slow torus whose fundamental rectangle contains the domain \eqref{pulse:band-chart} in the coordinates \((\sigma,R,Z)\).
For row \(j\), use \(\sigma _{j}= 2^{j}\lambda\). A broad band cutoff is supported in \((1/4,4)\) and equals one near supp \(\zeta _{j}\). Radial blocks also have a broad cutoff equal to one on the annulus. After short extension these cutoffs give a smooth periodization by zero. We apply them after differential reconstruction.

At common level \(c\), the largest domain is
\begin{equation}
\mathcal Q_c=\mathbb T^3_{\sigma,R,Z}\times\mathbb T^2_{Y_c}\times\mathbb T_\eta.\label{pulse:eq:2.30}
\end{equation}
The carrier has period \(2\pi\) and reality condition \(g_{-n}=\overline{g_n}\). The active variables are
\[
\begin{array}{l|l|c}
 \text{block} & \text{variables} & \text{dimension}\\ \hline
 \text{nonzero carrier} & \sigma,R,Z,Y_c,\eta & 6\\
 \text{angular mean, auxiliary fluctuation} & \sigma,R,Z,Y_c & 5\\
 \text{double-mean radial potentials} & \sigma,R,Z & 3\\
 P,J^{\theta z},J^{z\theta} & \sigma,Z & 2
\end{array}
\]
The first two blocks have zero carrier and auxiliary means, respectively. Inactive Fourier indices stay zero. In particular \(P, J\) receive no radial cutoff. Their lift supplies the interior compensators.

On each lifted rectangle, \(Y_{c}\) gives \((l_{r}, \tau _{\gamma })\). The clock factor is \(A_k\asymp k^2\). Local clock and torus norms therefore differ by polynomial factors on logarithmic windows. Endpoint cutoffs precede zero extension of forced amplitudes.

Let \(\mathfrak H_{\mathrm{phys}}\) be the physical free-source map. Let \(C_{j}\) convert row coordinates to physical units. Write \(R_j^{\mathrm{phys}}\) for normalized physical restriction, including the broad cutoff. Let \(\operatorname{res}_j\) restrict an ambient row to \(|Z|\le1\). Define
\begin{equation}
R_j^{\mathrm{ext}}=\operatorname{Ext}_j R_j^{\mathrm{phys}},\qquad A^{\mathrm{phys}}g=\sum_j\zeta_j^2 C_j\operatorname{res}_j g_j,\qquad H_j=R_j^{\mathrm{ext}}\mathfrak H_{\mathrm{phys}} A^{\mathrm{phys}}.\label{pulse:eq:2.31}
\end{equation}
Each row operator is fixed across windows. Only its physical restriction enters
synthesis, so we apply the short extension after normalization.

\noindent\textbf{Proposition \PulseResult{2.8}} (Physical synthesis)\textbf{.} \textit{Where every active band is included,}
\begin{equation}
A^{\mathrm{phys}}R^{\mathrm{ext}}=I.\label{pulse:eq:2.32}
\end{equation}
For the compatible block cutoff \(X\), the physical residual is the synthesis
of the row residuals, as in \eqref{pulse:eq:5.23}.

\textit{Extension preserves the derivative loss, scale gain, degree, band range and compatibility with pullback. Norms increase by polynomial factors and supports by the short extension distance.}

\textit{Proof.} We have \(\operatorname{res}_j\operatorname{Ext}_j=I\). The broad cutoff is one near \(\supp\zeta_j\). Thus \(\sum_j\zeta_j^2=1\) proves \eqref{pulse:eq:2.32}. Apply it to \(\mathfrak H_{\mathrm{phys}}(A^{\mathrm{phys}} Xg)\) and subtract. Proposition~\ref{pulse:thm:2.7} supplies the norm and support claims. Linearity preserves degree; exact \(\lambda\) support preserves band range. Cover pullback commutes with extension. Physical differentiation occurs before extension. \hfill\(\square\)

\subsubsection{Logarithmic covers and finite Fourier dimension}\label{pulse:sec:2.1.6}
\pulsebridge{2.1.6}

The index \(i(k)\) is nondecreasing with bounded increments for large \(k\), by the derivative computed above. Set
\begin{equation}
m_k=\lceil\log(2+k)\rceil,\qquad I_k=\{j:|j-k|\le m_k\},\qquad c_k=i(k-m_k-h_0).\label{pulse:eq:2.34}
\end{equation}
The integer \(h_{0}\) covers all native levels used in one reconstruction, and \(r_{0}\) is the rectangle width.

\noindent\textbf{Lemma \PulseResult{2.9}} (Common-cover costs)\textbf{.} \textit{All required native levels exceed} \(c_{k}\) \textit{by} \(O(\log k)\)\textit{. Fixed-order pullback costs, branch counts and inverse buffer widths are polynomial. For} \(j\in I_{k}\)\textit{,}
\begin{equation}
P(k)^{-1}\varepsilon_j^{-\delta_0}\le 2^{-jd_0}\Lambda_r^{c_k}\le P(k)\varepsilon_j^{-\delta_0},\qquad \Lambda_t^{-c_k}2^{(3+\delta)j}/2^{(1+\delta)j}\le P(k)\varepsilon_j.\label{pulse:eq:2.35}
\end{equation}
\textit{Bounded center distance gives bounded} \(|c_k-c_{k'}|\)\textit{.}

\textit{Proof.} Bounded index increments give the \(O(\log k)\) gap. Apply \eqref{pulse:eq:2.16} and the branch factors \(14^d,\Lambda_r^d,\Lambda_t^d\). Multiply \eqref{pulse:eq:2.13} by the corresponding level-change factors. Finally \(m_k-m_{k'}\) is bounded when \(k-k'\) is bounded. \hfill\(\square\)

Compatible pullback defines the global operators on \(Y_0\). Apply the common-cover estimates to pullbacks of the active native source families.

Choose a real even smooth symbol \(s= 1\) on \(|\xi|\le1\) and zero for \(|\xi|\ge2\). After fixed period rescaling, set
\begin{equation}
S_Lg=\sum_{n\in\mathbb Z^d}s(n/L)\widehat g(n)e^{in\cdot x},\qquad \mathcal E_L=\operatorname{range}S_L.\label{pulse:eq:2.36}
\end{equation}
The multiplier acts only on active variables and preserves reality and zero means. The range has dimension at most \(C_{d}(1 +L)^{d}\). Fourier weights give
\begin{equation}
\|g\|_{H^{m+n}}\le C_n L^n\|g\|_{H^m}\quad(g\in\mathcal E_L),\qquad \|(I-S_L)g\|_{H^m}\le C_n L^{-n}\|g\|_{H^{m+n}}.
\label{pulse:eq:2.37}
\end{equation}

Between separated label collars, \eqref{pulse:eq:2.36} has \(H^m\)-operator
norm \(O_{m,N}(L^{-N})\) times a polynomial in \(k\), also after short
extension. Appendix~\ref{app:pulse-tools} gives the kernel estimate.

\noindent\textbf{Corollary \PulseResult{2.10}} (Finite dimension)\textbf{.} \textit{For} \(L_j=\varepsilon_j^{-1/1320}\)\textit{, the real window dimension satisfies}
\begin{equation}
D_k\le C(1+k)^A L_k^6=C(1+k)^A\varepsilon_k^{-1/220}.\label{pulse:eq:2.38}
\end{equation}

\textit{Proof.} There are \(O(\log k)\) rows and \(O(k^{18})\) slow boxes per row. Signs and block types contribute fixed factors. Lower-cover branch counts and \(L_{j}/L_{k}\) are polynomial. Sum \(C(1 +L_{j})^{d}\) with \(d\le 6\). Reality and zero-mean conditions only reduce dimension. \hfill\(\square\)

We count only unknowns in the finite-dimensional range of the final cutoff \(S_L\). Since it need not commute with pullback, we compare the untruncated equations and their smoothing tails.

\subsection{Explicit phases, Gaussian pulses and positive covariance matching}\label{pulse:sec:2.2}
\pulsebridge{2.2}

\subsubsection{The background data used in this section}\label{pulse:sec:2.2.1}
\pulsebridge{2.2.1}

Fix a native band \(k\). Use the coordinates and normalized fields of
\eqref{pulse:eq:2.5}--\eqref{pulse:eq:2.8}, and write
\(x=(\widehat R,\widehat Z,\widehat T)\), \(\widehat r=\sqrt{2\widehat R}\),
\(\varepsilon=\varepsilon_k\) and \(U=U_k\). The carrier frequency and
the physical stress unit satisfy
\begin{equation}
 N=N_k=\lceil\varepsilon^{-1/2}\rceil,\qquad
 1\le\varepsilon N^2\le4,\qquad
 U^2\varepsilon\,\widehat T_B=\mathsf T_B .
 \label{pulse:eq:2.39}
\end{equation}

Proposition~\ref{fin:mixed-interface} (B2)--(B4) and the common positive
factor in \eqref{pulse:eq:2.8} give, on the closed stress annulus,
\begin{equation}
c\le\widehat F_B\le C,\qquad c\le-\widehat S_B^\theta\le C,\qquad |\widehat S_B|\le C,\qquad \kappa_B:=-\frac{|\widehat S_B|^2}{\widehat F_B\widehat S_B^\theta}\ge2+c.\label{pulse:eq:2.40}
\end{equation}

On its interior there is a band-independent \(c_{\rm cone}>0\) such that
\begin{equation}
-\widehat T_B\cdot\frac{\widehat S_B}{|\widehat S_B|}\ge c_{\mathrm{cone}}|\widehat T_B|,\qquad \sqrt{\frac{\kappa_B-2}{2}}\left|\frac{\widehat T_B\cdot\widehat S_B^\perp}{\widehat T_B\cdot\widehat S_B}\right|\le1-c_{\mathrm{cone}},\qquad \widehat S_B^\perp=(-\widehat S_B^z,\widehat S_B^\theta).\label{pulse:eq:2.41}
\end{equation}

The stress extends smoothly by zero across \(R=R_{\rm in}\) and \(R=R_e\).
With the weights in \eqref{intro:edge-weight}, for each \(m\) one has
\begin{equation}
 |\widehat T_B|\ge c\omega,\qquad
 |D_x^{\le m}\widehat T_B|\le C_m\omega d_e^{-q_m}.
 \label{pulse:eq:2.42}
\end{equation}
Here \(\partial_{\widehat R}=\sigma^{-2}\partial_R\) at fixed
\(\widehat Z,\widehat T\), by \eqref{pulse:eq:2.6}--\eqref{pulse:eq:2.7}.

The uniform \(C^1\) bounds in \eqref{pulse:eq:2.8} and the strict
margins extend \eqref{pulse:eq:2.40} to a fixed enlargement. On a native
label \(\gamma=(k,\alpha,\pm)\), the clock in
\eqref{pulse:eq:2.22}--\eqref{pulse:eq:2.23} depends only on physical
time after evaluation. Thus for \(f=f(x,\theta,\tau_\gamma)\),
\begin{equation}
 (D_t+u_B\cdot\nabla)f
  =(\partial_t+u_B\cdot\nabla)_{\tau}f
       +2^{(2+\delta)k}\partial_\tau f,\qquad
 0\le\tau_\gamma\le T_k\asymp k^2 .
 \label{pulse:eq:2.44}
\end{equation}
The subscript \(\tau\) holds the clock fixed, while \(\partial_\tau\)
holds the slow variables and angle fixed. Here \(D_\#^{\le m}\) denotes derivatives of order at most \(m\) in \((x,l_r,\tau_\gamma)\) at fixed carrier angle \(\eta\).

Let the subscript \(0\) denote evaluation at the reference point of
a slow box. By the uniform derivative bounds used in \eqref{pulse:eq:2.19}, on every outer box we have
\begin{equation}
|\widehat F_B-\widehat F_{B,0}|+|\widehat S_B-\widehat S_{B,0}|\le Ck^{-6}\label{pulse:eq:2.45}
\end{equation}
The stress in \eqref{pulse:eq:2.41} is evaluated pointwise throughout the box.

\subsubsection{Two scalar phases and their slowly varying gradients}\label{pulse:sec:2.2.2}
\pulsebridge{2.2.2}

On this box put \(T=T_k\), \(S_0=\widehat S_{B,0}\),
\(F_0=\widehat F_{B,0}\), \(s_0=S_0/|S_0|\),
\(s_0^\perp=(-s_0^z,s_0^\theta)\) and
\(\kappa_0=-|S_0|^2/(F_0S_0^\theta)\). By
\eqref{pulse:eq:2.40}, \(\kappa_0\asymp1\) and \(\kappa_0-2\ge c\).
Replacing the shear by \(S_0\) in the cone costs \(O(k^{-6})\)
relative to \(|\widehat T_B|\) by \eqref{pulse:eq:2.45}. Hence for
large \(k\), uniformly up to either edge,
\begin{equation}
-\widehat T_B\cdot s_0\ge\tfrac12c_{\mathrm{cone}}|\widehat T_B|,\qquad \sqrt{\frac{\kappa_0-2}{2}}\left|\frac{\widehat T_B\cdot s_0^\perp}{\widehat T_B\cdot s_0}\right|\le1-\tfrac12c_{\mathrm{cone}}.\label{pulse:eq:2.46}
\end{equation}

Choose \(M>0\) independent of the label so that
\begin{equation}
1-\tfrac14c_{\mathrm{cone}}<\frac{M}{\sqrt{1+M^2}}<1.\label{pulse:eq:2.47}
\end{equation}
Since \(2F_0s_0^\theta+|S_0|=|S_0|(1-2/\kappa_0)>0\), we may set
\begin{equation}
\begin{aligned}&\Lambda_0^2=-2F_0s_0^\theta(2F_0s_0^\theta+|S_0|)>0,\\&\Lambda_0\asymp1,\qquad B_p^2=\frac{\Lambda_0}{\varepsilon N^2(1+M^2)^{3/2}},\qquad B_p\asymp1.\end{aligned}\label{pulse:eq:2.48}
\end{equation}

For the two pulses, choose the unrounded tangential frequencies by
\begin{equation}
\left(\frac{\widetilde p_{\theta,\pm}}{\widehat r_0},p_{z,\pm}\right)
 = B_p\left(s_0^\perp\mp\frac{M}{T}\frac{S_0}{|S_0|^2}\right).\label{pulse:eq:2.49}
\end{equation}
Choose \(Np_{\theta,\pm}\in\mathbb Z\setminus\{0\}\) with
\(|p_{\theta,\pm}-\widetilde p_{\theta,\pm}|\le N^{-1}\), and set
\(p_{r,\pm}=\pm B_p M/2\). With these constant label parameters the phase is
\begin{equation}
\phi_\pm=p_{\theta,\pm}\theta+2^{\delta k}p_{z,\pm}\widehat Z
 +p_{r,\pm}\widehat r-\tau_\pm(p_{\theta,\pm}\widehat F_B+p_{z,\pm}\widehat u_B^z).\label{pulse:eq:2.50}
\end{equation}
The integer angular frequency makes \(e^{iN\phi_\pm}\)
\(2\pi\)-periodic. In the orthonormal cylindrical frame
\(D_r\tau_\pm=\partial_z\tau_\pm=0\), so
\begin{equation}
\begin{aligned}
 n_\pm:=2^{-k}\nabla\phi_\pm
 &=\left(p_{r,\pm},\frac{p_{\theta,\pm}}{\widehat r},p_{z,\pm}\right)
 -\tau_\pm\left[\left(\frac{p_{\theta,\pm}}{\widehat r},p_{z,\pm}\right)
                   \cdot\widehat S_B\right]e_r\\
 &\quad-\varepsilon\tau_\pm
   (p_{\theta,\pm}\partial_{\widehat Z}\widehat F_B
   +p_{z,\pm}\partial_{\widehat Z}\widehat u_B^z)e_z .
\end{aligned}\label{pulse:eq:2.51}
\end{equation}
The clock derivative in \eqref{pulse:eq:2.44} cancels the leading
angular and axial transport, leaving the phase defect
\begin{equation}
\begin{aligned}
 \mathfrak d_\pm:=2^{-(2+\delta)k}(D_t+u_B\cdot\nabla)\phi_\pm
 &=-\varepsilon\tau_\pm\bigl[p_{\theta,\pm}(-\partial_{\widehat T}
       +\widehat u_B^z\partial_{\widehat Z})\widehat F_B\\
 &\qquad+p_{z,\pm}(-\partial_{\widehat T}
       +\widehat u_B^z\partial_{\widehat Z})\widehat u_B^z\bigr]
       +\widehat u_B^r n_\pm^r .
\end{aligned}\label{pulse:eq:2.52}
\end{equation}
By \eqref{pulse:eq:2.8},
\(|\mathfrak d_\pm|\le C\varepsilon(1+\tau_\pm)\) and
\(|D_\#^{\le m}\mathfrak d_\pm|\le P_m(k)\varepsilon\).
This term enters the full force in \eqref{pulse:eq:3.38}.

\noindent\textbf{Lemma \PulseResult{2.11}} (Wavevector and principal coefficient bounds)\textbf{.} \textit{Put} \(s_\pm(\tau)=\pm M(1/2+\tau/T)\)\textit{. Throughout the outer slow box and} \(0\le \tau \le T\)\textit{,}
\begin{equation}
\begin{gathered}
 |n_\pm-B_p(s_\pm,s_0^\perp)|+|\partial_\tau n_\pm|\le Ck^{-2},\\
 c\le |n_\pm^{\theta,z}|\le |n_\pm|\le C.
\end{gathered}\label{pulse:eq:2.53}
\end{equation}
\textit{All fixed derivatives of} \(n_\pm\)\textit{,} \(n_\pm/|n_\pm|^2\) \textit{and the principal coefficients defined below are polynomially bounded.}

\textit{Proof.} Before rounding, the radial scalar product in \eqref{pulse:eq:2.51} is \(\mp B_p M/T\) at the reference point. Its error is \(O(N^{-1}+k^{-6})\), so the wavevector error is \(O(TN^{-1}+Tk^{-6}+T^{-1}+\varepsilon T)=O(k^{-2})\). Differentiation in \(\tau\) gives the same bound. Since \(B_p\asymp1\), the tangential component stays nonzero. All further derivatives cost only powers of \(T\asymp k^2\) and the polynomial cover factors of Lemma~\ref{pulse:thm:2.9}. Differentiating \(|n|^{-2}\) proves the rest. \hfill\(\square\)

\subsubsection{The exact constrained principal equation}\label{pulse:sec:2.2.3}
\pulsebridge{2.2.3}

Fix either pulse and suppress its sign where convenient. Set
\begin{equation}
K_B=\begin{pmatrix}
 0 & -2\widehat F_B & 0\\
 2\widehat F_B+\widehat S_B^\theta & 0 & 0\\
 \widehat S_B^z & 0 & 0
\end{pmatrix},\qquad d=\varepsilon N^2|n|^2.\label{pulse:eq:2.54}
\end{equation}

Take the real velocity \(Uh\cos(N\phi)\) and pressure \(U^2\operatorname{Re}(\pi e^{iN\phi})\). In units of \(2^{(3+2\delta )k}\), the principal equation is
\begin{equation}
\partial_\tau h+K_Bh+dh+iNn\pi=0,\qquad n\cdot h=0.\label{pulse:eq:2.55}
\end{equation}
In the common unit \(UD_t\tau=U^2 2^k=F_k/\varepsilon\), viscosity contributes
\(d\) and the pressure gradient contributes \(iNn\pi\). The entries of \(K_B\)
come from \(-2u_B^\theta/r,\ \partial_r u_B^\theta+u_B^\theta/r,\ \partial_r u_B^z\)
divided by \(D_t\tau\). The remaining terms are in \eqref{pulse:eq:3.38}.

Differentiating the constraint and contracting the equation with \(n\) gives the exact pressure and projected evolution
\begin{equation}
\pi=\frac iN\frac{n\cdot K_Bh-(\partial_\tau n)\cdot h}{|n|^2},\quad
\partial_\tau h=(A-dI)h,\quad A=-K_B+\frac{n(n^\top K_B-\partial_\tau n^\top)}{|n|^2}.\label{pulse:eq:2.56}
\end{equation}
Since \(n^\top A=-(\partial_\tau n)^\top\), the constraint is preserved. Real data give real \(h\), imaginary \(\pi\), and
\begin{equation}
\tfrac12\partial_\tau |h|^2=-h^r\widehat S_B\cdot(h^\theta,h^z)-d|h|^2.\label{pulse:eq:2.57}
\end{equation}

Write \(n_{t}= (n^{\theta }, n^{z})\) and \(n_t^\perp=(-n^z,n^\theta)\). An exact orthonormal frame for \(n^{\perp }\) is
\begin{equation}
e_x=\left(0,\frac{n_t^\perp}{|n_t|}\right),\qquad
 e_y=\frac{(|n_t|^2,-n^r n_t)}{|n|\,|n_t|},\qquad Q=[e_x\mid e_y].\label{pulse:eq:2.58}
\end{equation}
If \(h=QX\) with \(X=(\bar x,\bar y)^\top\), the equation becomes
\begin{equation}
\partial_\tau X_\pm=(B_\pm+E_\pm-d_\pm I)X_\pm,\qquad
 B_\pm=(1+s_\pm^2)^{-1/2}\begin{pmatrix}0&a_0\\ b_0&0\end{pmatrix},\label{pulse:eq:2.59}
\end{equation}
where
\[
a_0=2F_0s_0^\theta+|S_0|>0,\qquad b_0=-2F_0s_0^\theta>0,\qquad
 \Lambda_0=\sqrt{a_0b_0},\qquad r_*=\sqrt{\frac{b_0}{a_0}}=\sqrt{\frac{2}{\kappa_0-2}}.
\]
The remainder satisfies
\begin{equation}
\begin{gathered}|E_\pm|+|d_\pm-d_{\rm ref}|\le C/T,\qquad
 d_{\rm ref}(\tau)=\frac{\Lambda_0(1+s_\pm(\tau)^2)}{(1+M^2)^{3/2}},\\
 |D_\#^{\le m}E_\pm|+|D_\#^{\le m}d_\pm|\le P_m(k).
\end{gathered}\label{pulse:eq:2.60}
\end{equation}

Before damping, the moving-frame matrix is
\(-Q^\top K_BQ+(\partial_\tau Q)^\top Q\).
Its frozen part uses \(F_0,S_0\) and the reference frame
\[
 e_x^{\mathrm{ref}}=(0,-s_0),\qquad e_y^{\mathrm{ref}}=\frac{(1,-s_\pm s_0^\perp)}{\sqrt{1+s_\pm^2}}
\]
and equals \(B_\pm\). Lemma~\ref{pulse:thm:2.11} bounds the
frame error and its clock derivative by \(C/T\). For the damping,
\(|n_\pm|^2=B_p^2(1+s_\pm^2)+O(T^{-1})\), so
\eqref{pulse:eq:2.39} and \eqref{pulse:eq:2.48} give
\(|d_\pm-d_{\rm ref}|\le C/T\). Differentiating the frame and
damping formulas gives the fixed-order polynomial bounds.

\subsubsection{A Gaussian envelope and its propagator bound}\label{pulse:sec:2.2.4}

The two signs give the same reference growth rate because \(s_+^2=s_-^2\). Define
\begin{equation}
a_{\rm ref}(\tau)=\frac{\Lambda_0}{\sqrt{1+s_+(\tau)^2}}-
\frac{\Lambda_0(1+s_+(\tau)^2)}{(1+M^2)^{3/2}},\qquad
G(\tau)=\exp\!\left(\int_{T/2}^\tau a_{\rm ref}(u)\,du\right).\label{pulse:eq:2.61}
\end{equation}

\noindent\textbf{Proposition \PulseResult{2.12}} (Gaussian envelope and uniform principal propagator)\textbf{.} \textit{There are fixed} \(c, C >0\) \textit{such that}
\begin{equation}
e^{-C(\tau-T/2)^2/T}\le G(\tau)\le e^{-c(\tau-T/2)^2/T}\le1.\label{pulse:eq:2.62}
\end{equation}
\textit{Let} \(\mathcal U_1(\tau,s):n(s)^\perp\to n(\tau)^\perp\) \textit{be the propagator of} \eqref{pulse:eq:2.55}\textit{. For} \(0\le s\le \tau \le T\)\textit{,}
\begin{equation}
\|\mathcal U_1(\tau,s)\|\le C\frac{G(\tau)}{G(s)},\qquad d(\tau)\ge d_*>0.\label{pulse:eq:2.63}
\end{equation}
\textit{The constants are uniform in the native label. Every positive-order coefficient derivative of \(\log G\) is polynomially bounded.}

\begin{proof}
At \(\tau=T/2\), \(s_+=M\) and \(a_{\rm ref}=0\). On the full clock,
\[
 \partial_\tau a_{\rm ref}(\tau)
 =-\frac{\Lambda_0 M s_+(\tau)}{T}
  \left[(1+s_+(\tau)^2)^{-3/2}+2(1+M^2)^{-3/2}\right].
\]
Thus \(-C/T\le\partial_\tau a_{\rm ref}\le-c/T\). Integrating twice from \(T/2\) proves \eqref{pulse:eq:2.62}. Differentiation of \eqref{pulse:eq:2.61} gives the polynomial bounds for \(\log G\), also after pullback to a common cover.

The time-independent matrix \(D=\operatorname{diag}(r_*,1)\) is uniformly invertible and satisfies
\[
 DB_\pm D^{-1}=\frac{\Lambda_0}{\sqrt{1+s_\pm^2}}
 \begin{pmatrix}0&1\\1&0\end{pmatrix}.
\]
Its larger eigenvalue is \(a_{\rm ref}+d_{\rm ref}\). Taking the real part of the inner product
with \(DX\) and using \eqref{pulse:eq:2.60} gives
\[
 \tfrac12\partial_\tau |DX|^2\le\big(a_{\rm ref}(\tau)+C/T\big)|DX|^2.
\]
Since \(0\le\tau-s\le T\), integration gives the factor
\(e^CG(\tau)/G(s)\). Returning to the original frame proves the propagator
bound, and \eqref{pulse:eq:2.39}, \eqref{pulse:eq:2.53} give \(d\ge d_*>0\).
\end{proof}

\subsubsection{Real primary amplitudes and their flux directions}\label{pulse:sec:2.2.5}
\pulsebridge{2.2.5}

To normalize the initial radial component to \(G(0)\), prescribe the following coordinates in the frame \(Q_\pm\).
\begin{equation}
 \bar x_\pm(0)=\frac{\sqrt{1+s_\pm(0)^2}}{r_*}\,G(0),\qquad
 \bar y_\pm(0)=\frac{|n_\pm(0)|}{|n_\pm^{\theta,z}(0)|}\,G(0).
 \label{pulse:eq:2.64}
\end{equation}
By \eqref{pulse:eq:2.53} and \eqref{pulse:eq:2.58} these data give \(h^r_\pm(0)=G(0)\) and \(\bar y_\pm(0)/\bar x_\pm(0)=r_*+O(T^{-1})\). Solve \eqref{pulse:eq:2.59}, set \(h_\pm=Q_\pm X_\pm\), and recover \(\pi_\pm\) from \eqref{pulse:eq:2.56}.

\noindent\textbf{Proposition \PulseResult{2.13} (Primary amplitudes).}
\textit{The real velocity amplitudes \(h_\pm\) satisfy}
\begin{equation}
 \begin{aligned}
 cG(\tau)&\le h_\pm^r(\tau)\le |h_\pm(\tau)|\le CG(\tau),\\
 \frac{h_\pm^{\theta,z}}{h_\pm^r}
   &=-s_\pm(\tau)s_0^\perp
     -\sqrt{\frac{(\kappa_0-2)(1+s_\pm(\tau)^2)}2}\,s_0+O(T^{-1}).
 \end{aligned}
 \label{pulse:eq:2.65}
\end{equation}
\textit{The amplitudes are independent of \(l_r\). On the native cover and on every common logarithmic-window cover, for fixed \(m,j\),}
\begin{equation}
 |D_x^{\le m}\partial_\tau^j h_\pm|\le P_{m,j}(k)G,\qquad
 |D_x^{\le m}\partial_\tau^j\pi_\pm|\le N^{-1}P_{m,j}(k)G.
 \label{pulse:eq:2.66}
\end{equation}
\begin{proof}
Up to the first possible zero of \(\bar x\), put \(\rho=\bar y/\bar x-r_*\).
The scalar damping cancels in \eqref{pulse:eq:2.59}, giving
\[
 \partial_\tau\rho
 =-\frac{2\Lambda_0\rho+a_0\rho^2}{\sqrt{1+s_\pm^2}}
    +O(T^{-1})(1+|\rho|+\rho^2).
\]
By \eqref{pulse:eq:2.64}, \(|\rho(0)|\le C/T\). Choose \(K>C\) large enough
that the linear term points inward and dominates the error at \(\rho=\pm K/T\)
for large \(T\). Hence \(|\rho|\le K/T\). The first component gives
\(\partial_\tau\log(\bar x/G)=O(T^{-1})\), so \(cG\le\bar x\le CG\) on this
interval. Since \(G>0\), no zero occurs before or at \(T\).
The frame gives \(h^r=\bar y|n_t|/|n|\) and
\[
 \frac{h^{\theta,z}}{h^r}
 =\frac{\bar x}{\bar y}\frac{|n|}{|n_t|}\frac{n_t^\perp}{|n_t|}
       -\frac{n^r n_t}{|n_t|^2}.
\]
Substituting the estimates for \(n_\pm\) and \(\bar y/\bar x\) yields \eqref{pulse:eq:2.65}.

For a slow derivative of order \(m\), differentiating the frame
equation leaves the same principal operator. Inductively its source
is bounded by \(P_m(k)G(s)\) using \eqref{pulse:eq:2.60}. The
initial derivative is \(O(P_m(k)G(0))\) by \eqref{pulse:eq:2.64}
and the logarithmic derivative bounds for \(G\).
In variation of constants, \(G(\tau)/G(s)\) cancels the source weight
at \(s\). The remaining integration costs at most \(T\asymp k^2\).
Equation~\eqref{pulse:eq:2.59} gives the clock derivatives,
and \eqref{pulse:eq:2.56} retains \(N^{-1}\) for the pressure.
The data and coefficients are independent of \(l_r\).
Lemma~\ref{pulse:thm:2.9} gives the common-cover bounds.
\end{proof}

\paragraph{Time cutoff.}\label{pulse:sec:2.2.6}

Let \(\psi(\tau)=\chi_\tau(\tau/T)\), equal to one on \([3T/10,7T/10]\) and supported in \((T/6,5T/6)\). Multiplication of \eqref{pulse:eq:2.55} by \(\psi\) gives
\begin{equation}
 \partial_\tau(\psi h)+K_B(\psi h)+d\,\psi h+iNn\,\psi\pi
       =(\partial_\tau\psi)h .
 \label{pulse:eq:2.67}
\end{equation}
On \(\supp(\partial_\tau\psi)\), \(|\tau-T/2|\ge T/5\), so
\eqref{pulse:eq:2.62} and \eqref{pulse:eq:2.66} bound every fixed joint
derivative of the right-hand side by \(P_m(k)e^{-ck^2}\). The physical cost
\(e^{C_mk}\) leaves \(O_{m,A}(2^{-Ak})\) for every fixed \(m,A\), and the
cut-off fields extend smoothly by zero at both ends.

\subsubsection{The two averaged columns and their positive inverse}\label{pulse:sec:2.2.7}

On the rectangle of the label \(\gamma=(k,\alpha,\pm)\), define
\[
 b_\gamma=\chi_{\rm aux}(l_r)\chi_\tau(\tau_\gamma/T)\,h_\gamma\cos(N\phi_\gamma),
 \qquad H_\gamma=\big\langle\langle b_\gamma^r b_\gamma^{\theta,z}\rangle_\theta\big\rangle_Y,
 \qquad H_{k,\alpha}=[H_+\mid H_-] .
\]
We average first in the angle, then over the auxiliary torus with normalized Haar measure. The label rectangles are disjoint. We insert \(\zeta_k\zeta_{k,\alpha}\) when summing over boxes.

\noindent\textbf{Lemma \PulseResult{2.14} (Integrated covariance columns).}
\textit{There exist positive \(J_\pm\asymp k^{-1}\) such that}
\begin{equation}
 H_\pm=J_\pm[-A_cs_0\mp M s_0^\perp+e_\pm],\qquad
 A_c=\sqrt{\frac{(\kappa_0-2)(1+M^2)}2},\qquad |e_\pm|\le Ck^{-1}.
 \label{pulse:eq:2.68}
\end{equation}
\textit{On the outer slow box the matrix \(H_{k,\alpha}\) is invertible and}
\begin{equation}
 \begin{gathered}
 |\det H_{k,\alpha}|\ge ck^{-2},\qquad \|H_{k,\alpha}^{-1}\|\le Ck,\\
 |D_x^{\le m}H_{k,\alpha}|+|D_x^{\le m}H_{k,\alpha}^{-1}|\le P_m(k).
 \end{gathered}
 \label{pulse:eq:2.69}
\end{equation}
\begin{proof}
The integer \(Np_{\theta,\gamma}\) is nonzero, so the angular average of \(\cos^2(N\phi_\gamma)\) is exactly \(1/2\). Changing variables on the native rectangle to \((l_r,\tau)\) gives
\begin{equation}
 \begin{pmatrix}H_\gamma\\ J_\gamma\end{pmatrix}
 =\frac{r_0|\det(v_r,v_t)|}{T}
 \left(\int\chi_{\rm aux}(l_r)^2\,dl_r\right)
 \int_0^T\chi_\tau(\tau/T)^2 h_\gamma^r
 \begin{pmatrix} h_\gamma^{\theta,z}\\ h_\gamma^r\end{pmatrix}\,d\tau .
 \label{pulse:eq:2.70}
\end{equation}
Cover pullback preserves these averages. By the two-sided Gaussian
bound, \(h_\gamma^r\asymp G\), and \(\chi_\tau=1\) near the clock center,
the integral in \(J_\gamma\) is comparable to \(\sqrt T\). Inserting
\(|\tau-T/2|/T\) gives an integral of order \(O(1)\). Thus
\(J_\gamma\asymp T^{-1/2}\asymp k^{-1}\). Since \(s_\pm(T/2)=\pm M\),
\eqref{pulse:eq:2.65} bounds the difference from the center direction
\(-A_cs_0\mp M s_0^\perp\) by \(C(|\tau-T/2|/T+T^{-1})\).
Its weighted average is \(O(T^{-1/2})\), giving \eqref{pulse:eq:2.68}.
In the basis \(s_0,s_0^\perp\),
\[
 [s_0\mid s_0^\perp]^{\mathsf T} H_{k,\alpha}\operatorname{diag}(J_+^{-1},J_-^{-1})
   =\begin{pmatrix}-A_c&-A_c\\-M&M\end{pmatrix}+O(k^{-1}).
\]
The leading determinant is \(-2A_cM\), uniformly separated from zero. This gives the first line of \eqref{pulse:eq:2.69}. Differentiate \eqref{pulse:eq:2.70} and \(H_{k,\alpha}H_{k,\alpha}^{-1}=I\) to obtain the second.
\end{proof}

\noindent\textbf{Proposition \PulseResult{2.15} (Positive primary matching).}
\textit{Fix a sufficiently high band and a slow-box index \(\alpha\). On the intersection of the outer slow box for \(\alpha\) with the open stress annulus, both entries in}
\begin{equation}
 \begin{pmatrix}g_+^2\\g_-^2\end{pmatrix}=H_{k,\alpha}^{-1}\widehat T_B
 \label{pulse:eq:2.71}
\end{equation}
\textit{are strictly positive. For every fixed \(m\) there are polynomial bounds}
\begin{equation}
 \begin{aligned}
 ck|\widehat T_B|&\le g_\gamma^2\le Ck|\widehat T_B|,\qquad g_\gamma^2\ge ck\omega,\\
 |D_x^{\le m}g_\gamma^2|&\le P_m(k)\omega d_e^{-q'_m},\qquad
 |D_x^{\le m}g_\gamma|\le P_m(k)\sqrt{\omega}\,d_e^{-q'_m}.
 \end{aligned}
 \label{pulse:eq:2.72}
\end{equation}
\textit{The positive square roots are independent of the auxiliary variables and the carrier angle and extend smoothly by zero across both radial edges.}
\begin{proof}
Write \(\widehat T_B=T_n s_0+T_\perp s_0^\perp\). For the leading columns of
\eqref{pulse:eq:2.71} the two scalar equations and their solution are
\[
 -A_c(c_+^0+c_-^0)=T_n,\qquad
 -M(c_+^0-c_-^0)=T_\perp,\qquad
 c_\pm^0=\frac12\left(-\frac{T_n}{A_c}\mp\frac{T_\perp}{M}\right).
\]
The choice of \(M\) in \eqref{pulse:eq:2.46}--\eqref{pulse:eq:2.47} gives the strict cone bounds
\[
 T_n\le-c|\widehat T_B|,\qquad \frac{A_c|T_\perp|}{M(-T_n)}\le 1-c'
\]
with fixed \(c,c'>0\), so \(c_\pm^0\ge c''|\widehat T_B|>0\) and
\(c_\pm^0\le C|\widehat T_B|\). By \eqref{pulse:eq:2.68}--\eqref{pulse:eq:2.69},
\(J_\pm g_\pm^2=c_\pm^0+O(k^{-1}|\widehat T_B|)\). On sufficiently high bands,
\(J_\pm^{-1}\asymp k\) and the lower bound on
\(|\widehat T_B|\) in \eqref{pulse:eq:2.42} give the first line
of \eqref{pulse:eq:2.72}.

Differentiate \eqref{pulse:eq:2.71} using \eqref{pulse:eq:2.69} and
\eqref{pulse:eq:2.42}. Differentiating the positive square root of \(g_\gamma^2\)
gives terms with \(j\) differentiated factors of \(g_\gamma^2\), multiplied by
\(C(g_\gamma^2)^{1/2-j}\). The factors contribute \(\omega^j\), and \(g_\gamma^2\ge ck\omega\) gives the
last bound of \eqref{pulse:eq:2.72} and flatness at both edges.
\end{proof}

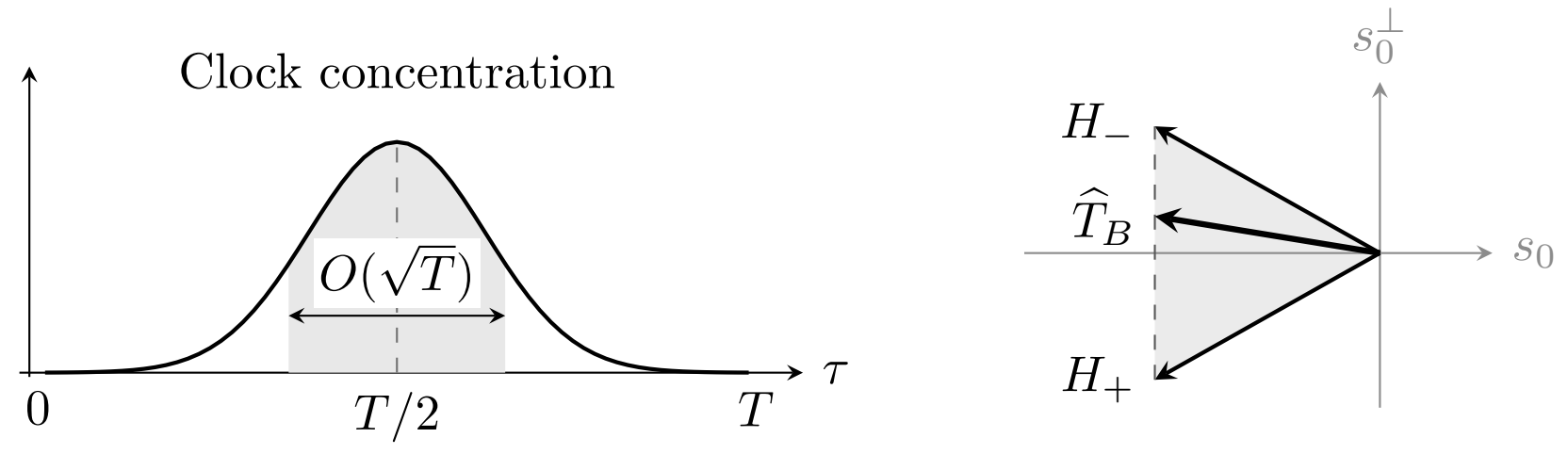
\begin{figure}[htbp]
\centering
\begin{tikzpicture}[font=\small,x=1cm,y=1cm,>=stealth]
 \begin{scope}
  \draw[->] (-2.65,0)--(2.85,0) node[right] {$\tau$};
  \draw[->] (-2.58,-.03)--(-2.58,2.15);
  \fill[black!9] plot[domain=-.76:.76,samples=35] (\x,{1.62*exp(-1.3*\x*\x)})
     --(.76,0)--(-.76,0)--cycle;
  \draw[thick] plot[domain=-2.47:2.47,samples=65] (\x,{1.62*exp(-1.3*\x*\x)});
  \draw[dashed,black!55] (0,0)--(0,1.62);
  \node[below] at (-2.52,0) {$0$};
  \node[below] at (0,0) {$T/2$};
  \node[below] at (2.52,0) {$T$};
  \draw[<->] (-.76,.40)--(.76,.40);
  \node[above,fill=white,inner sep=1pt] at (0,.45) {$O(\sqrt T)$};
  \node at (0,2.12) {Clock concentration};
 \end{scope}
 \begin{scope}[shift={(6.9,.84)},x=1.58cm,y=.89cm]
  \fill[black!8] (0,0)--(-1,-1)--(-1,1)--cycle;
  \draw[->,black!45] (-1.58,0)--(.50,0) node[right] {$s_0$};
  \draw[->,black!45] (0,-1.22)--(0,1.35) node[above] {$s_0^\perp$};
  \draw[thick,->] (0,0)--(-1,1);
  \draw[thick,->] (0,0)--(-1,-1);
  \draw[dashed,black!65] (-1,-1)--(-1,1);
  \draw[->,very thick] (0,0)--(-1,.29);
  \node[left] at (-1,1) {$H_-$};
  \node[left] at (-1,-1) {$H_+$};
  \node[left] at (-1,.29) {$\widehat T_B$};
 \end{scope}
\end{tikzpicture}
\caption{Gaussian concentration in the pulse clock (left) and the two leading flux directions (right). In units \(A_c,M\) along \(s_0,s_0^\perp\), the columns \(H_+,H_-\) point to leading order along \((-1,-1)\) and \((-1,1)\), respectively. The angular gap keeps both weights for \(\widehat T_B\) positive under \(O(k^{-1})\) direction errors.}
\label{fig:pulse-covariance}
\end{figure}

Insert these amplitudes into the physical leading waves by setting
\begin{equation}
 w_\gamma^{\rm pre}=U_k\varepsilon_k^{1/2}\zeta_k\zeta_{k,\alpha} g_\gamma
     \chi_{\rm aux}\chi_\tau h_\gamma\cos(N_k\phi_\gamma),\qquad
 W^{\rm pre}=\sum_\gamma w_\gamma^{\rm pre}.
 \label{pulse:eq:2.73}
\end{equation}
For each \(\lambda>0\) the sum is locally finite. Distinct active labels have vanishing cross products. By \eqref{pulse:eq:2.71} and the two squared partitions of unity,
\begin{equation}
 \mathcal C(W^{\rm pre})=
 \big\langle\langle(W^{\rm pre})^r(W^{\rm pre})^{\theta,z}\rangle_\theta\big\rangle_Y
 =\sum_{k,\alpha}\zeta_k^2\zeta_{k,\alpha}^2U_k^2\varepsilon_k\,\widehat T_{B,k}
 =\mathsf T_B .
 \label{pulse:eq:2.74}
\end{equation}
The normalization is \eqref{pulse:eq:2.39}. Figure~\ref{fig:pulse-covariance}
shows the normalized columns and the target direction.

\subsubsection{Exact curls and the coefficient bounds used later}\label{pulse:sec:2.2.8}
\pulsebridge{2.2.8}

Put \(a_\gamma=\varepsilon^{1/2}\zeta_k\zeta_{k,\alpha}g_\gamma\chi_{\rm aux}\chi_\tau h_\gamma\). We enforce incompressibility by using the physical velocity potential
\begin{equation}
 \mathcal A_\gamma=-\frac{U}{2^k N}\,
       \frac{n_\gamma\times a_\gamma}{|n_\gamma|^2}\sin(N\phi_\gamma),
 \qquad w_\gamma=\curl\mathcal A_\gamma .
 \label{pulse:eq:2.75}
\end{equation}
Evaluate the potential in physical variables. Let
\(\mathsf C=2^{-k}\curl_{\rm coef}\) differentiate the coefficients and the
cylindrical frame with the carrier angle held fixed. By the product rule,
\begin{equation}
 \begin{aligned}
 w_\gamma&=U\big(a_\gamma\cos(N\phi_\gamma)+b_\gamma^{\rm corr}\sin(N\phi_\gamma)\big),\\
 b_\gamma^{\rm corr}&=-N^{-1}\mathsf C\!\left(\frac{n_\gamma\times a_\gamma}{|n_\gamma|^2}\right).
 \end{aligned}
 \label{pulse:eq:2.76}
\end{equation}
Indeed, \(2^{-k}\nabla\phi_\gamma=n_\gamma\) and
\(-n_\gamma\times(n_\gamma\times a_\gamma)/|n_\gamma|^2=a_\gamma\)
by the constraint \(n_\gamma\cdot a_\gamma=0\). The smooth zero extension
of \(\mathcal A_\gamma\) gives a divergence-free pulse supported with its potential.

Take the corresponding real pressure
\[
 p_\gamma=U^2\operatorname{Re}\!\left(
 \varepsilon^{1/2}\zeta_k\zeta_{k,\alpha}g_\gamma\chi_{\rm aux}\chi_\tau\pi_\gamma
 e^{iN\phi_\gamma}\right).
\]
Equation~\eqref{pulse:eq:2.55} gives the principal pressure cancellation. Of the factors multiplying \(h_\gamma\), only \(\chi_\tau\) depends on \(\tau_\gamma\), so the remaining principal cutoff term is precisely \eqref{pulse:eq:2.67}.

\noindent\textbf{Corollary \PulseResult{2.16} (Size and covariance of an exact pulse).}
\textit{For every fixed coefficient order \(m\),}
\begin{equation}
 \begin{aligned}
 |D_\#^{\le m}a_\gamma|
   &\le P_m(k)\varepsilon^{1/2}\sqrt\omega\,d_e^{-q''_m}G_\gamma,\\
 |D_\#^{\le m}(\varepsilon^{1/2}\zeta_k\zeta_{k,\alpha}g_\gamma\chi_{\rm aux}\chi_\tau\pi_\gamma)|
   &\le P_m(k)\varepsilon\sqrt\omega\,d_e^{-q''_m}G_\gamma,\\
 |D_\#^{\le m} b_\gamma^{\rm corr}|
   &\le P_m(k)\varepsilon^{1-\delta_0}\sqrt\omega\,d_e^{-q''_m}G_\gamma .
 \end{aligned}
 \label{pulse:eq:2.77}
\end{equation}
\textit{On the retained charts \(d_e\ge c/k\), division by \(\sqrt\omega G_\gamma\) leaves bounds polynomial in \(k\) with the displayed powers of \(\varepsilon\). The exact angular covariance is}
\begin{equation}
 \langle w_\gamma\otimes w_\gamma\rangle_\theta=
   \frac{U^2}{2}\big(a_\gamma\otimes a_\gamma+
      b_\gamma^{\rm corr}\otimes b_\gamma^{\rm corr}\big).
 \label{pulse:eq:2.78}
\end{equation}
\textit{The resulting double-mean radial flux differs from \(\mathsf T_B\) by at most \(U^2P_m(k)\varepsilon^{2-2\delta_0}\omega\) on the retained charts, with the same bound at every fixed coefficient order.}
\begin{proof}
Combine \eqref{pulse:eq:2.66} and \eqref{pulse:eq:2.72}. The pressure carries
\(N^{-1}\le\varepsilon^{1/2}\). One coefficient curl in \eqref{pulse:eq:2.76}
costs \(P_m(k)\varepsilon^{-\delta_0}\), giving the last bound in
\eqref{pulse:eq:2.77}. The logarithmic derivatives of \(\omega,G_\gamma\)
on the retained charts are polynomial by \eqref{pulse:eq:2.42} and
\eqref{pulse:eq:2.61}; the original coefficients are flat at both edges.
The angular means of \(\cos^2,\sin^2,\cos\sin\) are \(1/2,1/2,0\).
Thus \eqref{pulse:eq:2.78}, disjoint label supports and \eqref{pulse:eq:2.74}
leave only the square of the curl correction in the averaged flux.
\end{proof}
\section{Compatible sources and initialization}\label{pulse:sec:3}
\pulsebridge{3}

\subsection{Compatible mean sources}\label{pulse:sec:3.1}
\pulsebridge{3.1}

In this subsection \(f=\langle\mathrm{NS}(u,p)\rangle\) since \(E_B=0\). Here \(\langle h\rangle=\int_{\mathbb T^2}\langle h\rangle_\theta\,dY\) uses normalized Haar measure before \(Y=Y_i(t,r)\) is evaluated. Subtracting an error of zero double mean preserves the identities. The radially compact velocity corrections preserve
\begin{equation}
\int_0^\infty r^2\langle u^\theta-u_B^\theta\rangle\,dr=0,\qquad \int_0^\infty r\langle u^z-u_B^z\rangle\,dr=0.
\label{pulse:eq:3.1}
\end{equation}

At fixed \((t,z)\) define the radial force defect \(P=\int_0^\infty f^r\,dr\) and the two flux defects
\[
\begin{aligned}
 J^{\theta z}&=\int_0^\infty r^2\langle u^\theta u^z-u_B^\theta u_B^z\rangle\,dr,\\
 J^{z\theta}&=\int_0^\infty r\left\langle (u^z)^2-(u_B^z)^2+p-p_B+
                 \tfrac r2(f^r-f_B^r)\right\rangle\,dr.
\end{aligned}
\]

The background has \(f_B^r=0\). Its compactly supported stress gives zero tangential force integrals with weights \(r^2,r\). Subtracting the background equations and integrating in \(r\), we obtain
\begin{equation}
\int r^2 f^\theta\,dr=\partial_z J^{\theta z},\qquad \int r f^z\,dr=\partial_z\left(J^{z\theta}-\frac12\int r^2f^r\,dr\right).
\label{pulse:eq:3.2}
\end{equation}
Compact support removes the boundary terms, and \eqref{pulse:eq:3.1} removes the time and axial viscosity terms.

\noindent\textbf{Proposition \NumberedPulseResult{3.1} (Free radial potentials).}
\textit{Choose smooth interior bumps with \(\int r^2b^\theta\,dr=\int rb^z\,dr=1\).
Put \(C^\theta=\partial_zJ^{\theta z}\) and \(C^z=\partial_z(J^{z\theta}-\tfrac12\int r^2f^r\,dr)\).
At fixed \(f^r,J^{\theta z},J^{z\theta}\), the map}
\begin{equation}
f^\theta=r^{-2}\partial_r g^\theta+b^\theta C^\theta,\qquad f^z=r^{-1}\partial_r g^z+b^z C^z
\label{pulse:eq:3.3}
\end{equation}
\textit{takes compact radial potentials bijectively onto the tangential forces satisfying \eqref{pulse:eq:3.2}. Their support can include the fixed bump window.}

\textit{Proof.} Integrating \eqref{pulse:eq:3.3} gives
\begin{equation}
\begin{aligned}g^\theta(r)&=\int_0^r s^2(f^\theta-b^\theta C^\theta)\,ds=-\int_r^\infty s^2(f^\theta-b^\theta C^\theta)\,ds,\\g^z(r)&=\int_0^r s(f^z-b^zC^z)\,ds=-\int_r^\infty s(f^z-b^zC^z)\,ds.\end{aligned}
\label{pulse:eq:3.4}
\end{equation}
The two tails agree by \eqref{pulse:eq:3.2} and vanish outside the radial hull of the source and bumps. Differentiation recovers the force. The integrals preserve smooth parameter dependence, and compact support gives uniqueness. \hfill\(\square\)

After the pressure adjustment the radial force equals \(\rho_r P\). Take \(\rho_r=(\sqrt2\,\lambda)^{-1}\rho_*(r/(\sqrt2\,\lambda))\), with \(\int \rho_*=1\), and set \(c_\rho=\int s^2\rho_*(s)\,ds\). Then \(\tfrac12\int r^2\rho_r\,dr=c_\rho\lambda^2\), and the axial coefficient becomes
\begin{equation}
C^z=\partial_zJ^{z\theta}-c_\rho\lambda^2\partial_z P-2c_\rho\lambda(\partial_z\lambda)P.\label{pulse:eq:3.5}
\end{equation}

For an intermediate source, use its full radial force in \eqref{pulse:eq:3.2}.
Since \(r^2dr=\lambda^3\sqrt{2R}\,dR\) and \(r\,dr=\lambda^2\,dR\), the other compensators are
\[
 b^\theta=\lambda^{-3}\beta^\theta(R),\quad b^z=\lambda^{-2}\beta^z(R),\qquad
 \int\sqrt{2R}\,\beta^\theta\,dR=\int\beta^z\,dR=1.
\]
Apply smoothing and radial cutoffs to the potentials before \eqref{pulse:eq:3.3}.
Smooth \(P,J^{\theta z},J^{z\theta}\) in their own variables before taking \(\partial_z\).
Reconstruction costs one normalized derivative. Its inverse has order zero by
\eqref{pulse:eq:3.4}, since \(C^\theta=\int r^2 f^\theta\,dr\) and \(C^z=\int r f^z\,dr\).

\subsection{Radial and pressure inverses}\label{pulse:sec:3.2}
\pulsebridge{3.2}

On a fixed band write \(U=U_k\), \(\varepsilon=\varepsilon_k\) and \(N=N_k\). The symbol \(O_m(\varepsilon^b\omega^pG^j)\) permits fixed powers of \(k\) and \(d_e^{-1}\) at order \(m\). On \(d_e\ge c/k\) the latter cost only powers of \(k\). Estimates with \(G\) refer to a single label. Supports of distinct labels are disjoint at a common physical point.

\subsubsection{Weighted radial integration}\label{pulse:sec:3.2.1}
\pulsebridge{3.2.1}

Use \(r_{\mathrm{sim}}=r/\lambda=\sqrt{2R}\) and \(w(r_{\mathrm{sim}})=\omega(r_{\mathrm{sim}}^2/2)\), with \(\omega\) and \(d_e\) as in \eqref{intro:edge-weight}.
The derivative of \(\log\omega\) with respect to \(\log R\) is
\(2c_a/\log^3(R/R_{\rm in})-2c_h/\log^3(R_e/R)\) and decreases strictly.
Choose a fixed interior interval \(K\) containing the maximizer of \(w\) and the supports of all moment profiles and compensators. On each complementary interval \(w\) increases toward \(K\), and \(\inf_K w>0\).

For \(\alpha\in\mathbb R\) fix a smooth \(b_\alpha\) supported in \(K\) with \(\int r_{\mathrm{sim}}^\alpha b_\alpha\,dr_{\mathrm{sim}}=1\). For a smooth radial input \(g=g(r_{\mathrm{sim}})\), with the other variables held fixed, set
\begin{equation}
m_\alpha(g)=\int_{\sqrt{2R_{\rm in}}}^{\sqrt{2R_e}}s^\alpha g(s)\,ds,\qquad
 (Q_\alpha g)(r_{\mathrm{sim}})=r_{\mathrm{sim}}^{-\alpha}\int_{\sqrt{2R_{\rm in}}}^{r_{\mathrm{sim}}}s^\alpha\bigl[g(s)-b_\alpha(s)m_\alpha(g)\bigr]\,ds.
 \label{pulse:eq:3.6}
\end{equation}
Then
\begin{equation}
 (\partial_{r_{\mathrm{sim}}}+\alpha/r_{\mathrm{sim}})Q_\alpha g=g-b_\alpha m_\alpha(g).
 \label{pulse:eq:3.7}
\end{equation}

\noindent\textbf{Lemma \NumberedPulseResult{3.2}\label{pulse:thm:3.3} (Weighted radial primitives).}
\textit{The support of \(Q_\alpha g\) lies in the radial hull of \(\operatorname{supp}g\) and \(K\). Outside \(K\),}
\begin{equation}
 (Q_\alpha g)(r_{\mathrm{sim}})=\begin{cases}
 r_{\mathrm{sim}}^{-\alpha}\int_{\sqrt{2R_{\rm in}}}^{r_{\mathrm{sim}}} s^\alpha g(s)\,ds,&r_{\mathrm{sim}}<\inf K,\\
 -r_{\mathrm{sim}}^{-\alpha}\int_{r_{\mathrm{sim}}}^{\sqrt{2R_e}}s^\alpha g(s)\,ds,&r_{\mathrm{sim}}>\sup K .
 \end{cases}\label{pulse:eq:3.8}
\end{equation}
\textit{For \(p>0\), if \(g=0\) on \(d_e<c/k\), then the normalized primitive extends smoothly by zero and}
\begin{equation}
 \|w^{-p}Q_\alpha(w^p g)\|_{H^m}\le P_m(k)\|g\|_{H^m}.
 \label{pulse:eq:3.9}
\end{equation}
\textit{The norms include the passive variables. With the physical compensator \(\lambda^{-\alpha-1}b_\alpha(r/\lambda)\), the physical primitive at \(r=\lambda r_{\mathrm{sim}}\) equals \(\lambda Q_\alpha[f(\lambda\,\cdot)](r_{\mathrm{sim}})\). Both primitives commute with auxiliary pullbacks and averages.}

\textit{For the same trimmed inputs, put \(\tau_B=|\lambda^{2+\delta}\mathsf T_B|\)
and \(a_B=\tau_B/w\). Then}
\begin{equation}
 \tau_B^{-p}Q_\alpha(\tau_B^p g)
   =a_B^{-p}\,w^{-p}Q_\alpha(w^p a_B^p g),
 \label{pulse:eq:3.10}
\end{equation}
\begin{equation}
 \|\tau_B^{-p}Q_\alpha(\tau_B^p g)\|_{H^m}
       \le P_m(k)\|g\|_{H^m}.
 \label{pulse:eq:3.11}
\end{equation}
\textit{On the full annulus, relative pointwise bounds allowing algebraic edge factors remain of that form. The unnormalized positive-weight primitive extends smoothly by zero when its normalized input and derivatives have only algebraic growth at the edges.}

\textit{Proof.} The compensated integral vanishes over the full annulus. Thus either endpoint gives \eqref{pulse:eq:3.8} and the support bound. Substitution of \(r=\lambda r_{\mathrm{sim}}\) gives the physical factor \(\lambda\).

For \(G=w^{-p}Q_\alpha(w^p g)\), the kernel weight on each tail satisfies
\((w(s)/w(r_{\mathrm{sim}}))^p\le1\). On \(K\), \(w\) is bounded below, so Schur's test on the
fixed annulus gives the \(L^2\) estimate. On the trimmed hull, derivatives of
\(\log w\) cost only powers of \(k\), and
\[
 \partial_{r_{\mathrm{sim}}}G=g-b_\alpha w^{-p} m_\alpha(w^p g)
                   -(\alpha/r_{\mathrm{sim}}+p\partial_{r_{\mathrm{sim}}}\log w)G .
\]
Passive differentiation and induction in \(r_{\mathrm{sim}}\) give \eqref{pulse:eq:3.9}.
By \eqref{pulse:eq:2.8} and \eqref{pulse:eq:2.42}, \(a_B\ge c>0\)
and every fixed derivative of \(a_B^{\pm p}\) is polynomially bounded on
the trimmed strip. Applying \eqref{pulse:eq:3.9} in \eqref{pulse:eq:3.10}
gives \eqref{pulse:eq:3.11}.

At either edge apply \eqref{fin:gaussian-primitive} to the tail in
\eqref{pulse:eq:3.8} on that side of \(K\). The radial identity above
costs only algebraic edge factors under differentiation. The Gaussian
then forces every derivative of the unnormalized primitive to vanish at the edge.
\hfill\(\square\)

\subsubsection{Translated radial integration}\label{pulse:sec:3.2.2}
\pulsebridge{3.2.2}

On a common cover of level \(c\), write \(Y=Y_c\), \(\Lambda=\Lambda_r^c\), and \(D_r=\partial_r+\Lambda d_0r^{d_0-1}v_r\cdot\nabla_{Y_c}\). In the physical variable \(r\), we also denote the rescaled interval \(\lambda K\) by \(K\). Choose \(\chi _{m}= 0\) below \(K\) and \(\chi _{m}= 1\) above \(K\), with \(\supp\chi_m'\subset K\). On \(K\) we have \(\omega \ge c_{K}>0\). For a compact radial source, set
\begin{gather}
\mathcal J H(r,Y)=\int_0^\infty H(s,Y+\Lambda(s^{d_0}-r^{d_0})v_r)\,ds,\label{pulse:eq:3.12} \\
I_0 H(r,Y)=\int_0^r H(s,Y+\Lambda(s^{d_0}-r^{d_0})v_r)\,ds
-\chi_m(r)\mathcal JH(r,Y).\label{pulse:eq:3.13}
\end{gather}
The slow parameters are held fixed.

\noindent\textbf{Proposition \NumberedPulseResult{3.4}} (Compact characteristic integration)\textbf{.} \textit{The real linear map} \(I_{0}\) \textit{has support in the radial hull of the source and} \(K\)\textit{. It satisfies}
\begin{equation}
D_r I_0H=H+R[H],\qquad R[H]=-\chi_m'\mathcal{J}H.\label{pulse:eq:3.14}
\end{equation}
\textit{Its auxiliary mean is} \(\langle I_0H\rangle_Y=\int_0^r\langle H\rangle_Yds-\chi_m\int_0^\infty\langle H\rangle_Yds\)\textit{. Thus} \(\langle R[H]\rangle_Y=0\) \textit{if} \(\int\langle H\rangle_Ydr=0\)\textit{. Both maps commute with permitted pullbacks to a lower auxiliary cover. In joint normalized norms,}
\begin{equation}
\|I_0H\|_{H^m}\le 2^{-k}P_m(k)\|H\|_{H^m}.\label{pulse:eq:3.15}
\end{equation}
\textit{For} \(p >0\) \textit{and sources zero on} \(d_{e}< c/k\)\textit{, one also has}
\begin{equation}
\|\omega^{-p}I_0(\omega^p\widehat H)\|_{H^m}\le 2^{-k}P_m(k)\|\widehat H\|_{H^m}.\label{pulse:eq:3.16}
\end{equation}
\textit{For untrimmed sources,} \(H=O_{m}(H_{0}\omega ^{p})\) \textit{implies} \(I_{0}H=O_{m}(2^{-k}H_{0}\omega ^{p})\)\textit{, allowing larger algebraic edge powers. All fixed derivatives extend by zero at the original edges. If} \(\int\langle H\rangle_Ydr=0\)\textit{, then}
\begin{equation}
\|R[H]\|_{H^m}\le P_{m,J}(k)M^{-J}\|H\|_{H^{m+2J}},\qquad M=2^{-k d_0}\Lambda,\quad J\ge 1.\label{pulse:eq:3.17}
\end{equation}
\textit{The remainder lies in} \(K\)\textit{. Since} \(M\ge P(k)^{-1}\varepsilon ^{-\delta_0 }\)\textit{, it is physically flat for families with polynomial norms at every fixed order.}

\textit{Proof.} In \eqref{pulse:eq:3.13} the auxiliary part of \(D_r\)
cancels the derivative of the translation \(-\Lambda r^{d_0}v_r\).
The upper endpoint and the cutoff derivative give \eqref{pulse:eq:3.14}.
Outside the radial hull the two integrals vanish or cancel. Translation
preserves Haar measure, and \(E^{d-c}v_r=\Lambda_r^{d-c}v_r\) carries the shift
at level \(c\) to that at level \(d\) under \(Y_d=E^{d-c}Y_c\).
This proves the mean and pullback identities.

For the estimates put \(u=(2^k r)^{d_0}\) and \(\sigma=(2^k s)^{d_0}-u\). Up to \(2^{-k}/d_{0}\), the primitive is
\begin{equation}
\begin{aligned}
 &(1-\chi_m)\int_{-\infty}^{0}(u+\sigma)^{1/d_0-1}H(2^{-k}(u+\sigma)^{1/d_0},Y+M\sigma v_r)\,d\sigma\\
 &\quad{}-\chi_m\int_{0}^{\infty}(u+\sigma)^{1/d_0-1}H(2^{-k}(u+\sigma)^{1/d_0},Y+M\sigma v_r)\,d\sigma .
\end{aligned}\label{pulse:eq:3.18}
\end{equation}
Extend the source by zero and the prefactor smoothly from the positive annulus.
For each value of the integration variable \(\sigma\), \(Y\mapsto Y+M\sigma v_r\) is an isometry independent of the slow variables.
Both integrals have bounded effective length, so Minkowski's inequality gives
\eqref{pulse:eq:3.15} without derivative loss. Using the tail on the same side
of \(K\), the weight argument of Lemma~\ref{pulse:thm:3.2} gives
\eqref{pulse:eq:3.16} and its untrimmed pointwise version.

For the remainder, the zero mode of \(\mathcal{J}H\) vanishes. A nonzero mode \(n\) has the form
\[
\frac{2^{-k}}{d_0}e^{-2\pi i M u(n\cdot v_r)}\int A_n(v)e^{2\pi i M v(n\cdot v_r)}\,dv.
\]
Compact support lets us integrate by parts \(J\) times. Using
\(\lvert n\cdot v_r\rvert^{-J}\le C_J\langle n\rangle^J\) costs \(J\) radial and \(J\) auxiliary derivatives.
For further output derivatives, first use \eqref{pulse:eq:3.18}; they then fall on \(A_{n}\), without creating powers of \(M\). Parseval gives \eqref{pulse:eq:3.17}; \(\chi_m'\) cancels \(2^{-k}\). Every fixed physical derivative costs a fixed exponential in \(k\). Choosing \(J\) afterward proves flatness. \hfill\(\square\)

\subsubsection{Pressure and stress primitives}\label{pulse:sec:3.2.3}
\pulsebridge{3.2.3}

Use a nonnegative bump \(\rho_r\) on \(K\) with \(\int\rho_r\,dr=1\) and size \(O(2^k)\). For an angular-mean radial force \(H\) on cover \(c\), set
\begin{equation}
 P_H=\int_0^\infty\langle H\rangle_Y\,dr,\qquad S_H=-H+\rho_rP_H,\qquad q_H=I_0S_H.
 \label{pulse:eq:3.19}
\end{equation}

\noindent\textbf{Proposition \NumberedPulseResult{3.5}\label{pulse:thm:3.6} (Compact pressure and uniform bounds).}
\textit{The pressure \(q_H\) is supported in the radial hull of \(\operatorname{supp}H\) and \(K\), and}
\begin{equation}
 D_rq_H=-H+\rho_rP_H+R_H,\qquad R_H=-\chi_m'\mathcal J S_H,\qquad
 \langle R_H\rangle_Y=0 .
 \label{pulse:eq:3.20}
\end{equation}
\textit{If \(H/F_k=O_m(\varepsilon^b\omega)\), then}
\begin{equation}
 q_H/U^2=O_m(\varepsilon^{b+1}\omega),\qquad
 \partial_z q_H/F_k=O_m(\varepsilon^{b+1}\omega).
 \label{pulse:eq:3.21}
\end{equation}
\textit{The axial bound uses one additional slow derivative of the input. Including \(R_He_r\) in the physical remainder leaves all three defect coordinates unchanged. In particular,}
\begin{equation}
 \int_0^\infty\left(r\langle q_H\rangle_Y+\tfrac12 r^2\partial_r\langle q_H\rangle_Y\right)\,dr=0 .
 \label{pulse:eq:3.22}
\end{equation}

\textit{With \(\|\cdot\|_m=\|\cdot\|_{H^{10+m}}\) in the coefficient variables, put}
\begin{equation}
 M=2^{-kd_0}\Lambda_r^c,\qquad \int_0^\infty\langle S_H\rangle_Y\,dr=0 .
 \label{pulse:eq:3.23}
\end{equation}
\textit{The pressure and the remainder in \eqref{pulse:eq:3.20}, with
\(\operatorname{supp}R_H\subset K\), satisfy}
\begin{equation}
 \|q_H\|_m\le 2^{-k}P_m(k)\|H\|_m,\qquad
 \|\partial_zq_H\|_m\le \varepsilon P_m(k)\|H\|_{m+1},
 \label{pulse:eq:3.25}
\end{equation}
\begin{equation}
 \|R_H\|_m\le P_{m,J}(k)M^{-J}\|H\|_{m+2J}\qquad (J\ge0).
 \label{pulse:eq:3.26}
\end{equation}
\textit{The relative estimates with \(H/\omega,q_H/\omega\) hold on fixed interior windows for interior-supported data, and on trimmed radial hulls as in Lemma~\ref{pulse:thm:3.2}.}

\textit{Proof.} The bump map \(H\mapsto\rho_r P_H\) has order zero in normalized
variables and gives \(\int\langle S_H\rangle_Y\,dr=0\). Proposition~\ref{pulse:thm:3.4}
applied to \(S_H\) gives the supports and \eqref{pulse:eq:3.20}.
In \eqref{pulse:eq:3.18} slow derivatives act on the input and cutoff at fixed
shift. Thus \eqref{pulse:eq:3.15} and \(2^{-k}\partial_z=\varepsilon\partial_{\widehat Z}\)
give \eqref{pulse:eq:3.25}, including the derivative of \(\lambda\) in the bump.
Equation~\eqref{pulse:eq:3.16} gives the trimmed bound. Its untrimmed
counterpart and \(2^{-k}F_k/U^2=\varepsilon\) give \eqref{pulse:eq:3.21}.
The zero integrated mean permits \eqref{pulse:eq:3.17} for \(J\ge1\), at
a cost of \(2J\) derivatives. For \(J=0\), the \(2^{-k}\) from the full
translated integral cancels \(\chi_m'=O(2^k)\), proving \eqref{pulse:eq:3.26}.

The pressure terms in \(P\) and \(J^{z\theta}\) are the endpoint differences
of \(\langle q_H\rangle_Y\) and \(r^2\langle q_H\rangle_Y/2\).
Compact support makes both zero, giving \eqref{pulse:eq:3.22}. The mixed flux
\(J^{\theta z}\) is independent of pressure, and \(R_H\) has zero Haar mean. \hfill\(\square\)

For a tangential double mean \(h_{d}\), set \(q_{\theta }= 2\), \(q_{z}= 1\) and
\begin{equation}
T^a=-r^{-q_a}\left[\int_0^r s^{q_a}h_d^a(s)ds-\chi_m(r)\int_0^\infty s^{q_a}h_d^a(s)ds\right].\label{pulse:eq:3.28}
\end{equation}
Then
\begin{equation}
(\partial_r+q_a/r)T^a=-h_d^a+\rho^a I_a,\qquad \rho^a=r^{-q_a}\chi_m',\quad I_a=\int r^{q_a}h_d^adr.\label{pulse:eq:3.29}
\end{equation}
Lemma~\ref{pulse:thm:3.2} gives compact support, outward-tail locality and weighted order-zero bounds, with factor \(2^{-k}\). The compensators are retained. Compatibility gives
\[
I_\theta=\partial_z J^{\theta z},\qquad I_z=\partial_z\left(J^{z\theta}-\tfrac12\int r^2 h_d^rdr\right).
\]
If these bracketed quantities have size \(\varepsilon ^{\beta }\) in units \(U^{2}2^{-3k}\) and \(U^{2}2^{-2k}\), respectively, then
\begin{equation}
\rho^a I_a/F_k=O_m(\varepsilon^\beta\omega).\label{pulse:eq:3.30}
\end{equation}
Indeed, \(\partial _{z}\) costs \(2^{(1-\delta )k}\), while \(\rho ^{\theta }=O(2^{3k})\) and \(\rho ^{z}=O(2^{2k})\) on \(K\).

\subsection{Exact moment repair by pure swirl}\label{pulse:sec:3.3}

Let \(V=\langle u\rangle_{\theta,Y}\). At fixed pressure add a pure swirl \(b(t,r,z)e_\theta\), supported in \(I_3\). This increment is divergence-free. Substitution in the definitions of Section~\ref{pulse:sec:3.1} gives
\begin{equation}
\begin{aligned}
 \Delta\langle f^r\rangle_{\theta,Y}&=-\frac{2V^\theta b+b^2}{r},&
 \Delta P&=-\int\frac{2V^\theta b+b^2}{r}\,dr,\\
 \Delta J^{\theta z}&=\int r^2V^zb\,dr,&
 \Delta J^{z\theta}&=-\frac12\int r(2V^\theta b+b^2)\,dr.
\end{aligned}\label{pulse:eq:3.32}
\end{equation}
In the first equality, the angular transport \((b/r)\partial_\theta u^r\) has zero angular mean. The same averaging removes \((b/r)\partial_\theta u^z\) from the axial force. Since the axial velocity is unchanged, the increment of \(J^{\theta z}\) is linear in \(b\).

\noindent\textbf{Proposition \NumberedPulseResult{3.8} (Four coefficients repair three defects).}
\textit{Suppose the normalized tangential means on \(I_3\) are sufficiently close to the background pair. For sufficiently small normalized defects, there is a small pure swirl satisfying}
\begin{equation}
\begin{aligned}
 \int r^2b\,dr&=0,&
 \int\frac{2V^\theta b+b^2}{r}\,dr&=P,\\
 \int r^2V^zb\,dr&=-J^{\theta z},&
 \int r(2V^\theta b+b^2)\,dr&=2J^{z\theta}.
\end{aligned}\label{pulse:eq:3.33}
\end{equation}
\textit{It sets \(P,J^{\theta z},J^{z\theta}\) equal to zero, preserves both linear velocity moments, and depends smoothly on the normalized defects and means with order-zero tame bounds.}

\textit{Proof.} Use \(x=R/R_b\) and \(\lambda^{1+\delta}(V^\theta,V^z)=A(\phi,\psi)\) with the normalization in \textup{(B5)}. For the four fixed profiles \(\xi_j\) of \eqref{fin:mixed-gram}, write
\[
 b=\lambda^{-1-\delta}Ah,\qquad h=\sum_{j=1}^4c_j\xi_j,\qquad
 B_{ij}=\int_{I_3/R_b}(\sqrt x,\phi/x,\sqrt x\,\psi,\phi)_i\,\xi_j\,dx .
\]
For the background, \(B=G_3\). Lemma~\ref{fin:mixed-inverse} gives \(\|B^{-1}\|\le 2/g_3\) for the stated neighborhood. Changing variables in the four equations yields
\begin{equation}
 Bc+N(c)=d,\qquad
 N(c)=\left(0,\frac12\int h^2/x\,dx,0,\frac12\int h^2\,dx\right)^T,\label{pulse:eq:3.34}
\end{equation}
where
\begin{equation}
 d=\left(0,\frac{\lambda^{2+2\delta}P}{A^2},
 -\frac{\lambda^{-1+2\delta}J^{\theta z}}{\sqrt2 R_b^{3/2}A^2},
 \frac{\lambda^{2\delta}J^{z\theta}}{R_b A^2}\right)^T .
 \label{pulse:eq:3.35}
\end{equation}
The three radial measures are \(dr/r=dx/(2x)\), \(r^2dr=\lambda^3\sqrt2 R_b^{3/2}\sqrt{x}\,dx\), and \(r\,dr=\lambda^2 R_b\,dx\).

Polarizing \(N\) gives a bilinear form \(\mathcal Q\) of bounded norm \(C_Q\)
on the fixed window. If \(32C_Q|d|\le g_3^2\), Lemma~\ref{pr:quadratic}
with defect \(-d\) gives a unique small solution with \(|c|\le4|d|/g_3\).
On the lemma's ball \(\|B^{-1}D_c N\|\le1/2\). Differentiating the solved
equation with respect to any parameter gives
\[
 (B+D_cN)Dc=Dd-(DB)c-(D\mathcal Q)(c,c).
\]
Thus \(d=O_m(\eta)\) gives \(c,h=O_m(\eta)\) by induction, with polynomial
band constants when the means and profiles have such bounds. This includes
the derivatives of \(\lambda(t,z)\). Equation~\eqref{pulse:eq:3.32} gives the
three cancellations. The first row of \eqref{pulse:eq:3.33} and the unchanged
axial velocity preserve \eqref{pulse:eq:3.1}. \hfill\(\square\)

The quadratic rows cancel the residuals
\(-\int b^2/r\,dr\) in \(P\) and \(-\tfrac12\int rb^2\,dr\) in
\(J^{z\theta}\) that a linear repair would leave.
The weight \(\omega\) is bounded below on the swirl window \(I_3\).

\subsection{Primary force and the signed stress inverse}\label{pulse:sec:3.4}

\subsubsection{The full-curl force identity}\label{pulse:sec:3.4.1}
\pulsebridge{3.4.1}

Fix one native wave label. Write \(n=2^{-k}\nabla\phi\) and \(d=\varepsilon N^2|n|^2\). The zero-mean waveform and its pressure solve
\[
 \partial_\tau h+Kh-d\partial_\eta^2h+Nn\partial_\eta\pi
       =-\varepsilon f_{\rm pr},\qquad n\cdot h=0.
\]
Taking the scalar product with \(n\) and using \(n\cdot\partial_\tau h=-(\partial_\tau n)\cdot h\) gives
\begin{equation}
 \pi=-N^{-1}\partial_\eta^{-1}
 \left(\frac{n\cdot Kh-(\partial_\tau n)\cdot h+\varepsilon n\cdot f_{\rm pr}}{|n|^2}\right).
 \label{pulse:eq:3.36}
\end{equation}
Let \(a=\psi h\) and \(\varpi=\psi\pi\), with \(\psi\) the clock cutoff. Brackets below denote evaluation at the carrier angle \(\eta=N\phi\). Define
\begin{equation}
\begin{aligned}
 A&=-\frac{U}{2^kN}\left[\frac{n\times\partial_\eta^{-1}a}{|n|^2}\right]_{\eta=N\phi},
 & w&=\nabla\times A=U[a+b]_{\eta=N\phi},\\
 b&=-N^{-1}\mathsf C\left(\frac{n\times\partial_\eta^{-1}a}{|n|^2}\right),
 & p&=U^2[\varpi]_{\eta=N\phi}.
\end{aligned}\label{pulse:eq:3.37}
\end{equation}
Here \(\mathsf C=2^{-k}\nabla_{\rm coef}\times\) differentiates the cylindrical frame and the auxiliary variables while holding \(\eta\) fixed. The carrier derivative of \(A\) gives \(Ua\), since \(-n\times(n\times a)/|n|^2=a\).

Write \(\widehat r=2^kr\), \(\widehat u_B=u_B/U\), \(\mathsf D_r=2^{-k}D_r\) and \(c=a+b\). The transport defect of the phase is \(\mathfrak d=2^{-(2+\delta)k}(D_t+u_B\cdot\nabla)\phi\), where \(D_t\) includes all physical time derivatives before evaluation. To expand the linearized force, set
\[
\begin{aligned}
 \mathsf T&=-\varepsilon\partial_{\widehat T}+\widehat u_B^r\mathsf D_r+
                 \varepsilon\widehat u_B^z\partial_{\widehat Z},\\
 \mathsf R c&=\varepsilon c^z\partial_{\widehat Z}\widehat u_B+
       (c^r\mathsf D_r\widehat u_B^r,\ \widehat r^{-1}\widehat u_B^r c^\theta,\ 0),\\
 \mathsf G\varpi&=(\mathsf D_r\varpi,\ 0,\ \varepsilon\partial_{\widehat Z}\varpi),\\
 \mathsf D_2c&=(\mathsf D_r^2+\widehat r^{-1}\mathsf D_r+
       \varepsilon^2\partial_{\widehat Z}^2)c-\widehat r^{-2}(c^r,c^\theta,0),\\
 \mathsf V_nc&=2(n^r\mathsf D_r+\varepsilon n^z\partial_{\widehat Z})c
       +(\mathsf D_r n^r+\widehat r^{-1}n^r+\varepsilon\partial_{\widehat Z}n^z)c\\
       &\hspace{1.5em}+2\widehat r^{-1}n^\theta(-c^\theta,c^r,0),\\
 \mathsf P&=\partial_\tau+K-d\partial_\eta^2.
\end{aligned}
\]
For \(L_B(w,p)=\partial_tw+(u_B\cdot\nabla)w+(w\cdot\nabla)u_B-\Delta w+\nabla p\), the exact identity is
\begin{equation}
\begin{aligned}
 F_k^{-1}\big[L_B(w,p)+F_k f_{\rm pr}-E_t\big]
  &=\varepsilon^{-1}\big[\mathsf P b+\mathsf T c+N\mathfrak d\,\partial_\eta c
                +\mathsf Rc+\mathsf G\varpi\\[-2pt]
  &\hspace{5em}{}-\varepsilon\mathsf D_2c-\varepsilon N\mathsf V_n\partial_\eta c
       \big]_{\eta=N\phi},\\
 E_t&=F_k\big[(1-\psi)f_{\rm pr}+\varepsilon^{-1}\psi'h\big]_{\eta=N\phi}.
\end{aligned}\label{pulse:eq:3.38}
\end{equation}
Here \(f_{\rm pr}\) is evaluated at \(\eta=N\phi\). Split each derivative into
its coefficient and carrier parts, using
\(\partial_t=2^{(2+\delta)k}\partial_\tau-2^{2k}\partial_{\widehat T}\)
and \(2^{-k}\partial_z=\varepsilon\partial_{\widehat Z}\).
Inserting the waveform equation for \(a=\psi h\) leaves \(\mathsf P b\) and the clock error \(E_t\).

\subsubsection{Primary and signed estimates}\label{pulse:sec:3.4.2}
\pulsebridge{3.4.2}

Use the coefficient bounds \eqref{pulse:eq:2.77} and the leading
covariance identity \eqref{pulse:eq:2.74}.

\noindent\textbf{Lemma \NumberedPulseResult{3.9} (Primary force).}
\textit{Write \(\mathrm{NS}(u_B+w_0,p_B+p_0)=E_0+f_1\). One can choose \(E_0-E_B\) physically flat with zero angular mean so that}
\begin{equation}
\begin{aligned}
 \frac{(f_1)_o}{F_k\sqrt\omega G}&=O_m(\varepsilon^{-\delta_0}),&
 \frac{(\langle f_1\rangle_\theta)_a^{\tan}}{F_k\omega}&=O_m(\varepsilon^{-\delta_0}),\\
 \frac{\langle f_1^{\tan}\rangle_{\theta,Y}}{F_k\omega}&=O_m(\varepsilon^{1-2\delta_0}),&
 \frac{\langle f_1^r\rangle_\theta}{F_k\omega}&=O_m(\varepsilon^{-\delta_0}).
\end{aligned}\label{pulse:eq:3.40}
\end{equation}
\textit{The subscripts \(o\) and \(a\) remove the angular and auxiliary means, respectively. The corrections lie in the stress annulus and preserve both moments in \eqref{pulse:eq:3.1}.}

\textit{Proof.} Apply \eqref{pulse:eq:3.38} to each homogeneous fundamental wave. The phase bounds give \(\mathfrak d,\widehat u_B^r=O_m(\varepsilon)\), with polynomial derivative bounds for \(n,K,d\). A radial coefficient derivative costs one source derivative and a factor \(P_m(k)\varepsilon^{-\delta_0}\). In the units \(F_k\sqrt\omega G\), the seven terms on the right of \eqref{pulse:eq:3.38}, in order, have sizes
\[
 \varepsilon^{-\delta_0},\quad \varepsilon^{1/2-\delta_0},\quad 1,\quad
 \varepsilon^{1/2},\quad\varepsilon^{-\delta_0},\quad
 \varepsilon^{1/2-2\delta_0},\quad \varepsilon^{-\delta_0}.
\]
On the support of the differentiated clock, \(G\le e^{-ck^2}\). The clock error has zero angular mean and is physically flat, since fixed physical derivatives cost at most \(e^{C_m k}\).

For the quadratic term use incompressibility in the form
\[
 \operatorname{div}_{\rm coef}c+Nn\cdot\partial_\eta c=0,\qquad
 Nn\cdot c=-\operatorname{div}_{\rm coef}\partial_\eta^{-1}c .
\]
The zero carrier mean gives the second identity. Thus \(Nn\cdot c\) costs one coefficient derivative and the quadratic force is \(F_k O_m(\varepsilon^{-\delta_0}\omega G^2)\).

Put \(Q^{ij}=\langle w_0^i w_0^j\rangle_\theta\). The linear oscillatory terms and the pulse pressure have zero angular mean. The remaining conservative cylindrical equations give
\begin{equation}
\begin{aligned}
 \langle f_1^r\rangle_\theta
    &=(D_r+1/r)Q^{rr}+\partial_z Q^{zr}-Q^{\theta\theta}/r,\\
 \langle f_1^j\rangle_\theta
    &=-(\partial_r+q_j/r)T_B^j+(D_r+q_j/r)Q^{rj}
                                  +\partial_z Q^{zj},\qquad j=\theta,z .
\end{aligned}\label{pulse:eq:3.41}
\end{equation}
Here \(q_\theta=2\) and \(q_z=1\). From \(Q/U^2=O_m(\varepsilon\omega)\), radial differentiation yields the auxiliary and radial estimates. Upon taking the auxiliary mean, the fast auxiliary derivative vanishes. By Corollary~\ref{pulse:thm:2.16},
\(\langle Q^{rj}\rangle_Y=T_B^j+U^2 O_m(\varepsilon^{2-2\delta_0}\omega)\).
Substituting in \eqref{pulse:eq:3.41} cancels the background stress divergence. The remaining radial and axial flux divergences have sizes \(F_k O_m(\varepsilon^{1-2\delta_0}\omega)\) and \(F_k O_m(\varepsilon\omega)\), respectively. The velocity is a compact curl with zero angular mean, which proves the support and moment assertions. \hfill\(\square\)

On the open stress annulus, let \(T\) be a real pair independent of \(\theta,Y\). Use the covariance matrix \(H_{k,\alpha}\) and its positive amplitudes \(g_\gamma\) from \eqref{pulse:eq:2.69} and \eqref{pulse:eq:2.72}. For \(\gamma=(k,\alpha,\pm)\), set
\begin{equation}
 \Sigma_k=2^{-(2+\delta)k}T,\qquad d_{k,\alpha}=H_{k,\alpha}^{-1}\Sigma_k,\qquad
 \delta g_{(k,\alpha,\pm)}=\frac{(d_{k,\alpha})_\pm}{2g_{(k,\alpha,\pm)}} .
 \label{pulse:eq:3.42}
\end{equation}
Insert these amplitudes into the homogeneous waves of \eqref{pulse:eq:3.37} and denote the resulting velocity and pressure by \(w_S[T],p_S[T]\).

\noindent\textbf{Proposition \NumberedPulseResult{3.10} (Signed stress inverse).}
\textit{For each label \(\gamma=(k,\alpha,\pm)\), on the intersection of its outer slow box with the open stress annulus, \eqref{pulse:eq:3.42} is real and linear in \(T\) for either sign of each component. With \(\mathcal B(U_1,V_1)=\langle U_1^r V_1^{\tan}+V_1^r U_1^{\tan}\rangle_{\theta,Y}\), we have}
\begin{equation}
 \mathcal B(W^{\rm pre},w_S^{\rm pre}[T])=T.\label{pulse:eq:3.43}
\end{equation}
\textit{Assume \(T/(U^2\omega)=O_m(\varepsilon^\beta)\), allowing algebraic edge factors as in Lemma~\ref{pulse:thm:3.3}. On trimmed radial supports these factors are polynomially bounded in \(k\). Then}
\begin{equation}
\begin{aligned}
 w_S^{\rm pre}/U&=O_m(\varepsilon^{\beta-1/2}\sqrt\omega G),&
 (w_S-w_S^{\rm pre})/U&=O_m(\varepsilon^{\beta-\delta_0}\sqrt\omega G),\\
 p_S/U^2&=O_m(\varepsilon^\beta\sqrt\omega G).
\end{aligned}\label{pulse:eq:3.44}
\end{equation}
\textit{Under the same hypothesis, the covariance of the full curls differs from \(T\) by \(U^2 O_m(\varepsilon^{\beta+1-2\delta_0}\omega)\).}

\textit{Proof.} By \eqref{pulse:eq:3.42},
\(H_{k,\alpha}(2g_+\delta g_+,2g_-\delta g_-)^{\mathsf T}=\Sigma_k\).
Polarizing \eqref{pulse:eq:2.74} gives \eqref{pulse:eq:3.43}, since distinct labels
have disjoint supports, the squared partitions sum to one, and \(U^2\varepsilon=2^{(2+\delta)k}\).
The lower bound \(g_\gamma^2\ge ck\omega\) gives
\(\delta g_\gamma=O_m(\varepsilon^{\beta-1}\sqrt\omega)\), with
derivatives obtained from the edge bounds for
\(\widetilde g_\gamma=g_\gamma/\sqrt\omega\) and its reciprocal.
The waveform and pressure bounds give the first and third estimates
in \eqref{pulse:eq:3.44}. One coefficient curl costs
\(N^{-1}\varepsilon^{-\delta_0}\) up to polynomial band factors. The mixed
cosine--sine terms average to zero. The two sine factors have exponents
\(1-\delta_0\) and \(\beta-\delta_0\), giving the stated full-curl error. \hfill\(\square\)

\subsection{Flat edge errors}\label{pulse:sec:3.5}
\pulsebridge{3.5}

For the pulses we use the same cutoff on overlapping bands by taking \(\kappa=1+\log(1/\lambda)\) as in \eqref{pulse:eq:2.29}. Choose \(d(R)>0\) smooth on \(R_{\rm in}<R<R_e\), equal to \(\log(R/R_{\rm in})\) near \(R_{\rm in}\) and to \(\log(R_e/R)\) near \(R_e\).

\noindent\textbf{Lemma \NumberedPulseResult{3.11}\label{pulse:thm:3.12} (Gaussian edge trimming).}
\textit{Fix \(0<a<b\) and a smooth cutoff \(\chi\) equal to zero on \(s\le a\)
and one on \(s\ge b\). If each fixed physical derivative of \(G\) is bounded on
\(d\le b/\kappa\) by \(C_m\lambda^{-C_m}\kappa^{B_m}d^{-q_m}e^{-c/d^2}\)
for some \(c>0\), then every fixed physical derivative of
\((1-\chi(\kappa d))G\) is \(O(\lambda^A)\) for every \(A>0\).
In particular, suppose the stress \(T_B=\lambda^{-2-\delta}\mathcal T(\kappa,R,Z)\)
extends smoothly by zero at both edges and satisfies}
\begin{equation}
 |D^m_{\kappa,R,Z}\mathcal T|\le C_m\kappa^{b_m} d^{-q_m} e^{-c/d^2},\qquad c>0.
 \label{pulse:eq:3.45}
\end{equation}
\textit{Then \(T_{\rm edge}=(1-\chi(\kappa d))T_B\) satisfies}
\begin{equation}
 |\partial_{\rm phys}^m T_{\rm edge}|+
 |\partial_{\rm phys}^m\mathcal D T_{\rm edge}|
 \le C_m\lambda^{-C_m}\kappa^{B_m}e^{-c_*\kappa^2}=O(\lambda^A).
 \label{pulse:eq:3.46}
\end{equation}

\textit{Proof.} Normalized derivatives of the cutoff cost powers of \(\kappa\).
In \eqref{pr:derivatives}, \(\lambda\partial_\lambda=-\partial_\kappa\),
\(1-\delta Z^2\ge1-\delta\), and the frame derivative costs \(1/r=O(\lambda^{-1})\).
Thus physical differentiation introduces only fixed scale powers, also for \(\mathcal D T_{\rm edge}\).
Half the Gaussian absorbs \(d^{-q_m}\). The other half is bounded by
\(e^{-c\kappa^2/(2b^2)}\), which absorbs all scale powers since \(\lambda=e^{1-\kappa}\).
\hfill\(\square\)

Choose \(0\le\chi_0\le1\) equal to zero on \(s\le1\) and one on \(s\ge2\). The cutoffs for the source and its reconstruction are
\begin{equation}
\begin{aligned}
 \chi_{\rm in}&=\chi_0\bigl((\kappa/a)\log(R/R_{\rm in})\bigr)
                 \chi_0\bigl((\kappa/a)\log(R_e/R)\bigr),\\
 \chi_r&=\chi_0\bigl((3\kappa/a)\log(R/R_{\rm in})\bigr)
          \chi_0\bigl((3\kappa/a)\log(R_e/R)\bigr).
\end{aligned}\label{pulse:eq:3.47}
\end{equation}
For large \(k\), \(\chi_r=1\) near both \(\operatorname{supp}\chi_{\rm in}\) and the compensator window \(K\). Fixed normalized derivatives cost powers of \(k\). These cutoffs act on the radial free coordinates in \eqref{pulse:eq:3.3}. The defects \(P,J^{\theta z},J^{z\theta}\) depend only on \(t,z\).

For \(G=\omega^p\widehat G\), \(p>0\), assume each fixed derivative of
\(\widehat G\) has an exponential scale bound and an algebraic edge bound.
On either discarded collar \(d\le2a/\kappa\). The factor \(\omega^p\) supplies
the Gaussian, and the auxiliary and carrier chain rules cost at most \(e^{C_m k}\).
Thus Lemma~\ref{pulse:thm:3.11} gives \(O(\lambda^A)\) for each fixed
physical derivative of \((1-\chi_{\rm in})G\), using finitely many input bounds.
Trimming preserves the background moments. Proposition~\ref{pulse:thm:3.8}
cancels the three defects while preserving both linear velocity moments.

\noindent\textbf{Proposition \NumberedPulseResult{3.13} (Cutoff through a pulse potential).}
\textit{Let \(A,q\) be defined on the retained region \(d\ge a/\kappa\), with \(v=\nabla\times A\). On the transition strip \(a/\kappa\le d\le b/\kappa\), assume scale-power bounds for the fixed derivatives of the background and}
\begin{equation}
 |\partial_{\rm phys}^m(A,v,q)|\le C_m\lambda^{-C_m}\kappa^{B_m}e^{-c_1/d^2},
 \qquad c_1>0.
 \label{pulse:eq:3.49}
\end{equation}
\textit{Choose a smooth cutoff \(\eta(\kappa d)\) that is zero for \(d\le a/\kappa\) and one for \(d\ge b/\kappa\). Then \(v_{\rm cut}=\nabla\times(\eta A)\) is divergence-free. Together with \(q_{\rm cut}=\eta q\), it agrees with the original pulse on \(d\ge b/\kappa\). All new equation errors on the transition strip are physically flat.}

\textit{Proof.} Since \(v_{\rm cut}=\eta v+\nabla\eta\times A\), each new
term contains a factor satisfying \eqref{pulse:eq:3.49}. All other factors
cost at most scale powers, so the Gaussian estimate in Lemma~\ref{pulse:thm:3.11} applies.
\hfill\(\square\)

Applying \(\chi_{\rm in}\) to the compatible coordinates preserves the polynomial inverse bounds at zeros of the cutoff.

\subsection{Initialization}\label{pulse:sec:3.6}

\noindent\textbf{Proposition \NumberedPulseResult{3.14} (Initialized state).}
\textit{On all sufficiently high bands there are real fields}
\[
 u_*=u_B+w_0+v_A+b_Me_\theta,\qquad
 p_*=p_B+p_0+q_*,\qquad \mathrm{NS}(u_*,p_*)=E_*+f_* .
\]
\textit{The velocity is divergence-free. Its corrections lie in the stress annulus, preserve the two moments \eqref{pulse:eq:3.1}, and set \(P,J^{\theta z},J^{z\theta}\) equal to zero. The force satisfies}
\begin{equation}
 \langle f_*^r\rangle_\theta=0,\qquad
 \frac{(f_*)_o}{F_k\sqrt\omega G}=O_m(\varepsilon^{-\delta_0}),\qquad
 \frac{\langle f_*\rangle_\theta}{F_k\omega}=O_m(\varepsilon^{1-3\delta_0}).
 \label{pulse:eq:3.50}
\end{equation}
Here \(O_m\) allows the algebraic factors in \(d_e^{-1}\) specified in Section~\ref{pulse:sec:3.2}. The ordinary joint Sobolev bounds for the normalized source follow after the free-coordinate cutoff in Section~\ref{pulse:sec:4.4.1}.
\textit{The correction fields satisfy}
\begin{equation}
\begin{aligned}
 v_A^{\theta,z}/U&=O_m(\varepsilon^{1-\delta_0}\omega),&
 v_A^r/U&=O_m(\varepsilon^{2-\delta_0}\omega),\\
 b_M/U&=O_m(\varepsilon^{1-\delta_0}),&
 q_*/U^2&=O_m(\varepsilon^{1-\delta_0}\omega).
\end{aligned}\label{pulse:eq:3.51}
\end{equation}
\textit{The swirl \(b_M\) has the same weighted bound because it is supported in \(I_3\), where \(\omega\) is bounded below. The error \(E_*-E_B\) is physically flat and has zero double mean.}

\textit{Proof. Auxiliary correction.} The zero auxiliary mean part
\(h=(\langle f_1\rangle_\theta)_a^{\tan}\) satisfies
\(h/F_k=O_m(\varepsilon^{-\delta_0}\omega)\) by \eqref{pulse:eq:3.40}.
On the common level-\(c\) cover set
\begin{equation}
\begin{aligned}
 A&=-\Lambda_t^{-c}(v_t\cdot\nabla_{Y_c})^{-1}h^\theta,&
 B&=-\Lambda_t^{-c}(v_t\cdot\nabla_{Y_c})^{-1}h^z,\\
 \Psi&=r^{-1}I_0(rB),&
 v_A&=Ae_\theta+\nabla\times(\Psi e_\theta).
\end{aligned}\label{pulse:eq:3.52}
\end{equation}
The directional inverse loses one derivative and contributes \(F_k/(U\Lambda_t^c)\le P(k)\varepsilon\) in velocity units. By Proposition~\ref{pulse:thm:3.4},
\[
\begin{aligned}
 \Psi&=O_m(2^{-k}U\varepsilon^{1-\delta_0}\omega),&
 v_A^r&=-\partial_z\Psi,\qquad v_A^\theta=A,\\
 v_A^z&=B+R_A,& R_A&=-\chi'_m\mathcal J(rB)/r .
\end{aligned}
\]
This gives \eqref{pulse:eq:3.51} for \(v_A\), a compact divergence-free swirl
and curl with zero auxiliary mean. Hence it preserves \eqref{pulse:eq:3.1}.
Its fast derivative cancels \(h\) up to
\begin{equation}
 E_A=\Lambda_t^c(v_t\cdot\nabla_{Y_c}R_A)e_z
     =\frac{\chi'_m}{r}\mathcal J(rh^z)e_z .\label{pulse:eq:3.53}
\end{equation}
This error has zero auxiliary mean and is flat by \eqref{pulse:eq:3.17}. Keep \(R_A\) in the velocity. Write \(u_2=u_B+w_0+v_A\), \(p_2=p_B+p_0\), \(E_2=E_0+E_A\), and \(f_2=\mathrm{NS}(u_2,p_2)-E_2\).

To estimate \(f_2\), put \(C^{ij}=u_B^i v_A^j+v_A^i u_B^j+v_A^i v_A^j\),
\(\Delta_0=D_r^2+r^{-1}D_r+\partial_z^2\), and
\(\partial_t^{\rm slow}=-2^{2k}\partial_{\widehat T}\).
Angular averaging gives the exact identities
\[
\begin{aligned}
 \langle f_2^\theta-f_1^\theta\rangle_\theta
 &= -h^\theta+\partial_t^{\rm slow}v_A^\theta-(\Delta_0-r^{-2})v_A^\theta
                       +(D_r+2/r)C^{r\theta}+\partial_z C^{z\theta},\\
 \langle f_2^z-f_1^z\rangle_\theta
 &= -h^z+\partial_t^{\rm slow}v_A^z-\Delta_0v_A^z
                       +(D_r+1/r)C^{rz}+\partial_z C^{zz},\\
 \langle f_2^r-f_1^r\rangle_\theta
 &= D_t v_A^r-(\Delta_0-r^{-2})v_A^r+(D_r+1/r)C^{rr}
                              +\partial_z C^{zr}-C^{\theta\theta}/r .
\end{aligned}
\]
In units of \(F_k\omega\), \eqref{pulse:eq:4.20} gives the tangential
slow-time, viscous, radial-transport and axial-transport exponents
\(1-\delta_0,1-3\delta_0,1-2\delta_0,1-\delta_0\), respectively.
The quadratic axisymmetric terms are smaller, but the radial centrifugal term
\(-2u_B^\theta v_A^\theta/r\) has exponent \(-\delta_0\).
The nonzero angular increment is \((v_A\cdot\nabla)w_0+(w_0\cdot\nabla)v_A\),
with normalized carrier derivative \(O_m(\sqrt\omega G)\) on \(w_0\).
Thus \eqref{pulse:eq:3.40} gives
\begin{equation}
\begin{aligned}
 \frac{(f_2)_o}{F_k\sqrt\omega G}&=O_m(\varepsilon^{-\delta_0}),&
 \frac{\langle f_2^{\tan}\rangle_\theta}{F_k\omega}&=O_m(\varepsilon^{1-3\delta_0}),\\
 \frac{\langle f_2^r\rangle_\theta}{F_k\omega}&=O_m(\varepsilon^{-\delta_0}).
\end{aligned}\label{pulse:eq:3.54}
\end{equation}

\textit{Moment repair.} We correct the three integrated defects at fixed pressure. The double mean of \(u_2\) is \(u_B\). Its defects satisfy
\begin{equation}
 \left(U^{-2}P_2,\ U^{-2}2^{3k}J_2^{\theta z},\
 U^{-2}2^{2k}J_2^{z\theta}\right)=O_m(\varepsilon^{1-\delta_0}).
 \label{pulse:eq:3.55}
\end{equation}
For \(P_2\), radial integration gains \(2^{-k}\), with \(2^{-k}F_k/U^2=\varepsilon\).
In the fluxes, the terms linear against \(u_B\) and the primary pressure
have zero double mean. The primary and auxiliary covariances are
\(O_m(U^2\varepsilon)\) and \(O_m(U^2\varepsilon^{2-2\delta_0})\).
The measures \(r^2\,dr\) and \(r\,dr\) give \(O(2^{-3k})\) and \(O(2^{-2k})\),
and the radial-force term in \(J_2^{z\theta}\) costs
\(F_k2^{-3k}\varepsilon^{-\delta_0}=U^2 2^{-2k}\varepsilon^{1-\delta_0}\).

With matrix \(G_3\) and data \eqref{pulse:eq:3.55}, Proposition~\ref{pulse:thm:3.8} gives
\(b_M/U=O_m(\varepsilon^{1-\delta_0})\) supported in \(I_3\),
with all three defects zero at fixed pressure. The radial angular-mean
increment is \(-[2(u_B^\theta+v_A^\theta)b_M+b_M^2]/r\).
Using \(\inf_{I_3}\omega>0\) for this term and
\eqref{pulse:eq:4.26}--\eqref{pulse:eq:4.27} for the other components
keeps the bounds \eqref{pulse:eq:3.54}.

\textit{Pressure projection.} Let \(H\) be the radial angular mean after moment repair.
Since \(P=0\), \(\int\langle H\rangle_Y\,dr=0\). Set
\(q_*=I_0(-H)\) and include \(R_H e_r\) in \(E_*\).
Proposition~\ref{pulse:thm:3.5} gives the pressure bound in
\eqref{pulse:eq:3.51}, sets \(\langle f_*^r\rangle_\theta=0\), and
adds \(F_k O_m(\varepsilon^{1-\delta_0}\omega)\) to the axial mean.
By \eqref{pulse:eq:3.22} the projection preserves all three cancellations. Each term of \(E_*-E_B\) is flat with zero double mean. The initialized tangential
double means are \((u_B^\theta+b_M,u_B^z)\). Their four-row matrix
differs from \(G_3\) by \(O_m(\varepsilon^{1-\delta_0})\) up to polynomial
band factors and has a uniform inverse on sufficiently high bands.
\hfill\(\square\)
\section{The correction operator}\label{pulse:sec:4}
\pulsebridge{4}

\subsection{Nonzero carrier modes}\label{pulse:sec:4.1}

Freeze \(u_*=u_B+v_*+w_0\) from Proposition~\ref{pulse:thm:3.14}. Here \(v_*\) is axisymmetric and \(w_0\) contains the two primary pulse families. A trial source can have any carrier harmonic. The forced inverse uses the additional viscous damping of higher harmonics. Its cross term with \(w_0\) supplies a source for the mean correction.

On a native label \(\gamma=(k,\alpha,\pm)\) we write \(U,\varepsilon,N\) for \(U_k,\varepsilon_k,N_k\). The force unit is \(F_k\).
The independent carrier angle \(\eta\) is evaluated at \(N\phi\) in physical space.

Every label retains its native scale and angular quantization. Overlapping scales differ by at most one. For the correction estimates in this section, \(\|\cdot\|_m\) is the joint \(H^{10+m}\) coefficient and carrier norm before physical pullback. The weighted estimates use the retained charts and their fixed extensions.

\textbf{The forced inverse.} Expand the source as \(f=\sum_{m\ne0} f_m e^{im\eta}\) in the independent \(2\pi\)-periodic carrier angle \(\eta\). At each native label and each \(m\ne0\) solve
\begin{equation}
\partial_\tau h_m+K_Bh_m+m^2 d\,h_m+imNn\pi_m=-f_m,\qquad n\cdot h_m=0,\qquad h_m(0)=0.
\label{pulse:eq:4.1}
\end{equation}
Differentiate \(n\cdot h_m=0\) in the fast time. The identity \(n\cdot\partial_\tau h_m=-(\partial_\tau n)\cdot h_m\) determines the pressure and gives the projected equation
\begin{equation}
\begin{aligned}
\partial_\tau h_m&=(A-m^2 dI)h_m-\Pi_n f_m,\qquad \Pi_n=I-\frac{nn^{\mathsf T}}{|n|^2},\\
\pi_m&=\frac{i}{mN}\,\frac{n\cdot K_Bh_m-(\partial_\tau n)\cdot h_m+n\cdot f_m}{|n|^2}.
\end{aligned}
\label{pulse:eq:4.2}
\end{equation}
The physical source is \((F_k/\varepsilon)f(x,l_r,\tau,N\phi)\). Summing \eqref{pulse:eq:4.1} over its carrier harmonics gives
\begin{equation}
\begin{aligned}
\partial_\tau h+K_Bh-d\,\partial_\eta^2 h+Nn\partial_\eta\pi&=-f,\qquad n\cdot h=0,\\
\pi&=-\frac1N\,\partial_\eta^{-1}\!\left(\frac{n\cdot K_Bh-(\partial_\tau n)\cdot h+n\cdot f}{|n|^2}\right).
\end{aligned}
\label{pulse:eq:4.3}
\end{equation}
Here \(\partial_\eta^{-1}\) is the zero-mean periodic primitive with multiplier \(1/(im)\).

\noindent\textbf{Proposition \PulseResult{4.1}} (All-harmonic forced inverse)\textbf{.} \textit{Suppose} \(f/G\) \textit{has all initial fast-time jets zero. Use joint Sobolev norms in} \((x, l_{r}, \tau, \eta )\) \textit{on a buffered native coefficient chart. The pair} \((l_r,\tau)\) \textit{parametrizes the auxiliary coordinates there. For every fixed integer} \(s\ge 0\)\textit{,}
\begin{equation}
\begin{aligned}
\left\|\frac{h}{G}\right\|_{H^s}+\left\|\frac{\partial_\tau h}{G}\right\|_{H^s}+\left\|\frac{\partial_\eta^2 h}{G}\right\|_{H^s}
 &\le P_s(k)\left\|\frac{f}{G}\right\|_{H^s},\\
\left\|\frac{\pi}{G}\right\|_{H^s}+\left\|\frac{\partial_\eta\pi}{G}\right\|_{H^s}
 &\le N^{-1} P_s(k)\left\|\frac{f}{G}\right\|_{H^s}.
\end{aligned}
\label{pulse:eq:4.4}
\end{equation}
\textit{The constants are independent of harmonic range. The estimates also hold with} \(G\) \textit{replaced by} \(\sqrt{\omega}\,G\) \textit{on retained strips. The inverse preserves slow and transverse support and commutes with lower-cover pullback retaining the native clock and initial section.}

\textit{Proof.} For either family, use the orthonormal frame \(Q\) for \(n^\perp\) in \eqref{pulse:eq:2.58}. Set \(z_m=Q^{\mathsf T}h_m/G\), so that \(Qz_m=h_m/G\). In this frame
\[
 \partial_\tau z_m=(B_\pm+E_\pm-a_{\rm ref}I-m^2dI)z_m-Q^{\mathsf T}f_m/G .
\]
Write \(\widetilde{\mathcal U}_m\) for its propagator. In the notation of \eqref{pulse:eq:2.63},
\(\widetilde{\mathcal U}_1(\tau,s)=G(s)Q(\tau)^{\mathsf T}\mathcal U_1(\tau,s)Q(s)/G(\tau)\).
The damping is scalar, so
\begin{equation}
\widetilde{\mathcal U}_m(\tau,s)=\exp\!\left[-(m^2-1)\int_s^\tau d(v)\,dv\right]\widetilde{\mathcal U}_1(\tau,s).
\label{pulse:eq:4.5}
\end{equation}
By \eqref{pulse:eq:2.63}, \(\|\widetilde{\mathcal U}_m(\tau,s)\|\le C e^{-(m^2-1)d_*(\tau-s)}\).
For \(|m|\ge2\) its time kernel has mass at most \(C/((m^2-1)d_*)\). For \(|m|=1\) the finite clock length gives \(CT\). Young's inequality yields
\[
 (1+m^2)\|z_m\|_{L^2_{\tau,x,l_r}}\le P(k)\|f_m/G\|_{L^2_{\tau,x,l_r}} .
\]
For a derivative \(D^\alpha\) in \((\tau,x,l_r)\), the equation for \(D^\alpha z_m\) has the same propagator and right-hand side
\[
 -D^\alpha(Q^{\mathsf T}f_m/G)+\sum_{0<\beta\le\alpha}\binom{\alpha}{\beta}
    (D^\beta(B_\pm+E_\pm-a_{\rm ref}I)-m^2D^\beta d\,I)D^{\alpha-\beta}z_m .
\]
The initial jets vanish by induction in the equation. Each coefficient derivative costs a polynomial in \(k\), and the preceding \(m^2\)-bound absorbs the differentiated damping at lower order. Thus
\[
 (1+m^2)\|D^\alpha z_m\|_{L^2}\le P_{|\alpha|}(k)
             \sum_{\beta\le\alpha}\|D^\beta(f_m/G)\|_{L^2}.
\]
Multiply by \(|m|^j\) for \(|\alpha|+j\le s\) and sum squares. The equation also bounds \(\partial_\tau z_m\). Reconstruction \(h_m=GQz_m\) uses the polynomial bounds on derivatives of \(Q\) and of \(\partial_\tau\log G\). The pressure follows from \eqref{pulse:eq:4.2}, since \(\partial_\eta^{-1}\) is bounded on zero-mean Sobolev spaces. Division by \(\sqrt\omega\) commutes with the clock equation. The solution is pointwise in the slow and transverse variables. Uniqueness gives conjugation symmetry for real data and commutation with any pullback preserving the native clock and its initial section. \hfill\(\square\)

A physical source \(F_k\varepsilon^{-3/10}\sqrt{\omega}\,Gg\) has fast-time coefficient \(f=\varepsilon^{7/10}\sqrt{\omega}\,Gg\). Thus \eqref{pulse:eq:4.4} gives weighted amplitude and pressure sizes \(P_s(k)\varepsilon^{7/10}\) and \(P_s(k)\varepsilon^{6/5}\), respectively, for a unit input.

\textbf{Time cutoffs and curl completion.} With the output cutoff \(\psi\) from \eqref{pulse:eq:2.67}, the forced equation gives
\begin{equation}
\partial_\tau(\psi h)+K_B\psi h-d\partial_\eta^2(\psi h)+Nn\partial_\eta(\psi\pi)+f=(\partial_\tau\psi)h+(1-\psi)f.\label{pulse:eq:4.6}
\end{equation}

After spectral smoothing use a fixed \(\chi_{\rm src}\in C_c^\infty(0,T)\) equal to
one near \([T/6,5T/6]\), so the source has zero initial jets. The cutoff \(\psi\)
extends the output smoothly by zero at both ends. The discarded source and both
errors in \eqref{pulse:eq:4.6} lie where \(|\tau-T/2|\ge T/5\). There the Gaussian
before division by \(G\) bounds every fixed physical derivative by
\(P_m(k)e^{-ck^2+C_mk}\) for inputs of polynomial weighted norms.

\textbf{Linear reconstruction.} Let \(g_{\rm raw}\) be real and smooth, with zero carrier mean and the assigned slow and transverse supports. Following spectral smoothing, put \(g=\chi_{\rm src}g_{\rm raw}\) and \(f=\varepsilon^{7/10}\sqrt{\omega}\,Gg\). The physical source is \(f_o[g]=(F_k/\varepsilon)[f]_{\eta=N\phi}\), with coefficient norm \(\|g\|_m\). With \(n=2^{-k}\nabla\phi\) and \(d=\varepsilon N^2|n|^2\), the velocity equation is
\begin{equation}
\partial_\tau h+K_Bh-d\,\partial_\eta^2 h+Nn\,\partial_\eta\pi=-f,\qquad n\cdot h=0,\quad h|_{\tau=0}=0.\label{pulse:eq:4.8}
\end{equation}

Use \eqref{pulse:eq:3.37} with \(a=\psi h\) for \(A[g],b,w[g]\).
The physical pressure is \(p[g]=U^2[\psi\pi]_{\eta=N\phi}\), with \(\pi\) from \eqref{pulse:eq:4.3}. The coefficient curl holds
\(\eta\) fixed and includes the auxiliary chain rule and cylindrical frame.
The pair is real and linear, has the assigned supports and has zero angular mean.
Its velocity is divergence-free and preserves both linear moments.

The cutoff force in \eqref{pulse:eq:4.6} is
\[
E_t[g]=\frac{F_k}{\varepsilon}\big[(1-\psi)f+(\partial_\tau\psi)h\big]_{\eta=N\phi}.
\]

\noindent\textbf{Proposition \PulseResult{4.2}} (The full nonzero-mode correction map)\textbf{.} \textit{Let} \(L_{\ast }\) \textit{be the exact linearization at} \(u_{\ast }\)\textit{. Then}
\begin{equation}
L_*(w[g],p[g])+f_o[g]-E_t[g]=R_o[g]+H[g],\qquad \langle R_o[g]\rangle_\theta=0,\label{pulse:eq:4.14}
\end{equation}
\textit{where}
\begin{equation}
\varepsilon^{3/10}\|R_o[g]/(F_k\sqrt{\omega}\,G)\|_m\le P_m(k)\varepsilon^{1/2-3\delta_0}\|g\|_{m+3}.\label{pulse:eq:4.15}
\end{equation}
\textit{All outputs are linear in} \(g\)\textit{, with loss three independent of harmonic range. The exact mean} \(H[g]\) \textit{is given by} \eqref{pulse:eq:4.16}\textit{.}

\textit{Proof.} The chain rules \eqref{pulse:eq:2.7} and \eqref{pulse:eq:2.14} give
\begin{equation}
\|\mathsf C z\|_m+\|2^{-k}D_rz\|_m\le P_m(k)\varepsilon^{-\delta_0}\|z\|_{m+1}.\label{pulse:eq:4.10}
\end{equation}
Coefficient differentiation does not add physical radial derivatives. Thus \(\delta_0\) is fixed, and conjugation by \(\sqrt{\omega}\,G\) only changes \(P_{m}\).

Proposition~\ref{pulse:thm:4.1} and the pressure formula imply
\begin{equation}
\begin{aligned}
 &\left\|\frac{(h,\partial_\tau h,\partial_\eta^2h)}{\sqrt{\omega}\,G}\right\|_m
 \le P_m(k)\varepsilon^{7/10}\|g\|_m,\qquad
 \left\|\frac{(\psi\pi,\partial_\eta(\psi\pi))}{\sqrt{\omega}\,G}\right\|_m
 \le P_m(k)\varepsilon^{6/5}\|g\|_m,\\
 &\left\|\frac{(b,\partial_\tau b,\partial_\eta^2 b)}{\sqrt{\omega}\,G}\right\|_m
 \le P_m(k)\varepsilon^{6/5-\delta_0}\|g\|_{m+1}.
\end{aligned}\label{pulse:eq:4.11}
\end{equation}
The cutoffs preserve these bounds. Apply \eqref{pulse:eq:4.10} once for the last line. The periodic primitive leaves only one carrier derivative on \(a\). Hence \(c=a+b\) has size \(\varepsilon ^{7/10}\) with one input derivative, independently of harmonic range.

\textbf{The full linear residual.} Apply \eqref{pulse:eq:3.38} at \(u_B\)
with \(K=K_B\), \(c=a+b\), \(\varpi=\psi\pi\), \(f_{\rm pr}=f/\varepsilon\), and the phase
defect \eqref{pulse:eq:2.52}. Its left side is \(F_k^{-1}[L_B(w[g],p[g])+f_o[g]-E_t[g]]\). The cutoff force \(E_t[g]\) has zero angular mean and is supported where \(G\le e^{-ck^2}\). It is physically flat for inputs with polynomial fixed-order norms. In the order of \eqref{pulse:eq:3.38}, the seven physical terms divided by \(F_k\sqrt{\omega}\,G\) have bounds
\[
\varepsilon^{1/5-\delta_0},\quad \varepsilon^{7/10-2\delta_0},\quad
 \varepsilon^{1/5},\quad \varepsilon^{7/10},\quad
 \varepsilon^{1/5-\delta_0},\quad \varepsilon^{7/10-3\delta_0},\quad
 \varepsilon^{1/5-2\delta_0},
\]
up to \(P_m(k)\lVert g\rVert_{m+3}\). The first bound uses all three bounds on \(b\) in \eqref{pulse:eq:4.11}. The fifth uses \(\varepsilon ^{-1}\varepsilon ^{-\delta_0 }\varepsilon ^{6/5}\). The sixth differentiates \(c\) twice, after its one curl derivative. This gives the common loss three. The other bounds use \(\widehat u_B^r,\mathsf R,\mathfrak d=O_m(\varepsilon)\). At most three physical coefficient differentiations occur, independently of the Sobolev order.

\textbf{Linearization at the initialized state.}
By \eqref{pulse:eq:3.51} and \eqref{pulse:eq:2.77},
\(v_*^{\tan}/U=O_m(\varepsilon^{1-\delta_0}\omega)\),
\(v_*^r/U=O_m(\varepsilon^{2-\delta_0}\omega)\), and
\(w_0/(U\sqrt{\omega}\,G)=O_m(\varepsilon^{1/2})\).
The cross terms with \(v_*\) have zero angular mean. The largest carrier term
has size \(\varepsilon^{-1}\varepsilon^{1-\delta_0}N\varepsilon^{7/10}
=O(\varepsilon^{1/5-\delta_0})\), and the remaining terms cost at most two input derivatives.
For two zero-mean divergence-free waveforms on the same label,
\begin{equation}
(Uc\cdot\nabla)(Ud)=U^2 2^k\big[c\cdot\nabla_{\rm coef}d-(\operatorname{div}_{\rm coef}\partial_\eta^{-1}c)\partial_\eta d\big]_{\eta=N\phi}.\label{pulse:eq:4.13}
\end{equation}
Here the coefficient derivatives are normalized by \(2^{-k}\).
Lemma~\ref{pulse:thm:3.9} gives \eqref{pulse:eq:4.13} by incompressibility of
the advecting waveform. Joint Sobolev multiplication gives a loss independent
of harmonic range, and distinct labels have disjoint supports.

The primary–new-wave cross terms have force size \(\varepsilon ^{-1}\varepsilon ^{1/2}\varepsilon ^{7/10}\varepsilon ^{-\delta_0 }=\varepsilon ^{1/5-\delta_0 }\). Their nonzero angular part retains \(\sqrt{\omega}\,G\) since \(G^{2}\le G\). This gives \eqref{pulse:eq:4.15}.

\textbf{Mean and defect outputs.} Define the physical tensor
\[
Q^{ab}[g]=\langle w_0^a w[g]^b+w[g]^a w_0^b\rangle_\theta .
\]
Only opposite harmonics survive. Thus only the new modes \(m=\pm 1\) enter this linear mean. Put \(Q_d=\langle Q\rangle_Y\). Tensor divergence gives
\begin{equation}
\begin{aligned}
 H^r&=(D_r+1/r)Q^{rr}+\partial_z Q^{zr}-Q^{\theta\theta}/r,\\
 H^\theta&=(D_r+2/r)Q^{r\theta}+\partial_z Q^{z\theta},\\
 H^z&=(D_r+1/r)Q^{rz}+\partial_z Q^{zz}.
\end{aligned}\label{pulse:eq:4.16}
\end{equation}
We have
\begin{equation}
\begin{aligned}
 \|Q/(U^2\omega)\|_m&\le P_m(k)\varepsilon^{6/5}\|g\|_{m+1},\\
 \|H/(F_k\omega)\|_m&\le P_m(k)\varepsilon^{1/5-\delta_0}\|g\|_{m+2},\\
 \|\langle H\rangle_Y/(F_k\omega)\|_m&\le P_m(k)\varepsilon^{1/5}\|g\|_{m+2}.
\end{aligned}\label{pulse:eq:4.17}
\end{equation}
The auxiliary part of \(D_{r}\) integrates to zero in the last bound. The double tangential mean is
\[
H_d^a=(\partial_r+q_a/r)Q_d^{ra}+\partial_z Q_d^{za},\qquad q_\theta=2,\quad q_z=1.
\]
Its radial term has order one in the mean-source grading. The axial term gains \(\varepsilon\). Adding \(H_d\) to the current double tangential source \(f_d\) in \eqref{pulse:eq:3.28} gives the next stress target
\begin{equation}
\widetilde T^a=-I_{q_a}(f_d^a+H_d^a).\label{pulse:eq:4.18}
\end{equation}
The part generated by this oscillatory block is explicitly \(-Q_d^{ra}-\mathcal I_{q_a}\partial_z Q_d^{za}\). Indeed integration by parts gives \(\int_0^\infty r^{q_a}(\partial_r+q_a/r)Q_d^{ra}\,dr=0\). Thus its compensated primitive is exactly \(Q_d^{ra}\). Its leading size is \(U^{2}\varepsilon ^{6/5}\omega\), or \(\varepsilon ^{1/5}\omega\) in the stress unit \(U^{2}\varepsilon = 2^{(2+\delta )k}\). Retain the other compensators in the force–moment coordinates. The physical first variations are
\begin{equation}
\begin{aligned}
 \delta P&=\int_0^\infty H_d^r\,dr,\qquad
 \delta J^{\theta z}=\int_0^\infty r^2 Q_d^{\theta z}\,dr,\\
 \delta J^{z\theta}&=\int_0^\infty r Q_d^{zz}\,dr+\frac12\int_0^\infty r^2H_d^r\,dr .
\end{aligned}\label{pulse:eq:4.19}
\end{equation}
The pressure and \(E_t[g]\) have zero angular mean. Each defect is at most \(P_m(k)\varepsilon^{6/5}\lVert g\rVert_{m+2}\) in normalized units, hence order one in defect units. Since \(\langle w[g]\rangle_\theta=0\), both linear velocity moments remain unchanged. Thus the output satisfies the linearized identities of Section~\ref{pulse:sec:3.1}. \hfill\(\square\)

\subsection{Moment corrections}\label{pulse:sec:4.2}
\pulsebridge{4.2}

For this section and Section~\ref{pulse:sec:4.3}, fix the phases, coefficient charts, envelopes, weights and cutoffs. Use the joint \(H^{10+m}\) norm \(\|\cdot\|_m\) on \eqref{pulse:eq:2.30}, including the carrier angle on labelled charts. Sum in the fixed label norm, with overlap factors included in \(P_m(k)\). Moment units are those in \eqref{pulse:eq:4.44}, equivalently the nonzero entries of \eqref{pulse:eq:3.35}.

\textbf{Component formulas and force gain.} Fixed reference-band coordinates retain the full \(\lambda (t, z)\) dependence and give
\begin{equation}
\begin{gathered}
 \partial_r=2^k O(1)D_{\widehat R},\qquad
 \partial_z=2^{(1-\delta)k}D_{\widehat Z},\qquad
 \partial_t=-2^{2k}D_{\widehat T},\\
 \frac{U_k^2 2^k}{F_k}=\varepsilon^{-1},\qquad
 \frac{U_k^2 2^{(1-\delta)k}}{F_k}=\frac{U_k2^{2k}}{F_k}=1.
\end{gathered}\label{pulse:eq:4.20}
\end{equation}
Here \(O_{m}\) allows polynomial factors in \(k\) for \(m\) coefficient derivatives, leaving the carrier undifferentiated. Before auxiliary evaluation let \(u=V+a+w\), where \(V=\langle u\rangle_{\theta,Y}\), \(a=\langle u\rangle_\theta-V\) and \(\langle w\rangle_\theta=0\). Assume the initialized bounds
\begin{equation}
\begin{gathered}
 V^{\theta,z}/U_k=O_m(1),\quad V^r/U_k=O_m(\varepsilon),\quad
 a^{\theta,z}/U_k=O_m(\varepsilon^{1-\delta_0}\omega),\quad
 a^r/U_k=O_m(\varepsilon^{2-\delta_0}\omega),\\
 w_\gamma/U_k=A_\gamma\cos(N_\ell\phi_\gamma)+B_\gamma\sin(N_\ell\phi_\gamma),\\
 A_\gamma=O_m(\varepsilon^{1/2}\sqrt{\omega}G_\gamma),\qquad
 B_\gamma=O_m(\varepsilon^{1-\delta_0}\sqrt{\omega}G_\gamma),
\end{gathered}\label{pulse:eq:4.21}
\end{equation}
with \(|\ell-k|\le1\), \(N_\ell=\lceil 2^{\delta\ell/2}\rceil\) and \(\partial_\theta\phi_\gamma=p_{\theta,\gamma}=O(1)\). The labelled coefficient bounds are summable. These are \eqref{pulse:eq:2.77} and \eqref{pulse:eq:3.51}. All moment corrections lie in a fixed interior window.

\noindent\textbf{Proposition \PulseResult{4.3}} (Moment force gain)\textbf{.} \textit{Assume} \eqref{pulse:eq:4.21} \textit{and the hypotheses of Proposition~\ref{pulse:thm:3.8}.} \textit{Suppose} \(\langle f^r\rangle_{\theta,Y}=\rho_r P\)\textit{,} \(\int\rho_r\,dr=1\)\textit{, and the radial hull of the supports of} \(\rho _{r}\) \textit{and the moment profiles lies in} \(K\)\textit{, where} \(\omega \ge c_{K}>0\)\textit{. If the target in} \eqref{pulse:eq:3.35} \textit{has size} \(d=O_m(\eta\varepsilon^{6/5})\) \textit{for a scalar factor} \(0\le\eta\le1\)\textit{, then the exact swirl repair and a compact slow pressure satisfy}
\begin{equation}
b/U_k=O_m(\eta\varepsilon^{6/5}),\qquad q/U_k^2=O_m(\eta\varepsilon^{6/5}).
\label{pulse:eq:4.22}
\end{equation}
\textit{For}
\begin{equation}
\mathcal R=\NS(u+b e_\theta,p+q)-\NS(u,p)+\rho_r P e_r,
\label{pulse:eq:4.23}
\end{equation}
\textit{let} \(\mathcal R_o,\mathcal R_a,\mathcal R_d\) \textit{denote the nonzero angular part, auxiliary fluctuation of the angular mean, and double mean. Then}
\begin{equation}
\begin{aligned}
 \mathcal R_{o,\gamma}/(F_k\sqrt{\omega}G_\gamma)&=O_m(\eta\varepsilon^{1/5}),&
 \mathcal R_a/(F_k\omega)&=O_m(\eta\varepsilon^{6/5-\delta_0}),\\
 \mathcal R_d/(F_k\omega)&=O_m(\eta\varepsilon^{6/5}),& \mathcal R_d^r&=0 .
\end{aligned}
\label{pulse:eq:4.24}
\end{equation}
\textit{In the grading} \eqref{pulse:eq:4.44}\textit{, these give gains} \(\varepsilon ^{1/2}, \varepsilon ^{1-\delta_0 }, \varepsilon\)\textit{. The repaired defects vanish exactly. The force estimate costs at most two slow derivatives of} \(d\)\textit{.}

\textbf{All carrier harmonics.} For a target \textbf{D} \(= (P, J^{\theta z}, J^{z\theta })\) of normalized size \(\varepsilon^{6/5}\mathbf d\), use the lift
\begin{equation}
\begin{aligned}
 \mathcal S(\mathbf D)&=\rho_r P e_r+b^\theta C^\theta e_\theta+b^z C^z e_z,\\
 C^\theta&=\partial_z J^{\theta z},\qquad C^z=\partial_z(J^{z\theta}-c_\rho\lambda^2P).
\end{aligned}
\label{pulse:eq:4.28}
\end{equation}
The bumps and \(c_\rho\) were fixed in Section~\ref{pulse:sec:3.1}. In \(C^z\) the axial derivative acts on \(\lambda^2\) as well as on \(P\). These tangential terms give the lift both compatibility identities. Their force has size \(F_k O_m(\varepsilon^{6/5}\|\mathbf d\|_{m+1})\).

Write \(u=V+a+w\), with \(\langle a\rangle_Y=0\), and let
\begin{equation}
w=\sum_\gamma U_\ell\varepsilon_\ell^{1/2}\sqrt{\omega}G_\gamma W_\gamma(x,\eta_\gamma),\qquad \eta_\gamma=N_\ell\phi_\gamma,\qquad \int_{\mathbb T} W_\gamma\,d\eta=0.\label{pulse:eq:4.29}
\end{equation}
The real smooth functions \(W_{\gamma }\) may contain every harmonic. Since \(N_\ell p_{\theta,\gamma}\in\mathbb Z\setminus\{0\}\), carrier and angular means agree. Also, \(\partial_\theta=N_\ell p_{\theta,\gamma}\partial_\eta\). Normalize the state by
\[
Z_{\rm st}=\bigl(V^{\theta,z}/U_k,V^r/(U_k\varepsilon),a^{\theta,z}/(U_k\varepsilon^{1-\delta_0}\omega),a^r/(U_k\varepsilon^{2-\delta_0}\omega),(W_\gamma)_\gamma\bigr).
\]

\noindent\textbf{Proposition \PulseResult{4.4}} (Canonical moment operator)\textbf{.} \textit{Assume the uniform inverse of Proposition~\ref{pulse:thm:3.8}} \textit{and polynomial coefficient bounds in the fixed charts. Assume that} \((\mathbf d,Z_{\mathrm{st}})\) \textit{has a bounded low joint Sobolev norm. Polynomial-radius bounds are allowed when} \(\varepsilon^{6/5}\mathbf d\) \textit{satisfies the smallness gate. The exact swirl repair has defect change} \(-\mathbf D\)\textit{. Apply Proposition~\ref{pulse:thm:3.6} to} \(H=\rho_r P-[2(V^\theta+a^\theta)b+b^2]/r\)\textit{. Its integrated double mean vanishes. Let} \(q, R_{H}\) \textit{be its pressure and tail. Then}
\begin{equation}
\mathcal M(\mathbf d,Z_{\mathrm{st}})=\NS(u+b e_\theta,p+q)-\NS(u,p)+\mathcal S(\mathbf D)-R_H e_r\label{pulse:eq:4.30}
\end{equation}
\textit{has zero radial angular mean and zero compatible defect coordinates. Its double tangential mean has zero weighted radial integrals. In the grading} \(\mathcal Y^m\) \textit{of} \eqref{pulse:eq:4.44}\textit{,}
\begin{equation}
\lVert\mathcal M(\mathbf d,Z_{\mathrm{st}})\rVert_{\mathcal Y^m}\le P_m(k)\varepsilon^{1/2}\bigl(\lVert\mathbf d\rVert_{m+2}+\lVert\mathbf d\rVert_0\lVert Z_{\mathrm{st}}\rVert_{m+2}\bigr),\label{pulse:eq:4.31}
\end{equation}
\textit{and}
\begin{equation}
\begin{aligned}&\lVert D\mathcal M(\mathbf d,Z_{\mathrm{st}})(\dot{\mathbf d},\dot Z_{\mathrm{st}})\rVert_{\mathcal Y^m}\\[-1pt]&\qquad\le P_m(k)\varepsilon^{1/2}\bigl[\lVert(\dot{\mathbf d},\dot Z_{\mathrm{st}})\rVert_{m+2}+(1+\lVert(\mathbf d,Z_{\mathrm{st}})\rVert_{m+2})\lVert(\dot{\mathbf d},\dot Z_{\mathrm{st}})\rVert_0\bigr].\end{aligned}\label{pulse:eq:4.32}
\end{equation}
\textit{Constants may depend polynomially on the low norms and inverse margin. Increase the fixed base order once if needed. The loss is two derivatives, independent of the harmonic range.}

\textit{Proof of Propositions~\ref{pulse:thm:4.3} and~\ref{pulse:thm:4.4}.}
The moment solve in Proposition~\ref{pulse:thm:3.8} gives the bound on \(b\) and \(\int[\rho_r P-(2V^\theta b+b^2)/r]\,dr=0\). Thus
\begin{equation}
q(r)=\int_r^\infty\!\left[\rho_s P-\frac{2V^\theta(s)b(s)+b(s)^2}{s}\right]ds
\label{pulse:eq:4.25}
\end{equation}
is compact in \(K\), removes the double radial mean, and satisfies \eqref{pulse:eq:4.22} and \(\partial _{z}q/F_{k}=O_{m}(\eta\varepsilon^{6/5})\). It preserves the repaired defects by \eqref{pulse:eq:3.22}. Direct substitution gives
\begin{equation}
\begin{aligned}
\langle\mathcal R^r\rangle_\theta&=-2a^\theta b/r,\\
\langle\mathcal R^\theta\rangle_\theta&=\partial_t b-(\partial_r^2+r^{-1}\partial_r+\partial_z^2-r^{-2})b\\
 &\quad +(V^r+a^r)(\partial_r b+b/r)+(V^z+a^z)\partial_z b,\\
\langle\mathcal R^z\rangle_\theta&=\partial_z q .
\end{aligned}
\label{pulse:eq:4.26}
\end{equation}
No derivative falls on \(a\). The double mean of its terms is zero. By \eqref{pulse:eq:4.20}, the terms involving \(V\), the time derivative and the viscous terms have size \(F_{k}O_{m}(\eta\varepsilon^{6/5})\); \(2a^{\theta }b/r\) has size \(F_{k}O_{m}(\eta\varepsilon^{6/5-\delta_0}\omega )\). The other \(a\) terms gain another \(\varepsilon\). This proves the mean estimates, since \(\omega\) is bounded below on \(K\).

Subtracting the angular mean leaves the interaction of the swirl repair with \(w\). Its three components are
\begin{equation}
\begin{aligned}
 \mathcal R_o^r&=(b/r)\partial_\theta w^r-2bw^\theta/r,\\
 \mathcal R_o^\theta&=(b/r)\partial_\theta w^\theta+w^r(\partial_r b+b/r)+w^z\partial_z b,\\
 \mathcal R_o^z&=(b/r)\partial_\theta w^z .
\end{aligned}
\label{pulse:eq:4.27}
\end{equation}
In \eqref{pulse:eq:4.27}, derivatives of \(w\) occur only in the angular direction. The largest such contribution, \((b/r)\partial_\theta w\), is \(F_k O_m(\eta\varepsilon^{1/5}\sqrt{\omega}G_\gamma)\) since \(N_\ell\varepsilon^{1/2}=O(1)\). For the curl correction the exponent is \(7/10-\delta_0\), and terms without a carrier derivative have exponent at least \(7/10\). This gives \eqref{pulse:eq:4.24} on each label. Viscosity uses two slow derivatives of \(b\) and the pressure gradient uses one, proving Proposition~\ref{pulse:thm:4.3}.

For the canonical operator, differentiating \(Bc+\mathcal N(c)=d\) gives \((B+D\mathcal N(c))\dot c=\dot d-\dot B c\). The uniform inverse proves order-zero tame dependence of \(b\) on the target and actual means.

For a general waveform \eqref{pulse:eq:4.29}, the same component formula
\eqref{pulse:eq:4.27} has leading term \((b/r)N_\ell p_{\theta,\gamma}\partial_\eta W_\gamma\). It gains
\(\varepsilon^{3/10-1+6/5}=\varepsilon^{1/2}\) in the grading and costs one carrier derivative.
The coefficient \(b\) is independent of the carrier, so multiplication preserves the
harmonics and labelled supports. Sobolev multiplication and label summation give the bound.

Split \(q=q_{0}+q_{a}\), with \(q_{0}\) given by \eqref{pulse:eq:4.25}. The source of \(q_{a}\) is \(-2a^{\theta }b/r\). Equations \eqref{pulse:eq:4.26} and \eqref{pulse:eq:3.25} give
\[
\begin{aligned}
 \mathcal M_d^{\rm tan}/(F_k\omega)&=O_m(\varepsilon^{6/5}),\\
 \mathcal M_a^\theta&=a^r(\partial_r b+b/r)+a^z\partial_z b=F_k O_m(\varepsilon^{11/5-\delta_0}\omega),\\
 \mathcal M_a^z&=\partial_z q_a=F_k O_m(\varepsilon^{11/5-\delta_0}\omega).
\end{aligned}
\]
The lift terms have the same \(\varepsilon^{6/5}\) scale. Thus the mean gains are \(\varepsilon\) and \(\varepsilon^{2-\delta_0}\).
Pressure and its zero-mean tail preserve the defect change \(-\mathbf D\).
Adding the lift cancels it, so \eqref{pulse:eq:3.2} gives zero weighted integrals.

Viscosity uses two slow derivatives of \(b\). Pressure and the lift each cost one slow derivative. The tame product estimates prove \eqref{pulse:eq:4.31}–\eqref{pulse:eq:4.32}. Uniform bounds on finite Fourier sums extend the result to every smooth carrier waveform. \hfill\(\square\)

\textbf{The frozen linear operator.} Freeze the initialized means and solve
\(B_0c=d\) for the linear moment correction. This is the differential of
the canonical moment correction at the initialized state. Estimate
\eqref{pulse:eq:4.32} gives its bounds in Section~\ref{pulse:sec:4.4}.
The pressure projection also acts on the full quadratic radial force.

\subsection{The auxiliary-mean correction}\label{pulse:sec:4.3}

Use the band coordinates \((\widehat R,\widehat Z,\widehat T,Y)\) and, on labelled charts,
the carrier angle \(\eta\) and native clock \eqref{pulse:eq:2.22}. Choose \(s\) above
half the largest chart dimension. We use the trimmed weights of Lemma~\ref{pulse:thm:3.2}
and the cutoff and compensator in \(K\), with polynomial derivative bounds.
On the common cover \(Y=Y_{i_0}\), evaluate a level-\(i\) coefficient at
\(E^{i-i_0}Y\) for \(i\ge i_0\). The physical derivatives are
\begin{equation}
\begin{aligned}
 D_t&=-2^{2k}\partial_{\widehat T}+\Lambda_t^{i_0}v_t\cdot\nabla_{Y_{i_0}},
 &\partial_z&=2^{(1-\delta)k}\partial_{\widehat Z},\\
 D_r&=2^k\bigl(c_1\partial_{\widehat R}+M c_2 v_r\cdot\nabla_{Y_{i_0}}\bigr),
 &M&=2^{-kd_0}\Lambda_r^{i_0}.
\end{aligned}\label{pulse:eq:4.33}
\end{equation}
Here \(c_{1}, c_{2}\) are fixed smooth functions of \(\widehat R\). By \eqref{pulse:eq:2.13} and \eqref{pulse:eq:2.14}, \(F_k/(U_k\Lambda_t^{i_0})\le P(k)\varepsilon\) and \(P(k)^{-1}\varepsilon^{-\delta_0}\le M\le P(k)\varepsilon^{-\delta_0}\).

Write the initialized field as \(u=V+a+w\), with \(\langle a\rangle_Y=\langle w\rangle_\theta=0\). Let \(A_{m}\) bound the joint \(H^{m}\) norms of
\begin{equation}
\frac{V^{\theta,z}}{U_k},\qquad \frac{V^r}{\varepsilon U_k},\qquad \frac{a^{\theta,z}}{U_k\varepsilon^{1-\delta_0}\omega},\qquad \frac{a^r}{U_k\varepsilon^{2-\delta_0}\omega},\qquad W_\gamma=\frac{w_\gamma}{U_k\varepsilon^{1/2}\sqrt{\omega}G_\gamma}.\label{pulse:eq:4.34}
\end{equation}
In the same label norm, physical reconstruction satisfies
\begin{equation}
\left\|\frac{\nabla w_\gamma}{2^k U_k\sqrt{\omega}G_\gamma}\right\|_{H^m}
\le P_m(k)\|W_\gamma\|_{H^{m+1}}.
\label{pulse:eq:4.35}
\end{equation}
This follows from \eqref{pulse:eq:2.51} and \eqref{pulse:eq:2.77}. A carrier derivative cancels \(\varepsilon ^{1/2}\) against \(N_\ell\asymp\varepsilon^{-1/2}\). Coefficient and frame derivatives are smaller. Thus the bound is independent of the harmonic range. Section~\ref{pulse:sec:4.5} gives the label summability.

For a trial source \(g\), we first correct the nonzero angular modes and the tangential double mean (Section~\ref{pulse:sec:4.4.3}). From the remaining tangential force we take the angular mean and subtract its mean in \(Y_{i_0}\). This gives the two components \(h_{\rm aux}^\theta,h_{\rm aux}^z\) to be canceled below.

\noindent\textbf{Proposition \PulseResult{4.5}} (Auxiliary correction)\textbf{.} \textit{Let \(h_{\rm aux}=F_k\varepsilon^{1/5}\omega\,\widehat h_{\rm aux}\), where the normalized input \(\widehat h_{\rm aux}\) is smooth, trimmed and has zero \(Y_{i_0}\)-mean. Set}
\begin{equation}
\begin{gathered}
 (A,B)=-\Lambda_t^{-i_0}(v_t\cdot\nabla_{Y_{i_0}})^{-1}(h_{\rm aux}^\theta,h_{\rm aux}^z),\\
 \Psi=r^{-1}I_0(rB),\qquad v=Ae_\theta+\nabla\times(\Psi e_\theta).
\end{gathered}
\label{pulse:eq:4.36}
\end{equation}
\textit{Here the inverse is the zero-mean inverse on the level-\(i_0\) torus. The velocity is real for real input, divergence-free and compact, with zero auxiliary mean. For} \(m\ge s\)\textit{,}
\begin{equation}
\left\|\frac{v^{\theta,z}}{U_k\varepsilon^{6/5}\omega}\right\|_{H^m}
+\left\|\frac{v^r}{U_k\varepsilon^{11/5}\omega}\right\|_{H^m}
+\left\|\frac{\Psi}{2^{-k}U_k\varepsilon^{6/5}\omega}\right\|_{H^m}
\le P_m(k)\|\widehat h_{\rm aux}\|_{H^{m+2}}.\label{pulse:eq:4.37}
\end{equation}
\textit{There are compact pressures} \(q=q_{a}+q_{d}\) \textit{and zero-mean tails} \(E_{a}, E_{p}\) \textit{such that}
\begin{equation}
\mathcal R=\mathrm{NS}(u+v,p+q)-\mathrm{NS}(u,p)
+h_{\rm aux}^\theta e_\theta+h_{\rm aux}^z e_z-E_a-E_p e_r
\label{pulse:eq:4.38}
\end{equation}
\textit{has} \(\langle\mathcal R^r\rangle_\theta=\rho_r\Delta P\)\textit{. Normalize its defect increment as} \(\Delta d=(U_k^{-2}\Delta P,2^{(1-2\delta)k}\Delta J^{\theta z},2^{-2\delta k}\Delta J^{z\theta})\)\textit{. Then}
\begin{gather}
\sum_\gamma
 \left\|\frac{\mathcal R_{o,\gamma}}{F_k\sqrt{\omega}\,G_\gamma}\right\|_{H^m}
 \le P_m(k)\varepsilon^{1/5}\mathcal N_m(\widehat h_{\rm aux}),\label{pulse:eq:4.39} \\
\left\|\frac{\langle\mathcal R\rangle_\theta}{F_k\omega}\right\|_{H^m}\le P_m(k)\varepsilon^{6/5-2\delta_0}\mathcal N_m(\widehat h_{\rm aux}),\label{pulse:eq:4.40} \\
\|\Delta d\|_{H^m}\le P_m(k)\varepsilon^{11/5-2\delta_0}\mathcal N_m(\widehat h_{\rm aux}),\label{pulse:eq:4.41}
\end{gather}
\textit{where}
\[
\mathcal N_m(\widehat h_{\rm aux})=
\bigl(1+A_{s+5}+\|\widehat h_{\rm aux}\|_{H^{s+5}}\bigr)
\|\widehat h_{\rm aux}\|_{H^{m+5}}+A_{m+5}\|\widehat h_{\rm aux}\|_{H^{s+5}}.
\]
\textit{The source grading gives a gain of at least} \(\varepsilon ^{1/2}\)\textit{. The projected force costs five derivatives of the given auxiliary input, or six if that input must first be decoded. The loss is independent of the harmonic range.}

\textit{Proof.} Proposition~\ref{pulse:thm:2.2} costs one auxiliary derivative.
Proposition~\ref{pulse:thm:3.4} gives the compact primitive at the same
Sobolev order. Applied to \eqref{pulse:eq:4.36}, as in
\eqref{pulse:eq:3.52}--\eqref{pulse:eq:3.53}, these give a divergence-free
field of zero Haar mean. Its tangential velocity has size
\(U_k\varepsilon^{6/5}\omega\), and
\(v^r=-\partial_z\Psi\) gains a further \(\varepsilon\).
This proves \eqref{pulse:eq:4.37}. Keep the translated remainder
\(-\chi_m'\mathcal J(rB)/r\) in \(v^z\). Its fast time derivative gives
\(E_a=(\chi_m'/r)\mathcal J(rh_{\rm aux}^z)e_z\) after cancellation of \(h_{\rm aux}\).

Let \(R_0\) denote the residual in \eqref{pulse:eq:4.38} before pressure projection
(\(q=E_p=0\)). Apply the component identities of Proposition~\ref{pulse:thm:3.14}
with \((u_B,v_A,w_0)=(V+a,v,w)\). Their tangential fast derivatives cancel
\(h_{\rm aux}\), and their tensor has entries \((V+a)^i v^j+v^i(V+a)^j+v^iv^j\).
The oscillatory terms are \((w\cdot\nabla)v+(v\cdot\nabla)w\).

Use \eqref{pulse:eq:4.20}. Tangential slow time and viscosity have exponents \(6/5\) and \(6/5-2\delta_0\) in force units. Background radial transport has exponent \(6/5-\delta_0\), because \(V^{r}\) and \(v^{r}\) each have an extra \(\varepsilon\). Axial transport has exponent \(6/5\). Old auxiliary and quadratic terms are smaller and are retained in \(\mathcal N_m\). The radial force is
\begin{equation}
\left\|\frac{\langle R_0^r\rangle_\theta+2V^\theta v^\theta/r}{F_k\omega}\right\|_{H^m}
 \le P_m(k)\varepsilon^{6/5-2\delta_0}\mathcal N_m(\widehat h_{\rm aux}).\label{pulse:eq:4.42}
\end{equation}
In force units the fast time derivative of the radial velocity has exponent \(6/5\); the two lower centrifugal contributions have exponents \(6/5-\delta_0\) and \(7/5\). The leading centrifugal term in \eqref{pulse:eq:4.42} has exponent only \(1/5\). Its auxiliary mean is zero, so the radial pressure projection can remove it without imposing another moment condition.

For \(S=\langle R_0^r\rangle_\theta\), set \(\Delta P=\int\langle S\rangle_Y\,dr\) and choose
\[
q_a=I_0(\langle S\rangle_Y-S),\qquad
 q_d(r)=\int_0^r[-\langle S(s)\rangle_Y+\rho_s\Delta P]\,ds.
\]
Proposition~\ref{pulse:thm:3.6} gives \(D_rq_a=\langle S\rangle_Y-S+E_p\) and zero auxiliary means for \(q_{a}, E_{p}\). The axial pressure derivative gains \(\varepsilon\) and costs one derivative. The compensated \(q_{d}\) has the same bound and no tail. Thus the complete leading radial force is projected to \(\rho_r\Delta P\). This proves \eqref{pulse:eq:4.40}.

By \eqref{pulse:eq:4.35}, \((v\cdot  \nabla )w_{\gamma }\) has size \(F_{k}\varepsilon ^{1/5}\omega ^{3/2}G_{\gamma }\) times its tame product. The term \((w_{\gamma }\cdot  \nabla )v\) has exponent \(7/10-\delta_0\). Since \(\omega^{3/2}\le\sqrt{\omega}\), \eqref{pulse:eq:4.39} follows. Each term retains its pulse support and Gaussian.

After double averaging, all terms linear against \(V\) vanish. The remaining covariance is \(\langle a^iv^j+v^ia^j+v^iv^j\rangle_Y\). Its tangential size is \(U_k^2\varepsilon^{11/5-\delta_0}\); each radial factor adds \(\varepsilon\). The averaged radial force is at most \(F_{k}\varepsilon ^{6/5-\delta_0 }\). Radial integration contributes \(2^{-k}F_k/U_k^2=\varepsilon\) to \(\Delta P\). The flux measures \(r^{2}dr\) and \(rdr\) give the same normalized bound. The radial-force part of \(J^{z\theta }\) also has this bound. Compact pressure preserves the defects by \eqref{pulse:eq:3.22}. This proves \eqref{pulse:eq:4.41}.

Viscosity contains at most two powers of \(M\), independent of \(m\). The torus inverse, curl, force and axial pressure gradient cost at most \(1 + 1 + 2 + 1\) derivatives. Sobolev products give \(\mathcal N_m\). The formulas have degree at most two in \(\widehat h_{\rm aux}\). Their first differences therefore have the same loss and gain. \hfill\(\square\)

Zero auxiliary mean preserves both velocity moments. The residual and its actual defect increments satisfy \eqref{pulse:eq:3.2}, with all tangential compensators retained.

Both tails contain only nonzero auxiliary modes. By \eqref{pulse:eq:3.17}, for each fixed \(J\),
\begin{equation}
\left\|\frac{E_a+E_p e_r}{F_k\varepsilon^{1/5}\omega}\right\|_{H^m}
\le P_{m,J}(k)M^{-J}\mathcal N_{m+2J}(\widehat h_{\rm aux}).
\label{pulse:eq:4.43}
\end{equation}
The tails lie in \(K\), where \(\inf_K\omega>0\). For each physical derivative and
target power of \(\lambda\), choose \(J\) in \eqref{pulse:eq:4.43} and use the
polynomial bound on the same input at order \(m+2J\). A Bernstein estimate
from a low norm would introduce powers of the Fourier cutoff, which may exceed \(M\).

\subsection{The complete source map}\label{pulse:sec:4.4}

All corrections are taken at the state \((u_*,p_*)\) of
Proposition~\ref{pulse:thm:3.14}.

\subsubsection{Source units and initialization}\label{pulse:sec:4.4.1}
\pulsebridge{4.4.1}

Use the units
\begin{equation}
\begin{array}{l|l}
 \text{source component}&\text{physical unit}\\ \hline
 \text{nonzero carrier }\gamma&F_k\varepsilon^{-3/10}\sqrt{\omega}\,G_\gamma\\
 \text{angular mean}&F_k\varepsilon^{1/5}\omega\\
 (P,J^{\theta z},J^{z\theta})&\varepsilon^{6/5}U^2(1,2^{-3k},2^{-2k})
\end{array}
\label{pulse:eq:4.44}
\end{equation}

Take the joint Sobolev norms with the label sums of Section~\ref{pulse:sec:4.5}. Fix an algebra order \(s_{\mathrm{alg}}\) and set \(s_{\mathrm{trial}}=s_{\mathrm{alg}}+ 33\). The radial potentials have units \(F_k \varepsilon ^{1/5}2^{-3k}\omega\) and \(F_k \varepsilon ^{1/5}2^{-2k}\omega\). Decoding costs one derivative. The inverse coordinate map has order zero by Proposition~\ref{pulse:thm:3.1}.

At the initialized state, let \(\mathbf f_*^{\rm raw}\) denote the free coordinates of
\(f_*=\NS(u_*,p_*)-E_*\) in Proposition~\ref{pulse:thm:3.14}.
All three defect coordinates and the radial angular mean are zero. Thus the two tangential
double means have zero weighted integrals by \eqref{pulse:eq:3.2}, and their potentials
are the compact primitives in \eqref{pulse:eq:3.4}. The raw estimates
\eqref{pulse:eq:3.50} permit algebraic factors \(d_e^{-q_m}\) at the two edges.

Use the cutoff \(\chi_{\rm in}\) from \eqref{pulse:eq:3.47}. Define
\(\mathbf f_*^{\rm in}\) by multiplying both tangential radial potentials,
the auxiliary fluctuation and each nonzero carrier component of
\(\mathbf f_*^{\rm raw}\) by \(\chi_{\rm in}\). Keep \(P,J^{\theta z},J^{z\theta}\)
unchanged. Write \(f_*^{\rm in}\) for the force obtained by decoding these coordinates.
The cutoff \(\chi_r\) equals one near their radial support and the middle window \(K\).
For every fixed \(m\), in the native charts with the label norm of
Section~\ref{pulse:sec:4.5},
\begin{equation}
\begin{aligned}
 \left\|\frac{(f_*^{\rm in})_{o,\gamma}}{F_k\sqrt\omega G_\gamma}\right\|_{H^m}
     &\le P_m(k)\varepsilon^{-\delta_0},&
 \left\|\frac{\langle f_*^{\rm in}\rangle_\theta}{F_k\omega}\right\|_{H^m}
     &\le P_m(k)\varepsilon^{1-3\delta_0},\\
 \langle (f_*^{\rm in})^r\rangle_\theta&=0,&
 P_*^{\rm in}=J_*^{\theta z,{\rm in}}=J_*^{z\theta,{\rm in}}&=0 .
\end{aligned}
\label{pulse:eq:4.45}
\end{equation}
Here the first inequality applies to each active label, and their polynomial count is included in \(P_m(k)\).
On the retained support \(d_e\ge c/k\). The logarithmic derivatives of \(\omega\)
and all fixed derivatives of \(\chi_{\rm in}\) in normalized coordinates cost only powers of \(k\).
Equation~\eqref{pulse:eq:3.50} therefore bounds the cutoff fluctuations in
\eqref{pulse:eq:4.45}. For the double mean, use the two tail formulas in
\eqref{pulse:eq:3.4} and the same-side tail estimate of Lemma~\ref{pulse:thm:3.3}
before cutting off the potentials. If \(g^\theta,g^z\) are the raw potentials,
the zero defects give
\[
 (f_*^{\rm in})_d^a
   =\chi_{\rm in}(f_*)_d^a+r^{-q_a}(\partial_r\chi_{\rm in})g^a,
 \qquad a=\theta,z,\qquad (q_\theta,q_z)=(2,1).
\]
The primitive supplies one radial length, which cancels the \(2^k\) in the physical derivative
of \(\chi_{\rm in}\). The remaining normalized cutoff derivatives are polynomial on the retained strip.
This proves the mean bound, including its two commutator terms. Weighted integration of
\(\partial_r(\chi_{\rm in} g^a)\) gives zero at both endpoints, so the two compatibility identities
hold after decoding.

Let \(\operatorname{Dec}\) denote the free-coordinate decoder, including the defect entries.
With the initialization error \(E_*\) embedded as a force with zero defect entries, define
\[
 \mathbf T_*=(E_*,0)+
      \operatorname{Dec}(\mathbf f_*^{\rm raw}-\mathbf f_*^{\rm in}).
\]
The decoded difference, including the derivatives of \(\chi_{\rm in}\), is supported in \(d_e\le C/k\).
Every fixed physical derivative is bounded by \(P_m(k)e^{-c_m k^2+C_m k}\)
by Lemma~\ref{pulse:thm:3.12}; both sources have identical defect entries.
Proposition~\ref{pulse:thm:3.14} controls \(E_*\).
Only \(\mathbf f_*^{\rm in}\) enters the window equation. Its normalized fixed-order norms
in the units of \eqref{pulse:eq:4.44} are polynomial in \(k\).

The inverse bounds are \(g_{\mathrm{match},\gamma}^2\ge ck\omega\) and the mixed Gram margin. After initialization the double means differ from the background by \(O_{m}(\varepsilon ^{1-\delta_0 })\). Thus the frozen matrix \(B_{\ast }\) satisfies \(\lVert B_*-G_3\rVert\le g_3/2\) on high bands.

\subsubsection{The stress inverse}\label{pulse:sec:4.4.2}
\pulsebridge{4.4.2}

Carry the physical defects of \eqref{pulse:eq:3.2}, including the full
compact intermediate radial force. Apply \eqref{pulse:eq:3.28}--\eqref{pulse:eq:3.30} to the full tangential
double mean \(h_d\), with \(\beta=6/5\). The order-zero primitive takes
\(h_d/(F_k\omega)=O_m(\varepsilon^{1/5})\) to \(T/(U^2\omega)=O_m(\varepsilon^{6/5})\).
The two compensators have size \(F_k O_m(\varepsilon^{6/5}\omega)\) by
\eqref{pulse:eq:3.2}, at the cost of one slow derivative. They remain in the source.

The bound on \(T/(U^2\omega)\) lets us apply Proposition~\ref{pulse:thm:3.10}
with \(\beta=6/5\). Thus \(\mathcal B(W^{\rm pre},w_S^{\rm pre}[T])=T\)
for either sign of each component.
With the full curl and pressure of \eqref{pulse:eq:3.37} and
\eqref{pulse:eq:4.3}, the linear pair \((w_S[T],p_S[T])\) satisfies
\begin{equation}
w_S^{\rm pre}/U=O_m(\varepsilon^{7/10}\sqrt{\omega}G),\qquad (w_S-w_S^{\rm pre})/U=O_m(\varepsilon^{6/5-\delta_0}\sqrt{\omega}G).\label{pulse:eq:4.48}
\end{equation}

The same seven-group estimate applies to homogeneous stress pulses. Their transverse amplitude and pressure have sizes \(\varepsilon^{7/10}\sqrt{\omega}G\) and \(\varepsilon^{6/5}\sqrt{\omega}G\). Their principal source consists only of the cutoff term.

The covariance error in Proposition~\ref{pulse:thm:3.10} has size
\(U^2 O_m(\varepsilon^{11/5-2\delta_0}\omega)\). Its double-mean force is
\(F_k O_m(\varepsilon^{6/5-2\delta_0}\omega)\), and the axial cross-flux
contributes \(F_k O_m(\varepsilon^{6/5}\omega)\).

Proposition~\ref{pulse:thm:4.2} bounds the nonzero force by \(F_k O_m(\varepsilon^{1/5-3\delta_0}\sqrt{\omega}G)\), with three derivatives lost. The auxiliary fluctuation and radial mean may be as large as \(F_k O_{m}(\varepsilon ^{1/5-\delta_0 }\omega )\). Pass both to the next operators. Also pass the three actual defect increments, of unit size \(\varepsilon ^{6/5}\) in physical moment units.

\subsubsection{One linear reconstruction}\label{pulse:sec:4.4.3}
\pulsebridge{4.4.3}

Let \(\mathscr L_*(v,q)\) denote the linearized force and its three
defect variations at \(u_*\). For a decoded source \(\mathbf G\), we carry
the full residual through the following five corrections.
\begin{enumerate}[leftmargin=*,nosep]
\item Cancel the nonzero angular modes by Proposition~\ref{pulse:thm:4.2},
retaining the covariance divergence and defect increments
\eqref{pulse:eq:4.16}--\eqref{pulse:eq:4.19}.
\item Apply \eqref{pulse:eq:3.28} and the full-curl stress inverse to the
complete tangential double mean.
\item Cancel \(h_{\rm aux}\), the auxiliary fluctuation of the remaining
tangential angular mean, by \eqref{pulse:eq:4.36}. Its radial force includes
\(-2V_*^\theta v_A^\theta/r\) and still enters the moment calculation.
\item For the resulting three defects \(D\), solve \(B_*c=d(D)\) with the
signs and units of \eqref{pulse:eq:3.35}. The swirl \(v_M=b_Me_\theta\)
has linearized defects \(-D\) and satisfies \(\int r^2 b_M\,dr=0\).
\item Apply Proposition~\ref{pulse:thm:3.6} to the full radial angular mean.
Since the repaired pressure defect is \(P=0\), only the translated pressure
tail remains in that component.
\end{enumerate}
Writing \(g\) for the free coordinates of the decoded source \(\mathbf G=\operatorname{Dec}(g)\),
these linear operations give the identity
\begin{equation}
\mathscr L_*(\mathcal V g,\mathcal Q_1 g)=\operatorname{Dec}(-g+\mathbf E g)+\mathbf T_1 g.\label{pulse:eq:4.49}
\end{equation}
In \eqref{pulse:eq:4.49}, \(\mathbf T_1\) contains the time-cutoff tail and the remainders of the translated primitives. The latter are flat by \eqref{pulse:eq:4.43} once every fixed-order norm of the input has a polynomial bound.

The pulses have zero angular mean, the auxiliary velocity has zero Haar mean,
and the pure swirl has zero weighted angular moment and no axial component.
Hence \eqref{pulse:eq:3.1}--\eqref{pulse:eq:3.2} hold throughout.
The defect lift is \eqref{pulse:eq:4.28}.

\subsubsection{Fixed gain and derivative loss}\label{pulse:sec:4.4.4}
\pulsebridge{4.4.4}

For a unit source in \eqref{pulse:eq:4.44}, the table gives the
powers of \(\varepsilon\) in the velocity and force bounds after their
weights are removed. A superscript \(+\) marks an input to the next step.
Here \(d,a,r\) denote the double mean, auxiliary fluctuation and radial
component, respectively. The last column counts derivatives.
\begin{equation}
\begin{array}{l|c|c|l|c}
 \text{correction}&v/U&\text{nonzero force}/F_k&\text{mean force}/F_k&\text{loss}\\ \hline
 \text{nonzero modes}&7/10&1/5-3\delta_0&(1/5-\delta_0)^+&3\\
 \text{stress}&7/10&1/5-3\delta_0&
 6/5-2\delta_0\ (d),\quad(1/5-\delta_0)^+\ (a,r)&3\\
 \text{auxiliary}&6/5-\delta_0&1/5-\delta_0&
 6/5-3\delta_0\ ({\rm tan}),\quad(1/5-\delta_0)^+\ (r)&4\\
 \text{moments}&6/5&1/5&
 6/5-\delta_0\ ({\rm tan}),\quad(1/5)^+\ (r)&2\\
 \text{pressure}&\text{zero}&\text{zero}&6/5-\delta_0\ (z),\quad\text{zero}\ (r)&1
\end{array}\label{pulse:eq:4.51}
\end{equation}
Axisymmetric radial velocities gain another \(\varepsilon\).
Before averaging, a radial derivative of the auxiliary input costs
\(\varepsilon^{-\delta_0}\). The radial mean and all three defects vanish
in the remaining free source. Including one derivative for decoding, we use
\begin{equation}
 \underbrace{1+3+3+4+2+1}_{\text{decoding and the five rows}}=14,\qquad q=33 .
 \label{pulse:eq:4.52}
\end{equation}
The auxiliary row uses the linear part of the unprojected \(R_0\) in
\eqref{pulse:eq:4.42}. The torus inverse and curl each cost one derivative,
and the force costs two.
The last row takes the axial pressure derivative after moment repair. These
two rows account for the five derivatives in Proposition~\ref{pulse:thm:4.5}.
The first summand counts decoding once by Section~\ref{pulse:sec:3.1}.
The primitives and matrix inverse have order zero. We compose the estimates
at a common Sobolev order. Using \(H^{m+10}\) throughout to bound \(H^m\)
therefore costs at most \(14+10\le q\). Chart changes cost only \(P_m(k)\).

\noindent\textbf{Proposition \PulseResult{4.6}.}  \textit{Let the radial input blocks have support in \(d_e\ge c/k\), with \(c>0\) fixed, and use the prescribed label and time cutoffs on carrier inputs. Lift the radially independent defect blocks into \(K\). In the coordinates and units of \eqref{pulse:eq:4.44}, for every \(m\ge s_{\rm trial}\),}
\begin{equation}
\|\mathbf E g\|_{H^m}\le P_m(k)\,\varepsilon^{1/2-4\delta_0}\|g\|_{H^{m+33}}.
\label{pulse:eq:4.53}
\end{equation}
\textit{The same bound holds with the radial and label exterior-tail cutoffs.}

\textit{Proof.} Divide the bounds in \eqref{pulse:eq:4.51} by the source units
in \eqref{pulse:eq:4.44}. For the nonzero modes the exponent becomes
\(1/5-3\delta_0+3/10=1/2-3\delta_0\). In the tangential mean the stress
compensators, curl covariance error and axial cross-flux have exponent at least
\(6/5-2\delta_0\), the auxiliary term \(6/5-3\delta_0\), and the axial pressure
gradient \(6/5-\delta_0\) by \eqref{pulse:eq:3.25}. Dividing by the mean unit
leaves at least \(1-3\delta_0\). The radial mean and the three defects have
been removed. The compact primitives encoding the tangential mean have order
zero by \eqref{pulse:eq:3.4} and Lemma~\ref{pulse:thm:3.3}.

The radial hull of each source and \(K\) stays in the retained strip,
where derivatives of the weight cost \(P_m(k)\). Common-cover differentiation
in \eqref{pulse:eq:4.10} costs the auxiliary frequency already counted in
\eqref{pulse:eq:4.51}. Summing \eqref{pulse:eq:4.52} gives the stated loss.
On a radial collar use the same-side primitive of \eqref{pulse:eq:3.8}
or \eqref{pulse:eq:3.18}. On a label collar use Lemma~\ref{pulse:thm:4.7}.
In both cases the global bound applies with the enlarged input cutoff.
\hfill\(\square\)

\subsubsection{The quadratic terms}\label{pulse:sec:4.4.5}
\pulsebridge{4.4.5}

For \(v=\mathcal Vg\), apply the canonical pressure projection to the angular mean of \((v\cdot\nabla)v\) and keep its quadratic physical defects. Write \(\mathcal Q_2(g,g)\) for the pressure, \(\mathbf T_2(g,g)\) for the flat tail and \(\mathbf B(g,g)\) for the remaining free source coordinates. Polarization defines symmetric bilinear maps.

Use the cutoffs and defect lifts of Proposition~\ref{pulse:thm:4.6} for both
inputs. In the following bounds we measure velocities in \(U\), forces in
\(F_k\), and moments in \(U^2(1,2^{-3k},2^{-2k})\).

\emph{Two waves.} Different labels have zero product. On the same label,
\eqref{pulse:eq:4.13} replaces \(Nn\cdot c\) by a coefficient divergence
of a zero-mean primitive. The force exponent is therefore
\(2(7/10)-1-\delta_0=2/5-\delta_0\), or \(7/10-\delta_0\) in the
nonzero-carrier unit. For the nonzero-mode inverse of Proposition~\ref{pulse:thm:4.2}
the product costs at most two derivatives and the pressure gradient one more.
In the double mean the integral of \(v_r\cdot\nabla_Y\) of the covariance
vanishes, so only the slow radial derivative contributes to its radial
divergence. The flux has exponent \(7/5\). Radial integration of the
double-mean radial force gives \(2^{-k}F_k/U^2=\varepsilon\), so the
moment entries also have exponent \(7/5\).

\emph{A wave and an axisymmetric field.} An angular derivative on the wave
costs \(\varepsilon^{-1/2}\) times one carrier derivative. A derivative
on the axisymmetric field costs at most \(\varepsilon^{-\delta_0}\).
The force exponent is at least
\((6/5-\delta_0)+7/10-1-1/2=2/5-\delta_0\), and its angular mean is zero.

\emph{Two axisymmetric fields.} The centrifugal force has exponent
\(2(6/5-\delta_0)-1=7/5-2\delta_0\). Radial advection gains an
\(\varepsilon\) in the radial velocity, which covers one auxiliary derivative.
Axial advection gains an \(\varepsilon\) in the axial derivative.
The fluxes and the integrated radial force therefore have at least the
wave moment exponent \(7/5\).

Pressure gains \(\varepsilon\) on the radial input and preserves the defects
by \eqref{pulse:eq:3.22}. Dividing by the units in \eqref{pulse:eq:4.44}
gives a force gain of at least \(2/5-\delta_0-1/5=1/5-\delta_0\) and a
moment gain of at least \(7/5-6/5=1/5\). Thus

\begin{equation}
\|\mathbf B(g,h)\|_{H^m}\le P_m(k)\,\varepsilon^{1/5-4\delta_0}\bigl(\|g\|_{H^{m+33}}\|h\|_{H^{s_{\rm trial}}}+\|h\|_{H^{m+33}}\|g\|_{H^{s_{\rm trial}}}\bigr).
\label{pulse:eq:4.55}
\end{equation}

Choose \(s_{\rm alg}>6/2\) large enough for the preceding low-order estimates.
The velocity costs at most \(1+3+3+4=11\) derivatives through the auxiliary
row, since the swirl depends on the resulting defects at order zero.
Allow two more for the wave product, one for the axial pressure gradient and
one for the moment coordinates. With the common shift above, each input
costs at most \(11+2+1+1+10=25\le q\). The tame Sobolev product puts one
factor at order \(m\) and the other at \(s_{\rm alg}\). Exchanging the
factors and summing the labels gives \eqref{pulse:eq:4.55}, since
\(s_{\rm alg}+q=s_{\rm trial}\). Only the higher input incurs the Bernstein
factor \(L^{33}\) on the \(H^{s_{\rm trial}}\) ball.

For the reconstruction
\begin{equation}
u=u_*+\mathcal V g,\qquad p=p_*+\mathcal Q_1g+\mathcal Q_2(g,g),
\label{pulse:eq:4.56}
\end{equation}
the exact force–defect identity is
\begin{equation}
\mathfrak F(u,p)=\operatorname{Dec}\bigl(\mathbf f_*^{\rm in}+\mathbf E g+\mathbf B(g,g)-g\bigr)
 +\mathbf T_1g+\mathbf T_2(g,g)+\mathbf T_* .
\label{pulse:eq:4.57}
\end{equation}
In \eqref{pulse:eq:4.57}, the addition of \(\mathbf T_*\) restores the discarded part of the initialized force exactly once. The linear and quadratic estimates both imply a gain \(\varepsilon^{1/20}\) in the units of \eqref{pulse:eq:4.44}.

\subsection{Locality of the reconstruction}\label{pulse:sec:4.5}
\pulsebridge{4.5}

Each label has one outer slow box. Its primary and stress pulses have the assigned inner transverse support, and their slow support stays at least \(ck^{-6}\) from the box boundary. Labels with overlapping slow supports have disjoint auxiliary rectangles by Proposition~\ref{pulse:thm:2.5}.

At the Sobolev order specified in each estimate, we take the \(\ell^2\) sum of the joint Sobolev norms of the labelled blocks in their active variables, with each unlabelled block included once.

In the coefficient frame of Proposition~\ref{pulse:thm:4.1}, each carrier harmonic satisfies
\[
 z_m(\tau)=-\int_0^\tau \widetilde{\mathcal U}_m(\tau,s)
                     Q(s)^{\mathsf T}f_m(s)/G(s)\,ds .
\]
This integral fixes \((x,l_r)\), so their cutoffs commute with the inverse.
Curl completion and pressure then use local derivatives. The clock cutoff
acts separately, and the background displacement fits inside the buffers
by the estimate following \eqref{pulse:eq:2.26}.

\noindent\textbf{Lemma \PulseResult{4.7}} (Exterior locality of labelled outputs)\textbf{.} \textit{Consider a finite composition of labelled principal inverses, curl and pressure reconstructions, local differential products, and averages. Also allow axisymmetric radial and torus inverses, moment repairs, and stress-generated pulses. Every} \textit{fixed carrier and synthesized stress pulse has the inner support above. Outside it, a carrier output requires an exterior same-label input. Mean and defect inputs alone create no exterior labelled} \textit{output.}

\textit{Proof.} The forced inverse fixes \((x,l_r)\), and nested cutoffs
absorb derivatives. Nonlocal radial and torus inverses act on angular
means. Such a mean can produce a carrier only through multiplication by a
labelled field or through a stress pulse. The fixed carriers and stress
pulses have inner support, and different labels have zero product.
Thus an exterior output requires an exterior input of the same label.
This implication persists under finite composition. \hfill\(\square\)

At fixed \(Y\), the short extension in Proposition~\ref{pulse:thm:2.7} moves \(Z\) by \(O(k^{-7})=o(k^{-6})\), so it fits between any fixed number of nested slow cutoffs inside the same outer box.

Let \(X_r,X_\ell,X_t\) denote the radial, slow/transverse and carrier-time cutoffs. The carrier-time multiplier \(X_t\) supplies \(\chi_{\rm src}\) for the forced equation. On radially independent moment coordinates all three act as the identity. For \(X=X_rX_\ell X_t\),
\[
I-X=(I-X_r)+X_r(I-X_\ell)+X_rX_\ell(I-X_t).
\]
The radial and label contributions in this decomposition use the exterior estimates above. For the time contribution we first reconstruct the physical source: on the discarded interval its factor \(G_\gamma\) is bounded by \(e^{-ck^2}\). This bound applies before division by the Gaussian. The physical cutoff is fixed; only the estimating collars vary with the required decay order.

The radial selections use \eqref{pulse:eq:3.8} and \eqref{pulse:eq:3.18}
on the same side. Only the translated middle remainder enters the flat tail
by \eqref{pulse:eq:3.17}. The ordinary compensators stay in the source.

\section{A finite construction on the scale windows}\label{pulse:sec:5}
\pulsebridge{5}
\subsection{A finite solve and its cutoff residual}\label{pulse:sec:5.1}
On one coefficient torus put \(\varepsilon=2^{-\delta k}\). Let \(K_k\) contain the source support and all fixed compensators in the interior of \(\{\chi=1\}\). An exterior radial cutoff selects the entire tail on either side. Assume that any fixed number of nested collars can be placed between \(K_k\) and \(\supp(1-\chi)\) with reciprocal gaps and fixed derivative norms bounded polynomially in \(k\). Assume also that, for each fixed integer \(r\ge0\), \(\|\chi\|_{C^r}\) is bounded by a polynomial in \(k\).

\begin{theorem}[Finite solve on a band]\label{pulse:thm:5.1}
Write \(H_k(v)=f_k+E_kv+B_k(v,v)\), with \(E_k\) linear, \(B_k\) symmetric bilinear, and \(f_k\) supported in \(K_k\) with polynomial norms at each fixed order. For the finite-window estimates in this section, write \(\|\cdot\|_m=\|\cdot\|_{H^m}\). For \(m\ge s\) suppose that
\begin{equation}\label{pulse:eq:5.1}
\begin{aligned}
 \|E_kv\|_m&\le P_m(k)\varepsilon^{1/20}\|v\|_{m+33},\\
 \|B_k(v,w)\|_m&\le P_m(k)\varepsilon^{1/20}
       (\|v\|_{m+33}\|w\|_s+\|w\|_{m+33}\|v\|_s).
\end{aligned}
\end{equation}
For nested exterior cutoffs with \(\widetilde\zeta=1\) near \(\supp\zeta\), assume
\begin{equation}\label{pulse:eq:5.2}
\begin{aligned}
 \|\zeta E_kv\|_m&\le P_m(k)\varepsilon^{1/20}\|\widetilde\zeta v\|_{m+33},\\
 \|\zeta B_k(v,v)\|_m&\le P_m(k)\varepsilon^{1/20}
  \big(\|\widetilde\zeta v\|_{m+33}\|v\|_s+\|\widetilde\zeta v\|_s\|v\|_{m+33}\big).
\end{aligned}
\end{equation}
Set \(L=\varepsilon^{-1/1320}\). For all sufficiently large \(k\), there is a unique \(g\) in the ball specified below such that
\begin{equation}\label{pulse:eq:5.3}
 g=S_L H_k(\chi g),\qquad g\in\operatorname{range}S_L .
\end{equation}
For each fixed \(m,A\) the same solution satisfies
\begin{align}
 \|g\|_m&\le P_m(k),\label{pulse:eq:5.4}\\
 \|(1-\chi)g\|_m+\|H_k(\chi g)-\chi g\|_m&\le C_{m,A}\varepsilon^A .\label{pulse:eq:5.5}
\end{align}
The same conclusion holds if the prescribed source has, in addition, tails bounded by \(C_{m,A}\varepsilon^A\) for every fixed \(m,A\).
\end{theorem}
\begin{proof}
On the range of \(S_L\), take the \(H^s\)-ball of radius \(R_k=1+4\max\{1,\sup_L\|S_L\|_{H^s\to H^s}\}\|f_k\|_s\). Multiplication by the cutoff and Bernstein's inequality give \(\|\chi g\|_{m+33}\le P_m(k)L^{33}\|g\|_m\). For \(g,h\) in this ball, polarization gives
\begin{equation}\label{pulse:eq:5.6}
 \|\Phi(g)-\Phi(h)\|_s\le P_s(k)\varepsilon^{1/20} L^{33}(1+R_k)\|g-h\|_s,
 \qquad \Phi(g)=S_L H_k(\chi g).
\end{equation}
Since \(\varepsilon^{1/20} L^{33}=\varepsilon^{1/40}\) and \(\|\Phi(0)\|_s\le(R_k-1)/4\), the map preserves the ball and has Lipschitz constant at most \(1/2\) for large \(k\). At any larger fixed \(m\), \eqref{pulse:eq:5.3} bounds \(\|g\|_m\) by \(\|S_Lf_k\|_m+P_m(k)\varepsilon^{1/40}(1+R_k)\|g\|_m\). The multiplier is bounded at each fixed order. Absorption proves \eqref{pulse:eq:5.4}. On the finitely many intervening bands Bernstein suffices. In particular \(H_k(\chi g)\) has polynomial norms at every fixed order.

For any fixed \(n\), integration by parts in the smooth Fourier kernel on separated supports gives
\begin{equation}\label{pulse:eq:5.7}
 \|\zeta S_L((1-\widetilde\zeta)h)\|_m\le P_{m,n}(k)L^{-n}\|h\|_s .
\end{equation}
The constants include inverse powers of the separation, which are polynomial for this family of cutoffs. Since \(g=S_{2L}g\), split \(g\) inside that operator using \(\widetilde\zeta\). Apply Bernstein on the near part and \eqref{pulse:eq:5.7} on the far part. This yields
\begin{equation}\label{pulse:eq:5.8}
 \|\zeta g\|_{m+33}\le P_m(k)L^{33}\|\widetilde\zeta g\|_m+P_{m,n}(k)L^{-n}\|g\|_s .
\end{equation}
Choose \(3J+1\) nested interior cutoffs \(\chi_i\) in the region where \(\chi=1\), starting near the support of \(f_k\), with \(\chi_{i+1}=1\) near \(\supp\chi_i\). Use \eqref{pulse:eq:5.7}, \eqref{pulse:eq:5.2} and \eqref{pulse:eq:5.8} successively on the three intervening collars. The remaining global factors are bounded by \eqref{pulse:eq:5.4}. Thus
\begin{equation}\label{pulse:eq:5.9}
 \|(1-\chi_{3\ell})g\|_m\le P_{m,J}(k)\varepsilon^{1/40}\|(1-\chi_{3(\ell-1)})g\|_m
                    +P_{m,J,n}(k)L^{-n}.
\end{equation}
After \(J\) substitutions, the right side is bounded by a polynomial times \(\varepsilon^{J/40}+\varepsilon^{n/1320}\). Choose \(J,n\) with strict slack for the desired power. The last exterior cutoff is one near \(\supp(1-\chi)\). Finally
\(H_k(\chi g)-\chi g=(I-S_L)H_k(\chi g)+(1-\chi)g\).
The polynomial high norms bound the first term, and the collar estimate bounds the second. A prescribed flat tail enters only finitely many substitutions at any fixed \(J\).
\end{proof}

\begin{remark}[Radially independent coordinates]\label{pulse:thm:5.2}
For \(g=(g_{\rm rad},g_{\rm par})\), use the product norm on the active
variables of each block and the block diagonal \(S_L\). Set
\(Xg=(\chi g_{\rm rad},g_{\rm par})\) and
\(Z_\zeta g=(\zeta g_{\rm rad},0)\). With \eqref{pulse:eq:5.1} in the product
norm and \eqref{pulse:eq:5.2} for \(Z_\zeta\), apply \eqref{pulse:eq:5.9}
to the radial block, which alone needs interior support. The parameter
block enters through the global norms and lifts to middle compensators.
Thus the cutoff error is \(((1-\chi)g_{\rm rad},0)\), and the decoder
\eqref{pulse:eq:3.3} retains \(P,J^{\theta z},J^{z\theta}\).
\end{remark}

\subsection{Synthesis on overlapping bands}\label{pulse:sec:5.2}
For cover levels \(k_1\le k_2\), a function \(g(Y_{k_2},\eta)\) pulls back to \(g(E^{k_2-k_1}Y_{k_1},\eta)\). We suppress the unchanged variables, including \(\eta\), below. By Lemma~\ref{pulse:thm:2.3} and Proposition~\ref{pulse:thm:3.4}, the directional inverse and the translated radial primitive commute with this substitution. The same holds for products, normalized Haar averages, moment matrices, physical derivatives and the exact curl. For a label \(\gamma\) native to band \(\ell\), retain its original clock, its zero initial section and its cutoff on every branch. On cover \(k_2\) the fast time derivative is
\[
 \partial_{\tau_\gamma}^{(k_2)}
   =2^{-(2+\delta)\ell}\Lambda_t^{k_2}v_t\cdot\nabla_{Y_{k_2}}.
\]
Since \(Ev_t=\Lambda_t v_t\), the chain rule gives the corresponding derivative on cover \(k_1\). Uniqueness for the constrained initial-value problem on each branch gives
\begin{equation}\label{pulse:eq:5.10}
  \bigl(w[g(E^{k_2-k_1}\cdot)],p[g(E^{k_2-k_1}\cdot)]\bigr)(Y_{k_1})
    =(w[g],p[g])(E^{k_2-k_1}Y_{k_1}).
\end{equation}
Thus the linear velocity, both pressure maps and their residuals commute with passage to the lower cover for data pulled back from above. On \(I_k\) use the common cover \(c_k\) of \eqref{pulse:eq:2.34}.
Lemma~\ref{pulse:thm:2.9} gives the polynomial transfer costs and the
physical scales \eqref{pulse:eq:2.35}.

For a band interval \(I\), let \(\mathcal A_I,\mathcal R_{I,j}\) denote
the synthesis and extension in \eqref{pulse:eq:2.31}, with the sum
restricted to \(I\). By Proposition~\ref{pulse:thm:2.8},
\(\mathcal A_I\mathcal R_I=\mathrm{Id}\) where all active bands are included.

Applying the decoder to \(\zeta_j^2D\) gives two additional terms.
\[
 \begin{aligned}
 C^\theta(\zeta_j^2D)&=\zeta_j^2 C^\theta(D)+(\partial_z\zeta_j^2)J^{\theta z},\\
 C^z(\zeta_j^2D)&=\zeta_j^2 C^z(D)+(\partial_z\zeta_j^2)(J^{z\theta}-c_\rho\lambda^2P).
 \end{aligned}
\]
Their physical lifts have size \(F_k O_m(\varepsilon_k^{6/5}\omega)\) and are included in the equation. Radial differentiation commutes with \(\zeta_j^2\) because \(\partial_r\lambda=0\). Physical radial integrals keep \(t,z\) fixed, and the fast-time integration keeps slow variables fixed. Extension preserves \(\lambda\). Thus each output row depends only on a fixed finite set of neighboring bands. On sufficiently high bands the extension changes \(Z\) by less than the reserved label width.

\begin{lemma}[Localized row estimates on a common cover]\label{pulse:thm:5.2a}
Take \(s=s_{\rm trial}\) and \(\varepsilon_j=2^{-\delta j}\).
On a common auxiliary cover write the canonical linear and polarized quadratic rows as
\[
 \mathsf L_j u=\mathcal R_{I,j} E(\mathcal A_I X u),\qquad
 \mathsf C_j(u,v)=\mathcal R_{I,j} B(\mathcal A_I X u,\mathcal A_I X v),
\]
where \(E\) and \(B\) are the graded linear and quadratic residual maps in \eqref{pulse:eq:4.53} and \eqref{pulse:eq:4.55}. The physical tails
\(T_1,T_2\) are accounted for separately. Let \(N_j\) contain the source
rows meeting row \(j\) through these compositions. Its radius and
cardinality are bounded by fixed stencil constants. Write
\(U_m=\sum_{i\in N_j}\|u_i\|_{H^m}\) in the native graded block norms.

Let \(Z_j\) select either an outward radial collar in the radially
dependent blocks or an exterior slow/transverse label collar in the
specified nonzero carrier blocks. Both selections annihilate independent
moment coordinates. The label selection also annihilates the mean.
Choose enlarged input selections \(\widetilde Z_i\) equal to one on the
corresponding exterior regions needed by the row and put
\(U_m^e=\sum_{i\in N_j}\|\widetilde Z_i u_i\|_{H^m}\).
For \(m\ge s\) the following estimates hold with a polynomial \(P_m\)
independent of the row data.
\[
\begin{aligned}
 \|Z_j \mathsf L_j u\|_{H^m}
 &\le P_m(j)\varepsilon_j^{1/20} U^e_{m+33},\\
 \|Z_j \mathsf C_j(u,u)\|_{H^m}
 &\le P_m(j)\varepsilon_j^{1/20}
       \big(U^e_{m+33}U_s+U^e_s U_{m+33}\big).
\end{aligned}
\]
Taking every selection equal to the identity gives the full-row estimates for the same cut-off input \(X u\).

\end{lemma}
\begin{proof}
Differentiate after physical synthesis. The two \(\partial_z\zeta_i^2\)
terms above have lifts in \(K\); radial integration commutes with the
weights. For each label, choose the nested slow cutoffs inside its fixed outer box, with gaps of order \(k^{-6}\) for any fixed number of cutoffs. At fixed \((\lambda,R,Y,\eta)\) the short extension moves
by \(O(k^{-7})=o(k^{-6})\), so chart pullbacks preserve the
buffers. Equation \eqref{pulse:eq:2.35} bounds the transfer costs
polynomially. Summing \eqref{pulse:eq:4.53} and \eqref{pulse:eq:4.55}
over the stencil gives the full linear and polarized quadratic estimates.

For the radial selection, first apply \(X_r\), which trims the ends
\(d_e<c/k\). The edge primitives in \eqref{pulse:eq:3.4}, \eqref{pulse:eq:3.18} and \eqref{pulse:eq:3.28} integrate
from the corresponding endpoint to the evaluation point. On an outward
collar the integration path stays between that collar and the same radial endpoint. The compensators
used to reconcile the two primitives are supported in \(K\).
Consequently the radial support stays inside the hull of the input
support and \(K\), and the edge estimates apply on the trimmed hull.
Indeed \eqref{pulse:eq:4.28} lifts \(P,J\) with the middle bumps. A tangential
source \(b^a C^a\) has a compensated primitive supported in \(K\).
The signed carrier amplitude in \eqref{pulse:eq:3.42} is pointwise in that primitive,
and the forced inverse holds the slow radius fixed. The swirl repair
and its derivatives also lie in \(K\). Each term in \eqref{pulse:eq:4.27} multiplies
the fixed primary or its derivative by such a factor. Thus a primary
extending across the annulus cannot propagate a middle lift beyond
the radial hull. The torus inverse also fixes the radius. In \eqref{pulse:eq:3.17}, \(\int\langle H\rangle_Y\,dr=0\).
Integration by parts controls the other modes of the translated tails in \(T_1,T_2\).
Ordinary compensators stay in the canonical row.

For the label selection, apply Lemma~\ref{pulse:thm:4.7} with the enlarged input cutoff still inside the same outer slow box and auxiliary rectangle. A mean can reenter a carrier only through a
primary or stress pulse with inner support, or a forced field of the same
label. Thus the selected output vanishes when its exterior labelled inputs vanish.

In either case let \(u_e=(\widetilde Z_i u_i)_{i\in N_j}\). The support statements give
\(Z_j\mathsf L_j(u-u_e)=Z_j\mathsf C_j(u-u_e,u-u_e)=0\). By symmetry the selected quadratic
row is \(Z_j \mathsf C_j(u_e,2u-u_e)\). Applying the full-row estimates and
the polynomial cutoff bounds proves both localized inequalities.

By polarization the derivative of the full row in direction \(h\) satisfies
\[
 \|\mathsf L_jh+2\mathsf C_j(u,h)\|_{H^s}
 \le P_s(j)\varepsilon_j^{1/20}\sum_{i\in N_j}
 \big((1+U_s)\|h_i\|_{H^{s+33}}+U_{s+33}\|h_i\|_{H^s}\big).
\]
This gives \eqref{pulse:eq:5.12}--\eqref{pulse:eq:5.13} on every fixed polynomial ball, including the
enlarged ball used under pullback.

The radial, label and discarded clock terms are those in the cutoff
decomposition at the end of Section~\ref{pulse:sec:4.5}.
\end{proof}

\subsection{Comparison of central components}\label{pulse:sec:5.3}
For the comparison statements, write \(a\) for the gain exponent and \(q\) for the derivative loss. Here \(a=1/20\) and \(q=33\). Fix \(N_j=\{i:|i-j|\le r\}\) containing every input row that can enter output \(j\). Write \(H_j(v)=f_j+R_j(v)\) with \(R_j(0)=0\), where \(R_j\) is the linear plus quadratic part and uses only inputs in \(N_j\). Take \(L_j=\varepsilon_j^{-a/(2q)}\). 

\begin{proposition}[Coupled equations on a window]\label{pulse:thm:5.3}
First assume each common row uses the same coefficient chart and the same Fourier cutoff in every window. Suppose \(f_j\) has polynomial fixed-order norms. On the polynomial-radius trial balls assume
\begin{align}
 \|R_j(Xv)\|_m&\le P_m(j)\varepsilon_j^a\sum_{i\in N_j}\|v_i\|_{m+q},\label{pulse:eq:5.12}\\
 \|D(R_j\circ X)(v)h\|_s&\le P_s(j)\varepsilon_j^a
 \left[\sum_{i\in N_j}\|h_i\|_{s+q}
       +\Big(\sum_{i\in N_j}\|v_i\|_{s+q}\Big)\Big(\sum_{i\in N_j}\|h_i\|_s\Big)\right].\label{pulse:eq:5.13}
\end{align}
The constants may depend polynomially on the low-norm radius. The charts, cutoffs and transfers are fixed by their physical indices, including under zero extension. The maps are continuously differentiable on the trial ball, and any admissibility or invertibility required to define them holds there. Then
\begin{equation}\label{pulse:eq:5.14}
 g_j^{[k]}=S_{L_j}H_j(Xg^{[k]})\ (j\in I_k),\qquad g_i^{[k]}=0\ (i\notin I_k)
\end{equation}
has a unique solution there for large \(k\), with joint polynomial norms at each fixed order. Its central components
\begin{equation}\label{pulse:eq:5.15} g_j=g_j^{[j]}\end{equation}
satisfy
\begin{equation}\label{pulse:eq:5.16}
 \|g_j-S_{L_j}H_j(Xg)\|_m\le C_m e^{-c_mj\log j}.
\end{equation}
\end{proposition}
\begin{proof}
Apply the contraction argument of Theorem~\ref{pulse:thm:5.1} to the product
\(\ell^2(H^s)\) ball of radius \(4\|\Phi(0)\|+1\) for
\(\Phi(v)_j=S_{L_j}H_j(Xv)\). Bounded row and column degree and
\(\varepsilon_j^aL_i^q\le C_r\varepsilon_j^{a/2}\) for \(i\in N_j\)
give contraction by \eqref{pulse:eq:5.13} and high-order absorption
by \eqref{pulse:eq:5.12}. The radius is polynomial since \(|I_k|=O(\log k)\).

Compare two windows whose centers differ by a bounded number of rows,
extending each array by zero. Write \(D_i\) for the \(H^s\) norm of the
row difference in common coordinates. We allow an additive error of norm
at most \(e_{\rm proj}\) in each row comparison, to cover different Fourier cutoffs below.
For identical Fourier cutoffs take \(e_{\rm proj}=0\). Applying \eqref{pulse:eq:5.13}
along the segment between the inputs and using Bernstein as above gives
\begin{equation}\label{pulse:eq:5.17}
 D_j\le e^{-bk}\max_{i\in N_j}D_i+e_{\rm proj},\qquad b>0 .
\end{equation}
Here the coefficient absorbs the polynomial factors and bounded neighbor count.
Each substitution moves at most \(r\) rows. Starting a bounded distance from
the center, we can substitute at least \(\lfloor c_r\log k\rfloor\) times
before leaving the intersection. Since \(D_i\le P(k)\) at its boundary,
the geometric sum gives
\begin{equation}\label{pulse:eq:5.18}
 D_j\le P(k)e^{-b'k\log k}+\frac{e_{\rm proj}}{1-e^{-bk}},\qquad b'>0 .
\end{equation}
With \(e_{\rm proj}=0\), compare the centers \(i\in N_j\) with the solution on
\(I_j\) and apply \eqref{pulse:eq:5.14}. This gives \eqref{pulse:eq:5.16}
at order \(s\). Bernstein to order \(m\) costs at most \(e^{C_mj}\),
which the first term in \eqref{pulse:eq:5.18} absorbs.
\end{proof}

\begin{proposition}[Comparison on a common cover]\label{pulse:thm:5.4}
\textit{Assume the hypotheses of Proposition~\ref{pulse:thm:5.3} on the common covers, with fixed \(a,q,r\). The untruncated rows and the physical cutoffs commute with lower-cover pullback and zero insertion. For bounded center gaps, require \eqref{pulse:eq:5.13} on a convex polynomial-radius ball containing the two pulled-back trial balls. With the smooth Fourier cutoffs \(S_{L_j}^{[k]}\) of \eqref{pulse:eq:2.36}, the window equation}
\[
 g_j^{[k]}=S_{L_j}^{[k]}H_j^{c_k}(Xg^{[k]})\quad(j\in I_k),
 \qquad g_j^{[k]}=0\quad(j\notin I_k)
\]
\textit{has a unique solution on the trial ball, with polynomial bounds at every fixed order. For bounded \(|k-\ell|+|j-k|\), \(c=\min(c_k,c_\ell)\), and fixed \(m,A\),}
\begin{equation}
 \bigl\|g_j^{[k]}(E^{c_k-c}Y_c)-g_j^{[\ell]}(E^{c_\ell-c}Y_c)\bigr\|_{H^m}
 \le C_{m,A}\varepsilon_k^A. \label{pulse:eq:5.19}
\end{equation}
\end{proposition}
\begin{proof} Apply the contraction and high-norm estimates of Proposition~\ref{pulse:thm:5.3} on each window separately. Lemma~\ref{pulse:thm:2.9} bounds the change of cover. Each solution has polynomial \(H^b\) norm for fixed \(b\). By \eqref{pulse:eq:5.12} the untruncated rows have the same bounds, including order \(s+n\). Thus \eqref{pulse:eq:2.37} gives
\[
 \bigl\|\big[(I-S_{L_i}^{[h]})H_i^{c_h}(Xg^{[h]})\big](E^{c_h-c}Y_c)\bigr\|_{H^s}
 \le P_n(k)L_k^{-n}
 \quad(h=k,\ell).
\]
Here \(L_i/L_k\) and its reciprocal grow at most polynomially on either
window. The untruncated rows agree as maps after pullback and zero insertion.
Subtracting them on the prescribed convex ball and using \eqref{pulse:eq:5.13}
therefore gives \eqref{pulse:eq:5.17} with \(e_{\rm proj}=P_n(k)L_k^{-n}\)
from the two Fourier tails. The bounded cover gap enlarges the Fourier
support only by a fixed factor, so Bernstein still gives the same gain.
Apply \eqref{pulse:eq:5.18} and raise the order by Bernstein at cost
\(C_mL_j^{\max(m-s,0)}\). Increasing \(n\) gives \eqref{pulse:eq:5.19}.
\end{proof}

\subsection{One physical correction}\label{pulse:sec:5.4}
Freeze \((u_*,p_*)\) at Proposition~\ref{pulse:thm:3.14}. On \(I_k\) take the common cover \(c_k=i(k-m_k-h_0)\), the cutoff \(X\) of Lemma~\ref{pulse:thm:5.2a}, and \(s=s_{\rm trial}\), \(L_j=\varepsilon_j^{-1/1320}\). Use the trimmed free source \(\mathbf f_*^{\rm in}\) from \eqref{pulse:eq:4.45}. Each overlapping band restricts the same physical cutoff \eqref{pulse:eq:3.47}.
On \(\supp\chi_{\rm in}\) both edge distances exceed or equal \(a/\kappa\),
while \(\chi_r=1\) for distances at least \(2a/(3\kappa)\), where \(a\)
is the edge cutoff parameter in \eqref{pulse:eq:3.47}.
This gap accommodates the nested radial collars of Theorem~\ref{pulse:thm:5.1}.
The finite map is
\begin{equation}\label{pulse:eq:5.20}
 \Phi_j^{[k]}(g)=S_{L_j}^{[k]}\mathcal R_{I_k,j}
  \big[\mathbf f_*^{\rm in}+E(\mathcal A_{I_k}Xg)+B(\mathcal A_{I_k}Xg,\mathcal A_{I_k}Xg)\big].
\end{equation}
Extend each input by zero outside \(I_k\). The inverse matrices are fixed
before choosing the input. Lemma~\ref{pulse:thm:2.9} bounds the two
pulled-back radii polynomially, and polarization in Lemma~\ref{pulse:thm:5.2a}
gives the row derivative bound on their common convex ball. Thus
Proposition~\ref{pulse:thm:5.4} gives the unique solution on the window
trial ball and polynomial norms at every fixed order. The equations are
quadratic in the real Fourier coefficients. The central components give
the physical free source
\begin{equation}\label{pulse:eq:5.21}
 g_j=g_j^{[j]},\qquad G=\sum_j\zeta_j^2 C_j\,\operatorname{res}_j(X_jg_j).
\end{equation}
We now apply \eqref{pulse:eq:4.56} to the source \(G\) defined in \eqref{pulse:eq:5.21}. The velocity and pressure increments are
\begin{equation}\label{pulse:eq:5.22}
 v=u_*-u_B+\mathcal V G,\qquad
 q=p_*-p_B+\mathcal Q_1G+\mathcal Q_2(G,G).
\end{equation}
At \(\lambda>0\) the sum is locally finite. Its chart representations
give identical physical fields on overlaps. Smooth data on the finitely
many low bands contribute a smooth force away from the vertex.

To find the residual, put \(H_j=\mathcal R_j H\mathcal A\) on any overlap containing the active bands. Proposition~\ref{pulse:thm:2.8} gives, in free coordinates before decoding,
\begin{equation}\label{pulse:eq:5.23}
 H(\mathcal A Xg)-\mathcal A Xg=
 \sum_j\zeta_j^2 C_j\operatorname{res}_j[H_j(Xg)-X_jg_j].
\end{equation}
On band \(k\), \eqref{pulse:eq:5.19} permits replacing the active
central components by those of \(g^{[k]}\) in \eqref{pulse:eq:5.23}
with error \(O_{m,A}(\varepsilon_k^A)\) for every fixed \(m,A\),
also at the higher input order required by the row estimates.
For this single window use the cutoff estimate of
Theorem~\ref{pulse:thm:5.1} on the radial and label selections of
Lemma~\ref{pulse:thm:5.2a}, leaving the independent moment coordinates
untouched by Remark~\ref{pulse:thm:5.2}. For each fixed depth \(J\)
all stencils involved lie in \(I_k\) for large \(k\). The collar,
Fourier and window-boundary errors are therefore bounded, up to
polynomials, by \(\varepsilon_k^{J/40}\), \(\varepsilon_k^{n/1320}\)
and \(e^{-c_m k\log(2+k)}\) by \eqref{pulse:eq:5.9} and
Proposition~\ref{pulse:thm:5.4}. The discarded clock blocks have the
physical Gaussian bound \(P_m(k)e^{-c_m k^2+C_m k}\) after restoring
the carrier weight and taking physical derivatives.

On the middle window, \eqref{pulse:eq:3.17} and \eqref{pulse:eq:4.43}
bound the translated-primitive tails from the initialization and from
\(\mathbf T_1,\mathbf T_2\) by \(P_{m,J'}(k)\varepsilon_k^{\delta_0 J'}\)
in their stated units. In \eqref{pulse:eq:4.57} add the physical
remainder \(\mathbf T_*\) once. It contains the initialization error
and the decoded difference \(\mathbf f_*^{\rm raw}-\mathbf f_*^{\rm in}\),
whose discarded radial edge has the physical Gaussian bound as well.
These are the coefficient and physical Gaussian estimates used in
Section~\ref{pulse:sec:6}. The fields are fixed before \(J,n,J'\) are chosen.
The velocity increments are curls or
axisymmetric swirls, so \(v\) is divergence free. The velocity and
pressure increments are supported in the annulus, preserving the exact
inner and exterior fields. The two moment identities
\eqref{pulse:eq:3.1} hold for the double mean before the physical
evaluation of \(Y\). For example, \(\chi(r)\cos(2\pi n\cdot Y)\),
\(0\ne n\in\mathbb Z^2\), has zero Haar mean. Choosing a nonzero
nonnegative \(\chi\) supported where the evaluated cosine has one
strict sign gives a nonzero physical radial moment.

\subsection{Dimension and scale threshold}\label{pulse:sec:5.5}
By Corollary~\ref{pulse:thm:2.10}, the real dimension satisfies
\begin{equation}\label{pulse:eq:5.25}
 D_k\le C(1+k)^A\varepsilon_k^{-1/220}
       \le C(1+|\log\lambda|)^A\lambda^{-\delta/220}.
\end{equation}
For the window trial radius \(R_k\) in Proposition~\ref{pulse:thm:5.4},
the contraction estimate has the form
\begin{equation}\label{pulse:eq:5.26}
 \sup_{\|g\|\le R_k}\|D\Phi^{[k]}(g)\|\le C_0(1+k)^B e^{-\mu k},
 \qquad \mu=\delta\log(2)/40 .
\end{equation}
Choose \(k_{\rm geom}\) so that the smallest native band in the halo, \(k-\lceil\log(2+k)\rceil-h_0\), exceeds all geometric thresholds when \(k\ge k_{\rm geom}\).
\begin{remark}[Explicit scale threshold]\label{pulse:thm:5.5}
For \(B\ge0\), with the power interpreted as \(1\) at \(B=0\), the coefficient in \eqref{pulse:eq:5.26} is at most \(1/2\) if
\begin{equation}\label{pulse:eq:5.27}
 k\ge\max\!\left\{k_{\rm geom},1,
 \left\lceil\frac{2}{\mu}\log\!\left(2\max\{1,C_0\}e^\mu\max\!\left\{1,\left(\frac{2B}{e\mu}\right)^B\right\}\right)\right\rceil\right\}.
\end{equation}
For \(B>0\), the maximum of \(x^B e^{-\mu x/2}\) on \(x>0\) is \((2B/(e\mu))^B\); for \(B=0\) the exponential is bounded by one. Take \(x=1+k\) in \eqref{pulse:eq:5.26}. The remaining decreasing exponential gives the stated ceiling. The threshold depends only on finite background data. Increasing the orders used in an estimate leaves the finite solution unchanged.
\end{remark}
\section{The physical force and finite-time blowup}\label{fin:blowup}\label{pulse:sec:6}
Let \(u_\infty,p_\infty\) be the fields of Section~\ref{pulse:sec:5}, and set
\(f_\infty=\NS(u_\infty,p_\infty)\). We estimate this force near \((1,0)\)
and localize the velocity by cutting off a vector potential.

Let \(f_{\mathrm G}\) be the sum of the clock and edge cutoff tails, and set \(f_{\mathrm{rem}}=f_\infty-f_{\mathrm G}\). Since \(E_B=0\), only the corrections contribute. On \(B_k\), let \(r_m(k)\) denote the sum of the order-\(m\) coefficient norms of \(f_{\mathrm{rem}}\) after decoding. We use the component units in \eqref{pulse:eq:4.44} and the largest force unit \(F_k\varepsilon_k^{-3/10}\) for the compact translated remainders. By Section~\ref{pulse:sec:5.4}, for fixed positive integers \(J,n,J'\), we have
\begin{equation}\label{pulse:eq:6.1}
 r_m(k)\le P_{m,J,n,J'}(k)
 \left(\varepsilon_k^{J/40}+\varepsilon_k^{n/1320}
       +\varepsilon_k^{\delta_0 J'}+e^{-c k\log(2+k)}\right).
\end{equation}
The four terms come from the iteration on the collars, Fourier truncation on the two covers being compared, integration by parts in the translated primitive, and propagation from the ends of the window, respectively.
We fix the fields before choosing \(J,n,J'\). On the middle window \(\inf\omega>0\).
The inputs to the two initialized tails in \eqref{pulse:eq:3.53} and
\eqref{pulse:eq:3.26} have unit \(F_k\varepsilon_k^{-\delta_0}\omega\)
by \eqref{pulse:eq:3.40} and \eqref{pulse:eq:3.54}. This and the unit in
\eqref{pulse:eq:4.43} are bounded by the largest force unit, so \(J'\)
integrations suffice. Restoring the physical weights gives
\[
 \|\partial^\alpha f_{\mathrm G}\|_{L^\infty(B_k)}
       \le P_\alpha(k)e^{-c_\alpha k^2+C_\alpha k}.
\]

\begin{proposition}[Flatness of the force]\label{pulse:thm:6.1}
For every integer \(a\ge0\), spatial multi-index \(\beta\), and \(N>0\), we have near \((1,0)\)
\begin{equation}\label{pulse:eq:6.2}\csname ltx@label\endcsname{force:pointwise}
 |\partial_t^a\partial_x^\beta f_\infty(t,x)|\le C_{a,\beta,N}\lambda^N .
\end{equation}
Consequently, all one-sided derivatives of \(f_\infty\) vanish at the vertex.
Moreover, on the cylinders \(Q_\rho\) of \eqref{intro:parabolic}, for every integer \(m\ge0\) and every \(N>0\), we have
\begin{equation}\label{pulse:eq:6.3}
 \max_{2a+|\beta|\le m}\rho^{3+2a+|\beta|}
 \|\partial_t^a\partial_x^\beta f_\infty\|_{L^\infty(Q_\rho)}
 \le C_{m,N}\rho^{N+3}\qquad(\rho\ \text{sufficiently small}).
\end{equation}
\end{proposition}
\begin{proof}
The largest force unit on \(B_k\) is
\(F_k\varepsilon_k^{-3/10}=2^{(3+13\delta/10)k}\).
By \eqref{pulse:eq:2.50} and the chain rule on the auxiliary torus, each spatial carrier derivative costs at most \(P(k)2^{(1+\delta/2)k}\).
Each time derivative costs at most \(P(k)2^{(2+3\delta/2)k}\).
Differentiating the cylindrical frame in Cartesian coordinates costs \(O(r^{-1})=O(2^k)\) on the annulus.
The derivatives of the normalized weights on the trimmed region and the number of active labels are bounded by polynomials in \(k\).
We sum over the labels and apply Sobolev embedding with \(m\) sufficiently large in terms of \(a,\beta\). We obtain
\begin{equation}\label{pulse:eq:6.4}
\begin{aligned}
 \|\partial_t^a\partial_x^\beta f_{\mathrm{rem}}\|_{L^\infty(B_k)}
    &\le P_{a,\beta}(k)\,2^{\Gamma_{a,\beta}k}\,r_m(k),\\
 \Gamma_{a,\beta}&=3+13\delta/10+(2+3\delta/2)a+(1+\delta/2)|\beta|.
\end{aligned}
\end{equation}
For the given \(N\), we choose \(J,n,J'\) so that
\begin{equation}\label{pulse:eq:6.5}
  \delta\min\{J/40,\ n/1320,\ \delta_0 J'\}>N+\Gamma_{a,\beta}+1 .
\end{equation}
Since \(\varepsilon_k=2^{-\delta k}\), the extra factor \(2^{-k}\) absorbs all polynomial factors.
The boundary and Gaussian errors decay faster than any power of \(2^{-k}\).
Since \(\lambda\asymp2^{-k}\) on the active annuli, we obtain \eqref{pulse:eq:6.2} there.
Off the annuli the exact background equation gives \(f_\infty=0\).

On \(Q_\rho\), we have \(\lambda^2\le1-t+z^2\le2\rho^2\) by
\eqref{intro:lambda} for small \(\rho\). Multiplying \eqref{pulse:eq:6.2}
by \(\rho^{3+2a+|\beta|}\) gives \eqref{pulse:eq:6.3}. For fixed
geometry, cutoffs, normalizations and lower bounds in \textup{(B3)--(B5)},
both \(C_{m,N}\) and the range \(0<\rho\le\rho_{m,N}\) are controlled by
finitely many background derivatives, including the weighted stress bounds.
We track their orders and the scale threshold in
\hyperref[app:force-constants]{Appendix~\ref*{app:pulse-tools}}.
The bound \(\lambda^2\le1-t+z^2\) also shows that every derivative is \(O((1-t+|x|)^M)\) near the vertex for every \(M>0\).
Hence each derivative extends continuously by zero at the vertex.
Applying the same estimates to successive difference quotients, one-sided in time, shows that all one-sided derivatives vanish there.
\end{proof}

\begin{lemma}[Velocity bounds near the vertex]\label{pulse:thm:6.2}\label{pulse:thm:6.3}
On the part of the active annulus in \(B_k\),
\begin{equation}\label{pulse:eq:6.9}
  |u_\infty|\le P(k)2^{(1+\delta)k},\qquad
  |\nabla u_\infty|\le P(k)2^{(2+\delta)k}.
\end{equation}
Consequently near the vertex there is a fixed \(B\) such that
\begin{equation}\label{pulse:eq:6.10}
 |u_\infty|+(r+\lambda)|\nabla u_\infty|
 \le C(1+|\log(r+\lambda)|)^B(r+\lambda)^{-1-\delta}.
\end{equation}
The exterior velocity and pressure extend smoothly to \(t=1\) from below at \(r>0\).
\end{lemma}
\begin{proof}
For the primary pulses, the carrier derivative costs \(2^{\delta k/2}\) after division
by \(2^k\), which the amplitude factor \(\varepsilon_k^{1/2}\) absorbs.
The polynomial bounds for the window solution, \eqref{pulse:eq:2.77},
\eqref{pulse:eq:3.51}, \eqref{pulse:eq:4.11} and
\eqref{pulse:eq:4.48}--\eqref{pulse:eq:4.51} give \eqref{pulse:eq:6.9}.
Here \(r\asymp\lambda\). For \(R\le R_{\rm in}\), use the bounded normalized
axis-regular quotients with scale \(\lambda\).

For \(R\ge R_e\), the heat integral \eqref{ue:heat-factor} gives
\(0<H_e\le1\) and \(|\xi\partial_\xi H_e|\le(1+\delta)H_e/2\). Thus
\eqref{ue:physical} gives \(|u_B|\le Cr^{-1-\delta}\) and
\(|\nabla u_B|\le Cr^{-2-\delta}\) for \(0<r,\lambda\le1\), including the
factor \(1/r\) from the cylindrical frame. Since \(r+\lambda\asymp r\) here,
this proves \eqref{pulse:eq:6.10}. For \(r\) bounded away from zero, the Gamma
weight has all polynomial moments and the differentiated pressure integrand
is \(O(s^{-3-2\delta})\) at infinity. Differentiation under the integral
gives the extension to \(t=1\).
\end{proof}

\begin{proposition}[Local energy and dissipation]\label{pulse:thm:6.4}
On \(Q_\rho=(1-\rho^2,1)\times B_\rho(0)\), the local kinetic energy and dissipation satisfy
\[
 \sup_{1-\rho^2<t<1}\int_{B_\rho(0)} |u_\infty(t,x)|^2\,dx
 +\int_{Q_\rho}|\nabla u_\infty|^2\,dx\,dt
 \le C\rho^{1-3\delta}(1+|\log\rho|)^{2B}
 \qquad(\rho\downarrow0),
\]
where \(B\) is the exponent in \eqref{pulse:eq:6.10}. In particular, both quantities tend to zero.
\end{proposition}
\begin{proof}
By \eqref{intro:coordinates}, \(r+\lambda\le\rho\) implies
\(r\le\rho\), \(|z|\le\rho^{1-\delta}\), and \(1-t\le\rho^2\).
The spatial volume at each time is therefore \(O(\rho^{3-\delta})\), and the space--time volume is \(O(\rho^{5-\delta})\).
On the shell \(\rho/2<r+\lambda\le\rho\), these volumes and \eqref{pulse:eq:6.10} bound each integral by \(C\rho^{1-3\delta}(1+|\log\rho|)^{2B}\).
On \(Q_\rho\) we have \(r+\lambda\le3\rho\) by \eqref{intro:coordinates}. Summing the shell bounds at radii \(3\rho,3\rho/2,3\rho/4,\ldots\) gives the claim since \(1-3\delta>0\).
\end{proof}

To pass from the local estimate to a divergence-free velocity on \(\mathbb R^3\), we multiply its potential by the cutoff before taking the curl.
This requires a potential smooth up to \(t=1\) away from the vertex. At \(t=1,z\ne0\), the implicit derivative in \eqref{intro:lambda} stays positive, so the profiles remain smooth, including on the axis. Near a point with \(t=1,z=0,r>0\), the fields agree with the exact exterior of Lemma~\ref{pulse:thm:6.2} since \(R=r^2/(2\lambda^2)\to\infty\).

Use the background streamfunction of Section~\ref{fin:section},
so \(r\Psi_B=\int_0^r s\,u_B^z\,ds\). Fix \(r_e\) in the exact exterior
throughout a small working rectangle in \(t,z\). There
\(\partial_z(r_e\Psi_B)=-r_e u_B^r=1\), hence
\begin{equation}\label{pulse:eq:6.11}\csname ltx@label\endcsname{fin:mixed-background-potential}
 \int_0^{r_e}s u_B^z(t,s,z)\,ds=z+C(t),\qquad
 \Psi_B=\frac{z+C(t)}r-\frac1r\int_r^{r_e}s u_B^z(t,s,z)\,ds .
\end{equation}
Evaluating the first identity in \eqref{pulse:eq:6.11} at a fixed \(z\ne0\) shows that \(C(t)\) is smooth up to \(t=1\). For \(r>0\) the second identity gives the same regularity for \(\Psi_B\), since the interval of integration stays away from the axis.
On the axis with \(z\ne0\), the defining integral gives
\(\Psi_B/r=\int_0^1 s\,u_B^z(t,sr,z)\,ds\), which is smooth and even in \(r\).
Thus \(\Psi_B e_\theta=(\Psi_B/r)(-x_2,x_1,0)\) is smooth in Cartesian coordinates.

Each velocity correction is a curl or an axisymmetric swirl. Adding the curl potentials to \(\Psi_B e_\theta\) and the swirls to the background swirl gives
\begin{equation}\label{pulse:eq:6.12}\csname ltx@label\endcsname{fin:full-potential}
 u_\infty=\curl A_\infty+b_\infty e_\theta,\qquad \partial_\theta b_\infty=0 .
\end{equation}
The sums are smooth wherever \(\lambda>0\) by local finiteness. Near \(t=1,z=0,r>0\) the corrections vanish, so \(A_\infty\) is given by the smooth potential in \eqref{pulse:eq:6.11} and \(b_\infty e_\theta\) by the exterior swirl.

\begin{proof}[Proof of Corollary~\ref{intro:blowup}]
\textbf{Blowup from rest.}
Choose a smooth axisymmetric cutoff \(\eta\) with compact spatial support in the working neighborhood, zero up to some \(t_0\in(0,1)\) and one near \((1,0)\), and for \(t<1\) set
\begin{equation}\label{pulse:eq:6.13}\csname ltx@label\endcsname{fin:localized-fields}
 u=\curl(\eta A_\infty)+\eta b_\infty e_\theta,\qquad
 p=\eta p_\infty,\qquad f=\NS(u,p).
\end{equation}
The curl and the axisymmetric swirl are divergence free. On the support of
derivatives of \(\eta\), the potentials and \(p_\infty\) are smooth up to \(t=1\).
Proposition~\ref{pulse:thm:6.4} therefore gives the energy and dissipation
bounds, and the axis trace \eqref{intro:blowup-trace} gives the Type~II rate.
Near \((1,0)\), \(f=f_\infty\) satisfies \eqref{intro:parabolic} and its
one-sided jet vanishes at that point by Proposition~\ref{pulse:thm:6.1}.
The force is thus smooth up to the entire terminal face. Extend it in time
at each fixed \(x\) by Section~\ref{pulse:sec:2.1.4} and apply a smooth
cutoff at later times. This preserves the spatial support and gives
\(f\in C_c^\infty((0,\infty)\times\mathbb R^3)\).

Finally, restart at \(1-\rho^2\) and rescale the remaining interval to \((0,1)\).
The initial velocity \(\rho u(1-\rho^2,\rho\,\cdot)\) is smooth and compactly
supported. Its axial component at the origin is \(j_0\rho^{-\delta}\) by
\eqref{intro:blowup-trace}. Each fixed rescaling preserves finiteness of
energy and dissipation and the Type~II exponent.
For \(f_\rho(s,y)=\rho^3 f(1-\rho^2+\rho^2s,\rho y)\),
\eqref{intro:parabolic} gives \(O_{m,N}(\rho^{N+3})\) in
\(C^m_{\rm par}((0,1)\times B_1)\). On the whole space we have
\[
\begin{aligned}
 \|f_\rho\|_{L^{4/3}(0,1;L^2)}
 &=\|f\|_{L^{4/3}(1-\rho^2,1;L^2)}
 \le \rho^{3/2}\sup_{0<t<1}\|f(t)\|_{L^2},\\
 \|f_\rho\|_{L^\infty((0,1)\times\mathbb R^3)}
 &=\rho^3\|f\|_{L^\infty((1-\rho^2,1)\times\mathbb R^3)} .
\end{aligned}
\]
The same bounds apply to the smooth compactly supported force in
\cite[Theorem~1.1]{OpenAI-NS}. Squaring the critical norm gives
\[
 \|f_\rho\|_{L^{4/3}(0,1;L^2)}^2
 +\|f_\rho\|_{L^\infty((0,1)\times\mathbb R^3)}=O(\rho^3).
\]
In rescaled variables, these solutions can start with arbitrarily large speed
at the origin and develop a Type~II singularity one unit of time later under
forcing arbitrarily small in both norms.
\end{proof}
\paragraph{AI assistance.}
The authors proposed the exact inner core, the mixed exterior, and the direct
connection between them, and posed the existence problem studied here.
Discussions with AI and reading OpenAI's manuscript helped rule out many
unsuccessful approaches and led to the construction presented in this paper.
AI tools also assisted in drafting and polishing the article.
\paragraph{Acknowledgement.}
ZL was supported in part by NSFC (No.~12431007), NSFC Excellent Research Group Project (No.~62588101), NSFC (No.~12494544) and New Cornerstone Science Foundation through the XPLORER PRIZE. XR was supported in part by NSFC (No.~12601444) and the National Key R\&D Program of China (No.~2023YFA1010700).
\appendix
\makeatletter
\renewcommand{\section@cntformat}{Appendix~\thesection.\quad}
\makeatother
\section{Moment corrections and Gaussian integrals}\label{app:tools}

\begin{lemma}[Gram inverse on a fixed window]\label{pr:gram}
Let $n=4$ or $5$, and let $\omega_1,\ldots,\omega_n$ be linearly independent
real-analytic functions on an open interval, with one or two components.
Fix a bounded open subinterval $I$ whose closure lies inside it, and
a smooth $\rho$, positive on $I$ and zero outside $\overline I$.
The matrix $G_{ij}=\int_I\rho \omega_i\cdot \omega_j\,dx$ is positive definite.
For any $d\in\mathbb R^n$, the coefficients $c=G^{-1}d$ therefore give
the unique correction of the form $\rho\sum_j c_j\omega_j$ with moments $d_i$
against the rows $\omega_i$.
Keep the test profiles $\rho \omega_j$ fixed and let the rows
$\widetilde\omega_i(\cdot,Z,\lambda)$ vary smoothly. Define
$B_{ij}(Z,\lambda)=\int_I\widetilde\omega_i\cdot(\rho \omega_j)\,dx$.
If $\sup_{Z,\lambda}\|G^{-1}(B-G)\|\le1/2$, then $B^{-1}$ is smooth
and $\sup_{Z,\lambda}\|B^{-1}\|\le2\|G^{-1}\|$. Each fixed composition
of $\partial_Z$ and $\lambda\partial_\lambda$ applied to $B^{-1}$ is
uniformly bounded if the corresponding row derivatives up to that order
are uniformly bounded on $I$.
\end{lemma}
\begin{proof}
If $c^TGc=\int_I\rho|\sum_j c_j\omega_j|^2\,dx=0$, then $\sum_j c_j\omega_j=0$ on
$I$. Analyticity extends this identity to the whole interval, so $c=0$.
The Neumann series for $(I+G^{-1}(B-G))^{-1}$ gives the uniform inverse bound.
Repeated differentiation of $B^{-1}B=I$ and the integral for $B$ gives the
claimed bounds for $\partial_Z$ and $\lambda\partial_\lambda$.
\end{proof}

\begin{lemma}[Small quadratic correction]\label{pr:quadratic}
Let $A$ be a bounded invertible operator on a Banach space,
$B=\|A^{-1}\|$, and let $Q$ be bilinear with
$\|Q(h,k)\|\le C\|h\|\|k\|$. If
$8B^2C\|d\|\le1$, then
$Ah+Q(h,h)=-d$ has a unique solution in $\|h\|\le2B\|d\|$.
It depends smoothly on any smooth auxiliary parameters for which
these bounds hold uniformly.
\end{lemma}
\begin{proof}
On this ball, $h\mapsto-A^{-1}(d+Q(h,h))$ has norm at most
$B\|d\|+4B^3C\|d\|^2\le3B\|d\|/2$ and Lipschitz constant at most
$4B^2C\|d\|\le1/2$. Thus it has a unique fixed point. The same bound gives
invertibility of $A+Q(h,\cdot)+Q(\cdot,h)$ and hence smooth dependence.
\end{proof}

\paragraph{Gaussian integrals at a flat edge.}
For $c>0$, $j\in\RR$ and smooth $a=a(s,\lambda,Z)$ with $\lambda>0$,
substituting $s^{-2}=d^{-2}+\tau$ gives
\begin{equation}\label{fin:gaussian-smooth-quotient}
\begin{aligned}
 &e^{c/d^2}d^{j-3}\int_0^d e^{-c/s^2}s^{-j}a(s,\lambda,Z)\,ds\\
 &\qquad=\frac12\int_0^\infty e^{-c\tau}
 (1+d^2\tau)^{(j-3)/2}
 a\!\left(\frac{d}{\sqrt{1+d^2\tau}},\lambda,Z\right)d\tau.
\end{aligned}
\end{equation}
Differentiated integrands are bounded by polynomials in $\tau$ times
$e^{-c\tau}$. Hence for $\lambda>0$ the right side is smooth in $(d,\lambda,Z)$ up to
$d=0$, where it equals $a(0,\lambda,Z)/(2c)$.
On $0\le d\le d_0$, $|Z|\le1$, bounds for fixed compositions of
$\partial_d,\lambda\partial_\lambda,\partial_Z$ follow from the
corresponding bounds for $a$ on $0\le s\le d_0$.
For $0<d\le d_0$,
\begin{equation}\label{fin:gaussian-primitive}
 \int_0^d e^{-c/s^2}s^{-j}\,ds
       \le C_j e^{-c/d^2}d^{3-j}.
\end{equation}
Here $C_j$ may also depend on $c$ and $d_0$. A smooth positive factor in the
integrand is included in $a$.

\begin{proof}[Proof of Lemma~\ref{rr:analytic}]
Extend $U$ smoothly with compact $R$ support and set
\begin{equation}\label{rr:gaussian-formula}
 \widetilde U(R,Z)=\frac1{\sqrt\pi\,\varepsilon\rho(R)}
       \int_{\mathbb R}
       \exp\left(-\frac{(R-s)^2}{\varepsilon^2\rho(R)^2}\right)
                         U(s,Z)\,ds.
\end{equation}
The kernel is holomorphic with nonzero denominator in a complex
neighborhood of each interior real $R$. The integrand has compact support.
We differentiate under the integral to prove radial analyticity.

Substitute $s=R+\varepsilon\rho(R)t$, differentiate, and subtract
$\partial_R^i\partial_Z^kU(R,Z)$. Each term is a
translation difference or contains \(\varepsilon\rho\) or one of its radial derivatives.
Since \(|\partial_R^k\rho|\le C_k\sqrt\rho\), the mean value theorem
and the finite Gaussian moments give \eqref{rr:analytic-error}. The difference and all its derivatives
extend by zero at both endpoints.
\end{proof}
\section{Short extension, Fourier localization and force bounds}\label{app:pulse-tools}

\paragraph{Proof of Lemma~\ref{pulse:thm:2.6}.}
Let \(b_j=2^j\), \(j\ge0\). Seeley's coefficients \cite{Seeley1964} are
\[
 a_j=\prod_{\ell\ne j}\frac{-1-b_\ell}{b_j-b_\ell}.
\]
In a product truncated at any \(N\ge j\), the factors with \(\ell<j\)
contribute \(2^{-j(j+1)/2}\) after extracting \(2^{\ell-j}\) from each.
The remaining product is bounded by
\(\prod_{\ell\ge0}(1+2^{-\ell})\prod_{n\ge1}(1-2^{-n})^{-2}\).
Thus \(|a_j|\le C2^{-j(j+1)/2}\). Interpolating \(x^m\) at
\(b_0,\ldots,b_N\) and evaluating at \(-1\) gives the last identity in
\eqref{pulse:eq:2.27} by dominated convergence.

Keep \(f(x)\) for \(x\ge0\). For \(x<0\), set
\[
 E_hf(x)=\sum_{j\ge0}a_j\chi(-b_jx/h)f(-b_jx),
 \qquad \chi\in C_c^\infty((-1,1)),\quad \chi=1\text{ near }0 .
\]
Only samples in \([0,h]\) occur, and the extension vanishes for \(x\le-h\).
By \eqref{pulse:eq:2.27} every differentiated series converges up to zero
and agrees there with the corresponding derivative of \(f\).
After \(\ell\le m\) normal derivatives, the substitution \(y=-b_jx\)
bounds the \(j\)th term in \(L^2\) by
\(C_m|a_j|b_j^{\ell-1/2}h^{-m}\|f\|_{H^m(0,h)}\).
Summation gives the bound in Lemma~\ref{pulse:thm:2.6}.
The formula acts pointwise in the passive variables, so passive derivatives
commute with extension and their supports are preserved. \hfill\(\square\)

\paragraph{Fourier cutoffs with separated supports.}
After fixed period rescaling, the kernel of the multiplier \(S_L\) in
\eqref{pulse:eq:2.36} is the periodization of \(L^d\check{s}(Lx)\).
Integration by parts in the Fourier integral gives
\(\left|\partial^\beta\check{s}(x)\right|\le C_{\beta,q}(1+|x|)^{-q}\)
for every \(q\). If two cutoffs have supports a torus distance \(b>0\)
apart, every term of the periodized kernel between them has argument
of length at least \(Lb\). Taking \(q>N+d+m\) and summing the translates
therefore bounds the kernel and its first \(m\) derivatives by
\(C_{m,N}L^{-N}\) times a polynomial in \(1+b^{-1}\), for \(L\ge1\).
Schur's test and the product rule give the same \(H^m\)-operator bound
with polynomial dependence on the cutoff \(C^m\) norms.
\paragraph{Constants in the local force estimate.}\phantomsection\label{app:force-constants}
Consider backgrounds satisfying Proposition~\ref{fin:mixed-interface} on a common domain. Fix \(\delta\), the annulus, windows, edge weight, normalizations
and cutoffs, including a common neighborhood on which the final cutoff
equals one. Fix also the lower bounds in \textup{(B3)--(B5)}. Fix \(m,N\)
and choose \(J,n,J'\) in \eqref{pulse:eq:6.5} for all \(2a+|\beta|\le m\).
Put \(s=s_{\rm trial}\) and \(b=\max\{s,m+5\}\), allowing one derivative
for decoding and four for Sobolev embedding in at most six variables.
To apply a row at order \(b\), compare its window inputs at order \(b+33\).
Bernstein then costs \(L^{b+33-s}\), so use the Fourier tail of order
\(n+b+33-s\) in the proof of Proposition~\ref{pulse:thm:5.4}.
The untruncated row needs \(b+n+33\) derivatives and its input needs \(b+n+66\).
The collar iteration \eqref{pulse:eq:5.9} stays at the same order at each step.
For the translated tails, \eqref{pulse:eq:3.26} and \eqref{pulse:eq:4.43}
use \(2J'\) extra derivatives of their inputs.

To bound the background order, start the row and initialization estimates at
\(b+n+66+2J'+10\), including the ten derivatives in the pressure norm
of \eqref{pulse:eq:3.25}. Follow their formulas backwards. Add one for
each first derivative, curl or directional inverse in Sobolev norm, two for
viscosity, and four at each use of \(H^{j+4}\hookrightarrow C^j\).
Products, ordinary primitives and matrix inverses require the same order
of their arguments. The differentiated clock equation \eqref{pulse:eq:2.59}
requires the same order of its coefficients, including the first derivatives
in its frame. The maximum order reached on the background factors suffices in
\textup{(B2)--(B3)} of Proposition~\ref{fin:mixed-interface}.
For instance, the square root in \eqref{pulse:eq:2.72} uses
\(g_\gamma^2\ge ck\omega\); the weighted bounds \eqref{intro:stress-bounds}
give only powers of \(d_e^{-1}\), hence powers of \(k\) on the trimmed region.
With these data fixed, the bounds in \textup{(B2)} and the finitely many
weighted constants \(C_j,q_j\) in \eqref{intro:stress-bounds} control
the coefficients and degrees of the finitely many polynomials above.

Start with \eqref{pulse:eq:5.27} and increase the band index until every
required high-order factor \(P_j(k)(1+R_k)\varepsilon_k^{1/40}\) is at most
\(1/2\) on all bands of the window. This gives the high norms by absorption
in \eqref{pulse:eq:5.3}. Increase the index once more to accommodate the
\(3J+1\) collars and to absorb the polynomials, Gaussian tails and
\(e^{-ck\log(2+k)}\) in \eqref{pulse:eq:6.4}--\eqref{pulse:eq:6.5}.
The clock Gaussian is uniform by \eqref{pulse:eq:2.62}. On the discarded edge
\(d_e\le C/k\), half the Gaussian in \eqref{intro:edge-weight} absorbs any
fixed power of \(d_e^{-1}\) and the rest is \(O(e^{-c k^2})\).
Thus the resulting index \(K\) and the constant
in \eqref{pulse:eq:6.2} are common to the stated class.

Choose \(\rho_{m,N}\le 2^{-K}/4\) inside the common neighborhood where
the final cutoff is one. Since \(\lambda\le\sqrt2\rho\) on \(Q_\rho\)
and \(2^k\lambda>1/2\) on \(B_k\), every active band meeting \(Q_\rho\)
has \(k\ge K\) for \(\rho\le\rho_{m,N}\). This gives a common
\(C_{m,N}\) and a common range of scales in \eqref{intro:parabolic}.

\end{document}